\input xy
\xyoption{all}


\magnification=1200
\hsize=13.50cm    
\vsize=18cm       
\hoffset=-2mm
\voffset=.8cm
\parindent=12pt   \parskip=0pt     
\hfuzz=1pt

\pretolerance=500 \tolerance=1000  \brokenpenalty=5000

\catcode`\@=11

\font\eightrm=cmr8
\font\eighti=cmmi8
\font\eightsy=cmsy8
\font\eightbf=cmbx8
\font\eighttt=cmtt8
\font\eightit=cmti8
\font\eightsl=cmsl8
\font\sevenrm=cmr7
\font\seveni=cmmi7
\font\sevensy=cmsy7
\font\sevenbf=cmbx7
\font\sevenit=cmti10 at 7pt
\font\sixrm=cmr6
\font\sixi=cmmi6
\font\sixsy=cmsy6
\font\sixbf=cmbx6

\font\douzebf=cmbx10 at 12pt

\font\twelvebf=cmbx10 at 12pt

\font\tencal=eusm10

\font\sevencal=eusm7

\font\fivecal=eusm5
\newfam\calfam
\textfont\calfam=\tencal
\scriptfont\calfam=\sevencal
\scriptscriptfont\calfam=\fivecal
\def\cal#1{{\fam\calfam\relax#1}}

\skewchar\eighti='177 \skewchar\sixi='177
\skewchar\eightsy='60 \skewchar\sixsy='60

\def\tenpoint{%
  \textfont0=\tenrm \scriptfont0=\sevenrm
  \scriptscriptfont0=\fiverm
  \def\rm{\fam\z@\tenrm}%
  \textfont1=\teni  \scriptfont1=\seveni
  \scriptscriptfont1=\fivei
  \def\oldstyle{\fam\@ne\teni}\let\old=\oldstyle
  \textfont2=\tensy \scriptfont2=\sevensy
  \scriptscriptfont2=\fivesy
  \textfont\itfam=\tenit
  \def\it{\fam\itfam\tenit}%
  \textfont\slfam=\tensl
  \def\sl{\fam\slfam\tensl}%
  \textfont\bffam=\tenbf
  \scriptfont\bffam=\sevenbf
  \scriptscriptfont\bffam=\fivebf
  \def\bf{\fam\bffam\tenbf}%
  \textfont\ttfam=\tentt
  \def\tt{\fam\ttfam\tentt}%
  \abovedisplayskip=12pt plus 3pt minus 9pt
  \belowdisplayskip=\abovedisplayskip
  \abovedisplayshortskip=0pt plus 3pt
  \belowdisplayshortskip=4pt plus 3pt 
  \smallskipamount=3pt plus 1pt minus 1pt
  \medskipamount=6pt plus 2pt minus 2pt
  \bigskipamount=12pt plus 4pt minus 4pt
  \normalbaselineskip=12pt
  \setbox\strutbox=\hbox{\vrule height8.5pt depth3.5pt width0pt}%
  \let\bigf@nt=\tenrm
  \let\smallf@nt=\sevenrm
  \normalbaselines\rm}
  
\def\eightpoint{%
  \textfont0=\eightrm \scriptfont0=\sixrm
  \scriptscriptfont0=\fiverm
  \def\rm{\fam\z@\eightrm}%
  \textfont1=\eighti  \scriptfont1=\sixi
  \scriptscriptfont1=\fivei
  \def\oldstyle{\fam\@ne\eighti}\let\old=\oldstyle
  \textfont2=\eightsy \scriptfont2=\sixsy
  \scriptscriptfont2=\fivesy
  \textfont\itfam=\eightit
  \def\it{\fam\itfam\eightit}%
  \textfont\slfam=\eightsl
  \def\sl{\fam\slfam\eightsl}%
  \textfont\bffam=\eightbf
  \scriptfont\bffam=\sixbf
  \scriptscriptfont\bffam=\fivebf
  \def\bf{\fam\bffam\eightbf}%
  \textfont\ttfam=\eighttt
  \def\tt{\fam\ttfam\eighttt}%
  \abovedisplayskip=9pt plus 3pt minus 9pt
  \belowdisplayskip=\abovedisplayskip
  \abovedisplayshortskip=0pt plus 3pt
  \belowdisplayshortskip=3pt plus 3pt 
  \smallskipamount=2pt plus 1pt minus 1pt
  \medskipamount=4pt plus 2pt minus 1pt
  \bigskipamount=9pt plus 3pt minus 3pt
  \normalbaselineskip=9pt
  \setbox\strutbox=\hbox{\vrule height7pt depth2pt width0pt}%
  \let\bigf@nt=\eightrm
  \let\smallf@nt=\sixrm
  \normalbaselines\rm}
\tenpoint


\def\pc#1{\bigf@nt#1\smallf@nt}

\catcode`\;=\active
\def;{\relax\ifhmode\ifdim\lastskip>\z@\unskip\fi \kern\fontdimen2 
 -1.2 \fontdimen3 \string;}

\catcode`\:=\active
\def:{\relax\ifhmode\ifdim\lastskip>\z@\unskip\fi\penalty\@M\ 
\fi\string:}

\catcode`\!=\active
\def!{\relax\ifhmode\ifdim\lastskip>\z@ \unskip\fi\kern\fontdimen2 
 -1.1 \fontdimen3 \string!}

\catcode`\?=\active
\def?{\relax\ifhmode\ifdim\lastskip>\z@ \unskip\fi\kern\fontdimen2 
 -1.1 \fontdimen3 \string?}

\frenchspacing

\newtoks\auteurcourant
\auteurcourant={\hfil}

\newtoks\titrecourant
\titrecourant={\hfil}

\newtoks\hautpagetitre
\hautpagetitre={\hfil}

\newtoks\baspagetitre
\baspagetitre={\hfil}

\newtoks\hautpagegauche     
\hautpagegauche={\eightpoint\rlap{\folio}\hfil\the
\auteurcourant\hfil}

\newtoks\hautpagedroite     
\hautpagedroite={\eightpoint\hfil\the
\titrecourant\hfil\llap{\folio}}

\newtoks\baspagegauche
\baspagegauche={\hfil} 

\newtoks\baspagedroite
\baspagedroite={\hfil}

\newif\ifpagetitre
\pagetitretrue  


\headline={\ifpagetitre\the\hautpagetitre
\else\ifodd\pageno\the\hautpagedroite\else\the
\hautpagegauche\fi\fi}

\footline={\ifpagetitre\the\baspagetitre\else
\ifodd\pageno\the\baspagedroite\else\the
\baspagegauche\fi\fi
\global\pagetitrefalse}

\def\raggedbottom{\topskip 10pt plus 36pt\r@ggedbottomtrue}

\def\point{\raise.2ex\hbox{\douzebf .}}
\def\pointir{\unskip . --- \ignorespaces}


\def\Medbreak{\vskip-\lastskip\medbreak}

\def\rem#1\endrem{%
\Medbreak {\it#1\unskip} : }


\long\def\th#1 #2\enonce#3\endth{%
\Medbreak {\pc#1} {#2\unskip}\pointir{\it #3}\medskip}


\long\def\thm#1 #2\enonce#3\endthm{%
\Medbreak {\pc#1} {#2\unskip}\pointir{\it #3}\medskip}

\def\decale#1{\smallbreak\hskip 28pt\llap{#1}\kern 5pt}

\def\decaledecale#1{\smallbreak\hskip 34pt\llap{#1}\kern 5pt}

\let\@ldmessage=\message

\def\message#1{{\def\pc{\string\pc\space}%
\def\'{\string'}\def\`{\string`}%
\def\^{\string^}\def\"{\string"}%
\@ldmessage{#1}}}

\def\({{\rm (}}

\def\){{\rm )}}


\def\up#1{\raise 1ex\hbox{\smallf@nt#1}}

\def\diagram#1{\def\normalbaselines{\baselineskip=0pt 
\lineskip=5pt}\matrix{#1}}


\def\longmaprightover#1#2{\smash{\mathop{\hbox to#2 
{\rightarrowfill}}\limits^{\scriptstyle#1}}}

\def\longmapleftover#1#2{\smash{\mathop{\hbox to#2 
{\leftarrowfill}}\limits^{\scriptstyle#1}}}

\def\longmaprightunder#1#2{\smash{\mathop{\hbox to#2 
{\rightarrowfill}}\limits_{\scriptstyle#1}}}

\def\longmapleftunder#1#2{\smash{\mathop{\hbox to#2 
{\leftarrowfill}}\limits_{\scriptstyle#1}}}

\def\longhookrightarrowover#1#2{\smash{\mathop{\lhook\joinrel 
\mathrel{\hbox to #2{\rightarrowfill}}}\limits^{\scriptstyle#1}}}

\def\longhookrightarrowunder#1#2{\smash{\mathop{\lhook\joinrel 
\mathrel{\hbox to #2{\rightarrowfill}}}\limits_{\scriptstyle#1}}}

\def\longhookleftarrowover#1#2{\smash{\mathop{{\hbox to #2 
{\leftarrowfill}}\joinrel\kern -0.9mm\mathrel\rhook} 
\limits^{\scriptstyle#1}}}

\def\longhookleftarrowunder#1#2{\smash{\mathop{{\hbox to #2 
{\leftarrowfill}}\joinrel\kern -0.9mm\mathrel\rhook} 
\limits_{\scriptstyle#1}}}

\def\longtwoheadrightarrowover#1#2{\smash{\mathop{{\hbox to #2 
{\rightarrowfill}}\kern -3.25mm\joinrel\mathrel\rightarrow} 
\limits^{\scriptstyle#1}}}

\def\longtwoheadrightarrowunder#1#2{\smash{\mathop{{\hbox to #2 
{\rightarrowfill}}\kern -3.25mm\joinrel\mathrel\rightarrow} 
\limits_{\scriptstyle#1}}}

\def\longtwoheadleftarrowover#1#2{\smash{\mathop{\joinrel\mathrel 
\leftarrow\kern -3.8mm{\hbox to #2{\leftarrowfill}}} 
\limits^{\scriptstyle#1}}}

\def\longtwoheadleftarrowunder#1#2{\smash{\mathop{\joinrel\mathrel 
\leftarrow\kern -3.8mm{\hbox to #2{\leftarrowfill}}} 
\limits_{\scriptstyle#1}}}


\def\longmapsto#1{\mapstochar\mathrel{\joinrel \kern-0.2mm\hbox to 
#1mm{\rightarrowfill}}}

\def\og{\leavevmode\raise.3ex\hbox{$\scriptscriptstyle 
\langle\!\langle\,$}}
\def\fg{\leavevmode\raise.3ex\hbox{$\scriptscriptstyle 
\,\rangle\!\rangle$}}

\def\section#1#2{\vskip 5mm {\bf {#1}. {#2}}\vskip 5mm}
\def\subsection#1#2{\vskip 3mm {\it #2}\vskip 3mm}

\catcode`\@=12

\showboxbreadth=-1  \showboxdepth=-1


\message{`lline' & `vector' macros from LaTeX}

\def\Grille{\setbox13=\vbox to 5\unitlength{\hrule width 109mm \vfill}
\setbox13=\vbox to 65mm
{\offinterlineskip\leaders\copy13\vfill\kern-1pt\hrule} 
\setbox14=\hbox to 5\unitlength{\vrule height 65mm\hfill} 
\setbox14=\hbox to 109mm{\leaders\copy14\hfill\kern-2mm \vrule height 
65mm}
\ht14=0pt\dp14=0pt\wd14=0pt \setbox13=\vbox to 0pt
{\vss\box13\offinterlineskip\box14} \wd13=0pt\box13}

\def\rule(#1,#2)\dir(#3,#4)\long#5{%
\noalign{\leftput(#1,#2){\lline(#3,#4){#5}}}}

\def\arrow(#1,#2)\dir(#3,#4)\length#5{%
\noalign{\leftput(#1,#2){\vector(#3,#4){#5}}}}

\def\put(#1,#2)#3{\noalign{\setbox1=\hbox{%
\kern #1\unitlength \raise #2\unitlength\hbox{$#3$}}%
\ht1=0pt \wd1=0pt \dp1=0pt\box1}}

\catcode`@=11

\def\{{\relax\ifmmode\lbrace\else$\lbrace$\fi}

\def\}{\relax\ifmmode\rbrace\else$\rbrace$\fi}

\def\newcount{\alloc@0\count\countdef\insc@unt}

\def\newdimen{\alloc@1\dimen\dimendef\insc@unt}

\def\newwrite{\alloc@7\write\chardef\sixt@@n}

\newwrite\@unused
\newcount\@tempcnta
\newcount\@tempcntb
\newdimen\@tempdima
\newdimen\@tempdimb
\newbox\@tempboxa

\def\@spaces{\space\space\space\space}

\def\@whilenoop#1{}

\def\@whiledim#1\do #2{\ifdim #1\relax#2\@iwhiledim{#1\relax#2}\fi}

\def\@iwhiledim#1{\ifdim #1\let\@nextwhile=\@iwhiledim 
\else\let\@nextwhile=\@whilenoop\fi\@nextwhile{#1}}

\def\@badlinearg{\@latexerr{Bad \string\line\space or \string\vector 
\space argument}}

\def\@latexerr#1#2{\begingroup 
\edef\@tempc{#2}\expandafter\errhelp\expandafter{\@tempc}%

\def\@eha{Your command was ignored.^^JType \space I <command> <return> 
\space to replace it with another command,^^Jor \space <return> \space 
to continue without it.}

\def\@ehb{You've lost some text. \space \@ehc}

\def\@ehc{Try typing \space <return> \space to proceed.^^JIf that 
doesn't work, type \space X <return> \space to quit.}

\def\@ehd{You're in trouble here.  \space\@ehc}

\typeout{LaTeX error.  \space See LaTeX manual for explanation.^^J 
\space\@spaces\@spaces\@spaces Type \space H <return> \space for 
immediate help.}\errmessage{#1}\endgroup}

\def\typeout#1{{\let\protect\string\immediate\write\@unused{#1}}}

\font\tenln = line10
\font\tenlnw = linew10

\newdimen\@wholewidth
\newdimen\@halfwidth
\newdimen\unitlength 

\unitlength =1pt

\def\thinlines{\let\@linefnt\tenln \let\@circlefnt\tencirc 
\@wholewidth\fontdimen8\tenln \@halfwidth .5\@wholewidth}

\def\thicklines{\let\@linefnt\tenlnw \let\@circlefnt\tencircw 
\@wholewidth\fontdimen8\tenlnw \@halfwidth .5\@wholewidth}

\def\linethickness#1{\@wholewidth #1\relax \@halfwidth .5
\@wholewidth}

\newif\if@negarg

\def\lline(#1,#2)#3{\@xarg #1\relax \@yarg #2\relax 
\@linelen=#3\unitlength \ifnum\@xarg =0 \@vline \else \ifnum\@yarg =0 
\@hline \else \@sline\fi \fi}

\def\@sline{\ifnum\@xarg< 0 \@negargtrue \@xarg -\@xarg \@yyarg 
-\@yarg \else \@negargfalse \@yyarg \@yarg \fi
\ifnum \@yyarg >0 \@tempcnta\@yyarg \else \@tempcnta - \@yyarg \fi
\ifnum\@tempcnta>6 \@badlinearg\@tempcnta0 \fi
\setbox\@linechar\hbox{\@linefnt\@getlinechar(\@xarg,\@yyarg)}%
\ifnum \@yarg >0 \let\@upordown\raise \@clnht\z@ 
\else\let\@upordown\lower \@clnht \ht\@linechar\fi
\@clnwd=\wd\@linechar
\if@negarg \hskip -\wd\@linechar \def\@tempa{\hskip -2\wd \@linechar} 
\else \let\@tempa\relax \fi
\@whiledim \@clnwd <\@linelen \do {\@upordown\@clnht\copy\@linechar 
\@tempa \advance\@clnht \ht\@linechar \advance\@clnwd \wd\@linechar}%
\advance\@clnht -\ht\@linechar \advance\@clnwd -\wd\@linechar 
\@tempdima\@linelen\advance\@tempdima -\@clnwd 
\@tempdimb\@tempdima\advance\@tempdimb -\wd\@linechar
\if@negarg \hskip -\@tempdimb \else \hskip \@tempdimb \fi
\multiply\@tempdima \@m\@tempcnta \@tempdima \@tempdima \wd\@linechar 
\divide\@tempcnta \@tempdima \@tempdima \ht\@linechar 
\multiply\@tempdima \@tempcnta \divide\@tempdima \@m \advance\@clnht 
\@tempdima
\ifdim \@linelen <\wd\@linechar \hskip \wd\@linechar 
\else\@upordown\@clnht\copy\@linechar\fi}

\def\@hline{\ifnum \@xarg <0 \hskip -\@linelen \fi
\vrule height \@halfwidth depth \@halfwidth width \@linelen
\ifnum \@xarg <0 \hskip -\@linelen \fi}

\def\@getlinechar(#1,#2){\@tempcnta#1\relax
\multiply\@tempcnta 8\advance\@tempcnta -9
\ifnum #2>0 \advance\@tempcnta #2\relax
 \else\advance\@tempcnta -#2\relax\advance\@tempcnta 64 \fi
\char\@tempcnta}

\def\vector(#1,#2)#3{\@xarg #1\relax \@yarg #2\relax 
\@linelen=#3\unitlength
\ifnum\@xarg =0 \@vvector \else \ifnum\@yarg =0 \@hvector \else 
\@svector\fi \fi}

\def\@hvector{\@hline\hbox to 0pt{\@linefnt \ifnum \@xarg <0 
\@getlarrow(1,0)\hss\else \hss\@getrarrow(1,0)\fi}}

\def\@vvector{\ifnum \@yarg <0 \@downvector \else \@upvector \fi}

\def\@svector{\@sline\@tempcnta\@yarg \ifnum\@tempcnta <0 
\@tempcnta=-\@tempcnta\fi \ifnum\@tempcnta <5 \hskip -\wd\@linechar 
\@upordown\@clnht \hbox{\@linefnt \if@negarg 
\@getlarrow(\@xarg,\@yyarg) \else \@getrarrow(\@xarg,\@yyarg) 
\fi}\else\@badlinearg\fi}

\def\@getlarrow(#1,#2){\ifnum #2 =\z@ \@tempcnta='33\else 
\@tempcnta=#1\relax\multiply\@tempcnta \sixt@@n \advance\@tempcnta -9 
\@tempcntb=#2\relax \multiply\@tempcntb \tw@ \ifnum \@tempcntb >0 
\advance\@tempcnta \@tempcntb\relax \else\advance\@tempcnta 
-\@tempcntb\advance\@tempcnta 64 \fi\fi \char\@tempcnta}

\def\@getrarrow(#1,#2){\@tempcntb=#2\relax \ifnum\@tempcntb < 0 
\@tempcntb=-\@tempcntb\relax\fi \ifcase \@tempcntb\relax 
\@tempcnta='55 \or \ifnum #1<3 \@tempcnta=#1\relax\multiply\@tempcnta 
24 \advance\@tempcnta -6 \else \ifnum #1=3 \@tempcnta=49 
\else\@tempcnta=58 \fi\fi\or \ifnum #1<3 
\@tempcnta=#1\relax\multiply\@tempcnta 24 \advance\@tempcnta -3 \else 
\@tempcnta=51\fi\or \@tempcnta=#1\relax\multiply\@tempcnta \sixt@@n 
\advance\@tempcnta -\tw@ \else \@tempcnta=#1\relax\multiply\@tempcnta 
\sixt@@n \advance\@tempcnta 7 \fi \ifnum #2<0 \advance\@tempcnta 64 
\fi \char\@tempcnta}

\def\@vline{\ifnum \@yarg <0 \@downline \else \@upline\fi}

\def\@upline{\hbox to \z@{\hskip -\@halfwidth \vrule width 
\@wholewidth height \@linelen depth \z@\hss}}

\def\@downline{\hbox to \z@{\hskip -\@halfwidth \vrule width 
\@wholewidth height \z@ depth \@linelen \hss}}

\def\@upvector{\@upline\setbox\@tempboxa 
\hbox{\@linefnt\char'66}\raise \@linelen \hbox to\z@{\lower 
\ht\@tempboxa \box\@tempboxa\hss}}

\def\@downvector{\@downline\lower \@linelen \hbox to 
\z@{\@linefnt\char'77\hss}}

\thinlines

\newcount\@xarg
\newcount\@yarg
\newcount\@yyarg
\newcount\@multicnt
\newdimen\@xdim
\newdimen\@ydim
\newbox\@linechar
\newdimen\@linelen
\newdimen\@clnwd
\newdimen\@clnht
\newdimen\@dashdim
\newbox\@dashbox
\newcount\@dashcnt
\catcode`@=12

\newbox\tbox
\newbox\tboxa

\def\leftzer#1{\setbox\tbox=\hbox to 0pt{#1\hss}%
\ht\tbox=0pt \dp\tbox=0pt \box\tbox}

\def\rightzer#1{\setbox\tbox=\hbox to 0pt{\hss #1}%
\ht\tbox=0pt \dp\tbox=0pt \box\tbox}

\def\centerzer#1{\setbox\tbox=\hbox to 0pt{\hss #1\hss}%
\ht\tbox=0pt \dp\tbox=0pt \box\tbox}

\def\leftput(#1,#2)#3{\setbox\tboxa=\hbox{%
\kern #1\unitlength \raise #2\unitlength\hbox{\leftzer{#3}}}%
\ht\tboxa=0pt \wd\tboxa=0pt \dp\tboxa=0pt\box\tboxa}

\def\rightput(#1,#2)#3{\setbox\tboxa=\hbox{%
\kern #1\unitlength \raise #2\unitlength\hbox{\rightzer{#3}}}%
\ht\tboxa=0pt \wd\tboxa=0pt \dp\tboxa=0pt\box\tboxa}

\def\centerput(#1,#2)#3{\setbox\tboxa=\hbox{%
\kern #1\unitlength \raise #2\unitlength\hbox{\centerzer{#3}}}%
\ht\tboxa=0pt \wd\tboxa=0pt \dp\tboxa=0pt\box\tboxa}

\unitlength=1mm

\expandafter\ifx\csname amssym.def\endcsname\relax \else
\endinput\fi
%
\expandafter\edef\csname amssym.def\endcsname{%
       \catcode`\noexpand\@=\the\catcode`\@\space}
\catcode`\@=11
%

\def\undefine#1{\let#1\undefined}
\def\newsymbol#1#2#3#4#5{\let\next@\relax
 \ifnum#2=\@ne\let\next@\msafam@\else
 \ifnum#2=\tw@\let\next@\msbfam@\fi\fi
 \mathchardef#1="#3\next@#4#5}
\def\mathhexbox@#1#2#3{\relax
 \ifmmode\mathpalette{}{\m@th\mathchar"#1#2#3}%
 \else\leavevmode\hbox{$\m@th\mathchar"#1#2#3$}\fi}
\def\hexnumber@#1{\ifcase#1 0\or 1\or 2\or 3
\or 4\or 5\or 6\or 7\or 8\or
 9\or A\or B\or C\or D\or E\or F\fi}

\font\tenmsa=msam10
\font\sevenmsa=msam7
\font\fivemsa=msam5
\newfam\msafam
\textfont\msafam=\tenmsa
\scriptfont\msafam=\sevenmsa
\scriptscriptfont\msafam=\fivemsa
\edef\msafam@{\hexnumber@\msafam}
\mathchardef\dabar@"0\msafam@39
\def\dashrightarrow{\mathrel{\dabar@\dabar@\mathchar"0
\msafam@4B}}
\def\dashleftarrow{\mathrel{\mathchar"0\msafam@4C
\dabar@\dabar@}}

\def\ulcorner{\delimiter"4\msafam@70\msafam@70 }
\def\urcorner{\delimiter"5\msafam@71\msafam@71 }
\def\llcorner{\delimiter"4\msafam@78\msafam@78 }
\def\lrcorner{\delimiter"5\msafam@79\msafam@79 }
\def\yen{{\mathhexbox@\msafam@55}}
\def\checkmark{{\mathhexbox@\msafam@58}}
\def\circledR{{\mathhexbox@\msafam@72}}
\def\maltese{{\mathhexbox@\msafam@7A}}

\font\tenmsb=msbm10
\font\sevenmsb=msbm7
\font\fivemsb=msbm5
\newfam\msbfam
\textfont\msbfam=\tenmsb
\scriptfont\msbfam=\sevenmsb
\scriptscriptfont\msbfam=\fivemsb
\edef\msbfam@{\hexnumber@\msbfam}
\def\Bbb#1{{\fam\msbfam\relax#1}}
\def\widehat#1{\setbox\z@\hbox{$\m@th#1$}%
 \ifdim\wd\z@>\tw@ em\mathaccent"0\msbfam@5B{#1}%
 \else\mathaccent"0362{#1}\fi}
\def\widetilde#1{\setbox\z@\hbox{$\m@th#1$}%
 \ifdim\wd\z@>\tw@ em\mathaccent"0\msbfam@5D{#1}%
 \else\mathaccent"0365{#1}\fi}
\font\teneufm=eufm10
\font\seveneufm=eufm7
\font\fiveeufm=eufm5
\newfam\eufmfam
\textfont\eufmfam=\teneufm
\scriptfont\eufmfam=\seveneufm
\scriptscriptfont\eufmfam=\fiveeufm
\def\frak#1{{\fam\eufmfam\relax#1}}

\csname amssym.def\endcsname

\expandafter\ifx\csname pre amssym.tex at\endcsname\relax \else 
\endinput\fi
\expandafter\chardef\csname pre amssym.tex at\endcsname=\the
\catcode`\@
\catcode`\@=11
\begingroup\ifx\undefined\newsymbol \else\def\input#1
{\endgroup}\fi
\input amssym.def \relax
\newsymbol\boxdot 1200
\newsymbol\boxplus 1201
\newsymbol\boxtimes 1202
\newsymbol\square 1003
\newsymbol\blacksquare 1004
\newsymbol\centerdot 1205
\newsymbol\lozenge 1006
\newsymbol\blacklozenge 1007
\newsymbol\circlearrowright 1308
\newsymbol\circlearrowleft 1309
\undefine\rightleftharpoons
\newsymbol\rightleftharpoons 130A
\newsymbol\leftrightharpoons 130B
\newsymbol\boxminus 120C
\newsymbol\Vdash 130D
\newsymbol\Vvdash 130E
\newsymbol\vDash 130F
\newsymbol\twoheadrightarrow 1310
\newsymbol\twoheadleftarrow 1311
\newsymbol\leftleftarrows 1312
\newsymbol\rightrightarrows 1313
\newsymbol\upuparrows 1314
\newsymbol\downdownarrows 1315
\newsymbol\upharpoonright 1316
 
\newsymbol\downharpoonright 1317
\newsymbol\upharpoonleft 1318
\newsymbol\downharpoonleft 1319
\newsymbol\rightarrowtail 131A
\newsymbol\leftarrowtail 131B
\newsymbol\leftrightarrows 131C
\newsymbol\rightleftarrows 131D
\newsymbol\Lsh 131E
\newsymbol\Rsh 131F
\newsymbol\rightsquigarrow 1320
\newsymbol\leftrightsquigarrow 1321
\newsymbol\looparrowleft 1322
\newsymbol\looparrowright 1323
\newsymbol\circeq 1324
\newsymbol\succsim 1325
\newsymbol\gtrsim 1326
\newsymbol\gtrapprox 1327
\newsymbol\multimap 1328
\newsymbol\therefore 1329
\newsymbol\because 132A
\newsymbol\doteqdot 132B
 
\newsymbol\triangleq 132C
\newsymbol\precsim 132D
\newsymbol\lesssim 132E
\newsymbol\lessapprox 132F
\newsymbol\eqslantless 1330
\newsymbol\eqslantgtr 1331
\newsymbol\curlyeqprec 1332
\newsymbol\curlyeqsucc 1333
\newsymbol\preccurlyeq 1334
\newsymbol\leqq 1335
\newsymbol\leqslant 1336
\newsymbol\lessgtr 1337
\newsymbol\backprime 1038
\newsymbol\risingdotseq 133A
\newsymbol\fallingdotseq 133B
\newsymbol\succcurlyeq 133C
\newsymbol\geqq 133D
\newsymbol\geqslant 133E
\newsymbol\gtrless 133F
\newsymbol\sqsubset 1340
\newsymbol\sqsupset 1341
\newsymbol\vartriangleright 1342
\newsymbol\vartriangleleft 1343
\newsymbol\trianglerighteq 1344
\newsymbol\trianglelefteq 1345
\newsymbol\bigstar 1046
\newsymbol\between 1347
\newsymbol\blacktriangledown 1048
\newsymbol\blacktriangleright 1349
\newsymbol\blacktriangleleft 134A
\newsymbol\vartriangle 134D
\newsymbol\blacktriangle 104E
\newsymbol\triangledown 104F
\newsymbol\eqcirc 1350
\newsymbol\lesseqgtr 1351
\newsymbol\gtreqless 1352
\newsymbol\lesseqqgtr 1353
\newsymbol\gtreqqless 1354
\newsymbol\Rrightarrow 1356
\newsymbol\Lleftarrow 1357
\newsymbol\veebar 1259
\newsymbol\barwedge 125A
\newsymbol\doublebarwedge 125B
\undefine\angle
\newsymbol\angle 105C
\newsymbol\measuredangle 105D
\newsymbol\sphericalangle 105E
\newsymbol\varpropto 135F
\newsymbol\smallsmile 1360
\newsymbol\smallfrown 1361
\newsymbol\Subset 1362
\newsymbol\Supset 1363
\newsymbol\Cup 1264
 
\newsymbol\Cap 1265
 
\newsymbol\curlywedge 1266
\newsymbol\curlyvee 1267
\newsymbol\leftthreetimes 1268
\newsymbol\rightthreetimes 1269
\newsymbol\subseteqq 136A
\newsymbol\supseteqq 136B
\newsymbol\bumpeq 136C
\newsymbol\Bumpeq 136D
\newsymbol\lll 136E
 
\newsymbol\ggg 136F
 
\newsymbol\circledS 1073
\newsymbol\pitchfork 1374
\newsymbol\dotplus 1275
\newsymbol\backsim 1376
\newsymbol\backsimeq 1377
\newsymbol\complement 107B
\newsymbol\intercal 127C
\newsymbol\circledcirc 127D
\newsymbol\circledast 127E
\newsymbol\circleddash 127F
\newsymbol\lvertneqq 2300
\newsymbol\gvertneqq 2301
\newsymbol\nleq 2302
\newsymbol\ngeq 2303
\newsymbol\nless 2304
\newsymbol\ngtr 2305
\newsymbol\nprec 2306
\newsymbol\nsucc 2307
\newsymbol\lneqq 2308
\newsymbol\gneqq 2309
\newsymbol\nleqslant 230A
\newsymbol\ngeqslant 230B
\newsymbol\lneq 230C
\newsymbol\gneq 230D
\newsymbol\npreceq 230E
\newsymbol\nsucceq 230F
\newsymbol\precnsim 2310
\newsymbol\succnsim 2311
\newsymbol\lnsim 2312
\newsymbol\gnsim 2313
\newsymbol\nleqq 2314
\newsymbol\ngeqq 2315
\newsymbol\precneqq 2316
\newsymbol\succneqq 2317
\newsymbol\precnapprox 2318
\newsymbol\succnapprox 2319
\newsymbol\lnapprox 231A
\newsymbol\gnapprox 231B
\newsymbol\nsim 231C
\newsymbol\ncong 231D
\newsymbol\diagup 201E
\newsymbol\diagdown 201F
\newsymbol\varsubsetneq 2320
\newsymbol\varsupsetneq 2321
\newsymbol\nsubseteqq 2322
\newsymbol\nsupseteqq 2323
\newsymbol\subsetneqq 2324
\newsymbol\supsetneqq 2325
\newsymbol\varsubsetneqq 2326
\newsymbol\varsupsetneqq 2327
\newsymbol\subsetneq 2328
\newsymbol\supsetneq 2329
\newsymbol\nsubseteq 232A
\newsymbol\nsupseteq 232B
\newsymbol\nparallel 232C
\newsymbol\nmid 232D
\newsymbol\nshortmid 232E
\newsymbol\nshortparallel 232F
\newsymbol\nvdash 2330
\newsymbol\nVdash 2331
\newsymbol\nvDash 2332
\newsymbol\nVDash 2333
\newsymbol\ntrianglerighteq 2334
\newsymbol\ntrianglelefteq 2335
\newsymbol\ntriangleleft 2336
\newsymbol\ntriangleright 2337
\newsymbol\nleftarrow 2338
\newsymbol\nrightarrow 2339
\newsymbol\nLeftarrow 233A
\newsymbol\nRightarrow 233B
\newsymbol\nLeftrightarrow 233C
\newsymbol\nleftrightarrow 233D
\newsymbol\divideontimes 223E
\newsymbol\varnothing 203F
\newsymbol\nexists 2040
\newsymbol\Finv 2060
\newsymbol\Game 2061
\newsymbol\mho 2066
\newsymbol\eth 2067
\newsymbol\eqsim 2368
\newsymbol\beth 2069
\newsymbol\gimel 206A
\newsymbol\daleth 206B
\newsymbol\lessdot 236C
\newsymbol\gtrdot 236D
\newsymbol\ltimes 226E
\newsymbol\rtimes 226F
\newsymbol\shortmid 2370
\newsymbol\shortparallel 2371
\newsymbol\smallsetminus 2272
\newsymbol\thicksim 2373
\newsymbol\thickapprox 2374
\newsymbol\approxeq 2375
\newsymbol\succapprox 2376
\newsymbol\precapprox 2377
\newsymbol\curvearrowleft 2378
\newsymbol\curvearrowright 2379
\newsymbol\digamma 207A
\newsymbol\varkappa 207B
\newsymbol\Bbbk 207C
\newsymbol\hslash 207D
\undefine\hbar
\newsymbol\hbar 207E
\newsymbol\backepsilon 237F
\catcode`\@=\csname pre amssym.tex at\endcsname


\def\subsection#1#2{\vskip 3mm #1. {\it #2}\vskip 3mm}
\def\abstract{\let\rm=\sevenrm \let\bf=\sevenbf \let\it=\sevenit
\rm}

\centerline{\twelvebf Le lemme fondamental pour les groupes unitaires}
\vskip 10mm
\centerline{G\'{e}rard LAUMON\footnote{${}^{\ast}$}{\sevenrm CNRS et
Universit\'{e} Paris-Sud, UMR8628, Math\'{e}matique, F-91405 Orsay
Cedex, gerard.laumon@math.u-psud.fr} et NG\^{O} Bao
Ch\^{a}u\footnote{${}^{\ast\ast}$}{\sevenrm Universit\'{e} Paris-Sud,
UMR8628, Math\'{e}matique, F-91405 Orsay Cedex,
bao-chau.ngo@math.u-psud.fr}}
\vskip 20mm

{\abstract
{\bf Abstract.} Let $\scriptstyle G$ be an unramified reductive group
over a non archimedian local field $\scriptstyle F$.  The so-called
{\it Langlands Fundamental Lemma} is a family of conjectural
identities between orbital integrals for $\scriptstyle G(F)$ and
orbital integrals for endoscopic groups of $\scriptstyle G$.  In this
paper we prove the Langlands fundamental lemma in the particular case
where $\scriptstyle F$ is a finite extension of $\scriptstyle{\Bbb
F}_{p}((t))$, $\scriptstyle G$ is a unitary group and $\scriptstyle
p>\,\hbox{rank}(G)$.  Waldspurger has shown that this particular case
implies the Langlands fundamental lemma for unitary groups of rank
$\scriptstyle <p$ when $\scriptstyle F$ is any finite extension of
$\scriptstyle {\Bbb Q}_{p}$.

We follow in part a strategy initiated by Goresky, Kottwitz and
MacPherson.  Our main new tool is a deformation of orbital integrals
which is constructed with the help of the Hitchin fibration for
unitary groups over projective curves.
}

\section{0}{Introduction}
\vskip - 3mm

\subsection{0.1}{Le lemme fondamental de Langlands et Shelstad}

Soient $F$ un corps local non archim\'{e}dien de caract\'{e}ristique
r\'{e}siduelle diff\'{e}rente de $2$, $F'$ {\og}{son}{\fg} extension
quadratique non ramifi\'{e}e et $\tau$ l'\'{e}l\'{e}ment non trivial
du groupe de Galois de $F'$ sur $F$.  On consid\`{e}re le groupe
unitaire quasi-d\'{e}ploy\'{e} $G={\rm U}(n)$ sur $F$ dont le groupe
des points rationnels sur $F$ est
$$
G(F)=\{g\in \mathop{\rm GL}(n,F')\mid \tau^{\ast}({}^{{\rm
t}}g)\Phi_{n}g=\Phi_{n}\}
$$
o\`{u} la matrice $\Phi_{n}$ a pour seules entr\'{e}es non nulles les
$(\Phi_{n})_{i,n+1-i}=1$.

Soient $n=n_{1}+n_{2}$ une partition non triviale et $H={\rm U}(n_{1})\times
{\rm U}(n_{2})$ le groupe endoscopique de $G$ correspondant.

Soient $\delta =(\delta_{1},\delta_{2})$ un \'{e}l\'{e}ment
semi-simple, r\'{e}gulier et elliptique de $H(F)$ et $T=T_{1}\times
T_{2}\subset {\rm U}(n_{1})\times {\rm U}(n_{2})=H$ son
centralisateur; $T_{1}$ et $T_{2}$ sont des tores maximaux 
de ${\rm U}(n_{1})$ et ${\rm U}(n_{2})$ qui sont anisotropes sur $F$.
Fixons un plongement de $T$ comme tore maximal
dans $G$ et notons $\gamma$ l'image de $\delta$ par ce plongement.
Supposons que l'\'{e}l\'{e}ment semi-simple et elliptique $\gamma$ est
r\'{e}gulier dans $G$.

L'ensemble des classes de conjugaison dans la classe de conjugaison
stable de $\gamma$ dans $G(F)$ est en bijection naturelle
$\lambda\mapsto\gamma^{\lambda}$ avec le groupe fini
$\Lambda=\Lambda_{r}=\{\lambda\in ({\Bbb Z}/2{\Bbb Z})^{r}\mid
\lambda_{1}+\ldots +\lambda_{r}=0\}$ o\`{u} $r$ est le rang du
$F'$-tore d\'{e}ploy\'{e} maximal contenu dans le centralisateur de
$\gamma$ dans $\mathop{\rm GL}(n,F')$.  De m\^{e}me l'ensemble des
classes de conjugaison dans la classe de conjugaison stable de
$\delta$ dans $H(F)$ est en bijection naturelle $\lambda\mapsto
\delta^{\lambda}=(\delta_{1}^{\lambda_{1}}, \delta_{2}^{\lambda_{2}})$
avec le sous-groupe $\Lambda^{H}=\Lambda_{r_{1}}\times\Lambda_{r_{2}}$
de $\Lambda$ o\`{u} $r_{1}$ et $r_{2}$ sont les rangs des $F'$-tores
d\'{e}ploy\'{e}s maximaux contenus dans les centralisateurs de
$\delta_{1}$ dans $\mathop{\rm GL}(n_{1},F')$ et $\delta_{2}$ dans
$\mathop{\rm GL}(n_{2},F')$.  On a bien s\^{u}r $r=r_{1}+r_{2}$.  Pour
chaque $\lambda\in\Lambda$, le centralisateur $T^{\lambda}$ de
$\gamma^{\lambda}$ est une forme int\'{e}rieure de $T$ et est donc
isomorphe \`{a} $T$.  De m\^{e}me, pour chaque $\lambda\in
\Lambda_{H}\subset\Lambda$, le centralisateur $S^{\lambda}$ de
$\delta^{\lambda}$ est une forme int\'{e}rieure de $T$ et est donc lui
aussi isomorphe \`{a} $T$.

Notons ${\cal O}_{F'}$ l'anneau des entiers de $F'$.  Soient
$K=K_{n}=G(F)\cap\mathop{\rm GL}(n,{\cal O}_{F'})$ et
$K^{H}=K_{n_{1}}\times K_{n_{2}}$ les sous-groupes maximaux standard
de $G(F)$ et $H(F)$.  On normalise les mesures de Haar ${\rm d}g$ et
${\rm d}h$ de $G(F)$ et $H(F)$ en demandant que $K$ et $K^{H}$ soient
de volume $1$.  On consid\`{e}re les int\'{e}grales orbitales
$$
\mathop{\rm O}\nolimits_{\gamma^{\lambda}}(1_{K})
=\int_{T^{\lambda}(F)\backslash G(F)}1_{K}(g^{-1}\gamma^{\lambda}g)
{{\rm d}g\over {\rm d}t^{\lambda}}
$$
pour $\lambda\in\Lambda$ et
$$
\mathop{\rm O}\nolimits_{\delta^{\lambda}}^{H}(1_{K^{H}})
=\int_{S^{\lambda}(F)\backslash H(F)}
1_{K^{H}}(h^{-1}\delta^{\lambda}h){{\rm d}h\over {\rm d}s^{\lambda}}
$$
pour $\lambda\in\Lambda^{H}$.  On a fix\'{e} une mesure de Haar sur
$T(F)$, par exemple celle qui donne le volume $1$ au sous-groupe
compact maximal, et on a transport\'{e}, par les isomorphismes entre
$T^{\lambda}$ et $T$ et entre $S^{\lambda}$ et $T$ signal\'{e}s plus
haut, cette mesure en la mesure de Haar ${\rm d}t^{\lambda}$ sur
$T^{\lambda}(F)$ pour chaque $\lambda\in\Lambda$ et en la mesure de
Haar ${\rm d}s^{\lambda}$ sur $S^{\lambda}(F)$ pour chaque
$\lambda\in\Lambda^{H}$.

Soit $\kappa :\Lambda\rightarrow\{\pm 1\}$ le caract\`{e}re dont le
noyau est exactement $\Lambda_{H}$.  On forme suivant Langlands et
Shelstad (cf.  [La-Sh]) les combinaisons lin\'{e}aires
d'int\'{e}grales orbitales suivantes: la {\it $\kappa$-int\'{e}grale
orbitale}
$$
\mathop{\rm O}\nolimits_{\gamma}^{\kappa}(1_{K})=
\sum_{\lambda\in\Lambda}\kappa (\lambda )\mathop{\rm
O}\nolimits_{\gamma^{\lambda}}(1_{K})
$$
et {\it l'int\'{e}grale orbitale stable endoscopique}
$$
\mathop{\rm SO}\nolimits_{\delta}^{H}(1_{K^{H}})=
\sum_{\lambda\in\Lambda^{H}}\mathop{\rm
O}\nolimits_{\delta^{\lambda}}^{H}(1_{K^{H}}).
$$

Langlands et Shelstad (cf.  [La-Sh]) ont d\'{e}fini un {\it facteur de
transfert} $\Delta (\gamma ,\delta )$, qui est le produit d'un signe et
de la puissance
$$
|D_{G/H}(\gamma)|^{{1\over 2}}
$$
du nombre d'\'{e}l\'{e}ments du corps r\'{e}siduel de $F$, et ils ont
conjectur\'{e}:

\thm LEMME {\pc FONDAMENTAL}
\enonce
On a l'identit\'{e}
$$
\mathop{\rm O}\nolimits_{\gamma}^{\kappa}(1_{K})=\Delta (\gamma
,\delta )\mathop{\rm SO}\nolimits_{\delta}^{H}(1_{K^{H}}).
$$
\endthm

Waldspurger a d\'{e}montr\'{e} que pour \'{e}tablir cette conjecture
pour $F$ une extension finie de ${\Bbb Q}_{p}$, il suffisait de le
faire lorsque $F$ est une extension finie de ${\Bbb F}_{p}((t))$ (cf.
[Wal~1]), et ce apr\`{e}s avoir remplac\'{e} les groupes $G$ et $H$
par leurs alg\`{e}bres de Lie (cf.  [Hal], [Wal~2], [Wal~3]).

L'objet de cet article est de terminer la d\'{e}monstration du lemme
fondamental pour les groupes unitaires de rang $n<p$ en traitant ce
dernier cas: voir le th\'{e}or\`{e}me 1.5.1 pour l'\'{e}nonc\'{e}
pr\'{e}cis.

\subsection{0.2}{Notre strat\'{e}gie}

Dans la preuve pr\'{e}sent\'{e}e ici, nous utilisons des id\'{e}es de
Goresky, Kottwitz et MacPherson, et du premier auteur, id\'{e}es qui
ont \'{e}t\'{e} introduites dans les travaux ant\'{e}rieurs [G-K-M] et
[Lau].  Comme dans [G-K-M] on exprime le facteur de transfert \`{a}
l'aide d'une fl\`{e}che en cohomologie \'{e}quivariante de sorte que
le lemme fondamental se d\'{e}duit d'un isomorphisme en cohomologie
\'{e}quivariante.  Comme dans [Lau] on utilise un argument de
d\'{e}formation, qui fait {\og}{glisser}{\fg} d'une situation
d'intersection tr\`{e}s compliqu\'{e}e vers une situation
d'intersection transversale.

Les r\'{e}sultats de [G-K-M] dans le cas non ramifi\'{e} pour un
groupe r\'{e}ductif quelconque, et de [Lau] dans le cas
\'{e}ventuellement ramifi\'{e}, mais pour le groupe unitaire
uniquement, supposent d\'{e}montr\'{e}e une conjecture de puret\'{e}
des fibres de Springer.  Une telle conjecture a \'{e}t\'{e}
formul\'{e}e par Goresky, Kottwitz et MacPherson.

Nous ne savons pas d\'{e}montrer cette conjecture, mais nous
contournons le probl\`{e}me en d\'{e}montrant en fait un autre
\'{e}nonc\'{e} de puret\'{e}, \`{a} savoir la puret\'{e} d'un faisceau
pervers li\'{e} \`{a} une famille {\og}{universelle}{\fg} de
$\kappa$-int\'{e}grales orbitales globales.  Pour cela nous nous
fondons sur une interpr\'{e}tation g\'{e}om\'{e}trique de la
th\'{e}orie de l'endoscopie de Langlands et Kottwitz (cf.  [Lan] et
[Kot~1]) \`{a} l'aide de la fibration de Hitchin ([Hit]).  Cette
interpr\'{e}tation, d\'{e}couverte par le second auteur et
pr\'{e}sent\'{e}e ici uniquement dans le cas des groupes unitaires,
vaut en fait en toute g\'{e}n\'{e}ralit\'{e} (cf.  [Ngo]).  Enfin, un
argument dans l'esprit de ([Lau]) permet de conclure.

\subsection{0.3}{Plan de l'article}

Passons bri\`{e}vement en revue l'organisation de cet article. Dans
le chapitre 1, nous explicitons l'\'{e}nonc\'{e} du lemme fondamental
pour les groupes unitaires, en termes de comptage des r\'{e}seaux
qui sont auto-duaux par rapport \`{a} une forme hermitienne et qui
sont stables par une transformation unitaire.

Dans le chapitre 2, nous explicitons la construction de la fibration
de Hitchin dans le cas du groupe unitaire.  Nous faisons le lien entre
la fibration Hitchin d'un groupe unitaire et la fibration de Hitchin
d'un de ces groupes endoscopiques.

Dans le chapitre 3, le c{\oe}ur de ce travail, nous d\'{e}montrons une
identit\'{e} globale, que l'on devrait pouvoir identifier \`{a} une
identit\'{e} globale qui appara\^{\i}t dans la stabilisation de la
formule des traces.  L'\'{e}nonc\'{e} principal de ce chapitre est le
th\'{e}or\`{e}me 3.9.3.  On le d\'{e}montre \`{a} l'aide
d'un isomorphisme en cohomologie \'{e}quivariante.  Comme nous l'avons
mentionn\'{e} plus haut, l'isomorphisme en cohomologie
\'{e}quivariante que nous construisons est analogue \`{a} celui
construit ant\'{e}rieurement dans [G-K-M].  Comme nous l'avons
d\'{e}j\`{a} dit notre construction s'appuie sur un \'{e}nonc\'{e} de
puret\'{e}, d\'{e}montr\'{e} dans le paragraphe 3.2, et d'un argument
de d\'{e}formation.

Dans le chapitre 4, nous expliquons comment passer d'une situation
locale donn\'{e}e, \`{a} une situation globale du type de celle
consid\'{e}r\'{e}e dans le chapitre 3.  Ici, l'outil de base est un
th\'{e}or\`{e}me de Bertini rationnel 4.4.1, d\'{e}montr\'{e} par
Gabber ([Gab]) et Poonen ([Poo]).  Le comptage de la section (4.6) est
analogue \`{a} celui du th\'{e}or\`{e}me (15.8) de [G-K-M].

Enfin, dans un appendice, nous d\'{e}montrons une variante A.1.2 du
th\'{e}or\`{e}me de localisation d'Atiyah-Borel-Segal.  Puis nous
pr\'{e}sentons le calcul de la cohomologie \'{e}quivariante d'un
fibr\'{e} en droites projective et d'un fibr\'{e} en droites
projectives pinc\'{e}es.  Nous d\'{e}montrons dans le dernier
appendice une formule de points fixes.

\subsection{0.4}{Pr\'{e}cautions d'emploi de nos r\'{e}sultats}

Dans ce travail nous avons admis certains r\'{e}sultats sur la
cohomologie $\ell$-adique des champs alg\'{e}briques.

\subsection{0.5}{Remerciements}

Nous remercions A.~Abbes, J.-B.~Bost, L.~Breen, M.~Brion, J.-F.~Dat,
O.~Gabber, D.~Gaitsgory, A.~Genestier, L.~Illusie, S.~Kleiman,
L.~Lafforgue, F.~Loeser, M.~Raynaud et J.-L.~Waldspurger pour l'aide
qu'ils nous ont apport\'{e}e durant la pr\'{e}paration de ce travail.

\section{1}{Int\'{e}grales orbitales et comptage de r\'{e}seaux}
\vskip - 3mm

\subsection{1.1}{Les donn\'{e}es}

Pour tout corps local non archim\'{e}dien $K$ on note ${\cal O}_{K}$
son anneau des entiers, $\varpi_{K}$ une uniformisante de $K$ et
$v_{K}:K^{\times}\rightarrow {\Bbb Z}$ la valuation discr\`{e}te
normalis\'{e}e par $v_{K}(\varpi_{K})=1$.

Pour toute extension finie $L$ de $K$, on note $\mathop{\rm
Tr}\nolimits_{L/K}:L\rightarrow K$ et $\mathop{\rm
Nr}\nolimits_{L/K}:L^{\times}\rightarrow K^{\times}$ les trace et
norme correspondantes.
\vskip 2mm

Soient $F$ un corps local non archim\'{e}dien d'\'{e}gales
caract\'{e}ristiques diff\'{e}rentes de $2$, $k={\Bbb F}_{q}$ son corps
r\'{e}siduel et $F'$ son extension quadratique non ramifi\'{e}e de $F$
(de corps r\'{e}siduel ${\Bbb F}_{q^{2}}$).

On se donne une famille finie $(E_{i})_{i\in I}$ d'extensions finies
s\'{e}parables de $F$ qui sont toutes disjointes de $F'$ et, pour
chaque $i\in I$, un \'{e}l\'{e}ment $\gamma_{i}$ de l'extension
compos\'{e}e $E_{i}'=E_{i}F'$.  On note $n_{i}$ le degr\'{e} de
$E_{i}$ sur $F$ et $\tau$ l'\'{e}l\'{e}ment non trivial des groupes de
Galois
$$
\mathop{\rm Gal}(F'/F)\cong \mathop{\rm Gal}(E_{i}'/E_{i}).
$$

On suppose que, pour chaque $i\in I$, $\gamma_{i}$ engendre $E_{i}'$ sur
$F'$, $\gamma_{i}\in {\cal O}_{E_{i}'}$ et
$$
\gamma_{i}^{\tau}+\gamma_{i}=0.
$$
On suppose de plus que pour tous $i\not=j$ dans $I$ les polyn\^{o}mes
minimaux $P_{i}(T)$ et $P_{j}(T)$ sur $F'$ de $\gamma_{i}$ et
$\gamma_{j}$ sont premiers entres eux.  On suppose enfin que la
caract\'{e}ristique de $k$ est $>n=\sum_{i\in I}n_{i}$.

Le tore $T$ de l'introduction est alors le tore anisotrope sur $F$
dont le groupe des $F$-points est
$$
T(F)=\prod_{i\in I}\{x\in E_{i}'\mid x^{\tau}x=1\}
$$
et $\gamma$ est vu comme un point sur $F$ de l'alg\`{e}bre de Lie
de ce tore.

\subsection{1.2}{Invariants num\'{e}riques}

Pour chaque $i\in I$, le polyn\^{o}me minimal $P_{i}(T)$ de
$\gamma_{i}\in {\cal O}_{E_{i}'}$ sur $F'$ est un polyn\^{o}me unitaire de
degr\'{e} $n_{i}$ \`{a} coefficients dans ${\cal O}_{F'}$.  Comme
$\gamma_{i}^{\tau}=-\gamma_{i}$, on a de plus
$P_{i}^{\tau}(T)=(-1)^{n_{i}}P_{i}(-T)$.

On note $\delta_{i}$ la dimension sur $k$ de ${\cal O}_{E_{i}'}/{\cal
O}_{F'}[\gamma_{i}]$, c'est-\`{a}-dire la co-longueur de ${\cal
O}_{F'}[\gamma_{i}]$ comme sous-${\cal O}_{F'}$-r\'{e}seau de ${\cal
O}_{E_{i}'}$.  D'apr\`{e}s Gorenstein, le conducteur ${\frak
a}_{i}\subset {\cal O}_{F'}[\gamma_{i}]\subset {\cal O}_{E_{i}'}$ de
${\cal O}_{E_{i}'}$ dans ${\cal O}_{F'}[\gamma_{i}]$ est de
co-longueur $\delta_{i}$ comme sous-${\cal O}_{F'}$-r\'{e}seau de
${\cal O}_{F'}[\gamma_{i}]$ et est donc \'{e}gal \`{a}
$$
{\frak a}_{i}=\varpi_{E_{i}}^{{2\delta_{i}e_{i}\over n_{i}}}{\cal O}_{E_{i}'}
$$
o\`{u} $e_{i}$ est l'indice de ramification de $E_{i}$ sur $F$.

Puisque l'extension $E_{i}/F$ est de degr\'{e} $n_{i}<p$, a fortiori
premier \`{a} $p$, la diff\'{e}rente ${\frak D}_{E_{i}/F}$ est
\'{e}gale \`{a} l'id\'{e}al $\varpi_{E_{i}}^{e_{i}-1}{\cal O}_{E_{i}}$
de ${\cal O}_{E_{i}}$ d'apr\`{e}s la proposition 13, \S6, ch.  III, de
[Ser].  De m\^{e}me, la diff\'{e}rente ${\frak D}_{E_{i}'/F'}$ est
\'{e}gale \`{a} l'id\'{e}al $\varpi_{E_{i}}^{e_{i}-1}{\cal
O}_{E_{i}'}$ de ${\cal O}_{E_{i}'}$.  En utilisant loc.  cit.  Cor.
1, on a donc
$$
v_{E_{i}'}\left({dP_{i}\over dT}(\gamma_{i})\right)=
{2\delta_{i}e_{i}\over n_{i}}+e_{i}-1.
$$
\vskip 2mm

Pour tous $i\not=j$ dans $I$, le r\'{e}sultant $\mathop{\rm
Res}(P_{i},P_{j})\in F'$ est non nul puisque les polyn\^{o}mes
$P_{i}(T)$ et $P_{j}(T)$ sont premiers entre eux.  C'est donc un
\'{e}l\'{e}ment non nul de ${\cal O}_{F'}$.  On a de plus $\mathop{\rm
Res} (P_{i},P_{j})^{\tau}=(-1)^{n_{i}n_{j}}\mathop{\rm
Res}(P_{i},P_{j})$.  La valuation $r_{ij}\geq 0$ de $\mathop{\rm
Res}(P_{i},P_{j})$ est \'{e}gale \`{a}
$$
r_{ij}={n_{i}v_{E_{i}'}(P_{j}(\gamma_{i}))\over e_{i}}=
{n_{j}v_{E_{j}'}(P_{i}(\gamma_{j}))\over e_{j}}.
$$
Elle est aussi \'{e}gale \`{a} l'indice du ${\cal O}_{F'}$-r\'{e}seau
${\cal O}_{F'}[\gamma_{i}\oplus\gamma_{j}] \subset E'$ comme
sous-r\'{e}seau de ${\cal O}_{F'}[\gamma_{i}]\oplus {\cal
O}_{F'}[\gamma_{j}]\subset E'$.  De plus, on a
$$
P_{j}(\gamma_{i}){\cal O}_{F'}[\gamma_{i}]\oplus
P_{i}(\gamma_{j}){\cal O}_{F'}[\gamma_{j}]\subset {\cal O}_{F'}
[\gamma_{i}\oplus\gamma_{j}]\subset {\cal O}_{F'}[\gamma_{i}]
\oplus {\cal O}_{F'}[\gamma_{j}]
$$
puisque $P_{i}(\gamma_{i}\oplus\gamma_{j} )=0\oplus P_{i}(\gamma_{j})$
et $P_{j}(\gamma_{i}\oplus\gamma_{j} )= P_{j}(\gamma_{i})\oplus 0$, et
l'indice de $P_{j}(\gamma_{i}){\cal O}_{F'}[\gamma_{i}]\oplus
P_{i}(\gamma_{j}){\cal O}_{F'}[\gamma_{j}]$ dans ${\cal O}_{F'}
[\gamma_{i}\oplus\gamma_{j}]$ est aussi \'{e}gal \`{a} $r_{ij}$.

Pour toute partie $J$ de $I$ on note $E_{J}=\bigoplus_{i\in J}E_{i}$
et $E_{J}'=\bigoplus_{i\in J}E_{i}'$.  Ce sont des espaces
vectoriels de dimension $n_{J}=\sum_{i\in J}n_{i}$ sur $F$ et $F'$
respectivement.  On note aussi ${\cal O}_{E_{J}}=\bigoplus_{i\in
J}{\cal O}_{E_{i}}$ et ${\cal O}_{E_{J}'}=\bigoplus_{i\in J}{\cal
O}_{E_{i}'}$.  Soit $\gamma_{J}\in {\cal O}_{E_{J}'}$
l'\'{e}l\'{e}ment $\gamma_{J}=\oplus_{i\in J}\gamma_{i}$.  La
sous-${\cal O}_{F'}$-alg\`{e}bre ${\cal O}_{F'}[\gamma_{J}]$ de ${\cal
O}_{E_{J}'}$ est de co-longueur
$$
\delta_{J}=\sum_{i\in J}\delta_{i}+{1\over 2}\sum_{i\not=j\in J}r_{ij}
$$
dans ${\cal O}_{E_{J}'}$.

Soit ${\frak a}_{J}\subset {\cal O}_{F'}[\gamma_{J}]\subset {\cal
O}_{E_{J}'}$ le conducteur de ${\cal O}_{E_{J}'}$ dans ${\cal
O}_{F'}[\gamma_{J}]$.  D'apr\`{e}s Gorenstein, ${\frak a}_{J}$ est de
co-longueur $\delta_{J}$ dans ${\cal O}_{F'}[\gamma_{J}]$ et est
\'{e}gal \`{a}
$$
{\frak a}_{J}=\varpi_{E_{J}}^{\underline{a}_{J}}{\cal O}_{E_{J}'}
$$
o\`{u} $\underline{a}_{J}=(a_{i})_{i\in J}$ est la famille des entiers
$$
a_{i}={\left(2\delta_{i}+\sum_{j\in J-\{i\}}r_{ij}\right)e_{i}\over n_{i}}
$$
et o\`{u} on a pos\'{e} $\varpi_{E_{J}'}^{\underline{a}_{J}}=\oplus_{i\in
J}\varpi_{E_{i}'}^{a_{i}}$.

\subsection{1.3}{Formes hermitiennes}

On rappelle que le groupe $F^{\times}/\mathop{\rm
Nr}\nolimits_{F'/F}(F'^{\times})$ est le groupe \`{a} deux
\'{e}l\'{e}ments engendr\'{e} par la classe de n'importe quelle
uniformisante $\varpi_{F}$ de $F$.  On l'identifie \`{a} ${\Bbb
Z}/2{\Bbb Z}$ dans la suite.

Pour chaque $i\in I$ et chaque $c_{i}\in E_{i}^{\times}$, on munit le
$F'$-espace vectoriel $E_{i}'$ de la forme hermitienne
$$
\Phi_{i ,c_{i}}:E_{i}'\times E_{i}'\rightarrow F',~(x,y)\mapsto
\mathop{\rm Tr}\nolimits_{E_{i}'/F'} (c_{i}x^{\tau}y).
$$
Si ${\frak d}_{E_{i}/F}$ est le discriminant de $E_{i}/F$,
c'est-\`{a}-dire l'id\'{e}al
$$
{\frak d}_{E_{i}/F}=\mathop{\rm Nr}\nolimits_{E_{i}/F}({\frak
D}_{E_{i}/F})
$$
de ${\cal O}_{F}$, le discriminant de $\Phi_{c_{i}}$ est la classe
$\lambda_{i}(c_{i})\in {\Bbb Z}/2{\Bbb Z}$ de l'\'{e}l\'{e}ment
$$
\mathop{\rm Nr}\nolimits_{E_{i}/F}(c_{i})x
\subset F^{\times}
$$
dans $F^{\times}/\mathop{\rm Nr}\nolimits_{F'/F}(F'^{\times})$ pour
n'importe quel g\'{e}n\'{e}rateur $x$ de l'id\'{e}al ${\frak
d}_{E_{i}/F}$.  Comme $n_{i}$ est premier \`{a} la caract\'{e}ristique
r\'{e}siduelle par hypoth\`{e}se, on a
$$
{\frak d}_{E_{i}/F}=\varpi_{F}^{n_{i}{e_{i}-1\over e_{i}}}{\cal O}_{F}
$$
et
$$
\lambda_{i}(c_{i})\equiv v_{F}(\mathop{\rm
Nr}\nolimits_{E_{i}/F}(c_{i}))+n_{i}-{n_{i}\over
e_{i}}={n_{i}v_{E_{i}}(c_{i})\over e_{i}}+n_{i}-{n_{i}\over
e_{i}}\hbox{ {\rm (mod} }2).
$$

Plus g\'{e}n\'{e}ralement, pour chaque partie $J$ de $I$ et chaque
$c_{J}=(c_{i})_{i\in J}\in E_{J}^{\times}$, on munit le $F'$-espace
vectoriel $E_{J}'$ de la forme hermitienne
$$
\Phi_{J,c_{J}}=\oplus_{i\in J}\Phi_{i,c_{i}}:E_{J}'\times E_{J}'\rightarrow F'.
$$
Le discriminant de $\Phi_{J,c_{J}}$ est la somme
$$
\sum_{i\in J}\lambda_{i}(c_{i})\in {\Bbb Z}/2{\Bbb Z}
$$
des discriminants des $\Phi_{i,c_{i}}$.

\subsection{1.4}{R\'{e}seaux auto-duaux}

Pour chaque partie $J$ de $I$ et chaque $c_{J}\in E_{J}^{\times}$, on
consid\`{e}re l'ensemble
$$
\{M_{J}\subset E_{J}'\mid M_{J}^{\perp_{c_{J}}}=M_{J}\hbox{ et
}\gamma_{J}M_{J}\subset M_{J}\}
$$
des ${\cal O}_{F'}$-r\'{e}seaux $M_{J}$ de $E_{J}'$ qui sont \`{a} la
fois auto-duaux pour $\Phi_{J,c_{J}}$ et stables par $\gamma_{J}$.
Ici on a not\'{e}
$$
M_{J}^{\perp_{c_{J}}}=\{x\in E_{J}'\mid \Phi_{J,c_{J}}(x,M_{J})\subset
{\cal O}_{F'}\}
$$
l'orthogonal de $M_{J}$ pour la forme hermitienne $\Phi_{J,c_{J}}$.

\thm LEMME 1.4.1
\enonce
Pour chaque partie $J$ de $I$ et chaque $c_{J}\in E_{J}^{\times}$,
l'ensemble de r\'{e}seaux ci-dessus est un ensemble fini.
\endthm

\rem D\'{e}monstration
\endrem
Pour tout r\'{e}seau dans cet ensemble, on a
$$
{\frak a}_{J}M_{J}\subset M_{J}\subset {\cal O}_{E_{J}'}M_{J}
$$
et
$$
({\cal O}_{E_{J}'}M_{J})^{\perp_{c_{J}}}\subset M_{J}\subset ({\frak
a}_{J}M_{J})^{\perp_{c_{J}}}.
$$
Or on a ${\cal
O}_{E_{J}'}M_{J}=\varpi_{E_{J}'}^{-\underline{m}_{J}}{\cal
O}_{E_{J}'}$ pour une famille d'entiers
$\underline{m}_{J}=(m_{i})_{i\in J}$ index\'{e}e par $J$, et donc
${\frak a}_{J}M_{J}=
\varpi_{E_{J}'}^{\underline{a}_{J}-\underline{m}_{J}}{\cal
O}_{E_{J}'}$, $({\cal O}_{E_{J}'}M_{J})^{\perp_{c_{J}}}=
\varpi_{E_{J}}^{\underline{m}_{J}}({\cal O}_{E_{J}'})^{\perp_{c_{J}}}
=\varpi_{E_{J}}^{-\underline{b}_{J}+\underline{m}_{J}}
{\cal O}_{E_{J}'}$ et $({\frak a}_{J}M_{J})^{\perp_{c_{J}}}=
\varpi_{E_{J}}^{-\underline{a}_{J}+\underline{m}_{J}}({\cal
O}_{E_{J}'})^{\perp_{c_{J}}}=\varpi_{E_{J}}^{-\underline{a}_{J}
-\underline{b}_{J}+\underline{m}_{J}}{\cal O}_{E_{J}'}$ o\`{u}
$\underline{b}_{J}=(b_{i})_{i\in J}$ est la famille d'entiers
d\'{e}finie par $c_{i}^{-1}{\frak D}_{E_{i}/F}^{-1}=
\varpi_{E_{i}}^{-b_{i}}{\cal O}_{E_{i}}$.  On en
d\'{e}duit que
$$
b_{i}\leq 2m_{i}\leq 2a_{i}+b_{i},~\forall i\in J,
$$
et que
$$
\varpi_{E_{J}'}^{\underline{a}_{J}-{\underline{b}_{J}\over 2}}{\cal
O}_{E_{J}'}\subset M_{J}\subset
\varpi_{E_{J}'}^{-\underline{a}_{J}-{\underline{b}_{J}\over 2}}{\cal
O}_{E_{J}'},
$$
d'o\`{u} le lemme.
\hfill\hfill$\square$
\vskip 3mm

L'ensemble des r\'{e}seaux $M_{J}$ de $E_{J}'$ qui sont \`{a} la fois
auto-duaux pour $\Phi_{J,c_{J}}$ et stables par $\gamma_{J}$, admet
encore la description suivante qui est le point de d\'{e}part de ce
travail.  Consid\'{e}rons la ${\cal O}_{F'}$-alg\`{e}bre $A_{J}={\cal
O}_{F'}[\gamma_{J}]={\cal O}_{F'}[T]/(P_{J}(T))$ munie de l'involution
qui induit $\tau$ sur ${\cal O}_{F'}$ et qui envoie $T$ sur $-T$.
Alors, $E_{J}'$ est l'anneau total des fractions de $A_{J}$ et les
${\cal O}_{F'}$-r\'{e}seaux $M_{J}\subset E_{J}'$ tels que
$\gamma_{J}M_{J}\subset M_{J}$ ne sont rien d'autre que les id\'{e}aux
fractionnaires de $A_{J}$.  Un tel id\'{e}al fractionnaire admet un
inverse
$$
M_{J}^{-1}=\{m\in E_{J}'\mid xM_{J}\subset A_{J}\}.
$$

\thm LEMME 1.4.2
\enonce
Pour tout $c_{J}\in E_{J}^{\times}$, l'orthogonal du ${\cal
O}_{F'}$-r\'{e}seau $A_{J}\subset E_{J}'$ relativement \`{a} la forme
hermitienne $\Phi_{J,c_{J}}$ est
$$
(A_{J})^{\perp_{c_{J}}}=\left(\oplus_{i\in J}{1\over c_{i}{dP_{i}\over
dT}(\gamma_{i})P_{J-\{i\}}(\gamma_{i})} \right)A_{J}\subset E_{J}'.
$$

Plus g\'{e}n\'{e}ralement, pour tout $c_{J}\in E_{J}^{\times}$ et tout
id\'{e}al fractionnaire $M_{J}$ de $A_{J}$ on a la relation
$$
M_{J}^{\perp_{c_{J}}}=\left(\oplus_{i\in J}{1\over c_{i}{dP_{i}\over
dT}(\gamma_{i})P_{J-\{i\}}(\gamma_{i})} \right)M_{J}^{-1}.
$$
\endthm

\rem D\'{e}monstration
\endrem
Voir la d\'{e}monstration de la proposition 11, \S6, ch.  III, de
[Ser].
\hfill\hfill$\square$
\vskip 3mm

En particulier, si on note
$$
c_{J,i}^{0}={\varepsilon^{n_{J}-1}\over {dP_{i}\over
dT}(\gamma_{i})P_{J-\{i\}}(\gamma_{i})}\in E_{i}^{\times},~\forall\,
i\in J,
$$
o\`{u} $\varepsilon$ est un \'{e}l\'{e}ment de ${\Bbb
F}_{q^{2}}\subset F'$ tel que $\varepsilon^{\tau}=-\varepsilon$, on a
$c_{J}^{0}=(c_{J,i}^{0})_{i\in J}\in E_{J}^{\times}$ et le ${\cal
O}_{F'}$-r\'{e}seau $A_{J}\subset E_{J}'$ est auto-dual pour la forme
hermitienne $\Phi_{J,c_{J}^{0}}$.

\thm LEMME 1.4.3
\enonce
On a
$$
\lambda_{J,i}^{0}:=\lambda_{i}(c_{J,i}^{0})\equiv \sum_{j\in J-\{i\}}
r_{ji}\hbox{ {\rm (mod }}2)
$$
pour toute partie $J$ de $I$ et tout $i\in J$.
\endthm

\rem D\'{e}monstration
\endrem
On a vu que
$$
\lambda_{i}(c_{J,i}^{0})\equiv {n_{i}v_{E_{i}}(c_{J,i}^{0})\over
e_{i}}+n_{i}-{n_{i}\over e_{i}}\hbox{ {\rm (mod} }2).
$$
Or
$$
v_{E_{i}'}\left({dP_{i}\over dT}(\gamma_{i})P_{J-\{i\}}(\gamma_{i})\right)
={2\delta_{i}e_{i}\over n_{i}}+e_{i}-1+\sum_{j\in J-\{i\}}{e_{i}r_{ji}\over
n_{i}}.
$$
d'o\`{u} le lemme.
\hfill\hfill$\square$
\vskip 3mm

\subsection{1.5}{\'{E}nonc\'{e} du lemme fondamental}

Pour chaque $J\subset I$ et pour chaque $c_{J}$ tel que $\sum_{i\in
J}\lambda_{i}(c_{i})=0$, le cardinal de l'ensemble fini de r\'{e}seaux
$$
\mathop{\rm O}\nolimits_{\gamma_{J}}^{c_{J}}=|\{M_{J}\subset
E_{J}'\mid M_{J}^{\perp_{c_{J}}}=M_{J}\hbox{ et
}\gamma_{J}M_{J}\subset M_{J}\}|
$$
ne d\'{e}pend que de $\lambda_{J}(c_{J})=(\lambda_{i}(c_{i}))_{i\in
J}$.  Soit $\Lambda_{J}=\{\lambda_{J}\in ({\Bbb Z}/2{\Bbb Z})^{J}\mid
\sum_{i\in J}\lambda_{i}=0\}$.  Pour chaque
$\lambda_{J}\in\Lambda_{J}$ on peut donc noter
$$
\mathop{\rm O}\nolimits_{\gamma_{J}}^{\lambda_{J}}
=\mathop{\rm O}\nolimits_{\gamma_{J}}^{c_{J}}
$$
pour n'importe quel $c_{J}$ tel que $\lambda_{J}(c_{J})=\lambda_{J}$.

En fait $\mathop{\rm O}\nolimits_{\gamma_{J}}^{\lambda_{J}}$ est une
int\'{e}grale orbitale.  Plus pr\'{e}cis\'{e}ment, choisissons $c_{J}\in
E_{J}^{\times}$ tel que $\lambda_{J}(c_{J})=\lambda_{J}$. Comme le
discriminant de $\Phi_{J,c_{J}}$ est \'{e}gal \`{a} $1$, le
$F'$-espace vectoriel hermitien $(E_{J}',\Phi_{J,c_{J}})$ est
isomorphe au $F'$-espace hermitien standard
$(F'^{n_{J}},\Phi_{n_{J}})$. Choisissons un tel isomorphisme et
consid\'{e}rons le plongement de l'alg\`{e}bre de Lie de
$$
T_{J}(F)=\prod_{i\in J}\{x\in E_{i}'\mid x^{\tau}x=1\}
$$
dans l'alg\`{e}bre de Lie ${\frak u}(n_{J})$ de ${\rm U}(n_{J})$ qu'il
induit.  La classe de $G(F)$-conjugaison de l'image
$\gamma_{J}^{\lambda_{J}}$ de $\gamma_{J}$ par ce dernier plongement
ne d\'{e}pend pas des choix que l'on vient de faire.  Alors,
$\mathop{\rm O}\nolimits_{\gamma_{J}}^{\lambda_{J}}$ est
pr\'{e}cis\'{e}ment l'int\'{e}grale orbitale de la fonction
caract\'{e}ristique du compact maximal standard ${\frak K}_{n_{J}}$ de
${\frak u}(n_{J})(F)$

Soit $I=I_{1}\amalg I_{2}$ une partition de $I$. La donn\'{e}e de
cette partition \'{e}quivaut \`{a} la donn\'{e}e d'un caract\`{e}re
$$
\kappa_{I_{1},I_{2}}:\Lambda_{I}\rightarrow \{\pm 1\}
$$
\`{a} savoir le caract\`{e}re $\kappa$ d\'{e}fini par
$$
\kappa (\lambda_{I})=(-1)^{\sum_{i\in I_{1}}\lambda_{i}}=
(-1)^{\sum_{i\in I_{2}}\lambda_{i}}.
$$

On a alors la {\og}{$\kappa$-int\'{e}grale orbitale}{\fg}
$$
\mathop{\rm O}\nolimits_{\gamma}^{\kappa}=
\sum_{\lambda_{I}\in\Lambda_{I}}\kappa
(\lambda_{I}-(\lambda_{I_{1}}^{0},\lambda_{I_{2}}^{0}))
\mathop{\rm O}\nolimits_{\gamma}^{\lambda_{I}}
$$
et l'{\og}{int\'{e}grale orbitale endoscopique stable}{\fg}
$$
\mathop{\rm SO}\nolimits_{\gamma}^{H}=
\sum_{{\scriptstyle\lambda_{I_{1}}\in\Lambda_{I_{1}}
\atop\scriptstyle\lambda_{I_{2}}\in\Lambda_{I_{2}}}} \mathop{\rm
O}\nolimits_{\gamma_{I_{1}}}^{\lambda_{I_{1}}}\times \mathop{\rm
O}\nolimits_{\gamma_{I_{2}}}^{\lambda_{I_{2}}}.
$$

L'objet de cet article est de d\'{e}montrer le th\'{e}or\`{e}me
suivant conjectur\'{e} par Langlands et Shelstad (cf.  [La-Sh]) et
appel\'{e} par eux le {\og}{lemme fondamental pour les groupes
unitaires}{\fg} (ou plut\^{o}t sa variante alg\`{e}bre de Lie).

\thm TH\'{E}OR\`{E}ME 1.5.1
\enonce
Sous les hypoth\`{e}ses pr\'{e}c\'{e}dentes, on a la relation
$$
\mathop{\rm O}\nolimits_{\gamma}^{\kappa}=(-1)^{r}q^{r}
\mathop{\rm SO}\nolimits_{\gamma}^{H}
$$
o\`{u} on a pos\'{e}
$$
r=r_{I_{1},I_{2}}=\sum_{{\scriptstyle i_{1}\in I_{1}\atop\scriptstyle
i_{2}\in I_{2}}}r_{i_{1},i_{2}}.
$$
\endthm

En faisant passer le terme $(-1)^{r}q^{r}$ de l'autre c\^{o}t\'{e} de
l'\'{e}galit\'{e}, on fait appara\^{\i}tre dans l'expression de la
$\kappa$-int\'{e}grale orbitale comme combinaison lin\'{e}aire 
d'int\'{e}grales orbitales, des coefficients
$$
\kappa (\lambda_{I}-(\lambda_{I_{1}}^{0},\lambda_{I_{2}}^{0}))
(-1)^{r}q^{-r}
$$
qui sont \'{e}gaux aux facteurs de transfert de Langlands-Shelstad
(cf. [La-Sh]), d'apr\`{e}s des calculs Waldspurger valables pour
tous les groupes classiques (cf. la proposition X.8 de [Wal 4]). Ces
facteurs sont aussi les m\^{e}mes que ceux d\'{e}finis par Kottwitz
\`{a} l'aide des sections de Kostant (cf. [Kot~2]). 

\section{2}{Fibration de Hitchin}
\vskip - 3mm

\subsection{2.1}{Sch\'{e}mas en groupes unitaires et fibr\'{e}s
hermitiens}

On fixe une courbe projective, lisse et g\'{e}om\'{e}triquement connexe
$X$ sur $k={\Bbb F}_{q}$, de genre g\'{e}om\'{e}trique $g\geq 1$, et
un rev\^{e}tement \'{e}tale, galoisien, de degr\'{e} $2$, $\pi
:X'\rightarrow X$, dont l'espace total $X'$ est donc une courbe projective
et lisse que l'on suppose aussi g\'{e}om\'{e}triquement connexe.  On
note $\tau$ l'\'{e}l\'{e}ment non trivial du groupe de Galois de $X'$
sur $X$.

On fixe un entier $n\geq 1$ et l'on suppose que la caract\'{e}ristique
de $k$ est $>n$, et impaire si $n=1$.

On munit le fibr\'{e} vectoriel trivial ${\cal O}_{X'}^{\oplus n}$ de rang
$n\geq 1$ sur $X'$ de la forme hermitienne
$$
\Phi_{n}: {\cal O}_{X'}^{n}\times {\cal O}_{X'}^{n}\rightarrow {\cal O}_{X'}
$$
dont la matrice $\Phi_{n}$ a pour seules entr\'{e}es non nulles les
$(\Phi_{n})_{i,n+1-i}=1$.

On d\'{e}finit alors le sch\'{e}ma en groupes unitaires $G$ sur $X$ par
$$
G(S)=\{g\in \mathop{\rm GL}\nolimits_{n}(H^{0}(X'_{S},{\cal
O}_{X'_{S}}))\mid \tau^{\ast}({}^{{\rm t}}g)\Phi_{n}g=\Phi_{n}\}
$$
o\`{u} pour tout $X$-sch\'{e}ma $S$, on a not\'{e} par un indice $S$ le
changement de base par le morphisme structural $S\rightarrow X$.

Le $X'$-sch\'{e}ma en groupes $G_{X'}=X'\times_{\pi ,X}G$ n'est autre que
$\mathop{\rm GL}\nolimits_{n,{\cal O}_{X'}}$ puisque $X'\times_{X}X'$ est
la somme disjointe de deux copies de $X'$, $\tau$ \'{e}changeant ces
deux copies.

Il s'en suit que la restriction de $G$ au compl\'{e}t\'{e} formel $\mathop{\rm
Spf}({\cal O}_{x})$ de $X$ en un point ferm\'{e} $x$ est isomorphe
\`{a} $\mathop{\rm GL}\nolimits_{n,{\cal O}_{x}}$ si $x$ est
d\'{e}compos\'{e} dans $X'$. Par contre, si $x$ est inerte dans $X'$,
cette restriction est un sch\'{e}ma en groupes unitaires non
ramifi\'{e}.

Le choix de la forme hermitienne $\Phi_{n}$ assure que $G$ est
quasi-d\'{e}ploy\'{e}: le drapeau standard
$$
(0)\subset {\cal O}_{X'}\subset {\cal O}_{X'}^{\oplus 2}\subset\cdots \subset
{\cal O}_{X'}^{\oplus n}
$$
est auto-dual et d\'{e}finit donc un $X$-sch\'{e}ma en sous-groupes de
Borel de $G$.
\vskip 2mm

Pour tout $X$-sch\'{e}ma $S$, un $G_{S}$-torseur peut se d\'{e}crire
concr\`{e}tement comme un fibr\'{e} hermitien $({\cal E},\Phi )$ de rang $n$
form\'{e} d'un fibr\'{e} vectoriel ${\cal E}$ de rang $n$ sur
$X'_{S}$ muni d'une forme hermitienne non d\'{e}g\'{e}n\'{e}r\'{e}e,
c'est-\`{a}-dire d'un isomorphisme
$$
\Phi :{\cal E}\buildrel\sim\over\longrightarrow \tau_{S}^{\ast}{\cal
E}^{\vee}
$$
dont le transpos\'{e} ${}^{{\rm t}}\Phi$ est \'{e}gal \`{a}
$\tau_{S}^{\ast}\Phi$.  On a not\'{e} ${\cal E}^{\vee}=\mathop{{\cal
H}{\it om}}\nolimits_{{\cal O}_{X_{S}'}}({\cal E},{\cal O}_{X_{S}'})$ le
fibr\'{e} vectoriel dual de ${\cal E}$.  Pour abr\'{e}ger nous
appellerons parfois la donn\'{e}e d'un tel isomorphisme $\Phi$ une
structure unitaire sur le fibr\'{e} vectoriel ${\cal E}$. Le $G_{S}$-torseur
correspondant est
$$
{\cal T}=\mathop{{\cal I}{\it som}}\nolimits_{{\cal O}_{X'_{S}}}(({\cal
O}_{X'}^{n},\Phi_{n}),({\cal E},\Phi ))
$$
muni de l'action \'{e}vidente \`{a} droite de $G_{S}$.

\subsection{2.2}{Paires de Hitchin}

\`{A} tout $G_{S}$-torseur ${\cal T}$ comme ci-dessus on peut associer
le fibr\'{e} vectoriel $\mathop{\rm ad}({\cal T})$ sur $S$ d\'{e}duit
de ${\cal T}$ par la repr\'{e}sentation adjointe de $G$.  Si ${\cal
T}$ correspond au fibr\'{e} hermitien $({\cal E},\Phi )$, $\mathop{\rm
ad}({\cal T})$ n'est autre que le sous-fibr\'{e} vectoriel de rang
$n^{2}$ de $\pi_{S,\ast}\mathop{{\cal E}{\it nd}}\nolimits_{{\cal
O}_{X'_{S}}}({\cal E})$ form\'{e} des endomorphismes hermitiens de
$({\cal E},\Phi )$.

Soit $D$ un diviseur effectif sur $X$ de degr\'{e} $\geq g+1$.

Pour chaque entier $i$, on note $(-)(iD)$ le foncteur
$(-)\otimes_{{\cal O}_{X}}{\cal O}_{X}(iD)$ de la cat\'{e}gorie des
${\cal O}_{S}$-Modules dans elle-m\^{e}me sur n'importe quel
$X$-sch\'{e}ma $S$.

Une {\it paire de Hitchin} sur un $k$-sch\'{e}ma $S$ (\`{a} valeurs
dans $2D$) est un couple $({\cal T},\theta )$ o\`{u} ${\cal T}$ est un
$G_{S\times_{k}X}$-torseur et o\`{u}
$$
\theta\in H^{0}(S\times_{k}X,\mathop{\rm ad}({\cal T})(2D)).
$$
En termes concrets, une paire de Hitchin sur un $k$-sch\'{e}ma $S$
\`{a} valeurs dans $2D$ est un triplet, dit aussi de Hitchin, $({\cal
E},\Phi ,\theta )$ o\`{u} $({\cal E},\Phi )$ est un fibr\'{e}
hermitien de rang $n$ sur $S\times_{k}X'$ et o\`{u}
$$
\theta :{\cal E}\rightarrow {\cal E}(2D)
$$
est un homomorphisme de fibr\'{e}s vectoriels sur $S\times_{k}X'$ tel
que
$$
\Phi (2D)\circ\theta +\tau^{\ast}({}^{{\rm t}}\!\theta )(2D)\circ\Phi =0.
$$

On consid\`{e}re le $k$-champ  ${\cal M}$ classifiant ces
paires.

\thm PROPOSITION 2.2.1
\enonce
Le $k$-champ ${\cal M}$ est alg\'{e}brique et localement de type fini
sur $k$.
\endthm

\rem D\'{e}monstration
\endrem
On sait que le champ des fibr\'{e}s vectoriels ${\cal E}$ de rang $n$
sur $X'$ est alg\'{e}brique localement de type fini sur $k$.  Pour
prouver la proposition, il suffit donc de montrer que les morphismes
d'oubli
$$
({\cal E},\Phi )\mapsto {\cal E}
$$
et
$$
({\cal E},\Phi ,\theta)\mapsto ({\cal E},\Phi )
$$
sont repr\'{e}sentables et de type fini, ce qui est \'{e}vident.
\hfill\hfill$\square$

\subsection{2.3}{Section de Kostant}

Soient $K'/K$ une extension quadratique de corps de caract\'{e}ristique
diff\'{e}rente de $2$ et $a\mapsto\overline{a}$ l'\'{e}l\'{e}ment non
trivial du groupe de Galois de $K'/K$.  Consid\'{e}rons le groupe
unitaire
$$
G=\{g\in \mathop{\rm GL}\nolimits_{n}(K')\mid {}^{{\rm
t}}\overline{g}\Phi_{n}g=\Phi_{n}\}
$$
et son alg\`{e}bre de Lie
$$
{\frak g}=\{\xi\in\mathop{\rm gl}\nolimits_{n}(K')\mid {}^{{\rm
t}}\overline{\xi}\Phi_{n} +\Phi_{n}\xi =0\}.
$$
Le polyn\^{o}me caract\'{e}ristique
$$
T^{n}+a_{1}T^{n_1}+\cdots +a_{n}\in K'[T]
$$
de tout $\xi\in {\frak g}$ v\'{e}rifie n\'{e}cessairement
$$
\overline{a}_{i}=(-1)^{i}a_{i},~\forall i=1,\ldots ,n.
$$
Kostant a d\'{e}fini une classe de section
$$
\bigoplus_{i=1}^{n}\{a\in K'\mid \overline{a}=(-1)^{i}a\}\rightarrow
{\frak g}
$$
du morphisme polyn\^{o}me caract\'{e}ristique.

En voici un exemple.  Dans l'anneau de polyn\^{o}mes ${\Bbb Z}[{1\over
2}][a_{1},\ldots ,a_{n}]$ muni de la graduation pour laquelle $a_{i}$
est de degr\'{e} $i$ quel que soit $i=1,\ldots ,n$, il existe des
\'{e}l\'{e}ments
$$
b_{1}={a_{1}\over 2},~b_{2}={a_{2}\over 2}-{a_{1}^{2}\over
8},~b_{3}={a_{3}\over 2}-{a_{1}a_{2}\over 4}+{a_{1}^{3}\over 16},~\ldots ,~
b_{i}={a_{i}\over 2}+c_{i}(a_{1},\ldots
,a_{i-1}),~\ldots
$$
tels que $b_{i}$ soit homog\`{e}ne de degr\'{e} $i$ quel que soit
$i=1,\ldots ,n$ et que la matrice
$$
\pmatrix{-b_{1} & -b_{2} & \cdots & \cdots & -b_{n-1} & -2b_{n}\cr
1 & 0 & \cdots & \cdots & 0 & -b_{n-1}\cr
0 & \ddots & \ddots &  & \vdots & \vdots\cr
\vdots & \ddots & \ddots & \ddots & \vdots & \vdots\cr
\vdots &  & \ddots & \ddots & 0 & -b_{2}\cr
0 & \cdots & \cdots & 0 & 1 & -b_{1}\cr}
$$
ait pour polyn\^{o}me caract\'{e}ristique
$$
T^{n}+a_{1}T^{n_1}+\cdots +a_{n}.
$$
L'application
$$
\bigoplus_{i=1}^{n}\{a\in K'\mid \overline{a}=(-1)^{i}a\}\rightarrow
{\frak g}
$$
qui envoie $(a_{1},\ldots ,a_{n})$ sur la matrice ci-dessus est une
section de Kostant.

\subsection{2.4}{Morphisme de Hitchin}

L'{\it espace de Hitchin} ${\Bbb A}$ est le $k$-sch\'{e}ma affine
naturellement associ\'{e} au sous-$k$-espace vectoriel
$$
\bigoplus_{i=1}^{n}H^{0}(X',{\cal O}_{X'}(2iD))^{\tau^{\ast}=(-1)^{i}}\subset
\bigoplus_{i=1}^{n}H^{0}(X',{\cal O}_{X'}(2iD)).
$$
Si on note
$$
\pi_{\ast}{\cal O}_{X'}={\cal O}_{X}\oplus {\cal L}
$$
la d\'{e}composition en sous-${\cal O}_{X}$-Modules propres pour
l'automorphisme $\tau^{\ast}$, ${\Bbb A}$ est encore le
$k$-sch\'{e}ma affine naturellement associ\'{e} au $k$-espace vectoriel
$$
\bigoplus_{i=1}^{n}H^{0}(X,({\cal L}_{D})^{\otimes i})
$$
o\`{u} on a pos\'{e} ${\cal L}_{D}={\cal L}(2D)$ (puisque ${\cal
L}^{\otimes 2}$ est canoniquement isomorphe \`{a} ${\cal O}_{X}$).

On d\'{e}finit la {\it caract\'{e}ristique} d'un triplet de Hitchin
$({\cal E},\Phi ,\theta )$ sur un $k$-sch\'{e}ma $S$ comme le
$S$-point de ${\Bbb A}$ suivant: pour chaque entier $i=1,\ldots ,n$,
on consid\`{e}re la trace de l'homomorphisme
$$
\wedge^{i}\theta :\bigwedge_{{\cal O}_{X_{S}'}}^{i}{\cal E}
\rightarrow \Bigl(\bigwedge_{{\cal O}_{X_{S}'}}^{i}{\cal
E}\Bigr)(2iD)
$$
qui est une section globale de ${\cal O}_{X_{S}'}(2iD)$;
on a
$$
\tau^{\ast}(\mathop{\rm tr}\wedge^{i}\theta )=(-1)^{i}\mathop{\rm
tr}\wedge^{i}\theta
$$
puisque $\Phi (2D)\circ\theta +(\tau^{\ast}\,{}^{{\rm t}}\!\theta
)(2D)\circ\Phi =0$, et $\oplus_{i=1}^{n}(-1)^{i}\mathop{\rm
tr}\wedge^{i}\theta$ est donc un $S$-point de ${\Bbb A}$. Ci-dessus
on a not\'{e} $X_{S}'=S\times_{k}X'$

Le {\it morphisme de Hitchin} est le morphisme de champs alg\'{e}briques
$$
f:{\cal M}\rightarrow {\Bbb A}
$$
qui associe \`{a} chaque triplet de Hitchin $({\cal E},\Phi
,\theta )$ sa caract\'{e}ristique.

Pour tout choix d'une section de Kostant comme dans (2.3), on a une
section correspondante du morphisme de Hitchin, que l'on appelle encore
{\it section de Kostant}.  Cette derni\`{e}re section associe \`{a} tout
point $a$ de ${\Bbb A}$ le triplet $({\cal E},\Phi ,\theta )$ o\`{u}
${\cal E}=\bigoplus_{i=1}^{n}{\cal O}_{X'}((n+1-2i)D)$, o\`{u} $\Phi$
a pour matrice $\Phi_{n}$ et o\`{u} $\theta$ est donn\'{e} par la
matrice de Kostant ci-dessus.

\subsection{2.5}{Courbes spectrales}

Rappelons la construction de la courbe spectrale associ\'{e}e \`{a}
un $S$-point $a$ de ${\Bbb A}$ (cf. [B-N-R]).

Soit $p:\Sigma={\Bbb P}({\cal O}_{X}\oplus ({\cal
L}_{D})^{\otimes -1})\rightarrow X$ la compl\'{e}tion projective du
fibr\'{e} en droites $p^{\circ}:\Sigma^{\circ}={\Bbb V}(({\cal
L}_{D})^{\otimes -1})\rightarrow X$ dont les sections sont celles de ${\cal
L}_{D}$; c'est un fibr\'{e} en droites projectives.  On a
$p_{\ast}{\cal O}_{\Sigma}(1)={\cal O}_{X}\oplus ({\cal
L}_{D})^{\otimes -1}$ et la section $(1,0)$ de cette image directe
d\'{e}finit donc une section globale $V$ de ${\cal O}_{\Sigma}(1)$; de
m\^{e}me, on a $p_{\ast}({\cal O}_{\Sigma}(1)\otimes_{{\cal
O}_{\Sigma}}p^{\ast}{\cal L}_{D})={\cal L}_{D}\oplus {\cal O}_{X}$ et
la section $(0,1)$ de cette image directe d\'{e}finit une section
globale $U$ de ${\cal O}_{\Sigma}(1)\otimes_{{\cal
O}_{\Sigma}}p^{\ast}{\cal L}_{D}$.  Le couple $(U;V)$ est un
syst\`{e}me de coordonn\'{e}es homog\`{e}nes relatives sur $\Sigma$;
le lieu des z\'{e}ros de $U$ (resp.  $V$) est la section nulle ${\Bbb
P}({\cal O}_{X})\subset \Sigma^{\circ}$ (resp.  la section \`{a}
l'infini ${\Bbb P}(({\cal L}_{D})^{\otimes -1})=\Sigma-\Sigma^{\circ}$) de
$\Sigma^{\circ}$.

Si $a=\oplus_{i=1}^{n}a_{i}$ est un point de ${\Bbb A}$ \`{a} valeurs
dans un $k$-sch\'{e}ma $S$, on a la
section
$$
U^{n}+(p^{\ast}a_{1})VU^{n-1}+\cdots +(p^{\ast}a_{n})V^{n}
$$
de ${\cal O}_{S}\boxtimes_{k}({\cal O}_{\Sigma}(n)\otimes_{{\cal
O}_{\Sigma}}p^{\ast}({\cal L}_{D})^{\otimes n})$ dont le lieu des
z\'{e}ros est une $S$-courbe projective $Y_{a}$ trac\'{e}e sur la
$S$-surface projective $\Sigma_{S}=S\times_{k}\Sigma$.  Cette courbe
est par construction un rev\^{e}tement ramifi\'{e} de degr\'{e} $n$ de
$X_{S}=S\times_{k}X$ par la restriction $p_{a}:Y_{a}\rightarrow X_{S}$
\`{a} $Y_{a}$ de la projection $p_{S}:\Sigma_{S}\rightarrow X_{S}$.

On remarquera que la courbe spectrale $Y_{a}$ ne coupe pas la section
infini et est par cons\'{e}quent enti\`{e}rement contenue dans la carte
$\Sigma_{S}^{\circ}=S\times_{k}\Sigma^{\circ}=\{V\not=0\}={\Bbb V}({\cal
O}_{S}\boxtimes_{k}({\cal L}_{D})^{\otimes -1})$.  C'est donc aussi le lieu
des z\'{e}ros dans $\Sigma_{S}^{\circ}$ de la section
$$
u^{n}+((p^{\circ})^{\ast}a_{1})u^{n-1}+\cdots +((p^{\circ})^{\ast}a_{n})
$$
de ${\cal O}_{S}\boxtimes_{k}(p^{\circ})^{\ast}({\cal L}_{D})^{\otimes
n}$ o\`{u} $u={U\over V}$.  En particulier, la ${\cal
O}_{X_{S}}$-Alg\`{e}bre $p_{a,\ast}{\cal O}_{Y_{a}}$ est
isomorphe \`{a}
$$
({\cal O}_{S}\boxtimes_{k}\mathop{\rm Sym}\nolimits_{{\cal
O}_{X}}(({\cal L}_{D})^{\otimes -1}))/{\cal I}_{a}
$$
o\`{u} ${\cal I}_{a}$ est l'Id\'{e}al engendr\'{e} par l'image de
l'homomorphisme
$$
{\cal O}_{S}\boxtimes_{k}({\cal L}_{D})^{\otimes -n}
\rightarrow\bigoplus_{i=0}^{n}{\cal O}_{S}\boxtimes_{k}({\cal
L}_{D})^{\otimes -i}\subset {\cal O}_{S}\boxtimes_{k}
\mathop{\rm Sym}\nolimits_{{\cal O}_{X}}(({\cal L}_{D})^{\otimes -1})
$$
de composantes $(a_{n},a_{n-1},\ldots ,a_{1},1)$.

On note $D(a)$ le discriminant de la caract\'{e}ristique $a$,
c'est-\`{a}-dire le r\'{e}sultant du polyn\^{o}me
$$
u^{n}+a_{1}u^{n-1}+\cdots +a_{n}
$$
et de sa d\'{e}riv\'{e}e
$$
nu^{n-1}+(n-1)a_{1}u^{n-2}+\cdots +a_{n-1};
$$
c'est une section globale de $({\cal L}_{D})^{\otimes n(n-1)}$.

En fait, on a une courbe spectrale universelle
$$\diagram{
Y&\kern -1mm\subset\kern -1mm&\Sigma\times_{k}{\Bbb A}\cr
\llap{$\scriptstyle $}\left\downarrow
\vbox to 4mm{}\right.\rlap{}&&\cr
{\Bbb A}&&\cr
\arrow(13,11)\dir(-1,-1)\length{9}}
$$
dont le changement de base par $a:S\rightarrow {\Bbb A}$ est $Y_{a}$.
Le morphisme $Y\rightarrow {\Bbb A}$ est projectif, plat,
localement d'intersection compl\`{e}te, purement de dimension relative
$1$.  Sa fibre la plus mauvaise est celle en $a=0\in {\Bbb A}$: c'est
le lieu des z\'{e}ros de $U^{n}$, c'est-\`{a}-dire la section nulle de
$\Sigma^{\circ}\rightarrow X$ compt\'{e}e avec multiplicit\'{e} $n$.

\thm LEMME 2.5.1
\enonce
Pour tout point g\'{e}om\'{e}trique $a$ de ${\Bbb A}$ les conditions
suivantes sont \'{e}quivalentes (on rappelle que $p>n$):
\decale{\rm (i)} la courbe spectrale $Y_{a}$ est r\'{e}duite,

\decale{\rm (ii)} le rev\^{e}tement $Y_{a}\rightarrow X$ est
\'{e}tale au-dessus du point g\'{e}n\'{e}rique de $X$,

\decale{\rm (iii)} le discriminant $D(a)$ n'est pas identiquement nul.
\hfill\hfill$\square$
\endthm

Le plus grand ouvert ${\Bbb A}^{{\rm red}}\subset {\Bbb A}$ au-dessus
duquel $Y\rightarrow {\Bbb A}$ est \`{a} fibres
g\'{e}om\'{e}triquement r\'{e}duites est donc l'ouvert des $a$ tels
que $D(a)$ n'est pas identiquement nul.

Ce lieu contient l'ouvert l'ouvert ${\Bbb A}^{{\rm lisse}}$ au-dessus
duquel $Y\rightarrow {\Bbb A}$ est lisse.

Ces lieux sont non vides d\`{e}s que $({\cal L}_{D})^{\otimes n}$
admet une section $a_{n}$ qui n'a que des z\'{e}ros simples
puisqu'alors la courbe spectrale $Y_{a}$ o\`{u} $a=0\oplus\cdots\oplus
0\oplus a_{n}$ est lisse.  D'apr\`{e}s le th\'{e}or\`{e}me de
Riemann-Roch et le th\'{e}or\`{e}me de Bertini, ces lieux sont donc
non vides d\`{e}s que $({\cal L}_{D})^{\otimes n}$ est tr\`{e}s ample,
c'est-\`{a}-dire d\`{e}s que $2n\mathop{\rm deg}(D)\geq 2g+1$.

\thm PROPOSITION 2.5.2
\enonce
La restriction ${\cal M}^{{\rm red}}$ de ${\cal M}$ au-dessus de
l'ouvert ${\Bbb A}^{{\rm red}}\subset {\Bbb A}$ est lisse sur ${\Bbb
F}_{q}$.
\endthm

\rem D\'{e}monstration 
\endrem
Soit $({\cal E},\Phi ,\theta )$ un triplet de Hitchin dont la
caract\'{e}ristique $a$ est dans l'ouvert ${\Bbb A}^{{\rm red}}$ de
${\Bbb A}$.  Pour all\'{e}ger les notations nous supposons dans la
suite qu'il s'agit d'un $k$-point de ${\cal M}$, mais l'argument est
g\'{e}n\'{e}ral.

La fibre en ce point du complexe tangent \`{a} ${\cal M}$ est le
complexe $R\Gamma (X,K)$ o\`{u}
$$
K=[(\pi_{\ast}\mathop{{\cal E}{\it nd}}\nolimits_{{\cal O}_{X'}}({\cal
E}))^{\tau^{\ast}=-1}\rightarrow (\pi_{\ast}\mathop{{\cal E}{\it
nd}}\nolimits_{{\cal O}_{X'}}({\cal E}))^{\tau^{\ast}=-1}\otimes_{{\cal
O}_{X}}{\cal L}_{D}]
$$
est un complexe parfait concentr\'{e} en degr\'{e}s $0$ et $1$, avec
pour diff\'{e}rentielle l'application $\xi\rightarrow [\theta ,\xi ]$.
Il s'agit de voir que $H^{2}(X,K)=(0)$.

En utilisant la forme de Killing, on peut identifier le dual du
complexe $K$ au complexe $K\otimes_{{\cal O}_{X}}{\cal L}_{D}^{\otimes
-1}$.  Par dualit\'{e} de Serre on est donc ramen\'{e} \`{a}
d\'{e}montrer que
$$
H^{0}(X,{\cal H}^{0}(K)\otimes_{{\cal O}_{X}}{\cal L}_{D}^{\otimes -1}
\otimes_{{\cal O}_{X}}\Omega_{X/k}^{1})=(0).
$$
Or
$$
{\cal H}^{0}(K)=(\pi_{\ast}p_{a,\ast}'\mathop{{\cal E}{\it nd}}
\nolimits_{{\cal O}_{Y_{a}'}}({\cal F}))^{\tau^{\ast}=-1}
$$
o\`{u} $({\cal F},\iota )\in \overline{P}_{a}(k)$ est le point
correspondant \`{a} $({\cal E},\Phi ,\theta )$.  Soit $\rho
:\widetilde{Y}_{a}\rightarrow Y_{a}$ la normalisation de la courbe
r\'{e}duite $Y_{a}$ et $\rho':\widetilde{Y}{}_{a}'=X'\times_{X}
\widetilde{Y}_{a}\rightarrow Y_{a}'$ son changement de base par $\pi$.
On a une injection naturelle
$$
\mathop{{\cal E}{\it nd}}\nolimits_{{\cal O}_{Y_{a}'}}({\cal F})
\hookrightarrow\rho_{\ast}'{\cal O}_{\widetilde{Y}{}_{a}'}.
$$
En effet, si $A$ est un anneau local de $Y_{a}'$ et si $\widetilde{A}$ est la 
normalisation de $A$ dans son anneau total des fractions $\mathop{\rm 
Frac}(A)$, pour tout $A$-module $M$ sans torsion de rangs 
g\'{e}n\'{e}riques $1$, on a 
$$
A\subset\mathop{\rm End}\nolimits_{A}(M)\subset \mathop{\rm
End}\nolimits_{\mathop{\rm Frac}(A)}(\mathop{\rm
Frac}(A)\otimes_{A}M)=\mathop{\rm Frac}(A)
$$
et donc $\mathop{\rm End}\nolimits_{A}(M)\subset \widetilde{A}$ 
puisque $\mathop{\rm End}\nolimits_{A}(M)$ est de type fini sur $A$. Par 
suite
$$
{\cal H}^{0}(K)\subset (\pi_{\ast}p_{a,\ast}'\rho_{\ast}' {\cal
O}_{\widetilde{Y}{}_{a}'})^{\tau^{\ast}=-1}\subset
\pi_{\ast}p_{a,\ast}'\rho_{\ast}' {\cal O}_{\widetilde{Y}{}_{a}'}
$$
et
$$
H^{0}(X,{\cal H}^{0}(K)\otimes_{{\cal O}_{X}}{\cal L}_{D}^{\otimes -1}
\otimes_{{\cal O}_{X}}\Omega_{X/k}^{1})\subset
H^{0}(\widetilde{Y}_{a}^{\,\prime},\rho'^{\ast}p_{a}'^{\ast}
\pi^{\ast}\Omega_{X/k}^{1}(-2D))
$$
est nul puisque le degr\'{e} de $\Omega_{X/k}^{1}(-2D)$ est
strictement n\'{e}gatif par hypoth\`{e}se sur $D$.
\hfill\hfill$\square$
\vskip 3mm

\rem Remarque 
\endrem
La proposition ci-dessus pourrait aussi se d\'{e}duire d'un r\'{e}sultat 
de Fantchi, G\"{o}ttsche et van Straten (cf. la section A de [F-G-S]).
\hfill\hfill$\square$
\vskip 3mm

\thm LEMME 2.5.3
\enonce
Soient $a$ un point g\'{e}om\'{e}trique de ${\Bbb A}^{{\rm red}}$ et
$Z$ une composante irr\'{e}ductible de $Y_{a}$.  Alors, il existe un
unique entier $m$ compris entre $1$ et $n$ et une famille unique de
sections $b_{j}\in \kappa (a)\otimes_{k} H^{0}(X,({\cal
L}_{D})^{\otimes j})$, $j=1,\ldots m$, tels que $Z$ soit le diviseur
de Cartier sur ${\Bbb V}(({\cal L}_{D})^{\otimes -1})$ d\'{e}fini
par l'\'{e}quation
$$
u^{m}+((p^{\circ})^{\ast}b_{1})u^{m-1}+\cdots +((p^{\circ})^{\ast}b_{m})=0.
$$
\endthm

\rem D\'{e}monstration
\endrem
Notons simplement $L$ la fibre de ${\cal L}_{D}$ au point
g\'{e}n\'{e}rique de $\kappa (a)\otimes_{k}X$.  Comme $Z$ est un
rev\^{e}tement fini g\'{e}n\'{e}riquement \'{e}tale de $\kappa
(a)\otimes_{k}X$ de degr\'{e} $m$ compris entre $1$ et $n$, la
th\'{e}orie de Galois assure qu'il existe des uniques $b_{j}\in
L^{\otimes j}$ et des uniques $c_{k}\in L^{\otimes k}$ tels que
$$
u^{n}+a_{1}u^{n-1}+\cdots +a_{n}=(u^{m}+b_{1}u^{m-1}+\cdots
+b_{n})(u^{n-m}+c_{1}u^{n-m-1}+\cdots +c_{n-m}).
$$
Il ne reste plus qu'\`{a} v\'{e}rifier que chaque $b_{j}$ est en fait
dans $\kappa (a)\otimes_{k}H^{0}(X,({\cal L}_{D})^{\otimes j})\subset
L^{\otimes j}$.

Il revient au m\^{e}me de se donner la section globale $a_{i}$ de $({\cal
L}_{D})^{\otimes i}$ ou de se donner une section
$$
a_{i}:\kappa (a)\otimes_{k}X\rightarrow \kappa (a)\otimes_{k}
{\Bbb V}(({\cal L}_{D})^{\otimes -i})
$$
de la projection canonique, ou encore de se donner un morphisme
$$
a_{i}:\kappa (a)\otimes_{k}{\Bbb V}({\cal L}_{D})\rightarrow {\Bbb
A}_{\kappa (a)}^{1}
$$
qui est ${\Bbb G}_{{\rm m},\kappa (a)}$-\'{e}quivariant au sens o\`{u}
$a_{i}(tv)=t^{i}a_{i}(v)$.  Par suite, la donn\'{e}e de $a$
\'{e}quivaut \`{a} celle d'un morphisme ${\Bbb G}_{{\rm m},\kappa
(a)}$-\'{e}quivariant
$$
\kappa (a)\otimes_{k}{\Bbb V}({\cal L}_{D})\rightarrow {\Bbb
A}_{\kappa (a)}^{n}
$$
et on cherche \`{a} factoriser ce morphisme en
$$
\kappa (a)\otimes_{k}{\Bbb V}({\cal L}_{D})\rightarrow {\Bbb
A}_{\kappa (a)}^{m}\times_{\kappa (a)}{\Bbb
A}_{\kappa (a)}^{n-m}\rightarrow {\Bbb
A}_{\kappa (a)}^{n}
$$
o\`{u} la seconde fl\`{e}che envoie $((y_{1},\ldots
,y_{m}),(z_{1},\ldots ,z_{n-m}))$ sur les coefficients
$(y_{1}+z_{1},y_{2}+z_{2}+y_{1}z_{1},\cdots ,y_{m}z_{n})$ du
polyn\^{o}me produit $(T^{m}+y_{1}T^{m-1}+\cdots
+y_{m})(T^{n-m}+z_{1}T^{n-m-1}+\cdots +z_{n-m})$.

Or cette seconde fl\`{e}che est un rev\^{e}tement (ramifi\'{e}) fini
qui est ${\Bbb G}_{{\rm m},\kappa (a)}$-\'{e}quivariant et on a d\'{e}j\`{a}
une factorisation au dessus du point g\'{e}n\'{e}rique de $\kappa
(a)\otimes_{k}X$ par l'argument de th\'{e}orie de Galois qui
pr\'{e}c\`{e}de. Par suite, on a par changement de base un
rev\^{e}tement fini ${\Bbb G}_{{\rm m},\kappa (a)}$-\'{e}quivariant de
$\kappa (a)\otimes_{k}{\Bbb V}({\cal L}_{D})$ est une section de ce
rev\^{e}tement au-dessus du point g\'{e}n\'{e}rique de $\kappa
(a)\otimes_{k}X$. En prenant l'adh\'{e}rence $Z$ de cette section, on
obtient un morphisme ${\Bbb G}_{{\rm m},\kappa (a)}$-\'{e}quivariant
$$
Z\rightarrow \kappa (a)\otimes_{k}{\Bbb V}({\cal L}_{D})
$$
qui est fini et un isomorphisme au-dessus du point g\'{e}n\'{e}rique de $\kappa
(a)\otimes_{k}X$. Un tel morphisme est n\'{e}cessairement un
isomorphisme puisque $\kappa (a)\otimes_{k}{\Bbb V}({\cal L}_{D})$ est
normal.
\hfill\hfill$\square$

\subsection{2.6}{Champs de Picard}

Pour tout $S$-point $a$ de ${\Bbb A}^{{\rm red}}$, on note
$\pi_{a}:Y_{a}'=X'\times_{X}Y_{a}\rightarrow Y_{a}$ le rev\^{e}tement
double \'{e}tale d\'{e}duit du rev\^{e}tement double \'{e}tale $\pi
:X'\rightarrow X$ et $p_{a}':Y_{a}'\rightarrow X'$ la projection
canonique.  L'involution $\tau$ de $X'$ au-dessus de $X$ induit une
involution not\'{e}e encore $\tau$ de $Y_{a}'$ au-dessus de $Y_{a}$.

Le {\it champ de Picard relatif} de la $S$-courbe $Y_{a}'$ (\`{a}
fibres g\'{e}om\'{e}triquement r\'{e}duites) est le champ des ${\cal
O}_{Y_{a}'}$-Modules inversibles.  C'est un champ alg\'{e}brique
localement de type fini sur $S$ que l'on note $\mathop{\rm
Pic}\nolimits_{Y_{a}'/S}$.  Ce champ est naturellement muni d'une
structure de groupe induite par le produit tensoriel. En fait, c'est
une ${\Bbb G}_{{\rm m}}$-gerbe sur le sch\'{e}ma en groupes de Picard
relatif de $Y_{a}'/S$ qui existe sous nos hypoth\`{e}se.  Le champ
$\mathop{\rm Pic}\nolimits_{Y_{a}'/S}$ est de plus muni d'une
involution compatible \`{a} sa structure de groupe, qui envoie ${\cal
F}$ sur ${\cal F}^{\otimes -1}=\mathop{{\cal H}{\it
om}}\nolimits_{{\cal O}_{Y_{a}'}}({\cal F},{\cal O}_{Y_{a}'})$.

Le {\it champ de Picard compactifi\'{e} relatif} de $Y_{a}'$ est le
champ des ${\cal O}_{Y_{a}'}$-Modules coh\'{e}rents ${\cal F}$ {\it
sans torsion de rang} $1$, c'est-\`{a}-dire $S$-plats et dont la
restriction \`{a} chaque fibre de $Y_{a}'\rightarrow S$ est partout
sans torsion et de rang $1$ en chaque point g\'{e}n\'{e}rique.  C'est
un champ alg\'{e}brique localement de type fini sur $S$ que l'on note
$\mathop{\overline{\rm Pic}}\nolimits_{Y_{a}'/S}$.  Il contient
$\mathop{\rm Pic}\nolimits_{Y_{a}'/S}$ comme un ouvert qui est dense
fibre \`{a} fibre de sa projection sur $S$ puisque $Y_{a}'$ est
plong\'{e}e dans une surface relative sur $S$ (cf.  [Reg] et [A-I-K]).
Le champ $\mathop{\overline{\rm Pic}}\nolimits_{Y_{a}'/S}$ est
naturellement muni d'une action de $\mathop{\rm
Pic}\nolimits_{Y_{a}'/S}$ qui prolonge l'action par translation de
$\mathop{\rm Pic}\nolimits_{Y_{a}'/S}$ sur lui-m\^{e}me et qui est
induite par le produit tensoriel.  Il est de plus muni d'une
involution qui envoie ${\cal F}$ sur
$$
{\cal F}^{\vee}=\mathop{{\cal H}{\it om}}\nolimits_{{\cal
O}_{Y_{a}'}}({\cal F},\omega_{Y_{a}'/X_{S}'})
$$
o\`{u} $\omega_{Y_{a}'/X_{S}'}=\omega_{Y_{a}'}\otimes_{{\cal
O}_{Y_{a}'}}(p_{a}'^{\ast}\omega_{X_{S}'/S})^{\otimes -1}$ est le
Module dualisant relatif de $Y_{a}'$ au-dessus de
$X_{S}'=S\times_{k}X'$.  L'action du champ de Picard sur le champ de
Picard compactifi\'{e} est compatible aux involutions que l'on vient
de d\'{e}finir.

\rem Remarque
\endrem
Pour tout $S$-point $a$ de ${\Bbb A}^{{\rm red}}$ on a
$$
\omega_{Y_{a}/S\times_{k}X}=p_{a}^{\ast}({\cal L}_{D})^{\otimes n-1}
$$
puisque $\Omega_{\Sigma^{\circ}/X}^{1}$ est une extension de
$(p^{\circ})^{\ast}\Omega_{X}^{1}$ par $(p^{\circ})^{\ast}{\cal
L}_{D}$ et que
$$
({\cal I}_{a}/{\cal I}_{a}^{2})|Y_{a}=p_{a}^{\ast}({\cal
L}_{D})^{\otimes -n}
$$
o\`{u} ${\cal I}_{a}$ est l'Id\'{e}al qui d\'{e}finit $Y_{a}$ dans
$S\times_{k}\Sigma^{\circ}$.  Comme $X'\rightarrow X$ est \'{e}tale,
on a aussi $\omega_{Y_{a}'/S\times_{k}X'}=p_{a}'^{\ast}{\cal
O}_{X'}(2(n-1)D)$.
\hfill\hfill$\square$
\vskip 3mm

L'involution $\tau$ de $Y_{a}'$ induit une autre involution
$\tau^{\ast}$ sur les champs $\mathop{\rm Pic}\nolimits_{Y_{a}'/S}$
et  $\mathop{\overline{\rm Pic}}\nolimits_{Y_{a}'/S}$. On
note
$$
P_{a}=\bigl(\mathop{\rm
Pic}\nolimits_{Y_{a}'/S}\bigr)^{\tau^{\ast}=(-)^{\otimes -1}}
$$
la partie primitive pour cette involution, c'est-\`{a}-dire le champ
des couples $({\cal F},\iota )$ o\`{u} ${\cal F}$ est un ${\cal
O}_{Y_{a}'}$-Module inversible et o\`{u} $\iota :{\cal
F}\buildrel\sim\over\longrightarrow \tau^{\ast}{\cal F}^{\otimes -1}$
est un isomorphisme de ${\cal O}_{Y_{a}'}$-Modules tel que $\iota
=\tau^{\ast}(\iota^{\otimes -1})$.  De m\^{e}me on note
$$
\overline{P}_{a}=\bigl(\mathop{\overline{\rm Pic}}
\nolimits_{Y_{a}'/S}\bigr)^{\tau^{\ast}=(-)^{\vee}}
$$
le champ des couples $({\cal F},\iota )$ o\`{u} ${\cal F}$ est un
${\cal O}_{Y_{a}'}$-Module coh\'{e}rent sans torsion de rang $1$ et
o\`{u} $\iota :{\cal F}\buildrel\sim\over\longrightarrow
\tau^{\ast}{\cal F}^{\vee}$ est un isomorphisme de ${\cal
O}_{Y_{a}'}$-Modules tel que $\iota =\tau^{\ast}(\iota^{\vee})$.  On a
encore une action de $P_{a}$ sur $\overline{P}_{a}$ mais on n'a plus a
priori de plongement de $P_{a}$ dans $\overline{P}_{a}$.  Un tel
plongement existe apr\`{e}s le choix d'une section, par exemple une
section de Kostant.

Bien s\^{u}r, pour $S={\Bbb A}^{{\rm red}}$ et $a$ l'identit\'{e} de
${\Bbb A}^{{\rm red}}$, on obtient des champs universels
$P\rightarrow {\Bbb A}^{{\rm red}}$ et $\overline{P}\rightarrow {\Bbb
A}^{{\rm red}}$.

On remarque que, pour tout ${\cal F}\in\overline{P}_{a}(S)$,
$p_{a,\ast}'{\cal F}$ est un fibr\'{e} vectoriel de rang $n$ sur
$X_{S}'=S\times_{k}X'$ et que
$$
p_{a,\ast}'({\cal F}^{\vee})=(p_{a,\ast}'{\cal F})^{\vee}
(:=\mathop{{\cal H}{\it om}}\nolimits_{{\cal O}_{X_{S}'}}
(p_{a,\ast}'{\cal F},{\cal O}_{X_{S}'})).
$$
On a donc un morphisme de ${\Bbb A}^{{\rm red}}$-champs
$$
\overline{P}\rightarrow {\Bbb A}^{{\rm red}}\times_{{\Bbb A}}{\cal M}
$$
qui envoie $({\cal F},\iota )\in \overline{P}_{a}(S)$ sur le triplet
de Hitchin $({\cal E},\Phi ,\theta )$ sur $S$ de
caract\'{e}ristique $a$ o\`{u} ${\cal E}=p_{a,\ast}{\cal F}$, o\`{u}
$$
\Phi =p_{a,\ast}'\iota :{\cal E}=p_{a,\ast}'{\cal F}\rightarrow
p_{a,\ast}'\tau^{\ast}{\cal F}^{\vee}=\tau^{\ast}(p_{a,\ast}'{\cal
F})^{\vee} =\tau^{\ast}{\cal E}^{\vee}
$$
v\'{e}rifie l'\'{e}quation $\Phi =\tau^{\ast}(\Phi^{\vee})$ et
o\`{u} $\theta$ est d\'{e}fini par l'homomorphisme
$$
{\cal O}_{X_{S}'}(-2D)\subset \mathop{\rm
Sym}\nolimits_{{\cal O}_{X_{S}'}}({\cal O}_{X_{S}'}(-2D))/{\cal
I}_{a}'=(p_{a}')_{\ast}{\cal O}_{Y_{a}'}\rightarrow \mathop{{\cal E}{\it
nd}}\nolimits_{{\cal O}_{X_{S}'}}({\cal E}),
$$
${\cal I}_{a}'$ \'{e}tant l'Id\'{e}al engendr\'{e} par l'image de
l'homomorphisme
$$
({\cal O}_{X_{S}'}(-2D))^{\otimes -n}\rightarrow
\bigoplus_{i=0}^{n}({\cal O}_{X_{S}'}(-2D))^{\otimes -i}\subset
\mathop{\rm Sym}\nolimits_{{\cal O}_{X_{S}'}}({\cal O}_{X_{S}'}(-2D))
$$
de composantes $(a_{n},a_{n-1},\ldots ,a_{1},1)$.

La proposition suivante est une variante d'un r\'{e}sultat de Beauville,
Narasimhan et Ramanan (cf. [B-N-R]).

\thm PROPOSITION 2.6.1
\enonce
Le morphisme de ${\Bbb A}^{{\rm red}}$-champs
$$
\overline{P}\rightarrow {\Bbb A}^{{\rm red}}\times_{{\Bbb A}}{\cal M}
$$
d\'{e}fini ci-dessus est un isomorphisme.
\endthm

\rem D\'{e}monstration
\endrem
Pour d\'{e}montrer que ce morphisme est un isomorphisme, nous allons
construire un inverse.  Soit $({\cal E},\Phi ,\theta )$ un triplet de
Hitchin sur $S$ de caract\'{e}ristique $a\in {\Bbb A}^{{\rm red}}(S)$.
Comme on l'a vu ci-dessus, la section globale $\theta\in
H^{0}(X_{S}',\mathop{{\cal E}{\it nd}}\nolimits_{{\cal
O}_{X_{S}'}}({\cal E})\otimes_{{\cal O}_{X'}}{\cal O}_{X'}(2D))$ munit
${\cal E}$ d'une structure de $({\cal O}_{S}\boxtimes_{k}\mathop{\rm
Sym}\nolimits_{{\cal O}_{X'}}({\cal O}_{X'}(-2D)))$-Module.  Puisque
ce triplet a pour caract\'{e}ristique $a$, ce $({\cal
O}_{S}\boxtimes_{k}\mathop{\rm Sym}\nolimits_{{\cal O}_{X'}}({\cal
O}_{X'}(-2D)))$-Module est en fait un $({\cal
O}_{S}\boxtimes_{k}\mathop{\rm Sym}\nolimits_{{\cal O}_{X'}}({\cal
O}_{X'}(-2D)))/{\cal I}_{a}'$-Module.  Le ${\cal O}_{Y_{a}'}$-Module
correspondant ${\cal F}$ est alors $S$-plat et fibres par fibres un
Module sans torsion de rang $1$ sur $Y_{a}'$ (voir [B-N-R]).  Comme on
l'a vu ci-dessus, par dualit\'{e}, la donn\'{e}e d'une structure
unitaire $\Phi$ sur ${\cal E}$ est \'{e}quivalente \`{a} la donn\'{e}e
d'un isomorphisme $\iota :{\cal F}\buildrel\sim\over\longrightarrow
\tau^{\ast}({\cal F}^{\vee})$ qui v\'{e}rifie $\iota
=\tau^{\ast}(\iota^{\vee})$.
\hfill\hfill$\square$
\vskip 3mm

En particulier, la restriction \`{a} ${\Bbb A}^{{\rm red}}$ d'une
section de Kostant (cf.  la fin de la section (2.4)) est l'image d'une
section de $\overline{P}$, qui est dite encore de Kostant.

Vu notre choix de ${\cal L}_{D}$, il y a une section de Kostant
$({\cal K},\iota_{{\cal K}})$ de $\overline{P}$ particuli\`{e}rement
jolie.  Elle est donn\'{e}e de la fa\c{c}on suivante.  Pour tout
$S$-point $a$ de ${\Bbb A}^{{\rm red}}$ le ${\cal O}_{Y_{a}'}$-Module
inversible $p_{a}'^{\ast}{\cal O}_{X'}((n-1)D)$ o\`{u}
$\pi_{Y_{a}'}:Y_{a}'\rightarrow Y_{a}$ est induit par $\pi$, est une
racine carr\'{e}e de $\omega_{Y_{a}'/S\times_{k}X'}$, et on pose
$$
a^{\ast}{\cal K}=p_{a}'^{\ast}{\cal O}_{X'}((n-1)D)=
\pi_{Y_{a}'}^{\ast}p_{a}^{\ast}{\cal L}((n-1)D)
$$
et on prend pour
$$
a^{\ast}\iota_{{\cal K}}:\tau^{\ast}(a^{\ast}{\cal K})
\buildrel\sim\over\longrightarrow \omega_{Y_{a}'/S\times_{k}X'}
\otimes_{{\cal O}_{Y_{a}'}}(a^{\ast}({\cal K})^{\otimes -1})=
a^{\ast}{\cal K}
$$
l'isomorphisme de descente de $a^{\ast}{\cal K}$ en $p_{a}^{\ast}
{\cal L}((n-1)D)$.

\subsection{2.7}{Variante endoscopique}

Soit $n_{1}+n_{2}=n$ une partition de $n$ en deux entiers $\geq 1$.
On peut consid\'{e}rer les $X$-sch\'{e}ma en groupes unitaires $G_{1}$
et $G_{2}$ d\'{e}finis comme $G$ mais apr\`{e}s avoir remplac\'{e} $n$
par $n_{1}$ et $n_{2}$ et le $X$-sch\'{e}ma en groupes produit
$$
H=G_{1}\times_{X}G_{2}.
$$

On a des morphismes de Hitchin $f_{1}:{\cal M}_{1}\rightarrow {\Bbb
A}_{1}$ et $f_{2}:{\cal M}_{2}\rightarrow {\Bbb
A}_{2}$ o\`{u} ${\Bbb A}_{\alpha}$ est le $k$-sch\'{e}ma affine
naturellement associ\'{e} au $k$-espace vectoriel
$$
\bigoplus_{i=1}^{n_{\alpha}}H^{0}(X,({\cal L}_{D})^{\otimes i})
$$
pour $\alpha =1,2$ et on peut consid\'{e}rer leur produit
$$
f_{H}={\cal M}_{H}={\cal M}_{1}\times_{k}{\cal M}_{2}\rightarrow {\Bbb
A}_{1}\times_{k}{\Bbb A}_{2}={\Bbb A}_{H}.
$$

On a un morphisme de $f_{H}$ dans $f$
$$\diagram{
{\cal M}_{H}&\kern -1mm\smash{\mathop{\hbox to 8mm{\rightarrowfill}}
\limits^{\scriptstyle i_{{\cal M}}}}\kern -1mm&{\cal M}\cr
\llap{$\scriptstyle f_{H}$}\left\downarrow
\vbox to 4mm{}\right.\rlap{}&+&\llap{}\left\downarrow
\vbox to 4mm{}\right.\rlap{$\scriptstyle f$}\cr
{\Bbb A}_{H}&\kern -1mm\smash{\mathop{\hbox to 8mm{\rightarrowfill}}
\limits_{\scriptstyle i}}\kern -1mm&{\Bbb A}\cr}
$$
qui envoie $(({\cal E}_{1},\Phi_{1},\theta_{1}),({\cal
E}_{2},\Phi_{2},\theta_{2}))$ sur
$$
({\cal E}_{1}\oplus {\cal E}_{2},
\Phi_{1}\oplus\Phi_{2},\theta_{1}\oplus \theta_{2})
$$
et $(a_{1},a_{2})$ sur
$$
a=(a_{1,1}+a_{2,1},a_{1,2}+a_{1,1}a_{2,1}
+a_{2,2},\ldots ,a_{1,n_{1}}a_{2,n_{2}})
$$
de sorte que
$$
(u^{n_{1}}+a_{1,1}u^{n_{1}-1}+\cdots +a_{1,n_{1}})
(u^{n_{2}}+a_{2,1}u^{n_{2}-1}+\cdots +a_{2,n_{2}})
=(u^{n}+a_{1}u^{n-1}+\cdots +a_{n}).
$$

En particulier, pour tous points $a_{1}$ de ${\Bbb A}_{1}$ et $a_{2}$
de ${\Bbb A}_{2}$ \`{a} valeurs dans un $k$-sch\'{e}ma $S$ la courbe
spectrale $Y_{i(a_{1},a_{2})}\subset \Sigma_{S}^{\circ}
=S\times_{k}\Sigma^{\circ}$ est le diviseur de Cartier relatif (\`{a}
$S$) somme des diviseurs de Cartier relatifs $Y_{1,a_{1}}$ et
$Y_{2,a_{2}}$ o\`{u} $Y_{\alpha ,a_{\alpha}}\subset
\Sigma_{S}^{\circ}$ est d\'{e}finie comme $Y_{a}\subset
\Sigma_{S}^{\circ}$ apr\`{e}s avoir remplac\'{e} $n$ par $n_{\alpha}$
et $a$ par $a_{\alpha}$.

La courbe spectrale endoscopique universelle $Y_{H}\rightarrow {\Bbb
A}_{H}$ est par d\'{e}finition la courbe relative somme disjointes
$$
Y_{H}=Y_{1}\times_{k}{\Bbb A}_{2}\amalg {\Bbb A}_{1}
\times_{k}Y_{2}
$$
o\`{u} $Y_{1}\subset {\Bbb A}_{1}\times_{k}\Sigma$ et $Y_{2}\subset
{\Bbb A}_{2}\times_{k}\Sigma$ sont les courbes spectrales universelles
pour $G_{1}$ et $G_{2}$.  D'apr\`{e}s ce qui pr\'{e}c\`{e}de, c'est
la normalisation partielle
$$
Y_{H}\rightarrow {\Bbb A}_{H}\times_{{\Bbb A}}Y=Y_{1}\times_{k}
{\Bbb A}_{2}+{\Bbb A}_{1}\times_{k}Y_{2}\subset {\Bbb A}_{H}
\times_{k}\Sigma
$$
qui s\'{e}pare les composantes $Y_{1}\times_{k}{\Bbb A}_{2}$ et ${\Bbb
A}_{1}\times_{k}Y_{2}$.

Avec des notations \'{e}videntes on a aussi un ${\Bbb A}_{H}$-champ de
Picard
$$
P_{H}=(\mathop{\rm Pic}\nolimits_{Y_{H}'/{\Bbb
A}_{H}})^{\tau^{\ast}=(-)^{\otimes -1}}=P_{1}\times_{k}P_{2}
$$
qui agit sur le ${\Bbb A}_{H}$-champ ${\cal M}_{H}$ par l'action
produit de celle de $P_{1}$ sur ${\Bbb A}_{1}$ et celle de $P_{2}$ sur
${\Bbb A}_{2}$.  Le morphisme de $f_{H}$ dans $f$ est
$P$-\'{e}quivariant o\`{u} $P$ agit sur ${\Bbb A}_{H}$ \`{a} travers
l'homomorphisme
$$
{\Bbb A}_{H}\times_{{\Bbb A}}P\twoheadrightarrow P_{H}
$$
d'image inverse pour la normalisation partielle $Y_{H}\rightarrow {\Bbb
A}_{H}\times_{{\Bbb A}}Y$.

On a comme pr\'{e}c\'{e}demment des ouverts
$$
{\Bbb A}_{H}^{{\rm lisse}}={\Bbb A}_{1}^{{\rm lisse}}\times_{k}{\Bbb
A}_{2}^{{\rm lisse}}\subset {\Bbb A}_{H}^{{\rm red}}={\Bbb A}_{1}^{{\rm red}}
\times_{k}{\Bbb A}_{2}^{{\rm red}}\subset
{\Bbb A}_{H}
$$
qui sont non vides d\`{e}s que $2\sup (n_{1},n_{2})\mathop{\rm
deg}(D)\geq 2g+1$. On a
$$
{\Bbb A}_{H}^{G-{\rm red}}:=i^{-1}({\Bbb A}^{{\rm red}})\subset
{\Bbb A}_{H}^{{\rm red}}.
$$

\thm LEMME 2.7.1
\enonce
Le morphisme $i:{\Bbb A}_{H}\rightarrow {\Bbb A}$ est fini.
\endthm

\rem D\'{e}monstration
\endrem
Le morphisme $i:{\Bbb A}_{H}\rightarrow {\Bbb A}$ est ${\Bbb G}_{{\rm
m},k}$-\'{e}quivariant pour les actions qui font de $a_{\alpha
,i_{\alpha}}$ une coordonn\'{e}e homog\`{e}ne de degr\'{e}
$i_{\alpha}$ pour chaque $i_{\alpha}=1,\ldots ,n_{\alpha}$ et chaque
$\alpha$, et de m\^{e}me de $a_{i}$ une coordonn\'{e}e homog\`{e}ne de
degr\'{e} $i$ pour chaque $i=1,\ldots ,n$.  De plus, le seul point de
${\Bbb A}_{H}$ d'image (identiquement) nulle dans ${\Bbb A}$ est $0$.
Par suite, cette application est finie.
\hfill\hfill$\square$
\vskip 3mm

\thm LEMME 2.7.2
\enonce
La restriction ${\Bbb A}_{H}^{G-{\rm red}}\rightarrow {\Bbb A}^{{\rm
red}}$ de $i$ \`{a} ${\Bbb A}_{H}^{G-{\rm red}}\subset {\Bbb
A}_{H}^{{\rm red}}$ est un morphisme net.

Plus pr\'{e}cis\'{e}ment, soit $a_{H}=(a_{1},a_{2})$ un point de
${\Bbb A}_{H}={\Bbb A}_{1}\times_{k}{\Bbb A}_{2}$ tel que les courbes
spectrales $Y_{a_{1}}$ et $Y_{a_{2}}$ soient g\'{e}om\'{e}triquement
int\`{e}gres et distinctes {\rm (}dans le cas o\`{u} $n_{1}=n_{2}${\rm
)}.  Alors, l'image du morphisme $i:{\Bbb A}_{H}\rightarrow {\Bbb A}$
est lisse en l'image $a$ de $a_{H}$ et le morphisme de ${\Bbb A}_{H}$
sur son image est un isomorphisme au dessus d'un voisinage de $a$ si
$n_{1}\not=n_{2}$ et est \'{e}tale de degr\'{e} $2$ au dessus d'un
voisinage de $a$ si $n_{1}=n_{2}$.
\endthm

\rem D\'{e}monstration
\endrem
Comme la source et le but de  $i:{\Bbb A}_{H}\rightarrow {\Bbb A}$
sont lisses sur $k$ et que ce morphisme est fini, il suffit de
consid\'{e}rer l'application tangente en $a_{H}$
$$\displaylines{
\qquad\Bigl(\kappa (a_{1})\otimes_{k}\bigoplus_{i_{1}=1}^{n_{1}}
H^{0}(X,({\cal L}_{D})^{\otimes i_{1}})\Bigr)\oplus \Bigl(\kappa
(a_{2})\otimes_{k} \bigoplus_{i_{2}=1}^{n_{2}}H^{0}(X,({\cal
L}_{D})^{\otimes i_{2}})\Bigr)
\hfill\cr\hfill
\rightarrow \kappa (a_{H})\otimes_{k}\bigoplus_{i=1}^{n}
H^{0}(X,({\cal L}_{D})^{\otimes i})\qquad}
$$
qui est donn\'{e}e par
$$
(\dot{P}_{1}(u),\dot{P}_{2}(u))\mapsto
P_{1}(u)\dot{P}_{2}(u)+P_{2}(u)\dot{P}_{1}(u)
$$
o\`{u} $P_{\alpha}(u)=u^{n_{\alpha}}+a_{\alpha
,1}u^{n_{\alpha}-1}+\cdots + a_{\alpha ,n_{\alpha}}$ et
$\dot{P}_{\alpha}(u)=\dot{a}_{\alpha ,1}u^{n_{\alpha}-1}+\cdots +
\dot{a}_{\alpha,n_{\alpha}}$ pour $\alpha =1,2$.  Mais cette
application est injective car, au point g\'{e}n\'{e}rique de
$\kappa_{a_{H}}\otimes_{k}X$, les polyn\^{o}mes $P_{1}(u)$ et
$P_{2}(u)$ sont premiers entre eux.
\hfill\hfill$\square$

\subsection{2.8}{L'ouvert elliptique}

On dira qu'un $S$-point $a$ de ${\Bbb A}^{{\rm red}}$ est une
caract\'{e}ristique {\it elliptique} si, pour tout point
g\'{e}om\'{e}trique $s$ de $S$ le rev\^{e}tement double
$Y_{a(s)}'\rightarrow Y_{a(s)}$ induit un isomorphisme de l'ensemble
des composantes irr\'{e}ductibles de $Y_{a(s)}'$ sur celui de
$Y_{a(s)}$.  Il revient au m\^{e}me de dire que $\tau$ agit
trivialement sur l'ensemble ${\rm Irr}(Y_{a(s)}')$ des composantes
irr\'{e}ductibles de $Y_{a(s)}'$. La notion d'ellipticit\'{e}
utilis\'{e}e ici est une notion g\'{e}om\'{e}trique, qui implique la
notion d'ellipticit\'{e} usuelle. De plus, un \'{e}l\'{e}ment
elliptique est pour nous automatiquement r\'{e}gulier.

\thm LEMME 2.8.1
\enonce
Les caract\'{e}ristiques elliptiques forment un ouvert dense ${\Bbb
A}^{{\rm ell}}$ de ${\Bbb A}^{\rm red}$.
\endthm

\rem D\'{e}monstration
\endrem
L'ensemble des caract\'{e}ristiques elliptiques \'{e}tant
constructible, il suffit de v\'{e}rifier que la propri\'{e}t\'{e}
elliptique est pr\'{e}serv\'{e}e par g\'{e}n\'{e}risation.  Soit
$a:S\rightarrow {\Bbb A}^{{\rm red}}$ un morphisme o\`{u} $S$ est le
spectre d'un anneau de valuation discr\`{e}te complet, de point
sp\'{e}cial g\'{e}om\'{e}trique $\overline{s}$ et de point
g\'{e}n\'{e}rique g\'{e}om\'{e}trique $\overline{\eta}$.  Supposons que
l'image de $\overline{s}$ est elliptique, il s'agit de d\'{e}montrer que
l'image de $\overline{\eta}$ est aussi elliptique.

Par d\'{e}finition, $a(\overline{s})$ est elliptique si et seulement si
$\tau$ agit trivialement sur l'ensemble des composantes
irr\'{e}ductibles ${\rm Irr}(Y_{a(\overline{s})}')$.  Pour
d\'{e}montrer que $a(\overline{\eta})$ est elliptique, il suffit donc
de d\'{e}montrer qu'il existe une application
$\tau$-\'{e}quivariante surjective
$$
{\rm Irr}(Y_{a(\overline{s})}')\twoheadrightarrow
{\rm Irr}(Y_{a(\overline{\eta})}').
$$
En effet, la surjectivit\'{e} de cette application force $\tau$ \`{a}
agir trivialement sur ${\rm Irr}(Y_{a(\overline{\eta})}')$.

Soit $Y_{S}'^{\circ}$ l'ouvert maximal de lissit\'{e} de
$Y_{S}'=S\times_{{\Bbb A}}Y'$ sur $S$.  Puisque $Y_{S}'$ est une
$S$-courbe plate \`{a} fibres g\'{e}om\'{e}triquement r\'{e}duites, on a
${\rm Irr}(Y_{a(s)}')= \pi_{0}(Y_{a(s)}'^{\circ})$ et ${\rm
Irr}(Y_{a(\eta)}')= \pi_{0}(Y_{a(\eta)}'^{\circ})$.  De plus,
d'apr\`{e}s le lemme 15.5.6 de [EGA IV], l'application
$\pi_{0}(Y_{a(\eta)}'^{\circ}) \rightarrow
\pi_{0}(Y_{S}'^{\circ})$ qui \`{a} une composante connexe associe
son adh\'{e}rence plate dans $Y_{S}'^{\circ}$, est bijective.

Consid\'{e}rons maintenant l'application $\pi_{0}(Y_{a(s)}'^{\circ})
\rightarrow \pi_{0}(Y_{S}'^{\circ})$ qui associe \`{a} une
composante connexe de $Y_{a(s)}'^{\circ}$ l'unique composante
connexe de $Y_{S}'^{\circ}$ qui la contient.  Compos\'{e}e avec
l'inverse de la bijection $\pi_{0}(Y_{a(\eta)}'^{\circ})
\buildrel\sim\over\longrightarrow \pi_{0}(Y_{S}'^{\circ})$, cette
application d\'{e}finit une application $\pi_{0}(Y_{a(s)}'^{\circ})
\rightarrow \pi_{0}(Y_{a(\eta)}'^{\circ})$.

Enfin, pour construire l'application cherch\'{e}e ${\rm
Irr}(Y_{a(\overline{s})}')\rightarrow {\rm Irr}(Y_{a(\overline{\eta})}')$,
il ne reste plus qu'\`{a} remplacer $S$ par le normalis\'{e} $S'$ de
$S$ dans une extension finie $\eta'$ de $\eta$ telle que ${\rm
Gal}(\overline{\eta}/\eta')$ agisse trivialement sur ${\rm
Irr}(Y_{a(\overline{\eta})}')$.  L'application ainsi construite est
clairement compatible \`{a} l'action de $\tau$.

La surjectivit\'{e} de ${\rm Irr}(Y_{a(\overline{s})}')\rightarrow {\rm
Irr}(Y_{a(\overline{\eta})}')$ r\'{e}sulte de la propret\'{e} de
$Y_{S}'$ sur $S$.
\hfill\hfill$\square$
\vskip 3mm

L'ouvert ${\Bbb A}^{{\rm ell}}$ est non vide car il contient l'ouvert
non vide ${\Bbb A}^{{\rm lisse}}$ o\`{u} la courbe spectrale $Y$ est
lisse.  En effet, au-dessus de ce lieu, les courbes $Y$ et $Y'$ qui
sont lisses et \`{a} fibres g\'{e}om\'{e}triquement connexes, ont
toutes leurs fibres g\'{e}om\'{e}triques irr\'{e}ductibles.

Remarquons cependant qu'il peut arriver que $Y_{a}$ soit
irr\'{e}ductible sans que le rev\^{e}tement \'{e}tale $Y_{a}'$ ne le
soit.

\thm LEMME 2.8.2
\enonce
Au-dessus de l'ouvert ${\Bbb A}^{\rm ell}$, le morphisme $P\rightarrow
{\Bbb A}^{{\rm red}}$ est de type fini.  De plus, pour tout point
g\'{e}om\'{e}trique $a$ de ${\Bbb A}^{\rm ell}$, le groupe des
composantes connexes $\pi_{0}(P_{a})$ de $P_{a}$ est canoniquement
isomorphe \`{a} $({\Bbb Z}/2{\Bbb Z})^{\mathop{\rm Irr}(Y_{a})}$
o\`{u} $\mathop{\rm Irr}(Y_{a})$ est l'ensemble des composantes
irr\'{e}ductibles de $Y_{a}$.

Soit $C$ un diviseur de Cartier dans l'ouvert de lissit\'{e} de
$Y_{a}'$ et consid\'{e}rons le fibr\'{e} en droites ${\cal
O}_{Y_{a}'}(C-\tau(C))$, muni de la structure unitaire \'{e}vidente.
Le point ainsi d\'{e}fini dans $P_{a}$ est dans la composante neutre
de $P_{a}$ si et seulement si le degr\'{e} de la restriction de $C$
\`{a} chaque composante irr\'{e}ductible de $Y_{a}'$ est pair.
\endthm

\rem D\'{e}monstration
\endrem
Soient $a$ un point g\'{e}om\'{e}trique de ${\Bbb A}^{\rm ell}$, $Y_{a}$
la courbe spectrale associ\'{e}e et $\pi_{a}:Y_{a}'\rightarrow Y_{a}$
son rev\^{e}tement \'{e}tale double.  Au-dessus de $Y_{a}$, on a une suite
exacte de faisceaux en groupes
$$
1\rightarrow {\Bbb G}_{{\rm m},Y_{a}}\rightarrow (\pi_{a})_{\ast}({\Bbb
G}_{m,Y_{a}'}) {\smash{\mathop{\hbox to 8mm{\rightarrowfill}}
\limits^{\scriptstyle \beta}}} (\pi_{a})_{\ast}({\Bbb G}_{{\rm m},
Y_{a}'})^{\tau^{\ast}=(-)^{-1}}\rightarrow 1
$$
o\`{u} $\beta$ est le morphisme {\og}{anti-norme}{\fg} donn\'{e} sur
les sections locales par $\beta (\xi)=\xi\tau(\xi)^{-1}$.

Par passage \`{a} la cohomologie on en d\'{e}duit une suite exacte 
longue
$$
\xymatrix @R=1mm @C=3mm {1\ar[r] & H^{0}(Y_{a},{\cal
O}_{Y_{a}}^{\times})\ar[r] & H^{0}(Y_{a}',{\cal
O}_{Y_{a}'}^{\times})\ar[r] & H^{0}(Y_{a}',{\cal
O}_{Y_{a}'}^{\times})^{\tau^{\ast}=(-)^{-1}}\ar[r] &\cr \ar[r] & {\rm
Pic}(Y_{a})\ar[r] & {\rm Pic}(Y_{a}')\ar[r] & P_{a}\ar[r] & 1}
$$
et donc une suite exacte \`{a} droite
$$
\pi_{0}({\rm Pic}(Y_{a}))\rightarrow \pi_{0}({\rm Pic}(Y_{a}'))
\rightarrow \pi_{0}(P_{a})\rightarrow 0.
$$

Puisque $a$ est elliptique, l'application $\mathop{\rm
Irr}(Y_{a}')\rightarrow \mathop{\rm Irr}(Y_{a})$ induite par $\pi_{a}$
est bijective.  On peut donc identifier $\pi_{0}({\rm Pic}(Y_{a}))$
et $\pi_{0}({\rm Pic}(Y_{a}'))$ \`{a} ${\Bbb Z}^{\mathop{\rm
Irr}(Y_{a})}$, et la fl\`{e}che $\pi_{0}({\rm Pic}(Y_{a}))\rightarrow
\pi_{0}({\rm Pic}(Y'_{a}))$ est alors la multiplication par $2$ dans
${\Bbb Z}^{\mathop{\rm Irr}(Y_{a})}$.  On en d\'{e}duit que
$\pi_{0}(P_{a})=({\Bbb Z}/2{\Bbb Z})^{\mathop{\rm Irr}(Y_{a})}$.

Soit $C$ un diviseur de Cartier comme dans l'\'{e}nonc\'{e}.
L'homomorphisme ${\rm Pic}(Y_{a}')\rightarrow P_{a}$ qui se d\'{e}duit
de $\beta$, envoie le fibr\'{e} inversible ${\cal O}_{Y_{a}'}(C)$ sur
le fibr\'{e} inversible ${\cal O}_{Y_{a}'}(C-\tau (C))$ muni de la
structure unitaire \'{e}vidente.  La description de la fl\`{e}che
induite sur les $\pi_{0}$ montre que ${\cal O}_{Y'_{a}}(C-\tau (C))$
est dans la composante neutre $P_{a}^{0}$ si et seulement si le
degr\'{e} de $C$ sur chaque composante irr\'{e}ductible de $Y_{a}'$
est pair.
\hfill\hfill$\square$

\rem Remarque
\endrem
L'ouvert ${\Bbb A}^{{\rm ell, int}}\subset {\Bbb A}^{{\rm ell}}\subset
{\Bbb A}^{{\rm red}}$ des caract\'{e}ristiques $a$ telles que $Y_{a}$
et $Y_{a}'$ soient int\`{e}gres est l'ouvert compl\'{e}mentaire dans
${\Bbb A}^{{\rm ell}}$ de la r\'{e}union des ${\Bbb A}^{{\rm ell}}\cap
i({\Bbb A}_{{\rm U}(n_{1})\times {\rm U}(n_{2})})$ pour toutes les
partitions non triviales $n=n_{1}+n_{2}$.
\hfill\hfill$\square$

\thm LEMME 2.8.3
\enonce
Il existe une application $\pi_{0}({\cal M}_{a})\rightarrow {\Bbb
Z}/2{\Bbb Z}$ \'{e}quivariante pour l'action de $\pi_{0}(P_{a})$ sur
$\pi_{0}({\cal M}_{a})$ induite par celle de $P_{a}$ sur ${\cal
M}_{a}$ et pour l'action de $\pi_{0}(P_{a})\cong ({\Bbb Z}/2{\Bbb
Z})^{\mathop{\rm Irr}(Y_{a})}$ sur ${\Bbb Z}/2{\Bbb Z}$ via la somme.
\endthm

\rem D\'{e}monstration
\endrem
On a une application canonique de l'espace des fibr\'{e}s unitaires
dans ${\Bbb Z}/2{\Bbb Z}$ d\'{e}finie comme suit.  \`{A} $({\cal
E},\Phi :{\cal E}\buildrel\sim\over\longrightarrow \tau^{\ast}{\cal
E})$ on associe le fibr\'{e} inversible $\bigwedge^{n}{\cal E}$ sur
$X'$ muni de la structure unitaire $\bigwedge_{}^{n}\Phi$.
D'apr\`{e}s le lemme pr\'{e}c\'{e}dent l'espace de module de ces
fibr\'{e}s inversibles unitaires a deux composantes connexes.
\hfill\hfill$\square$
\vskip 3mm

L'\'{e}nonc\'{e} suivant est crucial pour notre travail.  Il est pour
l'essentiel un cas particulier du th\'{e}or\`{e}me II.4 de [Fal].  Une
partie des arguments utilis\'{e}s figure aussi dans [Est].

\thm PROPOSITION 2.8.4
\enonce
{\rm (i)} La restriction ${\cal M}^{{\rm ell}}$ de ${\cal M}$ au-dessus de
l'ouvert ${\Bbb A}^{{\rm ell}}\subset {\Bbb A}$ est un champ de
Deligne-Mumford.  

\decale{\rm (ii)} Le morphisme de champs de Deligne-Mumford $f^{{\rm
ell}}:{\cal M}^{{\rm ell}}\rightarrow {\Bbb A}^{{\rm ell}}$ induit par
le morphisme de Hitchin est propre.
\endthm

\rem D\'{e}monstration 
\endrem
Soit $\overline{k}$ une cl\^{o}ture alg\'{e}brique de $k$.  Pour
l'assertion (i) on se contentera de d\'{e}montrer que les
automorphismes de tout $\overline{k}$-point de l'ouvert ${\cal
M}^{{\rm ell}}$ est fini.

Il r\'{e}sulte de la d\'{e}monstration de la proposition 2.5.2 que
l'alg\`{e}bre de Lie du groupe des automorphismes d'un objet $({\cal
E},\Phi ,\theta )$ de ${\cal M}(\overline{k})$ de caract\'{e}ristique
$a$ est
$$
H^{0}(\overline{k}\otimes_{k}X,{\cal H}^{0}(K))\subset H^{0}
(\overline{k}\otimes_{k}\widetilde{Y}_{a}^{\,\prime}, {\cal 
O}_{\overline{k}\otimes_{k}\widetilde{Y}_{a}^{\,\prime}})^{\tau^{\ast}=-1}
$$
avec les notations de cette d\'{e}monstration.  Mais si $a$ est dans
l'ouvert ${\Bbb A}^{{\rm ell}}$, $\tau^{\ast}$ agit trivialement sur
$H^{0}(\overline{k}\otimes_{k}\widetilde{Y}_{a}^{\,\prime}, {\cal
O}_{\overline{k}\otimes_{k}\widetilde{Y}_{a}^{\,\prime}})$. Par
cons\'{e}quent cette alg\`{e}bre de Lie est nulle et le groupe des
automorphismes de $({\cal E},\Phi ,\theta )$ est fini.
\vskip 2mm

Pour l'assertion (ii) nous proc\'{e}derons en trois temps en montrant
tout d'abord que le morphisme $f^{{\rm ell}}$ est de type fini, puis
qu'il satisfait la partie {\og}{existence}{\fg} du crit\`{e}re
valuatif de propret\'{e}, et enfin la partie {\og}{unicit\'{e}}{\fg}
de ce m\^{e}me crit\`{e}re.

Soit $a\in {\Bbb A}^{{\rm ell}}(\overline{k})$.  Le morphisme d'oubli
$f^{-1}(a)={\cal M}_{a}\rightarrow \mathop{\overline{\rm
Pic}}\nolimits_{Y_{a}'/\overline{k}}, ~({\cal F},\iota )\mapsto {\cal
F}$, est repr\'{e}sentable de type fini: sa fibre en un ${\cal
O}_{Y_{a}'}$-Module sans torsion ${\cal F}$ de rangs
g\'{e}n\'{e}riques $1$ est un ferm\'{e} de $\mathop{{\cal I}{\it
som}}\nolimits_{{\cal O}_{Y_{a}'}}({\cal F},\tau^{\ast}{\cal
F}^{\vee})$.  Pour d\'{e}montrer que ${\cal M}_{a}$ est de type fini
il suffit donc de voir que ce morphisme d'oubli se factorise \`{a}
travers un ouvert de type fini de $\mathop{\overline{\rm
Pic}}\nolimits_{Y_{a}'/\overline{k}}^{\,(n-1)\mathop{\rm
deg}(D)}\subset \mathop{\overline{\rm Pic}}
\nolimits_{Y_{a}'/\overline{k}}$.  On rappelle que le degr\'{e} de
$\omega_{Y_{a}/X}$ est \'{e}gal \`{a} $2(n-1)\mathop{\rm deg}(D)$ et 
donc que pour tout $({\cal F},\iota )\in {\cal M}_{a}$, ${\cal F}$ 
est de degr\'{e} $(n-1)\mathop{\rm deg}(D)$.

De tels ouverts sont obtenus en bornant les degr\'{e}s des
restrictions de ${\cal F}$ aux composantes irr\'{e}ductibles de
$Y_{a}'$.  Plus pr\'{e}cis\'{e}ment, soit
$$
\nu :Y_{a}^{\dagger}= \coprod_{C\in\mathop{\rm Irr}(Y_{a})}C
\rightarrow Y_{a}
$$
la normalisation partielle qui consiste \`{a} s\'{e}parer les
composantes irr\'{e}ductibles $C$ de $Y_{a}$ sans les modifier.
Alors, 
$$
\nu':Y_{a}^{\dagger\prime}= \coprod_{C\in\mathop{\rm Irr}(Y_{a})}
X'\times_{X}C\rightarrow Y_{a}'
$$
est aussi la normalisation partielle qui consiste \`{a} s\'{e}parer
les composantes irr\'{e}ductibles $C'=X'\times_{X}C$ de $Y_{a}'$ sans
les modifier, puisque $a$ est elliptique.  Pour tout
$\overline{k}$-point ${\cal F}$ de $\mathop{\overline{\rm Pic}}
\nolimits_{Y_{a}'/\overline{k}}^{\,(n-1)\mathop{\rm deg}(D)}$ notons
${\cal G}$ le plus grand quotient sans torsion de $\nu'^{\ast}{\cal
F}$ et $d_{C}({\cal F})+{1\over 2}\mathop{\rm deg}(\omega_{C/X})$ le
degr\'{e} de la restriction ${\cal G}_{C}$ de ${\cal G}$ \`{a} la
composante connexe $C'$ de $Y_{a}^{\dagger\prime}$.  On a une
fl\`{e}che injective d'adjonction
$$
{\cal F}\hookrightarrow\nu_{\ast}'{\cal G}
$$
dont le conoyau est annul\'{e} par le conducteur ${\frak a}$ de
$Y_{a}^{\dagger}/Y_{a}$, c'est-\`{a}-dire l'annulateur de
$\nu_{\ast}{\cal O}_{Y_{a}^{\dagger}}/{\cal O}_{Y_{a}}$.  On a donc
$$
{\frak a}\nu_{\ast}'{\cal G}\hookrightarrow {\cal F}
\hookrightarrow\nu_{\ast}'{\cal G}.
$$

On en d\'{e}duit que, pour toute famille d'entiers
$(e_{C})_{C\in\mathop{\rm Irr} (Y_{a})}$, les ${\cal F}$ de degr\'{e}
$(n-1)\mathop{\rm deg}(D)$ tels que $d_{C}({\cal F}) \geq e_{C}$ quel
que soit $C\in \mathop{\rm Irr}(Y_{a})$ forment une famille
limit\'{e}e. 

Maintenant, si $({\cal F},\iota )$ est un $\overline{k}$-point de 
${\cal M}_{a}$, on a
$$
\tau^{\ast}((\nu_{\ast}'{\cal G})^{\vee})\subset \tau^{\ast}({\cal F}^{\vee})
\cong {\cal F}\subset \nu_{\ast}'{\cal G}.
$$
On en d\'{e}duit que $d_{C}({\cal F})\geq 0$ pour chaque $C\in
\mathop{\rm Irr}(Y_{a})$.  Le champ de Deligne-Mumford ${\cal M}_{a}$
est donc bien de type fini.  Nous laissons au lecteur le soin de
g\'{e}n\'{e}raliser cet argument pour en d\'{e}duire que $f^{{\rm
ell}}$ est de type fini.
\vskip 2mm

Consid\'{e}rons maintenant le crit\`{e}re valuatif de propret\'{e}.
Soit $S$ un trait strictement hens\'{e}lien de point ferm\'{e} $s$ et
de point g\'{e}n\'{e}rique $\eta$.  Soit $a:S\rightarrow {\Bbb
A}^{{\rm ell}}$ un $S$-point de ${\Bbb A}^{{\rm ell}}$.  On a la $S$-courbe
spectrale $Y_{S}\rightarrow S$, et son rev\^{e}tement double \'{e}tale
$Y_{S}'\rightarrow Y_{S}$.  On a aussi le rev\^{e}tement fini
$p_{S}':Y_{S}'\rightarrow Y_{S}\rightarrow X_{S}$.  Soit $({\cal
F}_{\eta},\iota_{\eta})$ un $\kappa (\eta )$-point de
${\cal M}_{\eta}$, c'est-\`{a}-dire une ${\cal
O}_{Y_{\eta}'}$-Module coh\'{e}rent sans torsion de rang
g\'{e}n\'{e}rique $1$ muni d'une structure unitaire $\iota_{\eta}$.
On veut prolonger ce point en une section $({\cal F},\iota)$ de
${\cal M}_{S}\rightarrow S$ quitte \`{a} remplacer $S$ par un
rev\^{e}tement fini ramifi\'{e}.  Soit $U_{s}$ l'ouvert de
$X_{s}=s\times_{k}X$ au-dessus duquel $Y_{s}$ est \'{e}tale et
$V_{s}'\subset Y_{s}$ l'image r\'{e}ciproque de $U_{s}$ par $p_{s}'$.
Alors, $V_{s}'$ est r\'{e}union disjointe d'ouverts $V_{s,C}'$
index\'{e}s par les composantes irr\'{e}ductibles de $Y_{s}'$ ou ce
qui revient au m\^{e}me les composantes irr\'{e}ductibles $C$ de
$Y_{s}'$ puisque $a(s)$ est elliptique.  Les r\'{e}unions
$U=X_{\eta}\cup U_{s}\subset X_{S}$ et $V_{S}'=Y_{\eta}'\cup
V_{s}'\subset Y_{S}'$ sont des ouverts denses, et chaque $V_{s,C}'$
est un diviseur de Cartier sur $V'$.

On commence par prolonger $({\cal F}_{\eta},\iota_{\eta})$ \`{a} $V'$.
Pour cela on choisit un prolongement de ${\cal F}_{\eta}$ en un ${\cal
O}_{V'}$-Module coh\'{e}rent ${\cal G}$ sans torsion de rangs
g\'{e}n\'{e}riques $1$ tel que $\iota_{\eta}$ se prolonge en un
homomorphisme $\psi :{\cal G}\rightarrow \tau^{\ast}{\cal G}^{\vee}$
n\'{e}cessairement injectif.  On a alors
$$
\tau^{\ast}{\cal F}_{V'}^{\vee}={\cal G}\Bigl(\sum_{C}m_{C}
V_{s,C}'\Bigr)
$$
Quitte \`{a} ramifier $S$ en extrayant une racine carr\'{e}e de
l'uniformisante de $S$, on peut supposer que les $m_{C}$ sont tous
pairs et alors 
$$
{\cal F}_{V'}={\cal G}\Bigl(\sum_{C}{m_{C}\over 2}V_{s,C}'\Bigr)
$$
muni de la structure unitaire $\iota_{V'}$ induite par
$\psi$ r\'{e}pond \`{a} la question.

Maintenant on prend pour $({\cal F},\iota )$ l'image directe par
l'immersion ouverte $V'\hookrightarrow Y'$ de $({\cal
F}_{V'},\iota_{V'})$.  Pour voir que ${\cal F}$ est plat sur $S$ et
fibre \`{a} fibre sans torsion de rangs g\'{e}n\'{e}riques $1$, il
suffit de remarquer que $(p_{S}')_{\ast}{\cal F}$ est l'image directe
par l'immersion ouverte $U\hookrightarrow X_{S}$ de sa restriction
\`{a} $U$, et est donc un fibr\'{e} vectoriel puisque $X$ est lisse
sur $k$ et que $X_{S}-U$ est de codimension $2$ dans $X_{S}$.  On a
donc d\'{e}montr\'{e} la partie {\og}{existence}{\fg} du crit\`{e}re
valuatif de propret\'{e}.
\vskip 2mm
 
Pour conclure il ne reste plus qu'\`{a} traiter la partie
{\og}{unicit\'{e}}{\fg} de ce crit\`{e}re valuatif.  Soient donc
$({\cal F},\iota)$ et $({\cal F}^{1},\iota^{1}))$ deux sections de
${\cal M}_{S}\rightarrow S$ et $\varphi_{\eta}:({\cal
F}_{\eta},\iota_{\eta}) \buildrel\sim\over\longrightarrow ({\cal
F}_{\eta}^{1},\iota_{\eta}^{1})$ un isomorphisme entre leurs
restrictions \`{a} $Y_{\eta}'$.  Il s'agit de prolonger
$\varphi_{\eta}$ \`{a} $Y_{S}'$ tout entier.

Comme pr\'{e}c\'{e}demment, il suffit de prolonger $\varphi$ \`{a}
l'ouvert $V'\subset Y_{S}'$.  En se localisant au point
g\'{e}n\'{e}rique de chaque composante connexe $V_{s,C}'$ de $V_{s}'$,
on est ramen\'{e} \`{a} v\'{e}rifier l'assertion suivante.  Soient
$R'/R$ une extension \'{e}tale de degr\'{e} $2$ d'anneaux de
valuations discr\`{e}tes, $K'/K$ l'extension correspondantes entre les
corps des fractions, $N$ et $N^{1}$ deux $K'$-espace vectoriel de
dimension $1$ munis de structures unitaires relativement \`{a}
l'extension quadratique $K'/K$ et $M\subset N$ et $M^{1}\subset N^{1}$
deux $R'$-r\'{e}seaux auto-duaux relativement \`{a} ces structures
unitaires, alors tout isomorphisme unitaire $\psi
:N\buildrel\sim\over\longrightarrow N^{1}$ envoie $M$ sur $M^{1}$.

On peut supposer que $M=M^{1}=R'$, de sorte que $N=N^{1}=K'$ avec des
structures unitaires donn\'{e}es par les formes hermitiennes $\alpha
x^{\ast}y$ et $\alpha^{1}x^{\ast}y$ pour $\alpha ,\alpha^{1}\in
R^{\times}$.  Alors $\psi$ est donn\'{e} par $x\rightarrow \beta' x$
avec $\beta'\in K'^{\times}$ tel que $\alpha
=\beta'\beta'^{\ast}\alpha'$ et donc tel que $\beta'\in R^{\times}$,
d'o\`{u} l'assertion, et la partie (ii) de la proposition.
\hfill\hfill$\square$

\subsection{2.9}{La {\og}{glissade}{\fg}}

L'espace de Hitchin ${\Bbb A}_H$ du groupe endoscopique $H$ est un
produit ${\Bbb A}_H={\Bbb A}_{1}\times{\Bbb A}_{2}$ o\`{u} pour tous
$\alpha\in\{1,2\}$, ${\Bbb A}_{\alpha}$ est l'espace affine associ\'{e}
au $k$-espace vectoriel
$$
\bigoplus_{i=1}^{n_\alpha}H^{0}(X,({\cal L}_{D})^{\otimes i}).
$$
Le $k$-sch\'{e}ma en groupes vectoriels $\mathop{\rm Vect}(\Sigma /X)$
d\'{e}fini par le $k$-espace vectoriel $H^{0}(X,{\cal L}_{D})$ agit
par translation sur la surface r\'{e}gl\'{e}e $\Sigma ={\Bbb P}({\cal
O}_{X}\oplus ({\cal L}_{D})^{\otimes -1})\rightarrow X$ en
pr\'{e}servant la section infinie.  L'action de $\mathop{\rm
Vect}(\Sigma /X)$ se rel\`{e}ve au fibr\'{e} en droites ${\cal
O}_{\Sigma}(n_{\alpha})$ et donc induit une action de $\mathop{\rm
Vect}(\Sigma /X)$ sur $H^{0}(\Sigma ,{\cal O}_{\Sigma}(n_{\alpha}))=
\bigoplus_{i=0}^{n_{\alpha}}H^{0}(X,({\cal L}_{D})^{\otimes i})$ qui
est donn\'{e}e par
$$
v\cdot (a_{\alpha ,1},\ldots ,a_{\alpha ,n_{\alpha}})=(b_{\alpha ,1}(v),\ldots
,b_{\alpha ,n_{\alpha}}(v))
$$
o\`{u} $b_{\alpha ,1}(v),\ldots ,b_{\alpha ,n_{\alpha}}(v)$ sont d\'{e}finis par
$$
u^{n_{\alpha}}+b_{\alpha ,1}(v)u^{n_{\alpha}-1}+\cdots +b_{\alpha
,n_{\alpha}}(v)= (u+v)^{n_{\alpha}}+a_{\alpha
,1}(u+v)^{n_{\alpha}-1}+\cdots +a_{\alpha ,n_{\alpha}}
$$
pour tout $v\in \mathop{\rm Vect}(\Sigma /X)$.

Cette action de $\mathop{\rm Vect}(\Sigma /X)$ sur ${\Bbb A}_{\alpha}$
se rel\`{e}ve par construction \`{a} la courbe spectrale universelle
$Y_{\alpha}\rightarrow {\Bbb A}_{\alpha}$ et se rel\`{e}ve donc aussi
en une action sur la fibration de Hitchin $f_{\alpha}:{\cal
M}_{n_{\alpha}}\rightarrow {\Bbb A}_{\alpha}$ de ${\rm U}(n_{\alpha})$.

Par produit direct, on a donc une action de $\mathop{\rm Vect}(\Sigma
/X)\times_{k}\mathop{\rm Vect}(\Sigma /X)$ sur la fibration de Hitchin
$f_H:{\cal M}_{H}\rightarrow {\Bbb A}_{H}$ du groupe endoscopique.

Soit ${\Bbb A}_{H,\natural}^{G-{\rm red}}$ l'ouvert de ${\Bbb
A}_{H}^{G-{\rm red}}$ dont les points g\'{e}om\'{e}triques sont les
points g\'{e}om\'{e}\-tri\-ques $(a_{1},a_{2})$ de ${\Bbb A}_{H}^{G-{\rm
red}}$ tels que les deux courbes spectrales $Y_{a_{1}}$ et $Y_{a_{2}}$
trac\'{e}es sur $\Sigma$ se coupent transversalement et de plus, tels
qu'en tout point $z$ de leur intersection, $Y_{a_{1}}$ et $Y_{a_{2}}$
soient \'{e}tales sur $\kappa (a_{1},a_{2})\otimes_{k}X$.

\thm PROPOSITION 2.9.1
\enonce
L'ouvert
$$
\{((v_{1},v_{2}),(a_{1},a_{2}))\mid (v_{1}\cdot a_{1},v_{2}\cdot
a_{2})\in {\Bbb A}_{H,\natural}^{G-{\rm red}}\}\subset \mathop{\rm
Vect}(\Sigma /X)\times_{k}\mathop{\rm Vect}(\Sigma /X)\times_{k}{\Bbb
A}_{H}^{G-{\rm red}}
$$
image r\'{e}ciproque de l'ouvert ${\Bbb A}_{H,\natural}^{G-{\rm
red}}\subset {\Bbb A}_{H}$ par le morphisme d'action $\mathop{\rm
Vect}(\Sigma /X)\times_{k}\mathop{\rm Vect}(\Sigma /X)\times_{k}{\Bbb
A}_{H}\rightarrow {\Bbb A}_{H}$ s'envoie surjectivement sur ${\Bbb
A}_{H}^{G-{\rm red}}$ par la projection canonique $\mathop{\rm
Vect}(\Sigma /X)\times_{k}\mathop{\rm Vect}(\Sigma /X)\times_{k}{\Bbb
A}_{H}^{G-{\rm red}}\rightarrow {\Bbb A}_{H}^{G-{\rm red}}$.

En particulier ${\Bbb A}_{H,\natural}^{G-{\rm red}}$ est non vide.
\endthm

\rem D\'{e}monstration
\endrem
Soit $a=(a_{1},a_{2})$ un point g\'{e}om\'{e}trique de ${\Bbb
A}_{H}^{G-{\rm red}}$.  Comme $p>n$, le diviseur de Cartier effectif
$Y_{a_{1}}+Y_{a_{2}}$ de $\kappa (a)\otimes_{k}\Sigma$ est
g\'{e}n\'{e}riquement \'{e}tale sur $\kappa (a)
\otimes_{k}X$.  Il existe donc un ouvert dense $U$ de $\kappa
(a)\otimes_{k}X$ au-dessus duquel $Y_{a_{1}}$ et $Y_{a_{2}}$
sont \'{e}tales et ne se rencontrent pas.

Soit $\widetilde{U}\rightarrow U$ un rev\^{e}tement fini \'{e}tale qui
d\'{e}ploie compl\`{e}tement les restrictions finies \'{e}tales de
$Y_{1,a_{1}}$ et $Y_{2,a_{2}}$ \`{a} $U$.  Pour $\alpha =1,2$, on a
alors des sections
$$
b_{\alpha,1},\ldots , b_{\alpha,n_{\alpha}}\in H^{0}
(\widetilde{U},{\cal L}_{D})
$$
telles que
$$
u^{n_{\alpha}}+(a_{\alpha ,1}|\widetilde{U})u^{n_{\alpha}-1}+\cdots
+(a_{\alpha ,n_{\alpha}}|\widetilde{U})=
\prod_{i_{\alpha}=1}^{n_{\alpha}}(u-b_{\alpha ,i_{\alpha}}).
$$

Soit $v\in \kappa (a)\otimes_{k}H^{0}(X,{\cal L}_{D})$ qui
ne s'annule en aucun point de $\kappa (a)\otimes_{k}(X-U)$.
Alors, $v$ induit une base, not\'{e}e encore $v$, de la fibre ${\cal
L}_{D,K}$ de ${\cal L}_{D}$ au point g\'{e}n\'{e}rique $\mathop{\rm
Spec}(K)$ de $\kappa (a)\otimes_{k}X$ et les quotients
$b_{\alpha ,i_{\alpha}}/v$ sont des \'{e}l\'{e}ments bien d\'{e}finis
du corps des fonctions $\widetilde{K}=\kappa (a)
(\widetilde{U})$ de la courbe $\widetilde{U}$, \'{e}l\'{e}ments dont
on peut prendre les diff\'{e}rentielles
$$
{\rm d}(b_{\alpha ,i_{\alpha}}/v)\in\Omega_{\widetilde{K}/\kappa
(a)}^{1}.
$$

Choisissons arbitrairement $f\in K$ dont la diff\'{e}rentielle ${\rm
d}f\in \Omega_{K/\kappa (a)}^{1}$ est non nulle et tel que
$$
v'=fv\in \kappa (a)\otimes_{k}H^{0}(X,{\cal L}_{D}) \subset
{\cal L}_{D,K}.
$$
Il existe de tels $f=v'/v$ d'apr\`{e}s le th\'{e}or\`{e}me de
Riemann-Roch puisque $\mathop{\rm deg}({\cal L}_{D})\geq 2g+2$.

Montrons alors qu'il existe $c$ et $c'$ dans $\kappa (a)$
tels que $((cv+c'v')\cdot a_{1},a_{2})$ soit un point
g\'{e}om\'{e}trique de ${\Bbb A}_{H,\natural}^{G-{\rm red}}$, ce qui
terminera la d\'{e}monstration de la proposition.

Il existe $c'\in\kappa (a)$ tel que les expressions
$$
{\rm d}((b_{1 ,i_{1}}-c'v')/v)-{\rm d}(b_{2 ,i_{2}}/v)\in
\Omega_{\widetilde{K}/\kappa (a)}^{1}
$$
soient toutes non nulles et que $v$ ne s'annule en aucune des images
dans $\kappa (a)\otimes_{k}X$ des points d'intersection de
$Y_{c'v'\cdot a_{1}}$ et $Y_{a_{2}}$.  Pour un tel $c'$ notons
$U_{c'}$ un ouvert dense de $U$ tel que les fonctions rationnelles
$(b_{1,i_{1}}-c'v')/v$ et $b_{2,i_{2}}$ soient toutes
r\'{e}guli\`{e}res sur l'image inverse $\widetilde{U}_{c'}$ de
$U_{c'}$ dans $\widetilde{U}$ et que les diff\'{e}rences
$$
{\rm d}((b_{1 ,i_{1}}-c'v')/v)-{\rm d}(b_{2 ,i_{2}}/v)\in
\Omega_{\widetilde{U}_{c'}/\kappa (a)}^{1}
$$
ne s'annulent en aucun point de $\widetilde{U}_{c'}$.  Pour tout $c\in
\kappa (a)$ les courbes $Y_{(cv+c'v')\cdot a_{1}}$
et $Y_{a_{2}}$ sont \'{e}tales au-dessus de $U_{c'}$ et s'y coupent
transversalement.

Il ne reste plus qu'\`{a} choisir $c$ de telle sorte que les courbes
$Y_{(cv+c'v')\cdot a_{1}}$ et $Y_{a_{2}}$ ne se coupent pas au-dessus
de l'ensemble fini $\kappa (a) \otimes_{k}X-U_{c'}$.  C'est
possible puisqu'au-dessus de chaque point $x$ de $\kappa
(a)\otimes_{k}X-U_{c'}$, cette condition n'\'{e}carte qu'un
nombre fini de valeurs pour $c$, voire m\^{e}me aucune si $v(x)=0$,
$Y_{a_{1}}$ et $Y_{a_{2}}$ ne se coupant pas au-dessus d'un z\'{e}ro
de $v$.
\hfill\hfill$\square$

\section{3}{Un \'{e}nonc\'{e} global}
\vskip - 3mm

\subsection{3.1}{Un point particulier de ${\Bbb A}_{H}$}

On suppose dor\'{e}navant que la courbe $X$ admet un point
$x_{\infty}$ rationnel sur $k={\Bbb F}_{q}$ au-dessus duquel le
rev\^{e}tement double \'{e}tale $X'\rightarrow X$ est
d\'{e}compos\'{e}.

Comme dans le chapitre 2, on consid\`{e}re le sch\'{e}ma en groupes
unitaires $G$ sur $X$ \`{a} $n$ variables et son groupe endoscopique
$H=G_{1}\times_{X}G_{2}$ o\`{u} $G_{1}$ et $G_{2}$ sont les
sch\'{e}mas en groupes unitaires sur $X$ en $n_{1}$ et $n_{2}$
variables. Toujours comme dans le chapitre 2, une fois fix\'{e} le
diviseur effectif $D$ de degr\'{e} $\geq g+1$ sur $X$, on a la fibration
de Hitchin ${\cal M}\rightarrow {\Bbb A}$ et sa variante endoscopique
${\cal M}_{H}\rightarrow {\Bbb A}_{H}$. On a aussi les courbes
spectrales $Y\rightarrow {\Bbb A}$ et $Y_{H}\rightarrow {\Bbb A}_{H}$.

Fixons maintenant un point $a=(a_{1},a_{2})$, rationnel sur
$k={\Bbb F}_{q}$, de l'espace de Hitchin endoscopique ${\Bbb A}_{H}$.
On a donc des courbes spectrales $Y_{a_{1}}$, $Y_{a_{2}}$ et
$Y_{a}=Y_{a_{1}}+Y_{a_{2}}$ trac\'{e}es sur la surface $\Sigma ={\Bbb
P}({\cal O}_{X}\oplus ({\cal L}_{D})^{\otimes -1})$.

Faisons les hypoth\`{e}ses sur $a$ suivantes:
\vskip 1mm

\itemitem{-} $Y_{a_{1}}$ et $Y_{a_{2}}$ sont g\'{e}om\'{e}triquement 
irr\'{e}ductibles;
\vskip 1mm

\itemitem{-} les images inverses $Y_{a_{1}}'$ et $Y_{a_{2}}'$ de
$Y_{a_{1}}$ et $Y_{a_{2}}$ dans le rev\^{e}tement double \'{e}tale
$Y_{a}'=X'\times_{X}Y_{a}$ de $Y_{a}$ sont aussi
g\'{e}om\'{e}triquement irr\'{e}ductibles; en particulier, l'image
$i(a)$ de $a$ dans ${\Bbb A}$ est dans l'ouvert ${\Bbb A}^{{\rm ell}}$
(cf.  la section (2.8));
\vskip 1mm

\itemitem{-} le morphisme $Y_{a}\rightarrow X$ est \'{e}tale au-dessus
du point $x_{\infty}\in X(k)$ et pour chaque $\alpha\in\{1,2\}$, il existe
au moins un point de $Y_{a_{\alpha}}$ rationnel sur $k$ au-dessus de
$x_{\infty}$.
\vskip 1mm

\subsection{3.2}{Actions du groupe discret et puret\'{e}}

Notons $R$ l'hens\'{e}lis\'{e} de ${\Bbb A}$ en l'image de
$a$ par le morphisme canonique $i:{\Bbb A}_{H}\rightarrow
{\Bbb A}$ et $S$ l'hens\'{e}lis\'{e} de ${\Bbb A}_{H}$ en
$a$.  Puisque $a$ est elliptique, on a $R\subset
{\Bbb A}^{{\rm ell}}$.  Il r\'{e}sulte des lemmes 2.7.1 et 2.7.2 que le
morphisme induit par $i$ de $S$ dans $R$ est une immersion ferm\'{e}e.
On identifie dans la suite $S$ \`{a} son image par cette immersion
ferm\'{e}e, de sorte que $S\subset R$.  On note $s$ le point
ferm\'{e} commun de $S$ et $R$; son corps r\'{e}siduel $\kappa (s)$
est \'{e}gal \`{a} $k={\Bbb F}_{q}$.

Notons par un indice $R$ le changement de base par le morphisme
$R\rightarrow {\Bbb A}^{{\rm red}}\subset {\Bbb A}$ et par un indice
$S$ le changement de base par le morphisme $S\rightarrow {\Bbb
A}_{H}^{G-{\rm red}}\subset {\Bbb A}_{H}$. On note par un indice $s$ la
fibre en $s$ (fibre sp\'{e}ciale) des objets au-dessus de $R$ ou $S$.

On a donc la courbe spectrale relative $Y_{R}\rightarrow R$ et le
morphisme fini et plat $p_{R}:Y_{R}\rightarrow R\times_{k}X$ de
degr\'{e} $n$.  On a de plus le rev\^{e}tement double \'{e}tale
$Y_{R}'\rightarrow Y_{R}$.  La fibre sp\'{e}ciale $Y_{s}=Y_{a}$ se
d\'{e}compose en r\'{e}union de composantes g\'{e}om\'{e}triquement
irr\'{e}ductibles $Y_{s}=Y_{a_{1}}\cup Y_{a_{2}}$.  De m\^{e}me on a
$Y_{s}'=Y_{a_{1}}'\cup Y_{a_{2}}'$ o\`{u} les $Y'_{a_{\alpha}}$ sont aussi
g\'{e}om\'{e}triquement irr\'{e}ductibles.

La pr\'{e}-image $p_{R}^{-1}(x_{\infty})$ de
$R\times_{k}\{x_{\infty}\}$ par le morphisme fini plat
$p_{R}:Y_{R}\rightarrow R\times_{k}X$ est un $R$-sch\'{e}ma fini et
plat de degr\'{e} $n$.  Sa fibre sp\'{e}ciale \'{e}tant suppos\'{e}
r\'{e}duite et $k$ \'{e}tant parfait, $p_{R}^{-1}(x_{\infty})$ est
donc fini \'{e}tale de degr\'{e} $n$ au-dessus de $R$.  De m\^{e}me la
pr\'{e}-image $p_{R}'^{-1}(x_{\infty})$ de $R\times_{k}
\{x_{\infty}\}$ par le morphisme fini plat $p_{R}':Y_{R}'\rightarrow
R\times_{k}X$ est un $R$-sch\'{e}ma fini \'{e}tale de degr\'{e} $2n$.

Pour chaque $\alpha\in\{1,2\}$, il existe des $k$-points de
$Y_{a_{\alpha}}'$ au-dessus de $x_{\infty}$ et on en choisit
arbitrairement un que l'on note $y_{\alpha}'\in Y_{a_{\alpha}}'(k)$;
le point $y_{\alpha}'$ s'\'{e}tend de fa\c{c}on unique en une section
$$
y_{\alpha ,R}':R\rightarrow p_{R}'^{-1}(R).
$$

Les sections $y_{\alpha ,R}'$ induisent un homomorphisme de $R$-champs
de Picard
$$
R\times {\Bbb Z}^{2}\rightarrow \mathop{\rm Pic}\nolimits_{Y'_{R}/R}.
$$
Plus pr\'{e}cis\'{e}ment, ce morphisme est donn\'{e} par
$$
(d_{1},d_{2})\mapsto {\cal O}_{Y_{R}'}\bigl(d_{1}[y_{1 ,R}']+
d_{2}[y_{2,R}']\bigr).
$$
En le composant avec l'homomorphisme {\og}{anti-norme}{\fg} (voir la
preuve du lemme 2.8.2)
$$
\mathop{\rm Pic}\nolimits_{Y_{R}'/R}\rightarrow P_{R},~\xi\mapsto
\xi\tau (\xi )^{-1},
$$
on obtient un homomorphisme
$$
\rho: R\times{\Bbb Z}^{2}\rightarrow P_{R}
$$
qui est donn\'{e} concr\`{e}tement par
$$
(d_{1},d_{2})\mapsto {\cal O}_{Y_{R}'}(d_{1}[y_{1 ,R}']-d_{1}[\tau (y_{1
,R}')]+d_{2}[y_{2 ,R}']-d_{2}[\tau (y_{2 ,R}')]),
$$
le fibr\'{e} inversible ${\cal O}_{Y'_{R}}(d_{1}[y_{1 ,R}']-d_{1}[\tau
(y_{1 ,R}')]+d_{2}[y_{2 ,R}']-d_{2}[\tau (y_{2 ,R}')])$ \'{e}tant muni
de sa structure unitaire \'{e}vidente.

Soit $b$ un point g\'{e}om\'{e}trique de $R$ et $\rho_{b} :{\Bbb
Z}^{2}\rightarrow P_{b}$ la fibre de $\rho$ en $b$.  Consid\'{e}rons
l'homomorphisme compos\'{e}
$$
\pi_{0}(\rho_{b}):{\Bbb Z}^{2}\rightarrow P_{b}
\rightarrow\pi_{0}(P_{b})
$$
dont le but est $\pi_{0}(P_{b})=({\Bbb Z}/2{\Bbb Z})^{J_{b}}$ o\`{u}
$J_{b}:={\rm Irr}(Y_{b})$, d'apr\`{e}s le lemme 2.8.2.
Pour $\alpha\in\{1,2\}$ notons $j_{b}(\alpha)\in J_{b}$ l'indice de
l'unique composante irr\'{e}ductible $Y_{b,j_{b}(\alpha )}'$ de
$Y_{b}'$ qui contient le point $y_{\alpha ,R}'(b)$.  Alors
$\pi_{0}(\rho_{b})$ est donn\'{e} par
$$
(d_{1},d_{2})\mapsto \Bigl(\sum_{\alpha\in\{1,2\},~j_{b}(\alpha )=j}
\overline{d}_{\alpha }\Bigr)_{j\in J_{b}}
$$
o\`{u} $\overline{d}_{\alpha }$ est la classe de $d_{\alpha }$ modulo
$2$, de nouveau d'apr\`{e}s le lemme 2.8.2.  En particulier
$\pi_{0}(\rho_{b}):{\Bbb Z}^{2}\rightarrow \pi_{0}(P_{b})$ se
factorise \`{a} travers ${\Bbb Z}^{2}\twoheadrightarrow ({\Bbb
Z}/2{\Bbb Z})^{2}$.

\thm PROPOSITION 3.2.1
\enonce
Pour $\alpha =1,2$, consid\'{e}rons le caract\`{e}re
$\kappa_{\alpha}:{\Bbb Z}^{2}\rightarrow \{\pm 1\}$ d\'{e}fini par
$$
\kappa_{\alpha}(d_{1},d_{2})=(-1)^{d_{\alpha}}.
$$
Soit $b$ un point g\'{e}om\'{e}trique de l'hens\'{e}lis\'{e} $R$ de
$a$ dans ${\Bbb A}$.  Pour que $\kappa_{\alpha}$ se factorise \`{a}
travers l'homomorphisme $\pi_{0}(\rho_{b}):{\Bbb Z}^{2}\rightarrow
\pi_{0}(P_{b})$ d\'{e}fini plus haut, il est n\'{e}cessaire que $b\in
S\subset R$.
\endthm

\rem D\'{e}monstration
\endrem
Le caract\`{e}re $\kappa_{\alpha}$ se factorise \`{a} travers
$\pi_{0}(\rho_{b})$ si et seulement si l'application
$j_{b}:\{1,2\}\rightarrow J_{b}$ est injective.  Faisons donc cette
hypoth\`{e}se.

Pour d\'{e}montrer que $b\in S$ il suffit d'apr\`{e}s 2.5.3 de
d\'{e}montrer que, pour $\alpha\in\{1,2\}$ le rev\^{e}tement fini et
plat $p_{b,j_{b}(\alpha )}: Y_{b,j_{b}(\alpha )} \rightarrow
b\times_{k}X$ est de degr\'{e} $n_{\alpha}$, c'est-\`{a}-dire que le
nombre de points de la fibre $p_{b,j_{b}(\alpha )}^{-1}(b,x_{\infty})$
est \'{e}gal \`{a} $n_{\alpha}$.

Soit $\overline{R}$ l'hens\'{e}lis\'{e} strict de $R$ relatif \`{a} la
cl\^{o}ture alg\'{e}brique de $\kappa (s)=k$ dans $\kappa (b)$,
$\overline{s}$ son point ferm\'{e} et $\overline{b}$ le rel\`{e}vement
naturel de $b$ \`{a} $\overline{R}$.  On a $n$ sections distinctes
$\overline{R}\rightarrow Y_{\overline{R}}$ au-dessus du point
$x_{\infty}$.  Notons cet ensemble de sections
$p_{\overline{R}}^{-1}(x_{\infty})$.  On a donc une application
canonique
$$
\mu_{s}:p_{\overline{R}}^{-1}(x_{\infty})\rightarrow \{1,2\}
$$
qui associe \`{a} une section $y\in p_{\overline{R}}^{-1}(x_{\infty})$
l'unique indice $\alpha\in \{1,2\}$ tel que $y(\overline{s})$ soit un point
g\'{e}om\'{e}trique de $Y_{s,\alpha}$.  De m\^{e}me, on a l'application
$$
\mu_{b}:p_{\overline{R}}^{-1}(x_{\infty})\rightarrow J_{b}
$$
qui associe \`{a} une section $y\in p_{\overline{R}}^{-1}(x_{\infty})$
l'unique indice $j\in J_{b}$ tel que $y(\overline{b})$ soit un point
g\'{e}om\'{e}trique de $Y_{b,j}$.

Il suffit de d\'{e}montrer que $\mu_{b}=j_{b}\circ\mu_{s}$ car alors
$$
|p_{b,j_{b}(\alpha )}^{-1}(b,x_{\infty} )|=|\mu_{b}^{-1}(j_{b}(\alpha ))|
=|\mu_{s}^{-1}(\alpha )|=n_{\alpha}
$$
vu que $j_{b}$ est injective.

Montrons donc que $\mu_{b}=j_{b}\circ\mu_{s}$.  Par d\'{e}finition de
$j_{b}$ il revient au m\^{e}me de d\'{e}montrer que, pour tous $y,y'\in
p_{\overline{R}}^{-1}(x_{\infty})$ tels que $\mu_{s}(y)=\mu_{s} (y')$
on a $\mu_{b}(y)=\mu_{b}(y')$.

Cette assertion r\'{e}sulte du corollaire 15.6.7 de [EGA IV]
appliqu\'{e} \`{a} notre situation.  En effet, soit $Y_{R}^{\circ}$
l'ouvert de lissit\'{e} de la courbe $Y_{R}$ sur $R$.  Comme
$Y_{R}\rightarrow R$ est plat et \`{a} fibres g\'{e}om\'{e}triquement
r\'{e}duites, $Y_{R}-Y_{R}^{\circ}$ est fini sur $R$.  Les sections
$y\in p_{\overline{R}}^{-1}(x_{\infty})$ arrivent dans $Y_{R}^{\circ}$
puisque $Y_{R}\rightarrow R\times_{k}X$ est \'{e}tale au-dessus de
$R\times_{k}x_{\infty}$.  D'apr\`{e}s loc.  cit, pour chaque section $y\in
p_{\overline{R}}^{-1}(x_{\infty})$, il existe un {\it ouvert} $U_{y}$ de
$Y_{\overline{R}}^\circ$, contenant l'image de la section $y$ et tel
qu'en tout point g\'{e}om\'{e}trique $b$ de $S$, $U_{y,b}=U_{y}\cap
Y^\circ_{b}$ soit la composante connexe de $Y^\circ_{b}$ contenant le
point $y(b)$.  Si maintenant $y,y'\in p_{\overline{R}}^{-1}(x_{\infty})$
sont tels que $\mu_{s}(y)=\mu_{s}(y')$, alors $U_{y}\cap U_{y'}$ est
non vide puisque cette intersection est d\'{e}j\`{a} non vide dans la
fibre sp\'{e}ciale.  De plus, le morphisme $U_{y}\cap
U_{y'}\rightarrow \overline{R}$ \'{e}tant un morphisme lisse, a
fortiori universellement ouvert, et que l'image de ce morphisme
contient le point ferm\'{e} de $\overline{R}$, il est donc surjectif.
On a donc $U_{y,b}\cap U_{y',b}\not=\emptyset$ ce qui implique
$U_{y,b}= U_{y',b}$, c'est-\`{a}-dire $\mu_{b}(y)=\mu_{b}(y')$.
\hfill\hfill$\square$
\vskip 3mm

Suivant Deligne (cf. [Del]), nous dirons qu'un complexe de faisceaux
$\ell$-adiques sur l'hens\'{e}lis\'{e} en un point ferm\'{e} d'un
sch\'{e}ma de type fini sur ${\Bbb F}_{q}$ est potentiellement pur de
poids $w\in {\Bbb Z}$ s'il provient d'un complexe de faisceaux
$\ell$-adiques pur de poids $w$ sur un voisinage \'{e}tale de ce point
ferm\'{e}.  Un tel complexe $K$ potentiellement pur de poids $0$ est
automatiquement semi-simple, c'est \`{a} dire isomorphe \`{a} la somme
de ses faisceaux de cohomologie perverse d\'{e}cal\'{e}s
$$
K\cong \bigoplus_{n}{}^{{\rm p}}{\cal H}^{n}K[-n]
$$
o\`{u} chaque ${}^{{\rm p}}{\cal H}^{n}$ est pure de poids $n$  (cf.
la section (5.4) de [B-B-D]).

On appelle {\it endoscopiques} les deux caract\`{e}res
$\kappa_{1},\kappa_{2}:({\Bbb Z}/2{\Bbb Z})^{2} \rightarrow \{\pm 1\}$
d\'{e}j\`{a} apparus dans la proposition pr\'{e}c\'{e}dente et
d\'{e}finis par
$$
\kappa_{\alpha}(d_{1},d_{2})\mapsto (-1)^{d_{\alpha}}.
$$

\thm COROLLAIRE 3.2.2
\enonce
Soit $f_{R}:{\cal M}_{R}\rightarrow R$ le changement de base du
morphisme de Hitchin $f:{\cal M}\rightarrow {\Bbb A}$ par le morphisme
$R\rightarrow {\Bbb A}$ d'hens\'{e}lisation de ${\Bbb A}$ en
$a$.  Alors, l'action de ${\Bbb Z}^{2}$ sur chaque ${}^{{\rm
p}}{\cal H}^{n}(f_{R,\ast}{\Bbb Q}_{\ell})$ qui est induite par
l'homomorphisme ${\Bbb Z}^{2}\rightarrow P_{R}$ et par l'action de
$P_{R}$ sur ${\cal M}_{R}$, se factorise \`{a} travers le quotient
fini ${\Bbb Z}^{2}\twoheadrightarrow ({\Bbb Z}/2{\Bbb Z})^{2}$.

De plus, pour chaque entier $n$, dans la d\'{e}composition
$$
{}^{{\rm p}}{\cal H}^{n}(f_{R,\ast}{\Bbb Q}_{\ell})=
\bigoplus_{\kappa}{}^{{\rm p}}{\cal H}^{n}(f_{R,\ast}{\Bbb
Q}_{\ell})_{\kappa}
$$
suivant les caract\`{e}res $\kappa$ de $({\Bbb Z}/2{\Bbb Z})^{2}$,
tous les facteurs directs sont potentiellement purs de poids $n$ et,
pour $\kappa$ {\rm endoscopique}, le facteur ${}^{{\rm p}}{\cal
H}^{n}(f_{R,\ast}{\Bbb Q}_{\ell})_{\kappa}$ est \`{a} support dans le
ferm\'{e} $S$ de $R$.  En particulier, pour chaque caract\`{e}re
{\rm endoscopique} $\kappa$, la restriction \`{a} $S$ de ${}^{{\rm p}}{\cal
H}^{n}(f_{R,\ast}{\Bbb Q}_{\ell})_{\kappa}$ est potentiellement pure
de poids $n$.
\endthm

\rem D\'{e}monstration
\endrem
D'apr\`{e}s la proposition pr\'{e}c\'{e}dente, l'image de
$(d_{1},d_{2})\in {\Bbb Z}^{2}$ avec les $d_{\alpha }$ tous pairs, est
une section de $P_{S}$ sur $S$ dont la restriction \`{a} chaque fibre
g\'{e}om\'{e}trique de $P_{S}\rightarrow S$ est dans la composante
neutre de cette fibre.  D'apr\`{e}s le lemme d'homotopie ci-dessous,
$(d_{1},d_{2})$ agit donc trivialement sur chaque ${}^{{\rm p}}{\cal
H}^{n}(f_{R,\ast}{\Bbb Q}_{\ell})$.  En d'autres termes, ${\Bbb
Z}^{2}$ agit sur chaque ${}^{{\rm p}}{\cal H}^{n}(f_{R,\ast}{\Bbb
Q}_{\ell})$ \`{a} travers son quotient $({\Bbb Z}/2{\Bbb Z})^{2}$.

Comme l'ouvert ${\cal M}^{{\rm red}}\subset {\cal M}$ est lisse sur
${\Bbb F}_{q}$ (cf.  la proposition 2.5.2), le complexe ${\Bbb
Q}_{\ell ,{\cal M}^{{\rm red}}}[0]$ est pur de poids $0$.  Comme la
restriction $f^{{\rm ell}}:{\cal M}^{{\rm ell}}\rightarrow {\Bbb
A}^{{\rm ell}}$ de $f$ au-dessus de ${\Bbb A}^{{\rm ell}}$ est un
morphisme propre entre champs de Deligne-Mumford (cf.  la proposition
2.8.4), on en d\'{e}duit que $f_{\ast}^{{\rm ell}}{\Bbb Q}_{\ell}$ est
pur de poids $0$ sur ${\Bbb A}^{{\rm ell}}$ et donc que
$f_{R,\ast}{\Bbb Q}_{\ell}$ est potentiellement pur de poids $0$.

Pour d\'{e}montrer la derni\`{e}re assertion, il suffit de
v\'{e}rifier que la restriction de ${}^{{\rm p}}{\cal
H}^{n}(f_{R,\ast}{\Bbb Q}_{\ell})_{\kappa}$ \`{a} l'ouvert $R-S$ est
nulle pour $\kappa$ endoscopique.  Soit $b$ un point
g\'{e}om\'{e}trique de $R-S$.  D'apr\`{e}s la proposition
pr\'{e}c\'{e}dente, il existe $(d_{1},d_{2})\in{\Bbb Z}^{2}$ tel que
$\kappa (d_{1},d_{2})=-1$ et que $\rho_{b}(d_{1},d_{2})$ soit dans la
composante neutre de $P_{b}$.  Ceci implique que
$\rho_{b}(d_{1},d_{2})$ est dans la composante neutre de $P_{b'}$ pour
tout point g\'{e}om\'{e}trique $b'$ dans un voisinage \'{e}tale de
$b$.  En invoquant de nouveau le lemme d'homotopie, on voit que
$(d_{1},d_{2})$ agit trivialement sur la restriction de ${}^{{\rm
p}}{\cal H}^{n}(f_{R,\ast}{\Bbb Q}_{\ell})$ \`{a} un voisinage
\'{e}tale de $b$, de sorte que la restriction de la partie
$\kappa$-isotypique ${}^{{\rm p}}{\cal H}^{n}(f_{R,\ast}{\Bbb
Q}_{\ell})_{\kappa}$ \`{a} ce voisinage \'{e}tale est nulle.  Ceci
\'{e}tant vrai pour tous les points g\'{e}om\'{e}triques $b$ de $R-S$,
on conclut que la restriction de ${}^{{\rm p}}{\cal
H}^{n}(f_{R,\ast}{\Bbb Q}_{\ell})_{\kappa}$ \`{a} $R-S$ est nulle.
\hfill\hfill$\square$

\thm LEMME 3.2.3 (Lemme d'homotopie)
\enonce
Soit $f:X\rightarrow S$ un $S$-sch\'{e}ma et $\pi:G\rightarrow S$ un
$S$-sch\'{e}ma en groupes lisse \`{a} fibres g\'{e}om\'{e}triquement
connexes agissant sur $X$.  Alors le groupe des sections globales
$G(S)$ agit trivialement sur chaque faisceau de cohomologie perverse
${}^{\rm p}{\cal H}^{n}(f_{\ast}{\Bbb Q}_{\ell})$.
\endthm

\rem D\'{e}monstration
\endrem
Consid\'{e}rons le diagramme
$$\xymatrix{
G\times_{S} X \ar[dr]_{{\rm pr}_{G}} \ar[r]^{\alpha} & G\times_{S} X
\ar[d]_{{\rm pr}_{G}} \ar[r]^-{{\rm pr}_{X}}
& X \ar[d]^{f}\\
& G \ar[r]_{\pi} & S}
$$
o\`{u} la fl\`{e}che $\alpha$ est envoie $(g,x)$ sur $(g,gx)$.  Dans
ce diagramme, le triangle est commutatif et le carr\'{e} est
cart\'{e}sien.

Par le th\'{e}or\`{e}me de changement de base par un morphisme lisse,
on a $\pi^{\ast}({}^{\rm p}{\cal H}^{n}(f_{\ast}{\Bbb Q}_{\ell}))={}^{\rm
p}{\cal H}^{n}({\rm pr}_{G,\ast}{\Bbb Q }_{\ell})$. Le morphisme $\alpha$
induit un endomorphisme $[\alpha]$ de ${}^{\rm p}{\cal H}^{n}({\rm
pr}_{G,\ast}{\Bbb Q}_{\ell})$.

Puisque $\pi$ est un morphisme lisse \`{a} fibres
g\'{e}om\'{e}triquement connexe, $\pi^{\ast}$ convenablement
d\'{e}cal\'{e}, est un foncteur pleinement fid\`{e}le de la
cat\'{e}gorie des faisceaux pervers sur $S$ dans celle des faisceaux
pervers sur $G$ (cf.  la proposition 4.2.5 de [B-B-D]).  Il existe
donc un unique endomorphisme $\beta$ de ${}^{\rm p}{\cal
H}^{n}(f_{\ast}{\Bbb Q}_{\ell})$ tel que $\pi^{\ast}(\beta)=[\alpha]$.
Par fonctorialit\'{e}, pour toute section globale $g:S\rightarrow G$,
on a $g^{\ast}([\alpha])=\beta$.  En prenant la section neutre $g=1$,
on obtient que $\beta$ est l'identit\'{e}.  Pour toute autre section
$g:S\rightarrow G$, $g^{\ast}([\alpha])$ agit donc aussi comme
l'identit\'{e} sur ${}^{\rm p}{\cal H}^{n}({\rm pr}_{G,\ast}{\Bbb
Q}_{\ell})$.  C'est ce qu'on voulait d\'{e}montrer.
\hfill\hfill$\square$
\vskip 3mm

\thm COROLLAIRE 3.2.4
\enonce
Pour chaque caract\`{e}re {\rm endoscopique} $\kappa$, le faisceau
de cohomologie perverse ${}^{\rm p}{\cal H}^{n}(f_{S,\ast}{\Bbb
Q}_{\ell})_{\kappa}$ est pur de poids $n$ quel que soit l'entier $n$.
\endthm

\rem D\'{e}monstration
\endrem
D'apr\`{e}s le lemme ci-dessous, $f_{S,\ast}{\Bbb Q}_{\ell}$ se
d\'{e}compose en somme directe d'une partie $\kappa$ et d'une partie
hors $\kappa$.  De plus, la restriction \`{a} $S$ et les foncteurs de
cohomologie perverse commutent \`{a} cette d\'{e}composition.  Le
th\'{e}or\`{e}me r\'{e}sulte donc du th\'{e}or\`{e}me de changement de
base propre et du corollaire pr\'{e}c\'{e}dent.
\hfill\hfill$\square$

\rem Remarque
\endrem
Les faisceaux de cohomologie perverse ${}^{\rm p}{\cal
H}^{n}(g_{S,\ast}{\Bbb Q}_{\ell})$ et tous ses facteurs directs
${}^{\rm p}{\cal H}^{n}(g_{S,\ast}{\Bbb Q}_{\ell})_{\kappa}$ sont eux
aussi potentiellement purs de poids $n$ puisque ${\cal N}$ est lisse
sur $k$ et $g:{\cal N}\rightarrow {\cal A}_{H}$ est propre.
\hfill\hfill$\square$

\thm LEMME 3.2.5
\enonce
Soient $E$ un corps, ${\cal A}$ une cat\'{e}gorie ab\'{e}lienne
$E$-lin\'{e}aire, $K$ un objet de $D^{{\rm b}}({\cal A})$ et $\Gamma$
un groupe ab\'{e}lien op\'{e}rant de fa\c{c}on $E$-lin\'{e}aire sur
$K$.  On suppose que pour chaque entier $n$ l'objet de cohomologie
$H^{n}(K)$ admette dans ${\cal A}$ une d\'{e}composition
$\Gamma$-\'{e}quivariante
$$
H^{n}(K)=\bigoplus_{\chi}H^{n}(K)_{\chi}
$$
o\`{u} $\chi$ parcourt les caract\`{e}res de $\Gamma$ \`{a} valeurs
dans $E^{\times}$, o\`{u} pour chaque $\chi$ et chaque
$\gamma\in\Gamma$, $\gamma -\chi (\gamma)$ op\`{e}re de mani\`{e}re
nilpotente sur $H^{n}(K)_{\chi}$, et o\`{u} $H^{n}(K)_{\chi}=(0)$ pour
tous les $\chi$ sauf un nombre fini.

Alors, il existe une unique d\'{e}composition $\Gamma$-\'{e}quivariante
$$
K=\bigoplus_{\chi}K_{\chi}
$$
dans $D^{{\rm b}}({\cal A})$ o\`{u} $\chi$ parcourt les caract\`{e}res
de $\Gamma$ \`{a} valeurs dans $E^{\times}$, o\`{u} pour chaque $\chi$
et chaque $\gamma\in\Gamma$, $\gamma -\chi (\gamma)$ op\`{e}re de
mani\`{e}re nilpotente sur $K_{\chi}$, et o\`{u} $K_{\chi}=(0)$ pour
tous les $\chi$ sauf un nombre fini.  De plus, on a
$H^{n}(K_{\chi})=H^{n}(K)_{\chi}$ quels que soient l'entier $n$ et le
caract\`{e}re $\chi$.
\endthm

\rem D\'{e}monstration 
\endrem
Commen\c{c}ons par l'unicit\'{e}.  Soient $K'$ et $K''$ deux objets de
$D^{{\rm b}}({\cal A})$ munis d'actions de $\Gamma$ et soient
$\chi'$ et $\chi''$ deux caract\`{e}res distincts de $\Gamma$ \`{a}
valeurs dans $E^{\times}$ tels que, quel que soit $\gamma\in\Gamma$,
$\gamma -\chi'(\gamma )$ et $\gamma -\chi''(\gamma )$ op\`{e}rent
de mani\`{e}re nilpotente sur $K'$ et $K''$ respectivement. Il s'agit de 
v\'{e}rifier que tout morphisme $\Gamma$-\'{e}quivariant $f:K'\rightarrow 
K''$ est n\'{e}cessairement nul. Choisissons $\gamma\in\Gamma$ tel que 
$\chi'(\gamma )\not=\chi''(\gamma )$ et des entiers $n'$ et $n''$ tels que 
$(\gamma -\chi'(\gamma ))^{n'}$ et $(\gamma -\chi''(\gamma ))^{n''}$ 
annulent $K'$ et $K''$ respectivement. D'apr\`{e}s le 
th\'{e}or\`{e}me de Bezout, il existe des polyn\^{o}mes $P'(T)$ et 
$P''(T)$ dans $E[T]$ tels que $P'(T)(T-\chi'(\gamma 
))^{n'}+P''(T)(T-\chi''(\gamma ))^{n''}=1$. On a alors
$$
f=P''(\gamma )(\gamma -\chi''(\gamma ))^{n''}f+fP'(\gamma )(\gamma 
-\chi'(\gamma ))^{n'}=0.
$$

Passons \`{a} l'existence.  Si $K$ est concentr\'{e} en un seul
degr\'{e}, il n'y a rien \`{a} faire.  Sinon, soit $[a,b]\subset {\Bbb
Z}$ le plus petit intervalle tel que $K\in \mathop{\rm
ob}D^{[a,b]}({\cal A})$.  On raisonne par r\'{e}currence sur l'entier
$b-a\geq 1$.  On a le d\'{e}vissage
$$
\xymatrix{K'=\tau_{\leq b-1}K\ar[r]^-{u} & K\ar[r]^-{v} 
&K''=H^{b}(K)[-b]\ar[r]^-{\partial} & K'[1]}
$$
Par hypoth\`{e}se de r\'{e}currence on a les d\'{e}compositions voulues 
de $K'$ et $K''$. Consid\'{e}rons la fl\`{e}che de degr\'{e} $1$
$$
\partial :\bigoplus_{\chi}K_{\chi}''\rightarrow\bigoplus_{\chi}K_{\chi}'[1].
$$
D'apr\`{e}s notre argument pour l'unicit\'{e}, on a n\'{e}cessairement
$\partial =\bigoplus_{\chi}\partial_{\chi}$ pour des fl\`{e}ches
$\partial_{\chi}: K_{\chi}''\rightarrow K_{\chi}'[1]$.  On a donc une
d\'{e}composition de $K$ en $\bigoplus_{\chi}K_{\chi}$ o\`{u}
les $K_{\chi}[1]$ sont les c\^{o}nes des fl\`{e}ches $\partial{\chi}$.  Pour
n'importe quels $\gamma$ et $\chi$, on sait par hypoth\`{e}se de
r\'{e}currence qu'il existe des entiers $n'$ et $n''$ tels que
$(\gamma -\chi (\gamma ))^{n'}$ et $(\gamma -\chi (\gamma ))^{n''}$
annulent $K_{\chi}'$ et $K_{\chi}''$ respectivement.  Alors $(\gamma
-\chi (\gamma ))^{n'+n''}$ annule $K_{\chi}$ puisque $(\gamma -\chi
(\gamma ))^{n'}$ se factorise en $\xymatrix{K_{\chi}\ar[r] &
K_{\chi}'\ar[r]^-{u} & K_{\chi}}$, que $(\gamma -\chi (\gamma
))^{n''}$ se factorise en $\xymatrix{K_{\chi}\ar[r]^-{v} &
K_{\chi}''\ar[r] & K_{\chi}}$ et que $vu=0$.
\hfill\hfill$\square$

\subsection{3.3}{Actions du tore}

Au-dessus $S$ on a deux rev\^{e}tements
$$
p_{1}:Y_{1,S}\rightarrow S\times_{k}X\hbox{ et }
p_{2}:Y_{2,S}\rightarrow S\times_{k}X
$$
de degr\'{e} respectivement $n_{1}$ et $n_{2}$ qui, relativement \`{a}
$S$, sont des familles rev\^{e}tements de $X$.  On a aussi un
morphisme
$$
Y_{1,S}\amalg Y_{2,S}\rightarrow Y_{S}
$$
qui est une normalisation partielle en famille de la courbe relative
$Y_{S}\rightarrow S$.

Pour tout point g\'{e}om\'{e}trique $b$ de $S$, les courbes $Y_{1,b}$
et $Y_{2,b}$ sont trac\'{e}es sur la m\^{e}me surface r\'{e}gl\'{e}e
$\Sigma ={\Bbb P}({\cal O}_{X}\oplus ({\cal L}_{D})^{\otimes -1})$ au-dessus
de $X$.  Leur r\'{e}union est $Y_{b}$.  Leur intersection est un
sch\'{e}ma fini $Z_{b}$ de longueur
$$
r=2n_{1}n_{2}\mathop{\rm deg}(D)
$$
ind\'{e}pendante de $b$.  Les $Z_{b}$ s'organisent en une famille $Z$
finie et plate de degr\'{e} $r$ au-dessus de $S$.  Le sch\'{e}ma $Z$
se plonge naturellement dans chacune des courbes relatives $Y_{1,S}$
et $Y_{S}$ et on dispose en particulier d'un morphisme $q:
Z\rightarrow S\times_{k}X$ qui est fini, mais n'est pas plat.

Si on note par $(-)'$ le changement de base par le rev\^{e}tement
double \'{e}tale $X'\rightarrow X$ on a alors deux rev\^{e}tements
\'{e}tales doubles
$$
Y_{1,S}'\rightarrow Y_{1,S}\hbox{ et }Y_{2,S}'\rightarrow
Y_{2,S}
$$
et les rev\^{e}tement finis compos\'{e}s
$$
p_{1}':Y_{1,S}'\rightarrow S\times_{k}X\hbox{ et }
p_{2}':Y_{2,S}'\rightarrow S\times_{k}X.
$$
On a encore un morphisme
$$
Y_{1,S}'\amalg Y_{2,S}'\rightarrow Y_{S}'
$$
qui est une normalisation partielle en famille de la courbe relative
$Y_{S}'\rightarrow S$.

Le sch\'{e}ma
$$
Z'=Y_{1,S}'\cap Y_{2,S}'
$$
est un rev\^{e}tement double \'{e}tale de $Z$.  Il est donc fini et
plat sur $S$ de degr\'{e} $2r$ et l'involution $\tau$ agit sur $Z'$
sans point fixe.  On note $q': Z'\rightarrow S\times_{k}X$ le
morphisme compos\'{e} du rev\^{e}tement \'{e}tale double
$Z'\rightarrow Z$ et de $q$.

Consid\'{e}rons les faisceaux en groupes commutatifs
$$
J=\bigl(p_{\ast}'{\Bbb G}_{{\rm
m},Y_{S}'}\bigr)^{\tau^{\ast}=(-)^{-1}},
~ J_{\alpha}=\bigl(p_{\alpha ,\ast}'{\Bbb G}_{{\rm m},
Y_{\alpha ,S}'} \bigr)^{\tau^{\ast}=-1}\hbox{ et }
K=\bigl(q_{\ast}'{\Bbb G}_{{\rm m},Z'}\bigr)^{\tau^{\ast}=(-)^{-1}}
$$
pour la topologie {\rm fppf} sur $S\times_{k}X$.

On a une suite exacte de faisceaux
$$
0\rightarrow J\rightarrow J_{1}\times J_{2}\rightarrow K\rightarrow 0
$$
qui induit une fl\`{e}che co-bord
$$
\mathop{\rm pr}\nolimits_{S,\ast}K\rightarrow {\cal H}^{1}(
\mathop{\rm pr}\nolimits_{S,\ast}J)
$$
o\`{u} $\mathop{\rm pr}\nolimits_{S}:S\times_{k}X\rightarrow S$ est la
projection canonique.  Notons que par d\'{e}finition, $P_{S}$ est le
champ de Picard associ\'{e} au complexe $\tau_{\leq 1}(\mathop{\rm
pr}\nolimits_{S,\ast}J)$ de faisceaux en groupes ab\'{e}liens sur $S$.
En termes concrets, la donn\'{e}e d'une section globale de $K$ est
\'{e}quivalente \`{a} une donn\'{e}e de recollement des Modules
inversibles triviaux sur $Y_{1,S}'$ et $Y_{2,S}'$ avec
structures unitaires triviales le long de $Z'$, et la fl\`{e}che de
co-bord ci-dessus envoie cette derni\`{e}re donn\'{e}e sur le Module
inversible recoll\'{e}.

Consid\'{e}rons la fibre sp\'{e}ciale du morphisme fini et plat
$h=\mathop{\rm pr}\nolimits_{S}\circ q:Z\rightarrow S$.  C'est un
$k$-sch\'{e}ma artinien dont le r\'{e}duit $h_{s,{\rm red}}:Z_{s,{\rm
red}}\rightarrow s=\mathop{\rm Spec}(k)$ est le spectre d'un produit
fini d'extensions finies s\'{e}parables des $k$.  Comme $k$ est
parfait, on a une r\'{e}traction canonique $Z_{s}\rightarrow Z_{s,{\rm
red}}$ (cf.  le corollaire (19.6.2) du chapitre 0 de [EGA IV]) et
comme $S$ est hens\'{e}lien on a une factorisation canonique
$$
\xymatrix{h:Z\ar[r] &\overline{Z}\ar[r]^{\overline{h}} & S}
$$
de $h$ o\`{u} le morphisme $Z\rightarrow \overline{Z}$ est totalement
ramifi\'{e} au sens o\`{u} il induit un isomorphisme de $Z_{s,{\rm
red}}$ sur $\overline{Z}_{s}$ et o\`{u} $\overline{h}$ est
\'{e}tale (corollaire (18.5.12) de [EGA IV]). De m\^{e}me, on a une
factorisation
$$
\xymatrix{h'=\mathop{\rm pr}\nolimits_{S}\circ q':Z'\ar[r] &
\overline{Z}{}'\ar[r]^{\overline{h}{}'} & S}
$$
de $h'$ o\`{u} le morphisme $Z'\rightarrow \overline{Z}{}'$ est totalement
ramifi\'{e} au sens o\`{u} il induit un isomorphisme de $Z_{s,{\rm
red}}'$ sur $\overline{Z}{}_{s}'$ et o\`{u} $\overline{h}{}'$ est
\'{e}tale, factorisation qui s'ins\`{e}re dans un diagramme
commutatif \`{a} carr\'{e} cart\'{e}sien
$$
\xymatrix{Z'\ar[d]\ar[r] & \overline{Z}{}'\ar[d]\ar[dr]^{\overline{h}{}'} & \cr
Z\ar[r] & \overline{Z}\ar[r]_{\overline{h}} & S}
$$
o\`{u} les fl\`{e}ches verticales sont des rev\^{e}tements doubles
\'{e}tales.

On a donc un sous-$S$-tore
$$
\widetilde{T}=\overline{h}{}_{\ast}'{\Bbb G}_{{\rm
m},\overline{Z}{}'}\subset h_{\ast}'{\Bbb G}_{{\rm m},Z'}
$$
du $S$-sch\'{e}ma en groupes commutatifs lisse $h_{\ast}'{\Bbb
G}_{{\rm m},Z'}$, dont la fibre sp\'{e}ciale
$$
\widetilde{T}_{s}\subset h_{s,\ast}'{\Bbb G}_{{\rm m}Z_{s}'}
$$
est le sous-tore maximal et s'envoie isomorphiquement sur le tore
quotient maximal de $h_{s,\ast}'{\Bbb G}_{{\rm m}Z_{s}'}$ (comparer
avec le th\'{e}or\`{e}me 5.8 de l'expos\'{e} 4 de [SGA~3]).

L'involution $\tau^{\ast}$ de $h_{\ast}'{\Bbb G}_{{\rm m},Z'}$
pr\'{e}serve le sous-$S$-tore $\widetilde{T}$, et
$$
T=(\widetilde{T})^{\tau^{\ast}=(-)^{-1}}
$$
est le sous-$S$-tore de $\mathop{\rm
pr}\nolimits_{S,\ast}K=(h_{\ast}'{\Bbb G}_{{\rm
m},Z'})^{\tau^{\ast}=(-)^{-1}}$ qui rel\`{e}ve canoniquement le tore
maximal
$$
T_{s}\cong\bigl((h_{s,{\rm red}}')_{\ast}{\Bbb G}_{{\rm
m},Z_{s,{\rm red}}'}\bigr)^{\tau^{\ast}=(-)^{-1}}
$$
de $(\mathop{\rm pr}\nolimits_{S,\ast}K)_{s}$.

On a donc construit un $S$-tore $T$ et un $S$-homomorphisme
$T\rightarrow P_{S}$.

\thm PROPOSITION 3.3.1
\enonce
L'image de l'immersion ferm\'{e}e $i:{\cal N}_{S}\hookrightarrow {\cal
M}_{S}$ est pr\'{e}cis\'{e}ment le lieu des points fixes de $T$ agissant
sur ${\cal M}_{S}$ \`{a} travers le morphisme $T\rightarrow
P_{S}$.
\endthm

\rem D\'{e}monstration
\endrem
Comme $T$ agit sur ${\cal N}_{S}$ \`{a} travers la fl\`{e}che
compos\'{e}e \'{e}vidente
$$
T\subset\mathop{\rm pr}\nolimits_{S,\ast}K\rightarrow {\cal H}^{1}(
\mathop{\rm pr}\nolimits_{S,\ast}J)\rightarrow {\cal H}^{1}(
\mathop{\rm pr}\nolimits_{S,\ast}J_{1})\times_{S}{\cal H}^{1}(
\mathop{\rm pr}\nolimits_{S,\ast}J_{2})
$$
on a \'{e}videmment l'inclusion ${\cal N}_{S}\subset {\cal
M}_{S}^{T}$.  Il reste donc \`{a} montrer que tout point
g\'{e}om\'{e}trique de ${\cal M}_{S}$ qui est fixe par $T$ est
dans ${\cal N}_{S}$.

La donn\'{e}e d'un point g\'{e}om\'{e}trique $m$ de ${\cal M}_{S}$
est \'{e}quivalente aux donn\'{e}es suivantes:
\vskip 1mm

\itemitem{(a)} un point g\'{e}om\'{e}trique $b$ dans $S$,
\vskip 1mm

\itemitem{(b)} un Module coh\'{e}rent ${\cal F}$ sans torsion de rang
$1$ sur le rev\^{e}tement double \'{e}tale $Y_{b}'$ de la courbe
spectrale $Y_{b}$, ce Module \'{e}tant muni d'un isomorphisme $\iota
:{\cal F} \buildrel\sim\over\longrightarrow \tau^{\ast}{\cal
F}^{\vee}$ tel que ${}^{{\rm t}}\!\iota =\tau^{\ast}\iota$.
\vskip 1mm

De plus, il revient au m\^{e}me de se donner le couple $({\cal
F},\iota)$ ou de se donner:
\vskip 1mm

\itemitem{(1)} les restrictions ${\cal F}_{1}^{\circ}$ et ${\cal
F}_{2}^{\circ}$ de ${\cal F}$ aux composantes $Y_{1,b}'-Z_{b}'$ et
$Y_{2,b}-Z_{b}'$ de $Y_{b}'-Z_{b}'$, ces restrictions \'{e}tant munies
d'isomorphismes $\iota_{\alpha}: {\cal F}_{\alpha}^{\circ}
\buildrel\sim\over\longrightarrow \tau^{\ast}({\cal
F}_{\alpha}^{\circ})^{\vee}$ avec toujours ${}^{{\rm
t}}\!\iota_{\alpha} =\tau^{\ast}\iota_{\alpha}$ pour $\alpha =1,2$,
\vskip 1mm

\itemitem{(2)} pour chaque point (ferm\'{e}) $z$ de $Z_{b}$, les
restrictions ${\cal V}_{z'}$ et ${\cal V}_{z''}$ de ${\cal F}$ aux
compl\'{e}t\'{e}s formels de $Y_{b}'$ en les deux points $z'$ et $z''$
de $Z_{b}'$ au-dessus de $z$, ces restrictions \'{e}tant munies d'un
isomorphisme ${\cal V}_{z''}\cong {\cal V}_{z'}^{\vee}$,
\vskip 1mm

\itemitem{(3)} les donn\'{e}es de recollement \'{e}videntes.
\vskip 2mm

Pour chaque point $z$ de $Z_{b}$ choisissons arbitrairement un des
deux points de $Z_{b}'$ au-dessus de $z$, point que l'on note $z'$;
notons $A_{z'}$ l'anneau local compl\'{e}t\'{e} de $Y_{b}'$ en ce
point et $\mathop{\rm Frac}(A_{z'})$ l'anneau total des fractions de
$A_{z'}$. On a
$$
A_{z'}\subset A_{1,z'}\times A_{2,z'}\subset
\mathop{\rm Frac}(A_{1,z'})\oplus
\mathop{\rm Frac}(A_{2,z'})=\mathop{\rm Frac}(A_{z'})
$$
o\`{u}, pour $\alpha =1,2$, $A_{\alpha ,z'}$ est l'anneau local
compl\'{e}t\'{e} de $Y_{\alpha ,b}'$ en $z'$ et $\mathop{\rm
Frac}(A_{\alpha ,z'})$ est l'anneau total des fractions de $A_{\alpha
,z'}$.  Notons $V_{\alpha}$ la fibre de ${\cal F}_{\alpha}^{\circ}$ au
point $\mathop{\rm Spec}(\mathop{\rm Frac}(A_{\alpha ,z'}))$ de
$Y_{\alpha ,b}'-Z_{b}'$.

Les donn\'{e}es $(2)$ sont encore \'{e}quivalentes \`{a} la
donn\'{e}e, pour chaque point $z$ de $Z_{b}$ d'un $A_{z'}$-r\'{e}seau
${\cal V}_{z'}\subset V_{1}\oplus V_{2}$.

Si l'on fixe le point g\'{e}om\'{e}trique $b$ de $S$ on peut donc
\'{e}crire tout point g\'{e}om\'{e}trique $m$ de ${\cal M}_{b}$ sous
la forme
$$
m=\bigl(({\cal F}_{\alpha}^{\circ},\iota_{\alpha})_{\alpha =1,2},
({\cal V}_{z'})_{z\in Z_{b}},\hbox{\rm donn\'{e}es de recollement}\bigr).
$$

On a une description analogue pour les points g\'{e}om\'{e}trique de
${\cal N}_{S}\subset {\cal M}_{S}$, la seule diff\'{e}rence
\'{e}tant que l'on exige en plus que les r\'{e}seaux ${\cal V}_{z'}$ soient
d\'{e}compos\'{e}s au sens o\`{u}
$$
{\cal V}_{z'}={\cal V}_{z',1}\oplus {\cal V}_{z',2}\subset
V_{1}\oplus V_{2}
$$
o\`{u} ${\cal V}_{z',\alpha}$ est un $A_{\alpha ,z'}$-r\'{e}seau dans
$V_{\alpha}$, $\alpha =1,2$.
\vskip 2mm

Dans la suite, pour all\'{e}ger les notations, on privil\'{e}gie la
composante $Y_{1,b}$ de $Y_{b}$, ou ce qui revient au m\^{e}me on
identifie le quotient de $\kappa (b)^{\times}\times \kappa
(b)^{\times}$ par le $\kappa (b)^{\times}$ diagonal \`{a} $\kappa
(b)^{\times}$ via la premi\`{e}re projection.  Alors, le groupe des
$\kappa (b)$-points de la fibre $T_{b}$ de $T$ en $b$ admet la
description suivante:
$$
T_{b}(\kappa (b))=(\kappa (b)^{\times})^{Z_{b}}=\prod_{z\in Z_{b}}
\bigl(\kappa (b)^{\times}\times \kappa
(b)^{\times}\bigr)^{\tau^{\ast}=(-)^{-1}}\subset \prod_{z\in
Z_{b}}\bigl(A_{z'}^{\times}\times A_{z'}^{\times}\bigr)^{\tau^{\ast}=(-)^{-1}}
$$
o\`{u} $\tau^{\ast}$ \'{e}change les deux copies de $\kappa
(b)^{\times}$ et les deux copies de $A_{z'}$ (le premier facteur
$\kappa (b)^{\times}$ ou $A_{z'}$ correspond \`{a} $z'$ et le second
\`{a} $\tau (z')$).

De plus l'action de $t\in T_{b}(\kappa (b))$ sur les points
g\'{e}om\'{e}triques $m$ de ${\cal M}_{b}$ est donn\'{e}e par
l'action, pour chaque $z\in Z_{b}$, du facteur $t_{z}\in\kappa
(b)^{\times}$ correspondant sur le r\'{e}seau ${\cal V}_{z'}$ par
l'homoth\'{e}tie de rapport $(t_{z},1)\in A_{1,z'}^{\times}\times
A_{2,z'}^{\times}$
$$
t_{z}\cdot {\cal V}_{z}=(t_{z},1){\cal V}_{z}\subset V_{1}\oplus V_{2}.
$$
Par suite, $m$ est fixe sous l'action $T_{b}(\kappa (b))$ si et
seulement si, pour chaque $z\in Z_{b}$, ${\cal V}_{z'}$ est
d\'{e}compos\'{e}, d'o\`{u} la proposition.
\hfill\hfill$\square$

\subsection{3.4}{Un syst\`{e}me local de rang $1$ sur $S$}

Rappelons qu'au-dessus du sch\'{e}ma hens\'{e}lien $S$, nous avons construit
un sch\'{e}ma fini et plat $h:Z\rightarrow S$ de degr\'{e} $r$ et un
rev\^{e}tement double \'{e}tale $\pi_{Z}:Z'\rightarrow Z$. 

Notons $L_{Z'/Z}$ le syst\`{e}me local en ${\Bbb Z}$-modules libres de
rang $1$ et d'ordre $2$ qui est le conoyau de la fl\`{e}che
d'adjonction
$$
{\Bbb Z}_{Z}\rightarrow \pi_{Z,\ast}{\Bbb Z}_{Z'}.
$$
Comme $S$ et donc aussi $Z$ sont hens\'{e}liens, $L_{Z'/Z}$ est 
g\'{e}om\'{e}triquement constant.

Soit $S_{\natural}$ l'ouvert de $S$ dont les points
g\'{e}om\'{e}triques $b$ ont les deux propri\'{e}t\'{e}s suivantes:
\vskip 1mm

\item{-} les deux courbes $Y_{1,b}$ et $Y_{2,b}$, qui sont trac\'{e}es
sur le m\^{e}me surface r\'{e}gl\'{e}e $\kappa (b)\otimes_{k}\Sigma$,
se coupent transversalement,
\vskip 1mm

\item{-} les rev\^{e}tements $Y_{1,b}\rightarrow \kappa
(b)\otimes_{k}X$ et $Y_{2,b}\rightarrow \kappa (b)\otimes_{k}X$ sont
\'{e}tales en tout point d'intersection $z\in Z_{b}=Y_{1,b}\cap
Y_{2,b}$.
\vskip 2mm

L'ouvert $S_{\natural}$ est dense dans $S$ puisque l'ouvert ${\Bbb
A}_{H,\natural}^{G-{\rm red}}\subset {\Bbb A}_{H}$ est non vide (cf.
(2.9)).  On note par $(-)_{\natural}$ le changement de base par
l'immersion ouverte $S_{\natural}\hookrightarrow S$.  Les morphismes
$h_{\natural}:Z_{\natural}\rightarrow S_{\natural}$ et
$h_{\natural}':Z_{\natural}'\rightarrow S_{\natural}$ sont finis
\'{e}tales de degr\'{e} $r$ et $2r$ respectivement.  
\vskip 2mm

Le syst\`{e}me local en ${\Bbb Z}$-modules libres de rang $1$
$$
L_{Z_{\natural}'/Z_{\natural}/S_{\natural}}=\bigl(\bigwedge^{r}
h_{\natural ,\ast}L_{Z_{\natural}'/Z_{\natural}}\bigr)\otimes
\bigl(\bigwedge^{r}h_{\natural ,\ast}{\Bbb 
Z}_{Z_{\natural}}\bigr)^{\otimes -1}.
$$
sur $S_{\natural}$ est g\'{e}om\'{e}triquement constant.  Il se
prolonge donc trivialement en un syst\`{e}me local en ${\Bbb
Z}$-modules libres de rang $1$ sur $S$ tout entier.  Notons
$L_{Z'/Z/S}$ ce prolongement qui est bien s\^{u}r lui aussi
g\'{e}om\'{e}triquement constant.

\thm LEMME 3.4.1
\enonce
L'unique valeur propre de Frobenius agissant sur la fibre sp\'{e}ciale
de $L_{Z'/Z/S}$ est \'{e}gale \`{a} $(-1)^{m_{Z'/Z/S}}$ o\`{u}
$$
m_{Z'/Z/S}=\sum_{z}\mathop{\rm dim}\nolimits_{\kappa (z)}({\cal
O}_{Z_{s},z}),
$$
la somme portant sur les point ferm\'{e}s $z$ de $Z_{s}$ qui sont 
inertes dans $Z_{s}'$.
\endthm

\rem D\'{e}monstration 
\endrem
Comme $S$ est hens\'{e}lien, $Z$ est somme disjointe des ses 
localis\'{e}s en les points ferm\'{e}s de $Z_{s}$. On peut donc 
supposer qu'il n'y a qu'un seul point ferm\'{e} $z$ dans $Z_{s}$. 

Si ce point est d\'{e}compos\'{e} dans $Z_{s}'$, le rev\^{e}tement 
\'{e}tale $Z'$ de $Z$ est trivial et alors $L_{Z'/Z}={\Bbb Z}_{Z}$ 
et $L_{Z'/Z/S}={\Bbb Z}_{S}$, d'o\`{u} l'assertion dans ce cas.

Si au contraire $z$ est inerte dans $Z_{s}'$, notons $K$, $K_{Z}$ et
$K_{Z'}$ les corps des fonctions de $S$, $Z$ et $Z'$, et notons $k_{Z}$
et $k_{Z'}$ les corps r\'{e}siduels de $Z$ et $Z'$.  Fixons une
cl\^{o}ture s\'{e}parable $\overline{K}$ de $K$ et notons
$\overline{k}$ le corps r\'{e}siduel du normalis\'{e} de $S$ dans
$\overline{K}$.  On a un \'{e}pimorphisme de groupes pro-finis
$\mathop{\rm Gal}(\overline{K}/K)\twoheadrightarrow \mathop{\rm
Gal}(\overline{k}/k)$ et un carr\'{e} cart\'{e}sien de $\mathop{\rm
Gal}(\overline{K}/K)$-ensembles
$$
\xymatrix{\mathop{\rm Hom}\nolimits_{K}(K_{Z'},\overline{K})
\ar[d]_{}^{}\ar[r]_{}^{} &
\mathop{\rm Hom}\nolimits_{k}(k_{Z'},\overline{k})\ar[d]_{}^{}\cr
\mathop{\rm Hom}\nolimits_{K}(K_{Z},\overline{K})\ar[r]_{}^{} &
\mathop{\rm Hom}\nolimits_{k}(k_{Z},\overline{k})}
$$
o\`{u} toutes les fl\`{e}ches sont surjectives et o\`{u} les fibres des 
fl\`{e}ches verticales ont toutes $2$ \'{e}l\'{e}ments et celles des 
fl\`{e}ches horizontales ont toutes $m=m_{Z'/Z/S}$ \'{e}l\'{e}ments.  

Soit $\varphi\in\mathop{\rm Gal}(\overline{K}/K)$ est un
rel\`{e}vement arbitraire de l'\'{e}l\'{e}ment de Frobenius de
$\mathop{\rm Frob}\nolimits_{k}\in\mathop{\rm Gal}(\overline{k}/k)$.
Il s'agit de v\'{e}rifier que $\varphi$ agit par multiplication par
$(-1)^{m}$ sur
$$
\bigl(\bigwedge^{2r}{\Bbb Z}^{\mathop{\rm Hom}\nolimits_{K}(K_{Z'},
\overline{K})}\bigr)\otimes\bigl(\bigwedge^{r}{\Bbb Z}^{\mathop{\rm Hom}
\nolimits_{K}(K_{Z},\overline{K})}\bigr)^{\otimes -2}
$$
ou ce qui revient au m\^{e}me, que le d\'{e}terminant de l'action de
$\varphi$ sur ${\Bbb Z}^{\mathop{\rm Hom}\nolimits_{K}(K_{Z'},
\overline{K})}$ est $(-1)^{m}$.

Une fois fix\'{e} un point base $\iota_{0}'\in \mathop{\rm 
Hom}\nolimits_{k}(k_{Z'},\overline{k})$, on a 
$$
\mathop{\rm Hom}\nolimits_{k}(k_{Z'},\overline{k})
=\{\iota_{0}',\iota_{1}',\ldots ,\iota_{2r-1}'\}
$$
et
$$
\mathop{\rm Hom}\nolimits_{k}(k_{Z},\overline{k})
=\{\iota_{0}=\iota_{r},\iota_{1}=\iota_{r+1},\ldots 
,\iota_{r-1}=\iota_{2r-1}\}
$$
o\`{u} $\iota_{n}'=\mathop{\rm Frob}\nolimits_{k}^{n}\circ\iota_{0}'$
et $\iota_{n}=\iota_{n}'|k_{Z}$ pour $n=1,\ldots ,2r-1$. Par suite on a
$$
\mathop{\rm Hom}\nolimits_{K}(K_{Z},\overline{K})=A_{0}\amalg 
A_{1}\amalg\cdots\amalg A_{r-1}
$$
o\`{u} $A_{n}$ est la fibre de $\mathop{\rm
Hom}\nolimits_{K}(K_{Z},\overline{K}) \rightarrow \mathop{\rm
Hom}\nolimits_{k}(k_{Z},\overline{k})$ en $\iota_{n}$, et $\varphi$ induit des isomorphismes
$$
A_{0}\buildrel\sim\over\longrightarrow
A_{1}\buildrel\sim\over\longrightarrow\cdots
\buildrel\sim\over\longrightarrow
A_{r-1}\buildrel\sim\over\longrightarrow A_{0}.
$$
Si on note $\Phi$ le compos\'{e} de ces isomorphismes, on peut donc
identifier ${\Bbb Z}^{\mathop{\rm Hom}\nolimits_{K}(K_{Z},
\overline{K})}$ muni de l'action de $\varphi$ \`{a}
$({\Bbb Z}^{A_{0}})^{\{0,1,\ldots ,r-1\}}$ muni de l'automorphisme
$$
(x_{0},x_{1},\ldots 
,x_{r-1})\mapsto (\Phi (x_{r-1}),x_{0},\ldots ,x_{r-2}).
$$
De m\^{e}me, on peut identifier ${\Bbb Z}^{\mathop{\rm
Hom}\nolimits_{K}(K_{Z'}, \overline{K})}$ muni de l'action de
$\varphi$ \`{a} $({\Bbb Z}^{A_{0}})^{\{0,1,\ldots ,2r-1\}}$ muni de
l'automorphisme
$$
(x_{0},x_{1},\ldots 
,x_{2r-1})\mapsto (\Phi^{2}(x_{2r-1}),x_{0},\ldots ,x_{2r-2}).
$$

Le (signe du) d\'{e}terminant de ce dernier automorphisme ne
d\'{e}pend pas de $\Phi$, de sorte qu'on est ramen\'{e} au cas o\`{u}
$\Phi$ est l'identit\'{e}.  Dans ce cas, le d\'{e}terminant est
\'{e}gal \`{a} la signature de la permutation circulaire de
$\{0,1,\ldots ,2r-1\}$ \`{a} la puissance le nombre d'\'{e}l\'{e}ments
$A_{0}$, c'est-\`{a}-dire \`{a} $((-1)^{(2r-1)})^{m}=(-1)^{m}$.
\hfill\hfill$\square$

\subsection{3.5}{Le tore ${\cal T}$ sur l'ouvert o\`{u} $Y_{1}$ et
$Y_{2}$ se coupent transversalement}

Le $S_{\natural}$-sch\'{e}ma en groupes
$$
\widetilde{{\cal T}}=h_{\natural, \ast}'{\Bbb G}_{{\rm m},Z_{\natural}'}
$$
et le noyau
$$
{\cal T}=(h_{\natural, \ast}'{\Bbb G}_{{\rm
m},Z_{\natural}'})^{\tau^{\ast}=(-)^{-1}}
$$
sont des $S_{\natural}$-tores de rang $2r$ et $r$ respectivement qui
contiennent (en g\'{e}n\'{e}ral strictement) les
$S_{\natural}$-tores $\widetilde{T}_{\natural}$ et $T_{\natural}$.

Si $Y_{S_{\natural}}'\rightarrow S_{\natural}$ est la restriction
\`{a} $S_{\natural}$ de la courbe spectrale universelle, on a un
homomorphisme canonique de $S_{\natural}$-sch\'{e}mas en groupes
$$
\widetilde{{\cal T}}\rightarrow \mathop{\rm
Pic}\nolimits_{Y_{S_{\natural}}'/S_{\natural}}
$$
qui peut se construire comme un homomorphisme de co-bord d'une suite
exacte longue de cohomologie et qui se d\'{e}crit concr\`{e}tement de
la fa\c{c}on suivante: pour tout point $b$ de $S_{\natural}$ on
construit un fibr\'{e} en droites sur $Y_{b}'$ en recollant les
fibr\'{e}s triviaux ${\cal O}_{Y_{1,b}'}$ et ${\cal
O}_{Y_{2,b}'}$ \`{a} l'aide d'une fonction de recollement dans
$\widetilde{{\cal T}}_{b}=H^{0}(Z_{b}',{\cal O}_{Z_{b}'}^{\times})$.
L'homomorphisme pr\'{e}c\'{e}dent induit un homomorphisme de
$S_{\natural}$-sch\'{e}mas en groupes
$$
{\cal T}\rightarrow P_{S_{\natural}}=\bigl(\mathop{\rm
Pic}\nolimits_{Y_{S_{\natural}}'/S_{\natural}}\bigr)^{\tau^{\ast}=
(-)^{\otimes -1}}
$$
et donc une action de ${\cal T}$ sur la restriction ${\cal
M}_{S_{\natural}}$ de ${\cal M}$ \`{a} $S_{\natural}$. Le lieu des
points fixes pour cette derni\`{e}re action est le ferm\'{e} ${\cal
N}_{S_{\natural}}\subset {\cal M}_{S_{\natural}}$ induit par le
ferm\'{e} ${\cal N}$ de ${\cal M}$.
\vskip 2mm

Sur $Z_{\natural}'\times_{S_{\natural}}{\cal N}_{S_{\natural}}$ on a
un fibr\'{e} en droites muni d'une structure unitaire
$$
{\cal E}_{12}=(\widetilde{{\cal E}}_{12},\widetilde{{\cal E}}_{12}
\buildrel\sim\over\longrightarrow \tau^{\ast}(\widetilde{{\cal
E}}_{12})^{\otimes -1})
$$
dont la fibre en $(z',({\cal F}_{1},\iota_{1}),({\cal F}_{2},
\iota_{2}))$ est $({\cal F}_{1,z'}\otimes {\cal F}_{2,z'}^{\otimes
-1},\iota_{1,z'}\otimes\iota_{2,z'}^{\otimes -1})$ o\`{u}
$$
\iota_{1,z'}\otimes\iota_{2,z'}^{\otimes -1}:{\cal F}_{1,z'}\otimes
{\cal F}_{2,z'}^{\otimes -1} \buildrel\sim\over\longrightarrow {\cal
F}_{1,\tau (z')}^{\otimes -1} \otimes {\cal F}_{2,\tau (z')}=({\cal
F}_{1,\tau (z')}\otimes {\cal F}_{2,\tau (z')}^{\otimes -1})^{\otimes
-1}.
$$
A priori $\iota_{\alpha ,z'}$ est un isomorphisme entre ${\cal
F}_{\alpha ,z'}$ et ${\cal F}_{\alpha ,\tau (z')}^{\otimes
-1}\otimes\omega_{Y_{\alpha ,b}'/\kappa (b)\otimes_{k}X',\tau (z')}$,
mais comme $b$ est dans $S_{\natural}$, $Y_{\alpha ,b}'$ \'{e}tale sur
$\kappa (b)\otimes_{k}X'$ et la droite $\omega_{Y_{\alpha ,b}'/\kappa
(b)\otimes_{k}X',\tau (z')}$ est canoniquement trivialis\'{e}e.

La donn\'{e}e du fibr\'{e} en droites $\widetilde{{\cal E}}_{12}$ sur
$Z_{\natural}'\times_{S_{\natural}}{\cal N}_{S_{\natural}}$
\'{e}quivaut \`{a} celle d'un $\widetilde{{\cal T}}$-torseur sur
${\cal N}_{S_{\natural}}$.  Par suite, la donn\'{e}e de
${\cal E}_{12}$ \'{e}quivaut \`{a} la donn\'{e}e d'un ${\cal
T}$-torseur sur ${\cal N}_{S_{\natural}}$, torseur que l'on note dans
la suite encore par ${\cal E}_{12}$ (voir la section 6 de [La-Ra] pour
un autre point de vue sur ce torseur).
\vskip 2mm

Sur l'ouvert $S_{\natural}\subset S$ nous avons aussi un morphisme
injectif naturel de syst\`{e}mes locaux en ${\Bbb Z}$-modules libres
$$
L_{Z_{\natural}'/Z_{\natural}/S_{\natural}}\hookrightarrow\mathop{\rm
Sym}\nolimits_{{\Bbb Z}_{S_{\natural}}}^{r}(X^{\ast}({\cal T}))
$$
o\`{u} $L_{Z_{\natural}'/Z_{\natural}/S_{\natural}}$ est la restriction
\`{a} l'ouvert $S_{\natural}\subset S$ du syst\`{e}me local
$L_{Z'/Z/S}$ construit dans la section (3.4).

En effet, soit $a:S_{\natural,1}\rightarrow S_{\natural}$ un
rev\^{e}tement fini \'{e}tale qui d\'{e}ploie totalement
$h_{\natural}:Z_{\natural}\rightarrow S_{\natural}$ et
$h_{\natural}':Z_{\natural}'\rightarrow S_{\natural}$, de sorte que,
apr\`{e}s changement de base par $S_{\natural ,1}\rightarrow
S_{\natural}$, on ait les rev\^{e}tements triviaux
$$
S_{\natural ,1}\times_{S_{\natural}}Z_{\natural}'=S_{\natural ,1}^{\mathop{\rm
Mor}\nolimits_{S_{\natural}}(S_{\natural ,1},Z_{\natural}')}
\rightarrow S_{\natural ,1}\times_{S_{\natural}}Z_{\natural}=S_{\natural
,1}^{\mathop{\rm Mor}\nolimits_{S_{\natural}}(S_{\natural ,1},Z_{\natural})}
\rightarrow S_{\natural ,1}
$$
et les tores d\'{e}ploy\'{e}s
$$
S_{\natural ,1}\times_{S_{\natural}}\widetilde{{\cal T}}=
{\Bbb G}_{{\rm m},S_{\natural ,1}}^{\mathop{\rm Mor}
\nolimits_{S_{\natural}}(S_{\natural ,1},Z_{\natural}')}
$$
et
$$
S_{\natural ,1}\times_{S_{\natural}}{\cal T}=\mathop{\rm Ker}
({\Bbb G}_{{\rm m},S_{\natural ,1}}^{\mathop{\rm Mor}
\nolimits_{S_{\natural}}(S_{\natural ,1},Z_{\natural}')}\rightarrow
{\Bbb G}_{{\rm m},S_{\natural ,1}}^{\mathop{\rm Mor}
\nolimits_{S_{\natural}}(S_{\natural ,1},Z_{\natural})})
$$
o\`{u} l'homomorphisme norme est donn\'{e} par
$(x_{\varphi'})_{\varphi'}\rightarrow
(\prod_{\pi\circ\varphi'=\varphi}x_{\varphi'})_{\varphi}$.

Pour chaque $\varphi'\in \mathop{\rm Mor}
\nolimits_{S_{\natural}}(S_{\natural ,1},Z_{\natural}')$ on note alors
$$
\chi_{\varphi'}: S_{\natural ,1}\times_{S_{\natural}}{\cal T}\subset
{\Bbb G}_{{\rm m},S_{\natural ,1}}^{\mathop{\rm Mor}
\nolimits_{S_{\natural}}(S_{\natural ,1},Z_{\natural}')}\rightarrow
{\Bbb G}_{{\rm m},S_{\natural ,1}}.
$$
la projection canonique sur la composante d'indice $\varphi'$ et on
voit $\chi_{\varphi'}$ comme une section globale de
$X^{\ast}(S_{\natural ,1}\times_{S_{\natural}}{\cal T})$.

Choisissons arbitrairement un type {\og}{CM}{\fg}, c'est-\`{a}-dire un
ordre $(\varphi_{-}',\varphi_{+}')$ dans la fibre de $\mathop{\rm Mor}
\nolimits_{S_{\natural}}(S_{\natural ,1},Z_{\natural}')\rightarrow
\mathop{\rm Mor}\nolimits_{S_{\natural}}(S_{\natural ,1},
Z_{\natural})$ en chaque $\varphi\in \mathop{\rm Mor}
\nolimits_{S_{\natural}}(S_{\natural ,1},Z_{\natural})$, de sorte que
l'homomorphisme norme s'\'{e}crit encore $(x_{\varphi'})_{\varphi'}
\rightarrow (x_{\varphi_{-}'}x_{\varphi_{+}'})_{\varphi}$ et que
$\chi_{\varphi_{-}'}=-\chi_{\varphi_{+}'}$.  Alors, on peut former le
produit
$$
\prod_{\varphi\in \mathop{\rm Mor}\nolimits_{S_{\natural}}(S_{\natural
,1},Z_{\natural})}\chi_{\varphi_{+}'}\in H^{0}(S_{\natural
,1},\mathop{\rm Sym}\nolimits_{{\Bbb Z}_{S_{\natural
,1}}}^{r}(X^{\ast}(S_{\natural ,1}\times_{S_{\natural}}{\cal T})).
$$

\thm LEMME 3.5.1
\enonce
Il existe un unique morphisme injectif de syst\`{e}mes locaux sur
$S_{\natural}$
$$
L_{Z_{\natural}'/Z_{\natural}/S_{\natural}}\hookrightarrow\mathop{\rm
Sym}\nolimits_{{\Bbb Z}_{S_{\natural}}}^{r}(X^{\ast}({\cal T}))
$$
dont la restriction \`{a} $S_{\natural ,1}$ est l'injection
$$
{\Bbb Z}_{S_{\natural ,1}}\hookrightarrow \mathop{\rm Sym}
\nolimits_{{\Bbb  Z}_{S_{\natural ,1}}}^{r}(X^{\ast}(S_{\natural ,1}
\times_{S_{\natural}}{\cal T}))
$$
qui envoie la section globale $1$ sur la section globale
$\prod_{\varphi\in \mathop{\rm Mor}
\nolimits_{S_{\natural}}(S_{\natural
,1},Z_{\natural})}\chi_{\varphi_{+}'}$
\endthm

\rem D\'{e}monstration
\endrem
Comme $\chi_{\varphi_{-}'}=-\chi_{\varphi_{+}'}$, le produit
$\prod_{\varphi\in \mathop{\rm
Mor}\nolimits_{S_{\natural}}(S_{\natural ,1},Z_{\natural})}
\chi_{\varphi_{+}'}$ ne d\'{e}pend pas, au signe pr\`{e}s, du type CM
choisi.  De plus, le co-cycle de descente de ce produit est exactement
le m\^{e}me que celui du faisceau constant
$X^{\ast}(S_{\natural1}\times_{S_{\natural}}{\cal T})=
a^{\ast}X^{\ast}({\cal T})$ en le syst\`{e}me local $X^{\ast}({\cal
T})$, d'o\`{u} le lemme.
\hfill\hfill$\square$
\vskip 3mm

\subsection{3.6}{Le th\'{e}or\`{e}me g\'{e}om\'{e}trique sur l'ouvert
$S_{\natural}$}

Sur l'ouvert $S_{\natural}$ de $S$ on a le triangle commutatif de
morphismes de champs alg\'{e}briques
$$
\xymatrix{
{\cal N}_{S_{\natural}}\ar@{^{(}->}[rr]
\ar[rd]_{g_{S_{\natural}}} & & {\cal M}_{S_{\natural}}
\ar[ld]^{f_{S_{\natural}}}\cr
& S_{\natural} }
$$
o\`{u} l'image de l'immersion ferm\'{e}e horizontale est le lieu des
points fixes sous l'action du sous-$S_{\natural}$-tore $T_{\natural}\subset
{\cal T}$ et donc aussi du tore ${\cal T}$ puisque ce dernier agit
trivialement sur ${\cal N}_{S_{\natural}}$.

Le ${\cal T}$-torseur ${\cal E}_{12}$ sur ${\cal N}_{S_{\natural}}$ de
la section pr\'{e}c\'{e}dente induit un ${\cal T}$-torseur $[{\cal
E}_{12}/{\cal T}]$ sur le champ alg\'{e}brique $[{\cal
N}_{S_{\natural}}/{\cal T}]$.  Ici on a pris le quotient $[{\cal
E}_{12}/{\cal T}]$ pour l'action de ${\cal T}$ sur ${\cal E}_{12}$ de
la structure de torseur, alors qu'on a pris le quotient $[{\cal
N}_{S_{\natural}}/{\cal T}]$ pour l'action triviale de ${\cal T}$.

Si $\chi$ est une section de $X^{\ast}({\cal T})$ sur un ouvert
\'{e}tale $U\rightarrow S_{\natural}$ de $S_{\natural}$, on peut
pousser le ${\cal T}_{U}$-torseur $[{\cal E}_{12}/{\cal T}]_{U}$ par
$\chi$ et prendre la premi\`{e}re classe de Chern du fibr\'{e} en
droites $\chi ([{\cal E}_{12}/{\cal T}]_{U})$ sur $[{\cal
N}_{S_{\natural}}/{\cal T}]_{U}=[{\cal N}_{U}/{\cal T}_{U}]$ ainsi
obtenu.  Bien entendu, on a not\'{e} $(-)_{U}$ le changement de base
par le morphisme $U\rightarrow S_{\natural}$.  On voit cette classe de
Chern comme un morphisme de complexes $\ell$-adiques
$$
c_{1}(\chi ([{\cal E}_{12}/{\cal T}]_{U})):{\Bbb Q}_{\ell ,[{\cal
N}_{U}/{\cal T}_{U}]}\rightarrow {\Bbb Q}_{\ell ,[{\cal N}_{U}/{\cal
T}_{U}]}[2](1)
$$
sur $[{\cal N}_{U}/{\cal T}_{U}]$.  Cette classe de Chern induit par
image directe sur $U$ un morphisme de complexes $\ell$-adiques
$$
g_{U,\ast}^{{\cal T}_{U}}{\Bbb Q}_{\ell}\rightarrow g_{U,\ast}^{{\cal
T}_{U}}{\Bbb Q}_{\ell}[2](1)
$$
sur $U$, o\`{u} $g_{U}:{\cal N}_{U}:=U\times_{S_{\natural}}{\cal
N}_{S_{\natural}}\rightarrow U$ est la projection canonique, morphisme
qui \`{a} son tour induit un morphisme de faisceaux pervers
$\ell$-adiques gradu\'{e}s
$$
\bigoplus_{n}{}^{\rm p}{\cal H}^{n}(g_{U,\ast}^{{\cal T}_{U}}{\Bbb
Q}_{\ell})\rightarrow \bigoplus_{n}{}^{\rm p}{\cal
H}^{n+2}(g_{U,\ast}^{{\cal T}_{U}}{\Bbb Q}_{\ell})(1)
$$
sur $U$.

On a d\'{e}fini ainsi un morphisme de faisceaux pervers
$\ell$-adiques gradu\'{e}s
$$
e_{12}:X^{\ast}({\cal T})\otimes_{{\Bbb Z}_{S_{\natural}}}
\bigoplus_{n}{}^{\rm p}{\cal H}^{n} (g_{S_{\natural},\ast}^{{\cal
T}}{\Bbb Q}_{\ell})\rightarrow \bigoplus_{n}{}^{\rm p}{\cal
H}^{n+2}(g_{S_{\natural},\ast}^{{\cal T}}{\Bbb Q}_{\ell})(1)
$$
sur $S_{\natural}$.  En it\'{e}rant $r$ fois le morphisme $e_{12}$ on
obtient un morphisme de faisceaux pervers $\ell$-adiques gradu\'{e}s
$$
e_{12}^{r}:\mathop{\rm Sym}\nolimits_{}^{r}(X^{\ast}({\cal T}))
\otimes_{{\Bbb Z}_{S_{\natural}}}\bigoplus_{n}{}^{\rm p}{\cal H}^{n}
(g_{S_{\natural},\ast}^{{\cal T}}{\Bbb Q}_{\ell})\rightarrow
\bigoplus_{n}{}^{\rm p}{\cal H}^{n+2r}(g_{S_{\natural},\ast}^{{\cal T}}
{\Bbb Q}_{\ell})(r).
$$
sur $S_{\natural}$, et aussi un morphisme de faisceaux pervers
$\ell$-adiques gradu\'{e}s
$$
e_{12}^{r}:\mathop{\rm Sym}\nolimits_{}^{r}(X^{\ast}({\cal T}))
\otimes_{{\Bbb Z}_{S_{\natural}}}\bigoplus_{n}{}^{\rm p}{\cal H}^{n}
(g_{S_{\natural},\ast}^{{\cal T}}{\Bbb Q}_{\ell})_{\kappa}\rightarrow
\bigoplus_{n}{}^{\rm p}{\cal H}^{n+2r}(g_{S_{\natural},\ast}^{{\cal
T}}{\Bbb Q}_{\ell})_{\kappa}(r).
$$
sur $S_{\natural}$ pour chaque caract\`{e}re $\kappa$ de $({\Bbb
Z}/2{\Bbb Z})^{2}$.

On a le sous-syst\`{e}me local de rang $1$
$$
L_{Z_{\natural}'/Z_{\natural}/S_{\natural}}
\subset \mathop{\rm Sym}\nolimits_{}^{r}(X^{\ast}({\cal T}))
$$
construit dans la section pr\'{e}c\'{e}dente.

\thm LEMME 3.6.1
\enonce
Les restrictions
$$
L_{Z_{\natural}'/Z_{\natural}/S_{\natural}}\otimes_{{\Bbb
Z}_{S_{\natural}}}\bigoplus_{n}{}^{\rm p}{\cal H}^{n}
(g_{S_{\natural},\ast}^{{\cal T}}{\Bbb Q}_{\ell})\rightarrow
\bigoplus_{n}{}^{\rm p}{\cal H}^{n+2r}(g_{S_{\natural},\ast}^{{\cal T}}
{\Bbb Q}_{\ell})(r)
$$
et
$$
L_{Z_{\natural}'/Z_{\natural}/S_{\natural}}\otimes_{{\Bbb
Z}_{S_{\natural}}}\bigoplus_{n}{}^{\rm p}{\cal H}^{n}
(g_{S_{\natural},\ast}^{{\cal T}}{\Bbb Q}_{\ell})_{\kappa} \rightarrow
\bigoplus_{n}{}^{\rm p}{\cal H}^{n+2r}(g_{S_{\natural},\ast}^{{\cal T}}
{\Bbb Q}_{\ell})_{\kappa}(r)
$$
\`{a} $L_{Z_{\natural}'/Z_{\natural}/S_{\natural}}\hookrightarrow
\mathop{\rm Sym} \nolimits^{r}(X^{\ast}({\cal T}))$ des morphismes
$e_{12}^{r}$ ci-dessus sont injectives.
\endthm

Compte tenu de ce lemme on peut identifier et on identifiera les images
de ces morphismes avec leurs sources.

\rem D\'{e}monstration
\endrem
Il suffit de d\'{e}montrer l'injectivit\'{e} de la premi\`{e}re
fl\`{e}che apr\`{e}s le changement de base par le rev\^{e}tement fini
\'{e}tale $a:S_{\natural ,1}\rightarrow S_{\natural}$, qui d\'{e}ploie
totalement $h_{\natural}:Z_{\natural}\rightarrow S_{\natural}$ et
$h_{\natural}':Z_{\natural}'\rightarrow S_{\natural}$,
consid\'{e}r\'{e} dans la construction de la fl\`{e}che
$L_{Z_{\natural}'/Z_{\natural}/S_{\natural}}\subset \mathop{\rm
Sym}\nolimits^{r}(X^{\ast}({\cal T}))$, construction dont on reprend
les notations.

Ce changement de base d\'{e}ploie aussi le tore ${\cal T}$ et le choix
d'un type CM permet d'identifier ${\cal T}_{S_{\natural ,1}}$ \`{a}
${\Bbb G}_{{\rm m},S_{\natural ,1}}^{\mathop{\rm Mor}
\nolimits_{S_{\natural}}(S_{\natural ,1},Z_{\natural})}$.  Comme
${\cal T}$ agit trivialement sur ${\cal N}_{S_{\natural}}$, on a donc
$$
\bigoplus_{n}{}^{\rm p}{\cal H}^{n}(g_{S_{\natural ,1},\ast}^{{\cal
T}_{S_{\natural ,1}}}{\Bbb Q}_{\ell})=\Bigl(\bigoplus_{n}{}^{\rm
p}{\cal H}^{n} (g_{S_{\natural ,1},\ast}{\Bbb Q}_{\ell})\Bigr)
[(t_{\varphi_{+}'})_{\varphi\in \mathop{\rm Mor}
\nolimits_{S_{\natural}}(S_{\natural ,1},Z_{\natural})}]
$$
o\`{u} $t_{\varphi_{+}'}$ est la classe de Chern du fibr\'{e} en
droites sur le classifiant $[S_{\natural ,1}/{\cal T}_{S_{\natural
,1}}]$ obtenu en poussant le ${\cal T}_{S_{\natural ,1}}$-torseur
universel par le caract\`{e}re $\chi_{\varphi_{+}'}$.

Notons $\chi_{\varphi_{+}'} (\widetilde{{\cal E}}_{12})$ le fibr\'{e}
en droites sur ${\cal N}_{S_{\natural ,1}}$ obtenu en poussant la
restriction \`{a} ${\cal N}_{S_{\natural ,1}}$ du $\widetilde{{\cal
T}}$-torseur $\widetilde{{\cal E}}_{12}$ par le caract\`{e}re
$\chi_{\varphi_{+}'}$ et notons simplement
$$
c_{\varphi_{+}'}:\bigoplus_{n}{}^{\rm p}{\cal H}^{n} (g_{S_{\natural
,1},\ast}{\Bbb Q}_{\ell})\rightarrow \bigoplus_{n}{}^{\rm p}{\cal
H}^{n+2} (g_{S_{\natural ,1},\ast}{\Bbb Q}_{\ell})(1)
$$
la fl\`{e}che induite par sa premi\`{e}re classe de Chern
$$
c_{1}(\chi_{\varphi_{+}'} (\widetilde{{\cal E}}_{12})): {\Bbb Q}_{\ell
,{\cal N}_{S_{\natural ,1}}} \rightarrow {\Bbb Q}_{\ell ,{\cal
N}_{S_{\natural ,1}}}[2](1).
$$
Il s'agit alors de d\'{e}montrer que la fl\`{e}che de degr\'{e} $2r$
$$\displaylines{
\qquad\prod_{\varphi\in \mathop{\rm Mor}\nolimits_{S_{\natural}}
(S_{\natural ,1},Z_{\natural})}(t_{\varphi_{+}'}+c_{\varphi_{+}'}):
\Bigl(\bigoplus_{n}{}^{\rm p}{\cal H}^{n} (g_{S_{\natural
,1},\ast}{\Bbb Q}_{\ell})\Bigr) [(t_{\varphi_{+}'})_{\varphi\in
\mathop{\rm Mor} \nolimits_{S_{\natural}}(S_{\natural
,1},Z_{\natural})}]
\hfill\cr\hfill
\rightarrow \Bigl(\bigoplus_{n}{}^{\rm p}{\cal
H}^{n} (g_{S_{\natural ,1},\ast}{\Bbb Q}_{\ell})\Bigr)
[(t_{\varphi_{+}'})_{\varphi\in \mathop{\rm Mor}
\nolimits_{S_{\natural}}(S_{\natural ,1},Z_{\natural})}](r)\qquad}
$$
est injective (voir le lemme A.2.1), ce qui est
\'{e}vident.
\hfill\hfill$\square$
\vskip 3mm

Notre r\'{e}sultat principal concernant la cohomologie
\'{e}quivariante de ${\cal M}_{S}$ au-dessus de l'ouvert
$S_{\natural}$ de $S$ est alors le suivant:

\thm TH\'{E}OR\`{E}ME 3.6.2
\enonce
Pour chacun des deux caract\`{e}res endoscopiques $\kappa :{\Bbb
Z}^{2}\rightarrow \{\pm 1\}$, l'application de restriction
$$
\bigoplus_{n}{}^{\rm p}{\cal H}^{n}(f_{S_{\natural},\ast}^{{\cal
T}}{\Bbb Q}_{\ell})_{\kappa}\rightarrow \bigoplus_{n}{}^{\rm p} {\cal
H}^{n}(g_{S_{\natural},\ast}^{{\cal T}}{\Bbb
Q}_{\ell})_{\kappa}
$$
est injective et son image est pr\'{e}cis\'{e}ment le sous-${\Bbb
Q}_{\ell}[X^{\ast}({\cal T})(-1)]$-module en faisceaux pervers
gradu\'{e}s sur $S$
$$
L_{Z_{\natural}'/Z_{\natural}/S_{\natural}}\otimes_{{\Bbb
Z}_{S_{\natural}}} \bigoplus_{n}{}^{\rm p}{\cal H}^{n-2r}
(g_{S_{\natural},\ast}^{{\cal T}}{\Bbb Q}_{\ell})_{\kappa}(-r)\subset
\bigoplus_{n}{}^{\rm p}{\cal H}^{n}(g_{S_{\natural},\ast}^{{\cal
T}}{\Bbb Q}_{\ell})_{\kappa}.
$$
\endthm

\rem D\'{e}monstration
\endrem
Pour chaque entier $n$, ${}^{\rm p}{\cal H}^{n}(f_{S_{\natural},\ast}
{\Bbb Q}_{\ell})_{\kappa}$ est potentiellement pur de poids $n$
d'apr\`{e}s le corollaire 3.2.4.  L'injectivit\'{e} de la fl\`{e}che
de restriction est donc une cons\'{e}quence du corollaire A.1.3.

Le reste de la d\'{e}monstration du th\'{e}or\`{e}me fait l'objet des
deux paragraphes qui suivent.

\hfill\hfill$\square$
\vskip 3mm

\thm COROLLAIRE 3.6.3
\enonce
Pour chacun des deux caract\`{e}res endoscopiques $\kappa :{\Bbb
Z}^{2}\rightarrow \{\pm 1\}$ il existe un isomorphisme de faisceaux
pervers gradu\'{e}s sur $S_{\natural}$
$$
\bigoplus_{n}{}^{\rm p}{\cal H}^{n}(f_{S_{\natural},\ast} {\Bbb
Q}_{\ell})_{\kappa}\buildrel\sim\over\longrightarrow
L_{Z_{\natural}'/Z_{\natural}/S_{\natural}}
\otimes_{{\Bbb Z}_{S_{\natural}}}\bigoplus_{n}{}^{\rm p}
{\cal H}^{n-2r}(g_{S_{\natural},\ast} {\Bbb Q}_{\ell})_{\kappa}(-r).
$$
\endthm

\rem D\'{e}monstration
\endrem
D'apr\`{e}s le th\'{e}or\`{e}me on a un isomorphisme de ${\Bbb
Q}_{\ell}[X^{\ast}({\cal T})(-1)]$-modules en faisceaux pervers
gradu\'{e}s sur $S$
$$
\bigoplus_{n}{}^{\rm p}{\cal H}^{n}(f_{S_{\natural},\ast}^{{\cal
T}}{\Bbb Q}_{\ell})_{\kappa}\buildrel\sim\over\longrightarrow
L_{Z_{\natural}'/Z_{\natural}/S_{\natural}}\otimes_{{\Bbb
Z}_{S_{\natural}}} \bigoplus_{n}{}^{\rm p}{\cal H}^{n-2r}
(g_{S_{\natural},\ast}^{{\cal T}}{\Bbb Q}_{\ell})_{\kappa}(-r).
$$
De plus, compte tenu de la puret\'{e} l'isomorphisme canonique
$$
f_{S_{\natural},\ast}{\Bbb Q}_{\ell ,S_{\natural}}={\Bbb Q}_{\ell}
\otimes_{\varepsilon_{\ast}^{{\cal T}} {\Bbb Q}_{\ell}}^{{\rm L}}
f_{S_{\natural},\ast}^{{\cal T}} {\Bbb Q}_{\ell}
$$
induit une suite spectrale dont la partie $\kappa$
d\'{e}g\'{e}n\`{e}re en un isomorphisme
$$
\bigoplus_{n}{}^{\rm p}{\cal H}^{n}(f_{S_{\natural},\ast}{\Bbb
Q}_{\ell})_{\kappa}={\Bbb Q}_{\ell ,S_{\natural}}\otimes_{{\Bbb
Q}_{\ell}[X^{\ast}({\cal T})(-1)]}\bigoplus_{n}{}^{\rm p}{\cal
H}^{n}(f_{S_{\natural},\ast}^{{\cal T}}{\Bbb Q}_{\ell})_{\kappa}.
$$ 
De m\^{e}me, il existe un isomorphisme
$$
\bigoplus_{n}{}^{\rm p}{\cal H}^{n}(g_{S_{\natural},\ast}{\Bbb
Q}_{\ell})_{\kappa}={\Bbb Q}_{\ell ,S_{\natural}}\otimes_{{\Bbb
Q}_{\ell}[X^{\ast}({\cal T})(-1)]}\bigoplus_{n}{}^{\rm p}{\cal
H}^{n}(g_{S_{\natural},\ast}^{{\cal T}}{\Bbb Q}_{\ell})_{\kappa},
$$ 
d'o\`{u} le corollaire.
\hfill\hfill$\square$

\subsection{3.7}{Une construction \`{a} la Altman et Kleiman}

Nous aurons besoin d'une construction inspir\'{e}e des espaces de
modules de pr\'{e}\-sen\-ta\-tions d'Altman et Kleiman (cf. [Al-Kl]).

Notons simplement
$$
S_{\natural ,1}\subset \overbrace{Z_{\natural}\times_{S_{\natural}}
Z_{\natural}\cdots \times_{S_{\natural}}Z_{\natural}}^{r}
$$
le sous-sch\'{e}ma ouvert et ferm\'{e} du produit fibr\'{e} form\'{e}
des $\underline{z}=(z_{1},\ldots ,z_{r})$ tels que $z_{i}\not=z_{j}$
pour tous $i\not=j$.  C'est un $S_{\natural}$-sch\'{e}ma fini
\'{e}tale de degr\'{e} $r!$ qui, par changement de base, d\'{e}ploie
compl\`{e}tement le rev\^{e}tement fini \'{e}tale
$Z_{\natural}\rightarrow S_{\natural}$.

Pour tout point g\'{e}om\'{e}trique $\underline{z}=(z_{1},\ldots
,z_{r})$ de $S_{\natural ,1}$ d'image $b$ dans $S_{\natural}$ et
tout $i=0,1,\ldots ,r$, on consid\`{e}re la courbe
$Y_{\underline{z}}^{[i]}$ d\'{e}duite de $Y_{b}$ en normalisant cette
courbe en les points $z_{i+1},\ldots ,z_{r}\in Z_{b}=Y_{1,b}\cap
Y_{2,b}$.  On a donc $Y_{\underline{z}}^{[0]}=Y_{1,b}\amalg
Y_{2,b}$ et $Y_{\underline{z}}^{[r]}=Y_{b}$, et pour chaque
$i=0,1,\ldots ,r$, la courbe $Y_{\underline{z}}^{[i]}$ est une
normalisation partielle de $Y_{\underline{z}}^{[i+1]}$ qui a pour
composantes (non irr\'{e}ductibles en g\'{e}n\'{e}ral) $Y_{1,b}$ et
$Y_{2,b}$, composantes dont l'intersection est r\'{e}duite \`{a}
$\{z_{1},\ldots ,z_{i}\}$ et qui se coupent transversalement en ces
points.
$$
\xymatrix{Y_{1,b}\amalg Y_{2,b}=Y_{\underline{z}}^{[0]}\ar[r]
& Y_{\underline{z}}^{[1]}\ar[r] & \cdots\ar[r] &
Y_{\underline{z}}^{[r]}=Y_{b}\ar[d]\cr
& & & S_{\natural ,1}}
$$
On a aussi le rev\^{e}tement \'{e}tale double
$$
Y_{\underline{z}}'^{[i]}=X'\times_{X}Y_{\underline{z}}^{[i]}\rightarrow
Y_{\underline{z}}^{[i]}.
$$

Pour $i$ fix\'{e} et $\underline{z}$ variable, les courbes
$Y_{\underline{z}}^{[i]}$ et leurs rev\^{e}tements \'{e}tales doubles
ci-dessus se mettent naturellement en famille sur $S_{\natural ,1}$.
On peut donc former le $S_{\natural ,1}$-champ alg\'{e}brique ${\cal
M}^{[i]}$ dont les points g\'{e}om\'{e}triques sont les triplets
$(\underline{z},{\cal F}^{[i]},\iota^{[i]})$ o\`{u} $\underline{z}$
est un point g\'{e}om\'{e}trique de $S_{\natural ,1}$, ${\cal
F}^{[i]}$ est un Module sans torsion de rang $1$ sur
$Y_{\underline{z}}'^{[i]}$ et $\iota^{[i]}:{\cal F}^{[i]}
\buildrel\sim\over\longrightarrow \tau^{\ast}({\cal F}^{[i]})^\vee$
est une structure unitaire.

On a ${\cal M}^{[0]}=S_{\natural ,1}\times_{S_{\natural}}{\cal
N}_{S_{\natural}}$ et ${\cal M}^{[r]}=S_{\natural ,1}
\times_{S_{\natural}}{\cal M}_{S_{\natural}}$ et on a des immersions
ferm\'{e}es de $S_{\natural ,1}$-champs alg\'{e}briques
$$
\xymatrix{{\cal M}^{[0]}\ar@{^{(}->}[r]^-{i^{[0]}} & {\cal
M}^{[1]}\ar@{^{(}->}[r]^-{i^{[1]}} & {\cal
M}^{[2]}\cdots\ar@{^{(}->}[r]^-{i^{[r-1]}} & {\cal M}^{[r]}}
$$
induites par les foncteurs d'image directe pour les normalisations
partielles $Y_{\underline{z}}^{[i]}\rightarrow
Y_{\underline{z}}^{[i+1]}$.  Le compos\'{e} de ces immersions
ferm\'{e}es n'est autre que le changement de base par
$S_{\natural ,1}\rightarrow S_{\natural}\subset S$ de l'immersion
ferm\'{e}e ${\cal N}_{S}\hookrightarrow {\cal M}_{S}$
consid\'{e}r\'{e}e pr\'{e}c\'{e}demment.

Bien s\^{u}r, pour chaque $i=0,1,\ldots ,r$ on peut d\'{e}finir
pareillement le $S_{\natural ,1}$-sch\'{e}ma en groupes $P^{[i]}$
des Modules inversibles unitaires sur la famille des
$Y_{\underline{z}}'^{[i]}$ et on a une action par produit tensoriel de
$P^{[i]}$ sur ${\cal M}^{[i]}$.  On a des homomorphismes
$$
P^{[r]}\rightarrow\cdots\rightarrow P^{[1]}\rightarrow P^{[0]}
$$
et, compte tenu de ces homomorphismes, pour chaque $i=0,1,\ldots
,r-1$, l'immersion ferm\'{e}e $i^{[i]}$ ci-dessus est
$P^{[i+1]}$-\'{e}quivariante.

Le $S_{\natural ,1}$-tore ${\cal T}^{[r]}:=S_{\natural ,1}
\times_{S_{\natural}}{\cal T}$ se d\'{e}compose naturellement en un
produit
$$
{\cal T}^{[r]}={\cal T}_{1}\times_{S_{\natural ,1}}{\cal T}_{2}
\cdots \times_{S_{\natural ,1}}{\cal T}_{r},
$$
o\`{u} ${\cal T}_{i}$ est le noyau de l'homomorphisme
$P^{[i]}\rightarrow P^{[i-1]}$ et est le changement de base par la
$i$-\`{e}me projection $S_{\natural ,1}\rightarrow Z_{\natural}$ du
tore
$$
(\pi_{\ast}{\Bbb G}_{{\rm m},Z_{\natural}'})^{\tau^{\ast}=(-)^{-1}}.
$$

Notons
$$
S_{\natural ,1}'\subset \overbrace{Z_{\natural}'\times_{S_{\natural}}
Z_{\natural}'\cdots \times_{S_{\natural}}Z_{\natural}'}^{r}
$$
le sous-sch\'{e}ma ouvert et ferm\'{e} du produit fibr\'{e} form\'{e} des
$\underline{z}'=(z_{1}',\ldots ,z_{r}')$ dont l'image $\underline{z}$
dans $Z_{\natural}\times_{S_{\natural}} Z_{\natural}\cdots
\times_{S_{\natural}}Z_{\natural}$ est dans $S_{\natural ,1}$.  C'est
un rev\^{e}tement fini \'{e}tale de $S_{\natural ,1}$ de degr\'{e}
$2^{r}$.

Pour chaque $i=0,1,\ldots ,r-1$, consid\'{e}rons le
$S_{\natural ,1}'$-champ alg\'{e}brique ${\cal M}'^{[i,i+1]}$ des
$$
(\underline{z}',({\cal F}^{[i]},\iota^{[i]}),(\varepsilon :{\cal G}
\hookrightarrow {\cal F}^{[i,i+1]}))
$$
o\`{u}:
\vskip 1mm

\item{1)} $\underline{z}'$ est un point de
$S_{\natural ,1}'$, d'image $\underline{z}$ dans $S_{\natural ,1}$ et
$b$ dans $S_{\natural}$,
\vskip 1mm

\item{2)} le couple $(\underline{z}',({\cal F}^{[i]},\iota^{[i]}))$
est un point de ${\cal M}^{[i]}$,
\vskip 1mm

\item{3)} $\varepsilon$ est une modification \'{e}l\'{e}mentaire
inf\'{e}rieure en $z_{i+1}'$ de ${\cal F}^{[i,i+1]}$, c'est-\`{a}-dire
d'un homomorphisme injectif de ${\cal
O}_{Y_{\underline{z}}'^{[i+1]}}$-Modules
$$
\varepsilon :{\cal G}\hookrightarrow {\cal F}^{[i,i+1]}
$$
o\`{u}:
\vskip 1mm

\itemitem{-} $({\cal F}^{[i,i+1]},\iota^{[i,i+1]})$ est l'image directe
par le morphisme de normalisation partielle
$Y_{\underline{z}}'^{[i]}\rightarrow Y_{\underline{z}}'^{[i+1]}$ de
$({\cal F}^{[i]},\iota^{[i]})$,
\vskip 1mm

\itemitem{-} ${\cal G}$ est un Module sans torsion de rang $1$ sur
$Y_{\underline{z}}'^{[i+1]}$,
\vskip 1mm

\itemitem{-} $\mathop{\rm Coker}(\varepsilon )$ est de longueur $1$ et
\`{a} support dans $\{z_{i+1}'\}$.
\vskip 2mm

L'oubli de $\varepsilon$ est une fibration en droites projectives
$$
p'^{[i,i+1]}:{\cal M}'^{[i,i+1]}\rightarrow S_{\natural ,1}'
\times_{S_{\natural ,1}}{\cal M}^{[i]}
$$
et on a les deux sections
$$
\sigma_{1}^{[i]},\sigma_{2}^{[i]}:S_{\natural ,1}'
\times_{S_{\natural ,1}}{\cal M}^{[i]}\rightarrow {\cal M}'^{[i,i+1]}
$$
de $p'^{[i,i+1]}$ d\'{e}finies comme suit: soient $\alpha\in\{1,2\}$
et $(\underline{z}',({\cal F}^{[i]},\iota^{[i]}))$ un point de
$S_{\natural ,1}' \times_{S_{\natural ,1}}{\cal M}^{[i]}$, on
d\'{e}finit la modification inf\'{e}rieure \'{e}l\'{e}mentaire
$$
\varepsilon_{\alpha}:{\cal G}_{\alpha}\hookrightarrow {\cal F}^{[i,i+1]}
$$
qui intervient dans $\sigma_{1}^{[\alpha ]}(\underline{z}',({\cal
F}^{[i]},\iota^{[i]}))$ en prenant l'image directe par la
normalisation partielle $Y_{\underline{z}}'^{[i]}\rightarrow
Y_{\underline{z}}'^{[i+1]}$ de la modification \'{e}l\'{e}mentaire
$$
{\cal F}^{[i]}(-[z_{i+1,\alpha}'])\hookrightarrow {\cal F}^{[i]}
$$
o\`{u} $z_{i+1,\alpha}'$ est l'unique point lisse de
$Y_{\alpha ,b}'$ au-dessus du point double $z_{i+1}'$ de
$Y_{\underline{z}}'^{[i+1]}$.

On d\'{e}finit un morphisme de $S_{\natural ,1}'$-champs
alg\'{e}briques
$$
\rho'^{[i,i+1]}:{\cal M}'^{[i,i+1]}\rightarrow
S_{\natural ,1}'\times_{S_{\natural ,1}}{\cal M}^{[i+1]}
$$
en envoyant $\bigl(\underline{z}',({\cal
F}^{[i]},\iota^{[i]}),\varepsilon:{\cal G}\hookrightarrow {\cal
F}^{[i,i+1]}\bigr)$ sur $(\underline{z}',{\cal
F}^{[i+1]},\iota^{[i+1]})$ o\`{u} ${\cal F}^{[i+1]}$ est la somme
amalgam\'{e}e de la modification inf\'{e}rieure $\varepsilon :{\cal
G}\hookrightarrow {\cal F}^{[i,i+1]}$ en $z_{i+1}'$ et de la
modification sup\'{e}rieure $(\tau^{\ast}\varepsilon^{\vee})\circ
\iota^{[i,i+1]}: {\cal F}^{[i,i+1]}\hookrightarrow \tau^{\ast}{\cal
G}^{\vee}$ en $\tau (z_{i+1}')\not=z_{i+1}'$, et o\`{u}
$\iota^{[i+1]}$ est la structure unitaire \'{e}vidente.

Le lemme suivant est une variante du th\'{e}or\`{e}me 18 de [Al-Kl]
et on renvoie \`{a} cette r\'{e}f\'{e}rence pour les d\'{e}tails de sa
d\'{e}monstration.

\thm LEMME 3.7.1
\enonce
Le morphisme $\rho'^{[i,i+1]}:{\cal M}'^{[i,i+1]}\rightarrow
S_{\natural ,1}'\times_{S_{\natural ,1}}{\cal M}^{[i+1]}$ est un
pincement de ${\cal M}'^{[i,i+1]}$ identifiant les deux sections
$$
\sigma_{1}^{[i]},\sigma_{2}^{[i]}:S_{\natural ,1}'
\times_{S_{\natural ,1}}{\cal M}^{[i]}\rightarrow {\cal
M}'^{[i,i+1]}
$$
de $p'^{[i,i+1]}$. Plus pr\'{e}cis\'{e}ment,
\vskip 1mm

\item{-} $\rho'^{[i,i+1]}$ est fini,
\vskip 1mm

\item{-} l'image r\'{e}ciproque par $\rho'^{[i,i+1]}$ de l'image de
l'immersion ferm\'{e}e $S_{\natural ,1}'\times_{S_{\natural
,1}}i^{[i]}: S_{\natural ,1}'\times_{S_{\natural ,1}}{\cal M}^{[i]}
\hookrightarrow S_{\natural ,1}' \times_{S_{\natural ,1}}{\cal
M}^{[i+1]}$ est la r\'{e}union des images des deux sections
$\sigma_{1}^{[i]}$ et $\sigma_{2}^{[i]}$,
\vskip 1mm

\item{-} $\rho'^{[i,i+1]}$ est un isomorphisme en dehors des images de
ces deux sections,
\vskip 1mm

\item{-} il existe deux sections globales $h_{1}^{[i]}$ et $h_{2}^{[i]}$ du
$S_{\natural ,1}'$-sch\'{e}ma en groupes $S_{\natural ,1}'
\times_{S_{\natural ,1}}P^{[i]}$ induisant des automorphismes
de $S_{\natural ,1}'\times_{S_{\natural ,1}}{\cal
M}^{[i]}$ not\'{e}s aussi  $h_{1}^{[i]}$ et $h_{2}^{[i]}$, tels que
le carr\'{e}
$$
\xymatrix{S_{\natural ,1}'\times_{S_{\natural ,1}}{\cal
M}^{[i]}\ar[d]_{h_{\alpha}^{[i]}}\ar@{^{(}->}[r]^-{\sigma_{\alpha}^{[i]}} &
{\cal M}'^{[i,i+1]}\ar[d]^{\rho'^{[i,i+1]}}\cr
S_{\natural ,1}'\times_{S_{\natural ,1}}{\cal M}^{[i]}
\ar@{^{(}->}[r]_-{i^{[i]}} & S_{\natural ,1}'
\times_{S_{\natural ,1}} {\cal M}^{[i+1]}}
$$
soit commutatif pour $\alpha =1,2$.

\endthm

\rem D\'{e}monstration
\endrem
Pour $\alpha =1,2$, on a une structure unitaire \'{e}vidente induite
par $\iota^{[i]}$ sur le Module sans torsion de rang $1$, ${\cal
F}^{[i]}([\tau (z_{i+1,\alpha}')] -[z_{i+1,\alpha}'])$, structure que
l'on note $\iota^{[i]} ([\tau (z_{i+1,\alpha}')]-[z_{i+1,\alpha}'])$.
Avec cette notation on a
$$\displaylines{
\quad\rho'^{[i,i+1]}\circ\sigma_{\alpha}^{[i]}(\underline{z}',{\cal
F}^{[i]},\iota^{[i]})
\hfill\cr\hfill
=\Bigl(\underline{z}',(Y_{\underline{z}}'^{[i]}
\rightarrow Y_{\underline{z}}'^{[i+1]})_{\ast}\bigl({\cal F}^{[i]}
([\tau (z_{i+1,\alpha}')]-[z_{i+1,\alpha}']),\iota^{[i]}
([\tau (z_{i+1,\alpha}')]-[z_{i+1,\alpha}'])\bigr)\Bigr)\quad}
$$
pour $\alpha =1,2$ et
$$
i^{[i]}(\underline{z}',{\cal F}^{[i]},\iota^{[i]})=
\bigl(\underline{z}',(Y_{\underline{z}}'^{[i]}
\rightarrow Y_{\underline{z}}'^{[i+1]})_{\ast}({\cal
F}^{[i]},\iota^{[i]})\bigr)
$$
de sorte que la section $h_{\alpha}^{[i]}$ qui envoie $\underline{z}'$
sur ${\cal O}_{Y_{\underline{z}}'^{[i]}}([\tau
(z_{i+1,\alpha}')]-[z_{i+1,\alpha}'])$ muni de sa structure unitaire
\'{e}vidente, r\'{e}pond \`{a} la question.
\hfill\hfill$\square$
\vskip 3mm

Par dualit\'{e} et fonctorialit\'{e} on a une modification
\'{e}l\'{e}mentaire sup\'{e}rieure en $\tau (z_{i+1}')$
$$
(\tau^{\ast}\varepsilon^{\vee})\circ \iota^{[i,i+1]}:
{\cal F}^{[i,i+1]}\hookrightarrow \tau^{\ast}{\cal G}^{\vee}.
$$
Si on note (par abus) $[2z_{i+1}']$ le diviseur de Cartier sur
$Y_{\underline{z}}'^{[i+1]}$ induit par la fibre de
$Y_{\underline{z}}'^{[i+1]}\rightarrow \kappa (b)\otimes_{k}X'$
passant par le point $z_{i+1}'$, dans un voisinage suffisamment petit
de $z_{i+1}'$ (on veut \'{e}viter les autres $z_{j}'$ et $\tau
(z_{j}')$ qui pourrait \^{e}tre sur cette fibre), le ${\cal
O}_{Y_{\underline{z}}'^{[i+1]}}$-Module inversible ${\cal
O}_{Y_{\underline{z}}'^{[i+1]}}(-[2z_{i+1}'])$ est un Id\'{e}al de
co-longueur $2$ de ${\cal O}_{Y_{\underline{z}}'^{[i+1]}}$ et le
morphisme canonique
$$
\tau^{\ast}({\cal G}^{\vee}\otimes_{{\cal
O}_{Y_{\underline{z}}'^{[i+1]}}} {\cal
O}_{Y_{\underline{z}}'^{[i+1]}}(-[2z_{i+1}']))\hookrightarrow
\tau^{\ast}{\cal G}^{\vee}\otimes_{{\cal
O}_{Y_{\underline{z}}'^{[i+1]}}} {\cal
O}_{Y_{\underline{z}}'^{[i+1]}}=\tau^{\ast}{\cal
G}^{\vee}
$$
se factorise en une modification \'{e}l\'{e}mentaire inf\'{e}rieure en
$\tau (z_{i+1}')$
$$
\varepsilon^{\tau}:{\cal G}^{\tau}:=\tau^{\ast}({\cal G}^{\vee}
\otimes_{{\cal O}_{Y_{\underline{z}}'^{[i+1]}}}{\cal
O}_{Y_{\underline{z}}'^{[i+1]}}(-[2z_{i+1}'])) \hookrightarrow {\cal
F}^{[i,i+1]}
$$
suivi de $(\tau^{\ast}\varepsilon^{\vee})\circ\iota^{[i,i+1]}$.

On munit ${\cal M}'^{[i,i+1]}$ de la donn\'{e}es de descente,
relativement \`{a} $S_{\natural ,1}'\rightarrow
S_{\natural ,1}$, qui d\'{e}finie par l'involution $\tau$ qui
envoie $(\underline{z}',{\cal F}^{[i]},\iota^{[i]},\varepsilon :{\cal
G}\hookrightarrow {\cal F}^{[i,i+1]})$ sur
$$
\bigl(\tau (z'),{\cal F}^{[i]},\iota^{[i]},\varepsilon^{\tau}:{\cal
G}^{\tau}\hookrightarrow {\cal F}^{[i,i+1]}\bigr).
$$
Remarquons que les sections $\sigma_{1}^{[i]}$ et $\sigma_{2}^{[i]}$
ci-dessus sont \'{e}chang\'{e}es par cette involution.

On a donc un $S_{\natural ,1}$-champ alg\'{e}brique ${\cal
M}^{[i,i+1]}$ dont le changement de base par
$S_{\natural ,1}'\rightarrow S_{\natural ,1}$ est ${\cal
M}'^{[i,i+1]}$ et un morphisme d'oubli
$$
p^{[i,i+1]}:{\cal M}^{[i,i+1]}\rightarrow {\cal M}^{[i]}
$$
qui est une fibration en droites projectives {\og}{tordues}{\fg}.  Le
points g\'{e}om\'{e}triques de ${\cal M}^{[i,i+1]}$ sont les uplets
$$
\bigl(z,{\cal F}^{[i]},\iota^{[i]},\varepsilon':{\cal G}'
\hookrightarrow {\cal F}^{[i,i+1]}, \varepsilon'':{\cal
G}''\hookrightarrow {\cal F}^{[i,i+1]}\bigr)
$$
form\'{e}s d'un point g\'{e}om\'{e}trique $\underline{z}$ de
$S_{\natural ,1}$ d'image $b$ dans $S_{\natural}$, de $({\cal
F}^{[i]},\iota^{[i]})$ comme ci-dessus et de deux modifications
\'{e}l\'{e}mentaires inf\'{e}rieures $\varepsilon'$ et $\varepsilon''$
en $z_{i+1}'$ et $z_{i+1}''$ respectivement, o\`{u}
$\{z_{i+1}',z_{i+1}''\}$ est la fibre de $Z_{\natural}'\rightarrow
Z_{\natural}$ en $z_{i+1}$, modifications qui v\'{e}rifient la
compatibilit\'{e}
$$
\varepsilon'^{\tau}=\iota^{[i,i+1]}\circ\varepsilon''.
$$

Le morphisme $\rho'^{[i,i+1]}$ se descend un un morphisme de
$S_{\natural ,1}$-champs alg\'{e}briques
$$
\rho^{[i,i+1]}:{\cal M}^{[i,i+1]}\rightarrow {\cal M}^{[i+1]}.
$$

\subsection{3.8}{Preuve du th\'{e}or\`{e}me sur l'ouvert
$S_{\natural}$}

Pour terminer la d\'{e}monstration du th\'{e}or\`{e}me 3.6.2, il
suffit de le faire apr\`{e}s avoir remplac\'{e} $S_{\natural}$ par un
rev\^{e}tement fini \'{e}tale.

On peut donc se placer sur $S_{\natural ,1}$ et utiliser la
construction \`{a} la Altman et Kleiman de la section
pr\'{e}c\'{e}dente
$$\xymatrix{&{\cal
M}^{[i,i+1]}\ar[ld]_{p^{[i,i+1]}}^{}\ar[rd]^{\rho^{[i,i+1]}} &\cr
{\cal M}^{[i]}\ar[rd]_{f^{[i]}}  &  & {\cal M}^{[i+1]}\ar[ld]^{f^{[i+1]}}
\cr
& S_{\natural ,1} &}
$$
o\`{u} on a not\'{e} $f^{[i]}$ la projection canonique de ${\cal
M}^{[i]}$ sur $S_{\natural ,1}$.

On v\'{e}rifie que
$$
L_{Z_{\natural}'\times_{S_{\natural}}S_{\natural ,1}/
Z_{\natural}\times_{S_{\natural}}S_{\natural ,1}/
S_{\natural ,1}}=L_{1}\otimes_{S_{\natural ,1}}L_{2}
\otimes_{S_{\natural ,1}}\cdots \otimes_{S_{\natural ,1}}
L_{r}
$$
o\`{u} $L_{i}\subset X^{\ast}({\cal T}_{i})$ est l'image
r\'{e}ciproque par la $i$-\`{e}me projection
$S_{\natural ,1}\rightarrow Z_{\natural}$ de
$L_{Z_{\natural}'/Z_{\natural}/Z_{\natural}}$. On note
$$
L^{[i]}=L_{1}\otimes_{S_{\natural ,1}}L_{2}
\otimes_{S_{\natural ,1}}\cdots \otimes_{S_{\natural ,1}}L_{i}
\subset\mathop{\rm Sym}\nolimits_{{\Bbb Z}_{S_{\natural ,1}}}^{i}
(X^{\ast}({\cal T}^{[i]})).
$$
o\`{u} ${\cal T}^{[i]}={\cal T}_{1}\times_{S_{\natural ,1}}{\cal
T}_{2}\times_{S_{\natural ,1}}\cdots \times_{S_{\natural ,1}}
{\cal T}_{i}$.

On montre alors par r\'{e}currence sur $i=0,1,\ldots ,r$ que, pour
chacun des deux caract\`{e}res endoscopiques $\kappa :{\Bbb
Z}^{2}\rightarrow \{\pm 1\}$, la fl\`{e}che de restriction le long de
l'immersion ferm\'{e}e $i^{[i-1]}\circ\cdots i^{[1]}\circ
i^{[0]}:{\cal M}^{[0]}\hookrightarrow {\cal M}^{[i]}$,
$$
\bigoplus_{n}{}^{\rm p}{\cal H}^{n}((f^{[i]})_{\ast}^{{\cal T}^{[i]}}{\Bbb
Q}_{\ell})_{\kappa}\rightarrow \bigoplus_{n}{}^{\rm p} {\cal
H}^{n}((f^{[0]})_{\ast}^{{\cal T}^{[i]}}{\Bbb Q}_{\ell})_{\kappa}
$$
est injective et que son image est pr\'{e}cis\'{e}ment le sous-${\Bbb
Q}_{\ell}[X^{\ast}({\cal T}^{[i]})(-1)]$-module en faisceaux pervers
gradu\'{e}s sur $S$
$$
L^{[i]}\otimes_{{\Bbb Z}_{S_{\natural ,1}}}\bigoplus_{n}
{}^{\rm p}{\cal H}^{n-2i}((f^{[0]})_{\ast}^{{\cal T}^{[i]}}{\Bbb
Q}_{\ell})_{\kappa}(-i)\subset \bigoplus_{n}{}^{\rm p}{\cal
H}^{n}((f^{[0]})_{\ast}^{{\cal T}^{[i]}}{\Bbb Q}_{\ell})_{\kappa}.
$$

Le cran $i=0$ de la r\'{e}currence \'{e}tant tautologique, il suffit
de d\'{e}duire le cran $i$ du cran $i-1$.  Apr\`{e}s le changement de
base $S_{\natural ,1}'\rightarrow S_{\natural ,1}$ le diagramme
de champs ci-dessus devient un diagramme
$$\xymatrix{&{\cal
M}'^{[i-1,i]}\ar[ld]_{p'^{[i-1,i]}}^{}\ar[rd]^{\rho'^{[i-1,i]}} &\cr
{\cal M}'^{[i-1]}\ar[rd]_{f'^{[i-1]}} & & {\cal
M}'^{[i]}\ar[ld]^{f'^{[i]}} \cr
& S_{\natural ,1}' &}
$$
o\`{u} $p'^{[i-1,i]}$ est maintenant un fibr\'{e} en droites
projectives non tordues munies de deux sections $\sigma_{1}^{[i-1]}$
et $\sigma_{2}^{[i-1]}$, o\`{u} $S_{\natural ,1}'
\times_{S_{\natural ,1}}{\cal T}_{i}={\Bbb G}_{{\rm m},
S_{\natural ,1}'}$ agit par homoth\'{e}ties en fixant ces deux
sections et o\`{u} $\rho'^{[i-1,i]}$ est un pincement ${\Bbb G}_{{\rm
m},S_{\natural ,1}'}$-\'{e}quivariant de ce fibr\'{e} projectif le
long de ces deux sections.  On se retrouve en fait dans une situation
du type de celle consid\'{e}r\'{e}e dans la section (A.2).

Il r\'{e}sulte donc de la proposition A.2.5 que la fl\`{e}che de
restriction le long de l'immersion ferm\'{e}e $i'^{[i-1]}:{\cal
M}'^{[i-1]} \hookrightarrow {\cal M}'^{[i]}$
$$\eqalign{
\bigoplus_{n}{}^{\rm p}{\cal H}^{n}((f'^{[i]})_{\ast}^{{\Bbb G}_{{\rm
m},S_{\natural ,1}'}}{\Bbb Q}_{\ell})_{\kappa}\rightarrow
&\bigoplus_{n}{}^{\rm p} {\cal H}^{n}((f'^{[i-1]})_{\ast}^{{\Bbb
G}_{{\rm m},S_{\natural ,1}'}}{\Bbb Q}_{\ell})_{\kappa}\cr
&=\Bigl(\bigoplus_{n}{}^{\rm p} {\cal H}^{n}((f'^{[i-1]})_{\ast}
{\Bbb Q}_{\ell})_{\kappa}\Bigr)[t_{i}]\cr}
$$
est injective et que son image est pr\'{e}cis\'{e}ment le sous-${\Bbb
Q}_{\ell}[t_{i}]$-module en faisceaux pervers gradu\'{e}s sur $S$
$$
(t_{i}+c_{i})\Bigl( \bigoplus_{n}{}^{\rm p} {\cal H}^{n}((f'^{[i-1]})_{\ast}
{\Bbb Q}_{\ell})_{\kappa}\Bigr)[t_{i}]\subset
\Bigl(\bigoplus_{n}{}^{\rm p} {\cal H}^{n}((f'^{[i-1]})_{\ast}
{\Bbb Q}_{\ell})_{\kappa}\Bigr)[t_{i}]
$$
o\`{u} $c_{i}$ est la classe de Chern, not\'{e}e $c_{S}$ dans (A.2),
du fibr\'{e} vectoriel d\'{e}finissant le fibr\'{e} projectif
$p'^{[i-1,i]}$.  On rappelle que l'on a un isomorphisme canonique
$$
{\Bbb P}({\cal F}_{z_{i+1}'}^{[i,i+1]})={\Bbb P}({\cal F}_{z_{i+1,1}'}^{[i]}
\oplus {\cal F}_{z_{i+1,2}'}^{[i]})={\Bbb P}\bigl(({\cal F}_{z_{i+1,1}'}^{[i]}
\otimes ({\cal F}_{z_{i+1,2}'}^{[i]})^{\otimes -1})\oplus\kappa 
(z_{i+1}')\bigr).
$$
Par descente de $S_{\natural ,1}'$ \`{a} $S_{\natural
,1}$ on en d\'{e}duit que la fl\`{e}che de restriction le long de
l'immersion ferm\'{e}e $i^{[i-1]}:{\cal M}^{[i-1]} \hookrightarrow
{\cal M}^{[i]}$
$$
\bigoplus_{n}{}^{\rm p}{\cal H}^{n}((f^{[i]})_{\ast}^{{\cal
T}_{i}}{\Bbb Q}_{\ell})_{\kappa}\rightarrow \bigoplus_{n}{}^{\rm p}
{\cal H}^{n}((f^{[i-1]})_{\ast}^{{\cal T}_{i}}{\Bbb
Q}_{\ell})_{\kappa}
$$
est injective et que son image est pr\'{e}cis\'{e}ment le sous-${\Bbb
Q}_{\ell}[X^{\ast}({\cal T}_{i})(-1)]$-module en faisceaux pervers
gradu\'{e}s sur $S$
$$
L_{i}\otimes_{{\Bbb Z}_{S_{\natural ,1}}}\bigoplus_{n}{}^{\rm p}
{\cal H}^{n-2}((f^{[i-1]})_{\ast}^{{\cal T}_{i}}{\Bbb
Q}_{\ell})_{\kappa}(-1)\subset \bigoplus_{n}{}^{\rm p}{\cal
H}^{n}((f^{[i-1]})_{\ast}^{{\cal T}_{i}}{\Bbb Q}_{\ell})_{\kappa}.
$$

On conclut enfin la r\'{e}currence et donc la preuve du
th\'{e}or\`{e}me 3.6.2 par une formule de K\"{u}nneth.

\subsection{3.9}{L'argument de d\'{e}formation}

R\'{e}sumons la situation.  Au-dessus de l'hens\'{e}lis\'{e} $S$ de
${\Bbb A}_{H}$ en le point $a\in{\Bbb A}_H(k)$, on a deux
morphismes propres $f_{S}:{\cal M}_{S}\rightarrow S$ et $g_{S}:{\cal
N}_{S}\rightarrow S$.

Comme ${\cal N}_{S}$ est formellement lisse sur $k$, chaque faisceau
de cohomologie perverse ${}^{{\rm p}}{\cal H}^{n}(g_{S,\ast}{\Bbb
Q}_{\ell})$ est pur de poids $n$ et il en est de m\^{e}me de sa partie
$\kappa$, ${}^{{\rm p}}{\cal H}^{n}(g_{S,\ast}{\Bbb
Q}_{\ell})_{\kappa}$, pour n'importe quel caract\`{e}re $\kappa$ de
$({\Bbb Z}/2{\Bbb Z})^{2}$.  En particulier, les faisceaux pervers
${}^{{\rm p}}{\cal H}^{n}(g_{S,\ast}{\Bbb Q}_{\ell})$ et ${}^{{\rm
p}}{\cal H}^{n}(g_{S,\ast}{\Bbb Q}_{\ell})_{\kappa}$ sur $S$ sont
g\'{e}om\'{e}triquement semi-simples.

On a de plus l'action du groupe discret ${\Bbb Z}^{2}$ sur ${\cal
M}_{S}$ pour laquelle l'action induite de ${\Bbb Z}^{2}$ sur chaque
faisceau de cohomologie perverse ${}^{{\rm p}}{\cal
H}^{n}(f_{S,\ast}{\Bbb Q}_\ell)$ se factorise \`{a} travers le
quotient ${\Bbb Z}^{2}\rightarrow({\Bbb Z}/2{\Bbb Z})^{2}$.  Bien que
${\cal M}_{S}$ ne soit pas formellement lisse sur $k$, on sait que,
pour chacun des deux caract\`{e}res endoscopiques $\kappa$ de $({\Bbb
Z}/2{\Bbb Z})^{2}$ et chaque entier $n$ la partie $\kappa$-isotypique
${}^{{\rm p}}{\cal H}^{n}(f_{S,\ast}{\Bbb Q}_\ell)_\kappa$ de
${}^{{\rm p}}{\cal H}^{n}(f_{S,\ast}{\Bbb Q}_\ell)$ est pure de poids
$n$.  En particulier, cette partie $\kappa$-isotypique est un faisceau
g\'{e}om\'{e}triquement semi-simple.

On a enfin une action du $S$-tore $T$ sur ${\cal M}_{S}$ qui commute
\`{a} l'action de ${\Bbb Z}^{2}$ et dont le lieu des points fixes est
le ferm\'{e} image de ${\cal N}_{S}$ dans ${\cal M}_{S}$.   On peut
donc consid\'{e}rer les faisceaux pervers gradu\'{e}s de cohomologie
$T$-\'{e}quivariante
$$
\bigoplus_{n}{}^{{\rm p}}{\cal H}^{n}(f_{S,\ast}^{T}{\Bbb
Q}_{\ell})\hbox{ et }\bigoplus_{n}{}^{{\rm p}}{\cal
H}^{n}(g_{S,\ast}^{T} {\Bbb Q}_{\ell})
$$
et leur parties $\kappa$-isotypiques
$$
\bigoplus_{n}{}^{{\rm p}}{\cal H}^{n}(f_{S,\ast}^{T}{\Bbb
Q}_{\ell})_{\kappa}\hbox{ et }\bigoplus_{n}{}^{{\rm p}}{\cal
H}^{n}(g_{S,\ast}^{T} {\Bbb Q}_{\ell})_{\kappa}
$$
pour n'importe quel caract\`{e}re $\kappa$ de $({\Bbb Z}/2{\Bbb
Z})^{2}$.  Ce sont des faisceaux en modules gradu\'{e}s sur le
faisceau en ${\Bbb Q}_{\ell}$-alg\`{e}bres gradu\'{e}es ${\Bbb
Q}_{\ell}[X^{\ast}(T)(-1)]=\mathop{\rm Sym}(X^{\ast}(T)
\otimes_{{\Bbb Z}_{S}}{\Bbb Q}_{\ell ,S}(-1))$ sur $S$.

Fixons un caract\`{e}re endoscopique $\kappa$ de $({\Bbb Z}/2{\Bbb
Z})^{2}$.  Comme $T$ agit trivialement sur ${\cal N}_{S}$, on a
$$
N:=\bigoplus_{n}{}^{{\rm p}}{\cal H}^{n}(g_{S,\ast}^{T} {\Bbb
Q}_{\ell})_{\kappa}={\Bbb Q}_{\ell}[X^{\ast}(T)(-1)]
\otimes_{{\Bbb Q}_{\ell ,S}}\bigoplus_{n} {}^{{\rm p}}
{\cal H}^{n}(g_{S,\ast}{\Bbb Q}_{\ell})_{\kappa}
$$
et, tout comme $\bigoplus_{n} {}^{{\rm p}}{\cal H}^{n}(g_{S,\ast}{\Bbb
Q}_{\ell})_{\kappa}$, $N$ est un faisceau pervers gradu\'{e}
g\'{e}om\'{e}triquement semi-simple.

La puret\'{e} des ${}^{{\rm p}}{\cal H}^{n}(f_{S,\ast}{\Bbb
Q}_{\ell})_{\kappa}$ implique que la fl\`{e}che de
restriction
$$
M:=\bigoplus_{n}{}^{{\rm p}}{\cal H}^{n}(f_{S,\ast}^{T}{\Bbb
Q}_{\ell})_{\kappa} \rightarrow \bigoplus_{n}{}^{{\rm p}}{\cal
H}^{n}(g_{S,\ast}^{T} {\Bbb Q}_{\ell})_{\kappa}=N
$$
est injective, que tout comme $\bigoplus_{n} {}^{{\rm p}}{\cal
H}^{n}(f_{S,\ast}{\Bbb Q}_{\ell})_{\kappa}$, $M$ est un faisceau
pervers gradu\'{e} g\'{e}om\'{e}triquement semi-simple, et que
$$
\bigoplus_{n}{}^{{\rm p}}{\cal H}^{n}(f_{S,\ast}{\Bbb
Q}_{\ell})_{\kappa}={\Bbb Q}_{\ell}
\otimes_{{\Bbb Q}_{\ell}[X^{\ast}(T)(-1)]}M.
$$

\thm PROPOSITION 3.9.1
\enonce
Pour tout point g\'{e}om\'{e}trique $b$ de $S$ il existe un trait
strictement hens\'{e}lien $V$ de point ferm\'{e} $v$ et de point
g\'{e}n\'{e}rique $\eta$ et un morphisme de sch\'{e}mas $\varphi
:V\rightarrow S$ tel que
\vskip 1mm

\item{-} $\varphi (v)=b$ et $\varphi (\eta )\in S_{\natural}\subset
S$,
\vskip 1mm

\item{-} $\varphi^{\ast}N$ est \`{a} cohomologie ordinaire constante.
\endthm

\rem D\'{e}monstration
\endrem
C'est une cons\'{e}quence imm\'{e}diate de la proposition 2.9.1.

\hfill\hfill$\square$
\vskip 3mm

\thm COROLLAIRE 3.9.2
\enonce
Si $j:S_{\natural}\hookrightarrow S$ est l'inclusion, on a
$$
M=j_{!\ast}j^{\ast}M\hbox{ et }N=j_{!\ast}j^{\ast}N
$$
et aussi
$$
{}^{{\rm p}}{\cal H}^{n}(f_{S,\ast}{\Bbb
Q}_{\ell})_{\kappa}=j_{!\ast}j^{\ast}\,{}^{{\rm p}}{\cal
H}^{n}(f_{S,\ast}{\Bbb Q}_{\ell})_{\kappa}\hbox{ et }
{}^{{\rm p}}{\cal H}^{n}(g_{S,\ast}{\Bbb
Q}_{\ell})_{\kappa}=j_{!\ast}j^{\ast}\,{}^{{\rm p}}{\cal
H}^{n}(g_{S,\ast}{\Bbb Q}_{\ell})_{\kappa}
$$
quel que soit l'entier $n$.
\endthm

\rem D\'{e}monstration
\endrem
Comme $M$ (resp.  $N$) est pervers la fl\`{e}che d'oubli des supports
${}^{{\rm p}}{\cal H}^{0}(j_{!}j^{\ast}M)\rightarrow {}^{{\rm p}}{\cal
H}^{0}(j_{\ast}j^{\ast}M)$ (resp.  ${}^{{\rm p}}{\cal H}^{0}
(j_{!}j^{\ast}N)\rightarrow {}^{{\rm p}}{\cal H}^{0}
(j_{\ast}j^{\ast}N)$) se factorise en
$$
{}^{{\rm p}}{\cal H}^{0}(j_{!}j^{\ast}M)\rightarrow M\rightarrow
{}^{{\rm p}}{\cal H}^{0}(j_{\ast}j^{\ast}M)\hbox{ (resp.  }
{}^{{\rm p}}{\cal H}^{0}(j_{!}j^{\ast}N)\rightarrow N\rightarrow
{}^{{\rm p}}{\cal H}^{0}(j_{\ast}j^{\ast}N)\,)
$$
et il s'agit donc de d\'{e}montrer que les fl\`{e}ches de faisceaux
pervers ${}^{{\rm p}}{\cal H}^{0}(j_{!}j^{\ast}M)\rightarrow M$ et
${}^{{\rm p}}{\cal H}^{0}(j_{!}j^{\ast}N)\rightarrow N$ sont
surjectives et que les fl\`{e}ches de faisceaux pervers $M\rightarrow
{}^{{\rm p}}{\cal H}^{0}(j_{\ast}j^{\ast}M)$ et $N\rightarrow {}^{{\rm
p}}{\cal H}^{0}(j_{\ast}j^{\ast}N)$ sont injectives.  Le probl\`{e}me
est de nature g\'{e}om\'{e}trique et on peut donc remplacer $S$ par
son hens\'{e}lis\'{e} strict, ou ce qui revient au m\^{e}me supposer
que $M$ et $N$ qui sont a priori g\'{e}om\'{e}triquement semi-simples,
sont en fait semi-simples.

Pour tout ${\Bbb Q}_{\ell}$-faisceau pervers irr\'{e}ductible $K$
sur $S$ dont le support rencontre l'ouvert $S_{\natural}$, on a bien
s\^{u}r $K=j_{!\ast}j^{\ast}K$.  Il s'agit donc de d\'{e}montrer que
$N$, et donc a fortiori $M$, n'a pas de sous-objet simple $K$
dont le support dans $S$ ne rencontre pas $S_{\natural}$.

Raisonnons par l'absurde en supposant qu'il existe un tel $K$. Soit
$b$ un point g\'{e}om\'{e}trique de $S$ tel que la fibre de $K$ en $b$
soit non nulle. Prenons alors un morphisme $\varphi :V\rightarrow S$
comme dans la proposition. Alors, $\varphi^{\ast}K$ est un facteur
direct de $\varphi^{\ast}N$ et est donc \`{a} cohomologie ordinaire
constante. Par suite le support de $K$ rencontre $S_{\natural}$,
d'o\`{u} une contradiction.
\hfill\hfill$\square$
\vskip 3mm

Rassemblant les r\'{e}sultats de ce corollaire et ceux obtenus sur
l'ouvert $S_{\natural}$, on obtient finalement

\thm TH\'{E}OR\`{E}ME 3.9.3
\enonce
Pour $\kappa :({\Bbb Z}/2{\Bbb Z})^{2}\rightarrow\{\pm 1\}$
endoscopique on a un isomorphisme de faisceaux pervers gradu\'{e}s sur
$S$
$$
\bigoplus_{n}{}^{{\rm p}}{\cal H}^{n}(f_{S,\ast}{\Bbb
Q}_{\ell})_{\kappa}\buildrel\sim\over\longrightarrow
L_{Z'/Z/S}\otimes_{{\Bbb Z}_{S}} \bigoplus_{n}{}^{{\rm p}}{\cal
H}^{n-2r}(g_{S,\ast}{\Bbb Q}_{\ell})_{\kappa}(-r).\eqno{\square}
$$
\endthm

\subsection{3.10}{Comptage de points rationnels et le th\'{e}or\`{e}me
global}

Pour $\kappa :({\Bbb Z}/2{\Bbb Z})^{2}\rightarrow\{\pm 1\}$
endoscopique on vient de terminer la construction d'un isomorphisme de
faisceaux pervers gradu\'{e}s sur $S$
$$
\bigoplus_{n}{}^{{\rm p}}{\cal H}^{n}(f_{S,\ast}{\Bbb
Q}_{\ell})_{\kappa}\buildrel\sim\over\longrightarrow
L_{Z'/Z/S}\otimes_{{\Bbb Z}_{S}}\bigoplus_{n}{}^{{\rm p}}{\cal
H}^{n-2r}(g_{S,\ast}{\Bbb Q}_{\ell})_{\kappa}(-r).\leqno{(\ast )}
$$

Nous allons maintenant d\'{e}duire une cons\'{e}quence num\'{e}rique
de cet isomorphisme.

Soit $\overline{s}$ le point g\'{e}om\'{e}trique de $S$, localis\'{e}
au point ferm\'{e} $s$ de $S$ d\'{e}fini par une cl\^{o}ture
alg\'{e}brique $\overline{k}$ de $k={\Bbb F}_{q}=\kappa (s)$. On note
comme d'habitude $\mathop{\rm Frob}\nolimits_{s}$ l'endomorphisme de
Frobenius g\'{e}om\'{e}trique en $s$.

\thm LEMME 3.10.1
\enonce
On a les deux \'{e}galit\'{e}s
$$
\sum_{m,n}(-1)^{m+n}\mathop{\rm Tr}(\mathop{\rm Frob}\nolimits_{s},
H^{m}(({}^{{\rm p}}{\cal H}^{n}(f_{S,\ast}{\Bbb Q}_{\ell})_{\kappa}
)_{\overline{s}}))=\sum_{n}(-1)^{n}\mathop{\rm Tr}(\mathop{\rm Frob}
\nolimits_{s},H^{n}({\cal M}_{\overline{s}},{\Bbb Q}_{\ell})_{\kappa})
$$
et
$$
\sum_{m,n}(-1)^{m+n}\mathop{\rm Tr}(\mathop{\rm Frob}\nolimits_{s},
H^{m}(({}^{{\rm p}}{\cal H}^{n}(g_{S,\ast}{\Bbb Q}_{\ell})_{\kappa}
)_{\overline{s}}))=\sum_{n}(-1)^{n}\mathop{\rm Tr}(\mathop{\rm Frob}
\nolimits_{s},H^{n}({\cal N}_{\overline{s}},{\Bbb Q}_{\ell})_{\kappa}).
$$
\endthm

\rem D\'{e}monstration
\endrem
Consid\'{e}rons la cat\'{e}gorie ab\'{e}lienne ${\cal A}$ des ${\Bbb
Q}_{\ell}$-espaces vectoriels de dimension finie munis d'une action
continue de ${\rm Gal}(\overline{k}/k)$ et d'une action de ${\Bbb
Z}^{2}$ qui commute \`{a} l'action de ${\rm Gal}(\overline{k}/k)$ et
dont la restriction \`{a} $2{\Bbb Z}^{2}$ est unipotente.
Consid\'{e}rons aussi la cat\'{e}gorie ab\'{e}lienne ${\cal A}^{{\rm
ss}}$ des ${\Bbb Q}_{\ell}$-espaces vectoriels de dimension finie
munis d'une action continue de ${\rm Gal}(\overline{k}/k)$ et d'une
action de $({\Bbb Z}/2{\Bbb Z})^{2}$ qui commute \`{a} l'action de ${\rm
Gal}(\overline{k}/k)$.  La fl\`{e}che
$$
K_{0}({\cal A}^{{\rm ss}})\rightarrow K_{0}({\cal A})
$$
induite par l'inclusion de ${\cal A}^{{\rm ss}}$ dans ${\cal A}$ est
un isomorphisme.

Pour tout caract\`{e}re $\kappa:({\Bbb Z}/2{\Bbb Z})^{2}\rightarrow
\{\pm 1\}\subset {\Bbb Q}_{\ell}^{\times}$ et tout objet $V$ de ${\cal
A}$ on d\'{e}finit
$$
\mathop{\rm Tr}\nolimits_{\kappa}^{{\rm ss}}(\mathop{\rm Frob}
\nolimits_{q},V)\in {\Bbb Q}_{\ell}
$$
comme suit: $V$ d\'{e}finit un \'{e}l\'{e}ment $[V]$ de $K_{0}({\cal
A})$ donc un \'{e}l\'{e}ment $[V]^{{\rm ss}}$ de $K_{0}({\cal A}^{{\rm
ss}})$ et on pose
$$
\mathop{\rm Tr}\nolimits_{\kappa}^{{\rm ss}} (\mathop{\rm Frob}
\nolimits_{q},V)=\mathop{\rm Tr}\nolimits_{\kappa} (\mathop{\rm Frob}
\nolimits_{q},[V]^{{\rm ss}}).
$$

L'objet $C=f_{S,\ast}{\Bbb Q}_{\ell}$ de la cat\'{e}gorie
d\'{e}riv\'{e}e des complexes born\'{e}s de faisceaux $\ell$-adiques
sur $S$, est muni d'une action du groupe ${\Bbb Z}^{2}$ telle
que le sous-groupe $2{\Bbb Z}^{2}\subset{\Bbb Z}^{2}$ agit trivialement sur
$\bigoplus_{n}{}^{{\rm p}}{\cal H}^{n}(C)$.  En consid\'{e}rant les
tronqu\'{e}s successifs de $C$ pour la $t$-structure
interm\'{e}diaire, on voit que $2{\Bbb Z}^{2}$ agit de mani\`{e}re
unipotente sur $C$.

La fibre $C_{\overline{s}}$ de $C$ en $\overline{s}$ d\'{e}finit un
\'{e}l\'{e}ment
$$
[C_{\overline{s}}]:=\sum_{n}(-1)^{n}[H^{n}(C_{\overline{s}})]
\in K_{0}({\cal A})
$$
qui est \'{e}gal \`{a}
$$
[C_{\overline{s}}]=\sum_{n}(-1)^{n}[H^{n}({\cal M}_{\overline{s}},
{\Bbb Q}_{\ell})]
$$
d'apr\`{e}s le th\'{e}or\`{e}me de changement de base pour un morphisme
propre.  De m\^{e}me la fibre $({}^{{\rm p}}{\cal H}^{n}(C))_{\overline{s}}$
en $\overline{s}$  de chaque faisceau de cohomologie perverse
${}^{{\rm p}}{\cal H}^{n}(C)$ d\'{e}finit un
\'{e}l\'{e}ment
$$
[({}^{{\rm p}}{\cal H}^{n}(C))_{\overline{s}}] =
\sum_{m}(-1)^{m}[H^{m}(({}^{{\rm p}}{\cal H}^{n}(C))_{\overline{s}})]
\in K_{0}({\cal A}).
$$

Maintenant, on a l'\'{e}galit\'{e}
$$
[C_{\overline{s}}]=\sum_{n}(-1)^{n}[({}^{{\rm p}}{\cal H}^{n}
(C))_{\overline{s}}]
$$
dans $K_{0}({\cal A})$, d'o\`{u} l'\'{e}galit\'{e}
$$
\sum_{n}(-1)^{n}[({}^{{\rm p}}{\cal H}^{n}
(C))_{\overline{s}}]=
\sum_{n}(-1)^{n}[H^{n}({\cal M}_{\overline{s}},{\Bbb Q}_{\ell})],
$$
toujours dans $K_{0}({\cal A})$, et par cons\'{e}quent l'\'{e}galit\'{e}
$$
\sum_{n}(-1)^{n}\mathop{\rm Tr}\nolimits_{\kappa}^{{\rm ss}}
(\mathop{\rm Frob}\nolimits_{q},({}^{{\rm p}}{\cal H}^{n}
(C))_{\overline{s}})=\sum_{n}(-1)^{n} \mathop{\rm Tr}
\nolimits_{\kappa}^{{\rm ss}}(\mathop{\rm Frob}\nolimits_{q},
H^{n}({\cal M}_{\overline{s}},{\Bbb Q}_{\ell}))
$$
dans ${\Bbb Q}_{\ell}$. On a en fait deux \'{e}l\'{e}ments virtuels
de ${\cal A}^{{\rm ss}}$ qui, comme objets de ${\cal A}$, ont la
m\^{e}me classe dans $K_{0}({\cal A})$. Ils ont donc la m\^{e}me
classe dans $K_{0}({\cal A}^{{\rm ss}})$.

Enfin, puisque $2{\Bbb Z}^{2}$ agit trivialement aussi bien sur les
faisceaux pervers ${}^{{\rm p}}{\cal H}^{n}(C)$ que sur les groupes de
cohomologie ordinaires $H^{n}({\cal M}_{s},{\Bbb Q}_{\ell})$, et
puisque $({\Bbb Z}/2{\Bbb Z})^{2}$ est un groupe fini, de sorte que l'on a
$$
H^{m}(({}^{{\rm p}}{\cal H}^{n}(C)_{\kappa})_{\overline{s}})=
H^{m}(({}^{{\rm p}}{\cal H}^{n}(C))_{\overline{s}})_{\kappa},
$$
on obtient la premi\`{e}re \'{e}galit\'{e} cherch\'{e}e
$$
\sum_{m,n}(-1)^{m+n}\mathop{\rm Tr}(\mathop{\rm Frob}\nolimits_{q},
H^{m}(({}^{{\rm p}}{\cal H}^{n}(C)_{\kappa})_{\overline{s}})=\sum_{n}(-1)^{n}
\mathop{\rm Tr}(\mathop{\rm Frob}\nolimits_{q},H^{n}
({\cal M}_{\overline{s}},{\Bbb Q}_{\ell})_{\kappa}).
$$

La deuxi\`{e}me \'{e}galit\'{e} se d\'{e}montre de fa\c{c}on identique.
\hfill\hfill$\square$
\vskip 3mm

Compte tenu de ce lemme, on d\'{e}duit de l'isomorphisme $(\ast )$ la
proposition suivante.

\thm PROPOSITION 3.10.2
\enonce
Pour $\kappa :({\Bbb Z}/2{\Bbb Z})^{2}\rightarrow\{\pm 1\}$
endoscopique on a l'\'{e}galit\'{e}
$$
\sum_{n}(-1)^{n}\mathop{\rm Tr}(\mathop{\rm Frob}\nolimits_{s},
H^{n}({\cal M}_{\overline{s}},{\Bbb
Q}_{\ell})_{\kappa})=(-1)^{m_{Z'/Z/S}}q^{r}
\sum_{n}(-1)^{n}\mathop{\rm Tr}(\mathop{\rm Frob}\nolimits_{s},
H^{n}({\cal N}_{\overline{s}},{\Bbb Q}_{\ell})_{\kappa})
$$
o\`{u} $m_{Z'/Z/S}$ est l'entier de la section {\rm (3.4)}.
\hfill\hfill$\square$
\endthm

Nous allons maintenant r\'{e}crire cette \'{e}galit\'{e} \`{a} l'aide
de la formule des points fixes de la section (A.3), ou plut\^{o}t sa
variante champ\^{e}tre appliqu\'{e}e aux $k$-champs ${\cal M}_{a}$
muni de l'action du $k$-champ de Picard $P_{a}$, et ${\cal N}_{a}$
muni de l'action du $k$-champ de Picard $P_{H,a}$.

Comme dor\'{e}navant nous ne consid\'{e}rons plus d'autres fibres des
morphismes ${\cal M}\rightarrow {\Bbb A}$ et ${\cal N}\rightarrow
{\Bbb A}_{H}$ que celles en $a$, on peut all\'{e}ger les
notations en supprimant l'indice $a$.  On a donc une courbe
spectrale $Y\subset \Sigma$, not\'{e}e pr\'{e}c\'{e}demment
$Y_{a}$, r\'{e}union de composantes g\'{e}om\'{e}triquement
irr\'{e}ductibles $Y_{1}$ et $Y_{2}$.  On a de plus le
rev\^{e}tement double \'{e}tale $\pi_{Y}:Y'=X'\times_{X}Y\rightarrow
Y$ de composantes (g\'{e}om\'{e}triquement) irr\'{e}ductibles les
$Y_{\alpha}'=X'\times_{X}Y_{\alpha}$ pour $\alpha =1,2$.

La cat\'{e}gorie ${\cal M}(\overline{k})$ a pour objets les couples
$({\cal F},\iota_{{\cal F}})$ o\`{u} ${\cal F}$ est un ${\cal
O}_{\overline{k}\otimes_{k}Y'}$-Module coh\'{e}rent sans torsion de
rang $1$ et o\`{u} $\iota_{{\cal F}}:\tau^{\ast}{\cal
F}\buildrel\sim\over\longrightarrow {\cal F}^{\vee}$ est une structure
unitaire sur ${\cal F}$, les morphismes \'{e}tant les isomorphismes de
${\cal O}_{\overline{k}\otimes_{k}Y'}$-Modules qui respectent les
structures unitaires.  La cat\'{e}gorie de Picard $P(\overline{k})$ a
pour objet les couples $({\cal G},\iota_{{\cal G}})$ o\`{u} maintenant
${\cal G}$ est un ${\cal O}_{\overline{k}\otimes_{k}Y'}$-Module
inversible et o\`{u} $\iota :\tau^{\ast}{\cal G}
\buildrel\sim\over\longrightarrow {\cal G}^{\otimes -1}$ est une
structure unitaire sur ${\cal G}$, et $P(\overline{k})$ agit sur
${\cal M}(\overline{k})$ par produit tensoriel.

La cat\'{e}gorie ${\cal N}(\overline{k})$ est quant \`{a} elle la
sous-cat\'{e}gorie pleine de ${\cal M}(\overline{k})$ dont les objets
sont les couples $({\cal F},\iota_{{\cal F}})$ qui sont
d\'{e}compos\'{e}s relativement au d\'{e}coupage
$\overline{k}\otimes_{k}Y'= \overline{k}\otimes_{k} Y_{1}'\cup
\overline{k}\otimes_{k} Y_{2}'$.  On a donc
$$
{\cal N}(\overline{k})={\cal M}_{1}(\overline{k})\times {\cal
M}_{2}(\overline{k})
$$
o\`{u} la cat\'{e}gorie ${\cal M}_{\alpha}(\overline{k})$ est
d\'{e}finie comme ${\cal M}(\overline{k})$ mais apr\`{e}s avoir
remplac\'{e} $Y'$ par $Y_{\alpha}'$.  Si on pose $Q=
P_{H,a}$ pour all\'{e}ger les notations, la cat\'{e}gorie de
Picard $Q(\overline{k})=P_{1}(\overline{k})\times P_{2}
(\overline{k})$ o\`{u} $P_{\alpha} (\overline{k})$ est la
cat\'{e}gorie de Picard d\'{e}finie comme $P(\overline{k})$ mais
apr\`{e}s avoir remplac\'{e} $Y'$ par $Y_{\alpha}'$, agit sur
${\cal N}(\overline{k})$ par produit tensoriel.

Le point $a$ de ${\Bbb A}_{H}(k)$ \'{e}tant choisi comme en (3.1), on
sait d'apr\`{e}s le lemme 2.8.2, que $\pi_{0}(P)$ est le
$k$-sch\'{e}ma en groupes constant de valeur le groupe fini $({\Bbb
Z}/2{\Bbb Z})^{2}$.  On a de plus la suite exacte
$$
0\rightarrow P^{0}(k)\rightarrow P(k)\rightarrow ({\Bbb Z}/2{\Bbb Z})^{2}
\rightarrow 0.
$$

Rempla\c{c}ant successivement $Y'$ par $Y_{1}'$ et par
$Y_{2}'$, on voit que
$$
\pi_{0}(Q)(k)=\pi_{0}(P_{1})(k) \times \pi_{0}(P_{2})(k)
=({\Bbb Z}/2{\Bbb Z})^{2}
$$
et que l'on a aussi la suite exacte
$$
0\rightarrow Q^{0}(k)\rightarrow Q(k)\rightarrow \pi_{0}(Q)(k)
\rightarrow 0.
$$

La cat\'{e}gorie quotient $[{\cal M}/P](\overline{k})$ a pour objets
les $({\cal F},\iota_{{\cal F}})\in \mathop{\rm ob}{\cal
M}(\overline{k})$ et pour fl\`{e}ches de $({\cal F},\iota_{{\cal F}})$
vers un autre objet $({\cal F}',\iota_{{\cal F}'})$ de $[{\cal
M}/P](\overline{k})$, les classes d'isomorphie d'objets $({\cal
G},\iota_{{\cal G}})$ de $P(\overline{k})$ tels que $({\cal
G}\otimes_{{\cal O}_{\overline{k}\otimes_{k}Y'}}{\cal F}, \iota_{{\cal
G}}\otimes\iota_{{\cal F}})$ soit isomorphe \`{a} $({\cal
F}',\iota_{{\cal F}'})$ dans ${\cal M}(\overline{k})$.  On a un
foncteur de Frobenius naturel
$$
\mathop{\rm Frob}\nolimits_{q}:[{\cal M}/P](\overline{k})
\rightarrow [{\cal M}/P](\overline{k})
$$
\`{a} l'aide duquel on peut retrouver la cat\'{e}gorie $[{\cal
M}/P](k)$ comme dans la section (A.3).  Comme dans loc.
cit.  on note $[{\cal M}/P](k)_{\sharp}$ l'ensemble des classes
d'isomorphie des objets de $[{\cal M}/P](k)$ et on d\'{e}finit une
application
$$
\mathop{\rm cl}\nolimits_{{\cal M}}:[{\cal M}/P](k)_{\sharp}
\rightarrow\pi_{0}(P)=({\Bbb Z}/2{\Bbb Z})^{2}.
$$

La cat\'{e}gorie quotient $[{\cal N}/Q](\overline{k})$ est le produit
de cat\'{e}gories
$$
[{\cal N}/Q](\overline{k})=[{\cal M}_{1}/P_{1}](\overline{k})
\times [{\cal M}_{2}/P_{2}](\overline{k}).
$$
On a un foncteur de Frobenius naturel
$$
\mathop{\rm Frob}\nolimits_{q}:[{\cal N}/Q](\overline{k})
\rightarrow [{\cal N}/Q](\overline{k})
$$
produit des foncteurs de Frobenius pour $[{\cal M}_{1}/P_{1}]
(\overline{k})$ et $[{\cal M}_{2}/P_{2}](\overline{k})$, \`{a}
l'aide duquel on retrouve la cat\'{e}gorie $[{\cal N}/Q](k)=[{\cal
M}_{1}/P_{1}](k) \times [{\cal M}_{2}/P_{2}](k)$.  On
note $[{\cal N}/Q](k)_{\sharp}=[{\cal M}_{1}/P_{1}]
(k)_{\sharp} \times [{\cal M}_{2}/P_{2}](k)_{\sharp}$
l'ensemble des classes d'isomorphie des objets de $[{\cal N}/Q](k)$ et
on a alors une application
$$
\mathop{\rm cl}\nolimits_{{\cal N}}:[{\cal N}/Q](k)_{\sharp}
\rightarrow\pi_{0}(Q)=({\Bbb Z}/2{\Bbb Z})^{2}
$$
qui n'est autre que $\mathop{\rm cl}\nolimits_{{\cal N}}= \mathop{\rm
cl}\nolimits_{{\cal M}_{1}}\times \mathop{\rm cl}\nolimits_{{\cal
M}_{2}}$ o\`{u} $\mathop{\rm cl}\nolimits_{{\cal M}_{\alpha}}:
[{\cal M}_{\alpha}/P_{\alpha}](k)_{\sharp}
\rightarrow\pi_{0}(P_{\alpha})={\Bbb Z}/2{\Bbb Z}$
est d\'{e}finie comme $\mathop{\rm cl}\nolimits_{{\cal M}}$ mais
apr\`{e}s avoir remplac\'{e} $Y'$ par $Y_{\alpha}'$.

\thm LEMME 3.10.3
\enonce
Pour chacun des deux caract\`{e}res endoscopiques $\kappa$ on a
$\kappa\circ\mathop{\rm cl}\nolimits_{{\cal N}}\equiv 1$.
\endthm

\rem D\'{e}monstration
\endrem
Si on a un isomorphisme $\mathop{\rm Frob}\nolimits_{q}^{\ast}({\cal
F}_{\alpha},\iota_{{\cal F}_{\alpha}})\cong ({\cal F}_{\alpha},
\iota_{{\cal F}_{\alpha}})\otimes ({\cal G}_{\alpha}, \iota_{{\cal
G}_{\alpha}})$ o\`{u} $({\cal F}_{\alpha}, \iota_{{\cal F}_{\alpha}})
\in \mathop{\rm ob}{\cal M}_{\alpha}(\overline{k})$ et $({\cal
G}_{\alpha}, \iota_{{\cal G}_{\alpha}})\in P_{\alpha}(\overline{k})$,
alors la composante connexe $\lambda_{\alpha}\in\pi_{0}(P_{\alpha})
(\overline{k})={\Bbb Z}/2{\Bbb Z}$ contenant ${\cal G}_{\alpha}$ est
n\'{e}cessairement nulle.

En effet, d'apr\`{e}s le lemme 2.8.3, il existe une
application de $\pi_{0}({\cal M}_{\alpha})(\overline{k})$ dans ${\Bbb
Z}/2{\Bbb Z}$ \'{e}quivariante par rapport \`{a} l'action \'{e}vidente
de $\pi_{0}(P_{\alpha})(\overline{k})$ sur $\pi_{0}({\cal
M}_{\alpha})(\overline{k})$ et \`{a} l'action de
$\pi_{0}(P_{\alpha})(\overline{k})={\Bbb Z}/2{\Bbb Z}$ sur
lui-m\^{e}me par translation.

Or $\mathop{\rm Frob}\nolimits_{k}^{\ast}({\cal F}_{\alpha},
\iota_{{\cal F}_{\alpha}})$ et $({\cal F}_{\alpha}, \iota_{{\cal
F}_{\alpha}})$ sont dans la m\^{e}me composante connexe de ${\cal
M}_{\alpha}$, d'o\`{u} l'assertion.
\hfill\hfill$\square$
\vskip 3mm

Appliquant alors la proposition A.3.1 on obtient les formules de
points fixes suivantes.

\thm PROPOSITION 3.10.4
\enonce
On a les \'{e}galit\'{e}s
$$
\sum_{n}(-1)^{n}\mathop{\rm Tr}(\mathop{\rm Frob}
\nolimits_{q},H^{n}(\overline{k}\otimes_{k}{\cal M}, {\Bbb
Q}_{\ell})_{\kappa})
=|P^{0}(k)|\sum_{m\in [{\cal M}/P](k)_{\sharp}}{\kappa (\mathop{\rm
cl}\nolimits_{{\cal M}}(m))\over |\mathop{\rm Aut}(m)|}
$$
et
$$
\qquad\sum_{n}(-1)^{n}\mathop{\rm Tr}(\mathop{\rm Frob}
\nolimits_{q},H^{n}(\overline{k}\otimes_{k}{\cal N}, {\Bbb
Q}_{\ell})_{\kappa})
=|Q^{0}(k)|\sum_{n\in [{\cal N}/Q](k)_{\sharp}}{1\over
|\mathop{\rm Aut}(n)|}.
$$
\hfill\hfill$\square$
\endthm

\rem Remarque
\endrem
Les champs ${\cal M}$ et $P$ sont de Deligne-Mumford et n'importe quel
objet sur $\overline{k}$ de l'un de ces champs a pour seuls
automorphismes 
$$
\{\pm 1\}=\overline{k}^{\tau^{\ast}=(-)^{-1}}=
(H^{0}(\overline{k}\otimes_{k}Y', {\cal O}_{\overline{k}
\otimes_{k}Y'}^{\times}))^{\tau^{\ast} =(-)^{-1}}.
$$
Ces automorphismes sont donc les m\^{e}mes pour ${\cal M}$ et $P$.
Ils disparaissent donc dans le champ quotient $[{\cal M}/P]$ et dans la
suite on peut donc les n\'{e}gliger.  Le m\^{e}me raisonnement vaut pour
l'autre quotient $[{\cal N}/Q]$.
\hfill\hfill$\square$
\vskip 3mm

Nous allons finalement {\og}{localiser}{\fg} les seconds membres des
formules des points fixes ci-dessus.  Pour cela nous aurons besoin
d'un point base de ${\cal M}(k)$ et d'un point base de ${\cal N}(k)$.
\vskip 2mm

Commen\c{c}ons par pr\'{e}ciser ses points base.  \`{A} la fin de la
section (2.6) nous avons construit une section de Kostant de la
fibration de Hitchin au-dessus de ${\Bbb A}^{{\rm red}}$.  C'est un
${\cal O}_{Y'}$-Module inversible racine carr\'{e}e de
$\omega_{Y'/X'}$ muni d'une structure unitaire.  Rempla\c{c}ant
successivement $Y'$ par $Y_{1}'$ et par $Y_{2}'$ nous obtenons
donc un objet de Kostant $({\cal K},\iota_{{\cal K}})=(({\cal
K}_{1},\iota_{{\cal K}_{1}}), ({\cal K}_{2},\iota_{{\cal
K}_{2}}))$ de ${\cal N}(k)$.  Par image directe par la
normalisation partielle $Y_{1}'\amalg Y_{2}'\rightarrow Y'$,
c'est-\`{a}-dire par image par l'immersion ferm\'{e}e $i:{\cal
N}\hookrightarrow {\cal M}$, cet objet de Kostant de ${\cal N}$
d\'{e}finit un objet qui est diff\'{e}rent de l'objet de Kostant de
${\cal M}(k)$.  Remarquons en particulier que l'image directe de
${\cal K}$ est sans torsion de rang $1$ sur $Y'$ mais n'est pas
inversible comme l'est l'objet de Kostant.

Nous prendrons pour point base de ${\cal N}(k)$ l'objet de Kostant
$({\cal K},\iota_{{\cal K}})$ et pour point base de ${\cal M}(k)$
l'image de cet objet de Kostant de ${\cal N}(k)$ par $i$, objet que
nous noterons aussi $({\cal K},\iota_{{\cal K}})$.

Notons au passage que le signe $(-1)^{r}$ qui appara\^{\i}t dans
l'\'{e}nonc\'{e} du th\'{e}or\`{e}me 1.5.1 provient bien
du fait que nous avons pris comme point base l'objet de Kostant de
${\cal N}$ et non l'objet de Kostant de ${\cal M}$ comme dans [Kot~2].
\vskip 2mm

Proc\'{e}dons maintenant \`{a} la localisation pour ${\cal M}$.  Pour
chaque point ferm\'{e} $y$ de $Y$, on a un analogue local de la
cat\'{e}gorie $[{\cal M}/P](\overline{k})$.  Soit $A_{y}$ le
compl\'{e}t\'{e} de l'anneau semi-local de $Y'$ le long de
$\pi_{Y}^{-1}(y)$ et $\mathop{\rm Frac}(A_{y})$ l'anneau total des
fractions de $A_{y}$.  L'involution $\tau$ de $Y'$ induit une
involution not\'{e}e encore $\tau$ sur $A_{y}$.

Pour chaque point ferm\'{e} $y$ de $Y$, notons $K_{y}$ le
compl\'{e}t\'{e} de ${\cal K}$ le long de $\pi_{Y}^{-1}(y)$, de sorte
que $K_{y}$ est un $A_{y}$-module sans torsion de rang $1$ et que l'on
a une structure unitaire $\iota_{K_{y}}:\tau^{\ast}K_{y}\buildrel\sim
\over\longrightarrow K_{y}^{\vee}=\omega_{X/Y,y} \otimes_{{\cal
O}_{Y,y}}K_{y}^{-1}$.  On en d\'{e}duit un $\mathop{\rm Frac}
(A_{y})$-module
$$
V_{y}:=\mathop{\rm Frac}(A_{y})\otimes_{A_{y}}K_{y}
$$
de dimension $1$ et une structure unitaire
$$
\iota_{V_{y}}:\tau^{\ast}V_{y}\buildrel\sim\over\longrightarrow
V_{y}^{\vee}=\mathop{\rm Hom}\nolimits_{\mathop{\rm Frac}
(A_{y})}(V_{y},\mathop{\rm Frac}(A_{y}))
$$
sur cette droite vectorielle.  Ici le produit tensoriel par
$\omega_{X/Y,y}$ est inutile puisque $\omega_{X/Y}$ est \`{a} support
dans le ferm\'{e} fini de $Y$ o\`{u} le rev\^{e}tement fini
g\'{e}n\'{e}riquement \'{e}tale $Y\rightarrow X$ est ramifi\'{e}.

Soit $\overline{A}_{y}$ le compl\'{e}t\'{e} de la
$\overline{k}$-alg\`{e}bre semi-locale $\overline{k}
\otimes_{k}A_{y}$.  Pour tout $A_{y}$-module $M_{y}$ on note
$\overline{M}_{y}=\overline{A}_{y}\otimes_{A_{y}}M_{y}$.
Consid\'{e}rons l'ensemble de sous-$\overline{A}_{y}$-modules de
$\overline{V}_{y}$
$$
{\cal M}(y)(\overline{k})=\{{\cal V}_{y}\subset \overline{V}_{y}\mid
\mathop{\rm Frac}(\overline{A}_{y}){\cal V}_{y}=\overline{V}_{y}\hbox{
et } \iota_{\overline{V}_{y}}(\tau^{\ast}{\cal V}_{y})={\cal
V}_{y}^{\vee}\}
$$
o\`{u} $\iota_{\overline{V}_{y}}$ est induit par $\iota_{V_{y}}$ et
o\`{u} ${\cal V}_{y}^{\vee}=\omega_{X/Y,y}\otimes_{{\cal O}_{Y,y}}
{\cal V}_{y}^{-1}$ et
$$
{\cal V}_{y}^{-1}=\{\nu\in \mathop{\rm Hom}\nolimits_{\mathop{\rm
Frac}(\overline{A}_{y})}(\overline{V}_{y}, \mathop{\rm Frac}
(\overline{A}_{y})) \mid \nu ({\cal V}_{y})\subset \overline{K}_{y}\}.
$$
Sur cet ensemble on a l'action par homoth\'{e}ties du groupe
$$
P(y)(\overline{k})=\{a\in\mathop{\rm Frac}
(\overline{A}_{y})^{\times}/\overline{A}_{y}^{\times}\mid
\tau (a)=a^{-1}\}.
$$

Consid\'{e}rons alors la cat\'{e}gorie quotient $[{\cal M}(y)/P(y)]
(\overline{k})$ de ${\cal M}(y)(\overline{k})$ pour cette action de
$P(y)(\overline{k})$.  Si $y$ est un point lisse de $Y$,
$P(y)(\overline{k})$ agit librement et transitivement sur ${\cal
M}(y)(\overline{k})$, de sorte que la cat\'{e}gorie $[{\cal M}(y)/P(y)]
(\overline{k})$ n'a qu'un seul objet \`{a} isomorphisme pr\`{e}s et
cet objet n'a pas d'automorphisme non trivial.  Le produit de
cat\'{e}gories
$$
\prod_{y\in Y}[{\cal M}(y)/P(y)](\overline{k})
=\prod_{y\in Y^{{\rm sing}}}[{\cal M}(y)/P(y)](\overline{k})
$$
o\`{u} $Y^{{\rm sing}}$ est ferm\'{e} des points singuliers de $Y$ et
o\`{u} les $y$ sont des points ferm\'{e}s, est donc un produit fini.

On construit un foncteur
$$
\prod_{y\in Y^{{\rm sing}}}{\cal M}(y)(\overline{k})\rightarrow {\cal
M} (\overline{k})
$$
comme suit.  \`{A} toute collection de r\'{e}seaux $({\cal
V}_{y})_{y}$ on associe le ${\cal O}_{\overline{k}
\otimes_{k}Y'}$-Module coh\'{e}rent sans torsion de rang $1$ obtenu en
recollant la restriction du Module ${\cal K}$ \`{a} l'ouvert de
lissit\'{e} de $\overline{k}\otimes_{k}Y'$ et les ${\cal V}_{y}\subset
\overline{V}_{y}$.  On munit ce ${\cal O}_{\overline{k}
\otimes_{k}Y'}$-Module de la structure unitaire qui est la structure
$\iota_{{\cal K}}$ sur l'ouvert de lissit\'{e} de
$\overline{k}\otimes_{k}Y'$ et celle qui est donn\'{e}e sur chaque
${\cal V}_{y}$.  De mani\`{e}re similaire on a un morphisme de
cat\'{e}gories de Picard
$$
\prod_{y\in Y^{{\rm sing}}}P(y)(\overline{k})\rightarrow P
(\overline{k})
$$
qui, \`{a} toute collection de scalaires $(a_{y})_{y}$ associe le
${\cal O}_{\overline{k}\otimes_{k}Y'}$-Module inversible obtenu en
recollant le Module trivial sur l'ouvert de lissit\'{e} de
$\overline{k}\otimes_{k}Y'$ et les r\'{e}seaux $a_{y}
\overline{A}_{y}\subset \mathop{\rm Frac}(\overline{A}_{y})$, ce
recollement \'{e}tant muni de la structure unitaire \'{e}vidente.

Le morphisme entre les ${\cal M}$ est \'{e}quivariant relativement au
morphisme entre les $P$ et induit par passage au quotient un foncteur
$$
\prod_{y\in Y^{{\rm sing}}}[{\cal M}(y)/P(y)](\overline{k})\rightarrow
[{\cal M}/P](\overline{k}).
$$

\thm LEMME 3.10.5
\enonce
Ce foncteur $\prod_{y\in Y^{{\rm sing}}}[{\cal M}(y)/P(y)](\overline{k})
\rightarrow [{\cal M}/P](\overline{k})$ ci-dessus est une
\'{e}qui\-va\-len\-ce de cat\'{e}gories.
\endthm

\rem D\'{e}monstration
\endrem
Montrons tout d'abord que le foncteur est essentiellement surjectif.

Soit $({\cal F},\iota_{{\cal F}})$ un objet de ${\cal M}
(\overline{k})$.  La restriction de ${\cal F}$ \`{a} l'ouvert de
lissit\'{e} de $\overline{k}\otimes_{k}Y'$ est inversible.  On peut
donc trouver un Module inversible unitaire sur
$\overline{k}\otimes_{k}Y'$ qui prolonge la restriction de $({\cal
F},\iota_{{\cal F}})$ \`{a} cet ouvert.  En effet, pour avoir un tel
prolongement, il suffit de choisir, pour chaque $y\in Y^{{\rm sing}}$,
un vecteur de base unitaire $e_{y}$ de la $\mathop{\rm
Frac}(\overline{A}_{y})$-droite vectorielle $\mathop{\rm
Frac}(\overline{A}_{y}) \otimes_{\overline{A}_{y}}{\cal V}_{y}$ et de
recoller la restriction de $({\cal F},\iota_{{\cal F}})$ \`{a}
l'ouvert de lissit\'{e} et les sous-$\overline{A}_{y}$-modules
inversibles unitaires $\overline{A}_{y}e_{y}\subset \mathop{\rm
Frac}(\overline{A}_{y})\otimes_{\overline{A}_{y}}{\cal V}_{y}$.

Par suite, quitte \`{a} remplacer notre objet de d\'{e}part par un
objet qui lui est isomorphe dans $[{\cal M}/P](\overline{k})$, on peut
supposer que $({\cal F},\iota_{{\cal F}})$ et l'objet $({\cal
K},\iota_{{\cal K}})$ ont la m\^{e}me restriction \`{a} l'ouvert de
lissit\'{e} de $\overline{k}\otimes_{k}Y'$.  La surjectivit\'{e} est
maintenant triviale puisque tout Module sur
$\overline{k}\otimes_{k}Y'$ peut s'obtenir par recollement de Modules
sur l'ouvert de lissit\'{e} et de $\overline{A}_{y}$-modules.

Pour terminer la preuve du lemme, montrons que le foncteur est
pleinement fid\`{e}le.

Soient $({\cal V}_{y})_{y\in Y^{{\rm sing}}}$ et $({\cal V}_{y}')_{y\in Y^{{\rm
sing}}}$ deux objets de $\prod_{y\in Y^{{\rm sing}}}{\cal
M}(y)(\overline{k})$, d'images $({\cal F},\iota_{{\cal F}})$ et
$({\cal F}',\iota_{{\cal F}'})$ par le foncteur de recollement.  Il
s'agit de voir qu'il revient au m\^{e}me de se donner un
\'{e}l\'{e}ment $(a_{y})_{y\in Y^{{\rm sing}}}\in
P(y)(\overline{k})$ tel que $a_{y}{\cal V}_{y}={\cal V}_{y}'$ pour chaque $y$,
ou de se donner un objet $({\cal G},\iota_{{\cal G}})$ de
$P(\overline{k})$ et un isomorphisme $({\cal G},\iota_{{\cal
G}})\otimes ({\cal F},\iota_{{\cal F}})\buildrel\sim\over
\longrightarrow ({\cal F}',\iota_{{\cal F}'})$.  Or, par construction
les restrictions de $({\cal F},\iota_{{\cal F}})$ et $({\cal
F}',\iota_{{\cal F}'})$ \`{a} l'ouvert de lissit\'{e} de
$\overline{k}\otimes_{k}Y'$ sont toutes les deux isomorphes \`{a} la
restriction de l'objet $({\cal K},\iota_{{\cal K}})$, de sorte que la
restriction de $({\cal G},\iota_{{\cal G}})$ \`{a} cet ouvert de
lissit\'{e} est n\'{e}cessairement trivialis\'{e}e.  L'objet $({\cal
G},\iota_{{\cal G}})$ s'obtient donc de mani\`{e}re unique par
recollement de sous-$\overline{A}_{y}$-modules inversibles unitaires
$\overline{A}_{y}a_{y}\subset\mathop{\rm Frac}
(\overline{A}_{y})$
\hfill\hfill$\square$
\vskip 3mm

Pour chaque $y\in Y^{{\rm sing}}$ on a un foncteur de Frobenius naturel
$$
\mathop{\rm Frob}\nolimits_{q}=[{\cal M}(y)/P(y)](\overline{k})
\rightarrow [{\cal M}(y)/P(y)](\overline{k})
$$
\`{a} l'aide duquel ont peut d\'{e}finir une cat\'{e}gorie $[{\cal
M}(y)/P(y)](k)$ comme dans la section (A.3).  Soit
$[{\cal M}(y)/P(y)](k)_{\sharp}$ l'ensemble des classes d'isomorphie
d'objets de cette derni\`{e}re cat\'{e}gorie.

Pour chaque $y\in Y$, on a un plongement
$$
[{\cal M}(y)/P(y)](k)_{\sharp}\hookrightarrow
\prod_{y'\in Y^{{\rm sing}}}[{\cal M}(y')/P(y')](k)_{\sharp}
$$
dont la composante en $y$ est l'identit\'{e} et la composante en
n'importe quel $y'\not=y$ est constante de valeur la section
$K_{y'}\subset V_{y'}$.  On note
$$
\mathop{\rm cl}\nolimits_{{\cal M}(y)}:[{\cal M}(y)/P(y)](k)_{\sharp}
\rightarrow ({\Bbb Z}/2{\Bbb Z})^{2}
$$
la restriction de l'application $\mathop{\rm cl}\nolimits_{{\cal M}}$
via ce plongement.  L'application compos\'{e}e
$$
\prod_{y\in Y^{{\rm sing}}}[{\cal M}(y)/P(y)](k)_{\sharp}
\buildrel\sim\over\longrightarrow [{\cal M}/P](k)_{\sharp}
\,\smash{\mathop{\hbox to 8mm{\rightarrowfill}}
\limits^{\scriptstyle \mathop{\rm cl}\nolimits_{{\cal M}}}}\,
({\Bbb Z}/2{\Bbb Z})^{2}
$$
est alors la somme des applications $\mathop{\rm cl}\nolimits_{{\cal
M}(y)}$ pour $y\in Y^{{\rm sing}}$.

On a donc
$$
\sum_{m\in [{\cal M}/P](k)_{\sharp}}{\kappa (\mathop{\rm cl}
\nolimits_{{\cal M}} (m))\over |\mathop{\rm Aut}(m)|}=
\prod_{y\in Y^{{\rm sing}}}\left(\sum_{m_{y}\in 
[{\cal M}(y)/P(y)](k)_{\sharp}} {\kappa (\mathop{\rm cl}
\nolimits_{{\cal M}(y)}(m_{y}))\over |\mathop{\rm
Aut}(m_{y})|}\right).\leqno{(\ast\ast )_{{\cal M}}}
$$

En rempla\c{c}ant successivement $Y'$ par $Y_{1}'$ et par
$Y_{2}'$ dans les d\'{e}finitions et \'{e}nonc\'{e}s ci-dessus et
en faisant le produit des r\'{e}sultats obtenus, on obtient une
\'{e}quivalence de cat\'{e}gories
$$
\prod_{y\in Y^{{\rm sing}}}[{\cal N}(y)/Q(y)](k)_{\sharp}
\buildrel\sim\over\longrightarrow [{\cal N}/Q](k)_{\sharp},
$$
et donc l'\'{e}galit\'{e}
$$
\sum_{n\in [{\cal N}/Q] (k)_{\sharp}}{1\over |\mathop{\rm Aut}(n)|}=
\prod_{y\in Y^{{\rm sing}}} \left(\sum_{n_{y}\in [{\cal N}(y)/Q(y)]
(k)_{\sharp}}{1\over |\mathop{\rm
Aut}(n_{y})|}\right).\leqno{(\ast\ast )_{{\cal N}}}
$$

Notons $Z^{{\rm inerte}}$ l'ensemble des points ferm\'{e}s de
$Z=Y_{1}\cap Y_{2}$ qui sont inertes dans le rev\^{e}tement
double \'{e}tale $Z'\rightarrow Z$. Pour chaque $z\in Z^{{\rm
inerte}}$ notons $r_{z}$ la longueur de la $\kappa (z)$-alg\`{e}bre
artinienne ${\cal O}_{Z,z}$ (anneau local de $Z$ en $z$).

Notre th\'{e}or\`{e}me global est alors le suivant.

\thm TH\'{E}OR\`{E}ME 3.10.6
\enonce
Pour chacun des deux caract\`{e}res endoscopiques $\kappa :({\Bbb
Z}/2{\Bbb Z})^{2}\rightarrow\{\pm 1\}$ on a l'\'{e}galit\'{e}
$$\displaylines{
\quad\prod_{z\in Z^{{\rm inerte}}}\left|
{(\widetilde{A}_{z}^{\times})^{\tau^{\ast}=(-)^{-1}}
\over (A_{z}^{\times})^{\tau^{\ast}=(-)^{-1}}}\right|
\left(\sum_{m_{z}\in [{\cal M}(z)/P(z)](k)_{\sharp}}{\kappa
(\mathop{\rm cl}\nolimits_{{\cal M}(z)}(m_{z}))\over
|\mathop{\rm Aut}(m_{z})|}\right)
\hfill\cr\hfill
=\prod_{z\in Z^{{\rm inerte}}}(-|\kappa (z)|)^{r_{z}}
\left|{(\widetilde{A}_{z}^{\times})^{\tau^{\ast} =(-)^{-1}}
\over (A_{1,z}^{\times}\times A_{2,z}^{\times})^{\tau^{\ast}
=(-)^{-1}}}\right|\left(\sum_{n_{z}\in [{\cal N}(z)/Q(z)] (k)_{\sharp}}
{1\over |\mathop{\rm Aut}(n_{z})|}\right).\quad}
$$
\endthm

\rem D\'{e}monstration
\endrem
Commen\c{c}ons pas trois remarques.
\vskip 1mm

\itemitem{1)} Comme on a la suite exacte de cat\'{e}gories de Picard
$$
1\rightarrow\prod_{y}{(A_{1,y}^{\times}
\times A_{2,y}^{\times})^{\tau^{\ast} =(-)^{-1}}
\over (A_{y}^{\times})^{\tau^{\ast} =(-)^{-1}}}
\rightarrow P^{0}(k)
\rightarrow Q^{0}(k)\rightarrow 0,
$$
on a l'\'{e}galit\'{e}
$$
{|Q^{0}(k)|\over |P^{0}(k)|}=\prod_{y}\left|{(A_{1,y}^{\times}
\times A_{2,y}^{\times})^{\tau^{\ast} =(-)^{-1}}
\over (A_{y}^{\times})^{\tau^{\ast} =(-)^{-1}}}\right|.
$$
Pour chaque point ferm\'{e} $y$ de $Y^{{\rm sing}}$ notons
$\widetilde{A}_{y}$ le normalis\'{e} de l'anneau $A_{y}$ et
$\widetilde{A}_{\alpha ,y}$ celui de $A_{\alpha ,y}$ pour
$\alpha =1,2$, o\`{u} $A_{\alpha ,y}$ est d\'{e}fini comme $A_{y}$
mais apr\`{e}s avoir remplac\'{e} $Y'$ par $Y_{\alpha}'$.  On a
$$
A_{y}\subset A_{1,y}\times A_{2,y}\subset
\widetilde{A}_{1,y}\times \widetilde{A}_{2,y}=
\widetilde{A}_{y}.
$$
\`{A} l'aide ce ces notations on peut r\'{e}crire chaque facteur du
second membre de l'\'{e}galit\'{e} ci-dessus sous la forme
$$
\left|{(A_{1,y}^{\times}
\times A_{2,y}^{\times})^{\tau^{\ast} =(-)^{-1}}
\over (A_{y}^{\times})^{\tau^{\ast} =(-)^{-1}}}\right|
=\left|{(\widetilde{A}_{y}^{\times})^{\tau^{\ast}=(-)^{-1}}
\over (A_{y}^{\times})^{\tau^{\ast}=(-)^{-1}}}\right|
\left|{(\widetilde{A}_{y}^{\times})^{\tau^{\ast} =(-)^{-1}}
\over (A_{1,y}^{\times}\times A_{2,y}^{\times})^{\tau^{\ast}
=(-)^{-1}}}\right|^{-1}.
$$

\vskip 1mm

\itemitem{2)} L'entier $m_{Z'/Z/S}$ dans la proposition 3.10.2 est par
d\'{e}finition (cf. la section (3.4)) \'{e}gal \`{a}
$$
m_{Z'/Z/S}=\sum_{z\in Z^{{\rm inerte}}}r_{z}.
$$

\vskip 1mm

\itemitem{3)} L'entier $r_{z}$ a un sens pour tout point ferm\'{e} $z$ de
$Z$ et on a \'{e}videmment
$$
r=\sum_{z\in Z}r_{z}.
$$
\vskip 1mm

Compte tenu de ces remarques, des deux derni\`{e}res propositions et
des formules $(\ast\ast )_{{\cal M}}$ et $(\ast\ast )_{{\cal N}}$, il
ne reste plus qu'\`{a} v\'{e}rifier que
$$\displaylines{
\qquad\left|{(\widetilde{A}_{y}^{\times})^{\tau^{\ast}=(-)^{-1}}
\over (A_{y}^{\times})^{\tau^{\ast}=(-)^{-1}}}\right|
\left(\sum_{m_{y}\in
[{\cal M}(y)/P(y)](k)_{\sharp}}{\kappa (\mathop{\rm
cl}\nolimits_{{\cal M}(y)}(m_{y}))\over |\mathop{\rm
Aut}(m_{y})|}\right)
\hfill\cr\hfill
=|\kappa (y)|^{r_{y}}
\left|{(\widetilde{A}_{y}^{\times})^{\tau^{\ast} =(-)^{-1}}
\over (A_{1,y}^{\times}\times A_{2,y}^{\times})^{\tau^{\ast}
=(-)^{-1}}}\right|
\left(\sum_{n_{y}\in [{\cal N}(y)/Q(y)] (k)_{\sharp}}{1\over
|\mathop{\rm Aut}(n_{y})|}\right)\qquad}
$$
quel que soit $y\in Y^{{\rm sing}}-Z^{{\rm inerte}}$, pour achever la
preuve du th\'{e}or\`{e}me.

Si $y$ n'est pas dans $Z$, cette \'{e}galit\'{e} est tautologique
puisqu'on a $A_{y}=A_{\alpha ,y}$ pour l'unique $\alpha\in\{1,2\}$ tel
que $y\in Y_{\alpha}$.

Supposons donc que $y\in Z-Z^{{\rm inerte}}$.  Notons $B_{y}$ le
compl\'{e}t\'{e} de l'anneau local de $Y$.  On a $A_{y}=B_{y}\times
B_{y}$, l'involution $\tau$ \'{e}changeant les deux copies de $B_{y}$.
De fa\c{c}on compatible on a $K_{y}=L_{y}\oplus (\omega_{Y/X,y}
\otimes_{{\cal O}_{Y,y}}L_{y}^{\otimes -1})$ o\`{u} $L_{y}$ est un
$B_{y}$-module sans torsion de rang $1$.  Posons $W_{y}=\mathop{\rm
Frac} (B_{y})\otimes_{B_{y}}L_{y}$.  Notons $\overline{B}_{y}$ le
compl\'{e}t\'{e} de $\overline{k}\otimes_{k}B_{y}$ et $\overline{(-)}$
le foncteur $\overline{B}_{y}\otimes_{B_{y}}(-)$.

Avec ces notations, la cat\'{e}gorie ${\cal M}(y)(\overline{k})$ est
\'{e}quivalente \`{a} la cat\'{e}gorie des sous-$B_{y}$-modules
${\cal W}_{y}\subset\overline{W}_{y}$ tels que $\mathop{\rm
Frac}(\overline{B}_{y}){\cal W}_{y}=\overline{W}_{y}$.  Le groupe
$P(y)(\overline{k})$ est isomorphe \`{a} $\mathop{\rm
Frac}(\overline{B}_{y})^{\times}/\overline{B}_{y}^{\times}$ et agit
sur les ${\cal W}_{y}$ par homoth\'{e}ties.  Par cons\'{e}quent la
cat\'{e}gorie $[{\cal M}(y)/P(y)](k)$ est \'{e}quivalente \`{a} la
cat\'{e}gorie des couples $({\cal W}_{y},b)\in \mathop{\rm ob}{\cal M}
(\overline{k})\times P(y)(\overline{k})$ tels que
$$
\mathop{\rm Frob}\nolimits_{q}{\cal W}_{y}=b{\cal W}_{y}\subset
\overline{W}_{y}.
$$
Un tel $b$ est automatiquement dans la {\og}{composante neutre}{\fg}
$$
\overline{\widetilde{B}}_{y}^{\times}/\overline{B}_{y}^{\times}\subset
\mathop{\rm Frac}(\overline{B}_{y})^{\times}/\overline{B}_{y}^{\times}
$$
o\`{u} $\widetilde{B}_{y}$ est le normalis\'{e} de $B_{y}$ dans
$\mathop{\rm Frac}(B_{y})$. On a donc d'une part $\kappa (\mathop{\rm
cl}\nolimits_{{\cal M}(y)}({\cal W}_{y},b))=1$ et d'autre part
$$
\mathop{\rm Frob}\nolimits_{q}(c^{-1}{\cal W}_{y})=c^{-1}{\cal W}_{y}
$$
o\`{u} $c\in \overline{\widetilde{B}}_{y}^{\times}/
\overline{B}_{y}^{\times}$ est n'importe quelle solution de
l'\'{e}quation $\mathop{\rm Frob}\nolimits_{q}(c)c^{-1}
=b$.  On en d\'{e}duit que $[{\cal M}(y)/P(y)] (k)_{\natural}$ est
l'ensemble des classes de sous-$B_{y}$-modules ${\cal W}_{y}\subset W_{y}$
tels que $\mathop{\rm Frac}(B_{y}){\cal W}_{y}=W_{y}$, modulo l'action par
homoth\'{e}ties de $\mathop{\rm Frac}(B_{y})^{\times}/B_{y}^{\times}$
et que
$$
\mathop{\rm Aut}({\cal W}_{y})=\{b\in \mathop{\rm Frac}
(B_{y})^{\times}\mid b{\cal W}_{y}={\cal W}_{y}\}/B_{y}^{\times}.
$$
\vskip 2mm

De m\^{e}me, notons $B_{\alpha ,y}$ le compl\'{e}t\'{e} de l'anneau
local de $Y_{\alpha}$ en $y$.  On a $A_{\alpha ,y}=
B_{\alpha ,y}\times B_{\alpha ,y}$, l'involution $\tau$
\'{e}changeant les deux copies de $B_{\alpha ,y}$.  De fa\c{c}on
compatible, le compl\'{e}t\'{e} de ${\cal K}_{\alpha}$ le long de
$\pi_{Y_{\alpha}}^{-1}(y)$ se d\'{e}compose en
$K_{\alpha ,y}=L_{\alpha ,y}\oplus
(\omega_{Y_{\alpha}/X,y}\otimes_{{\cal O}_{Y_{\alpha},y}}
L_{\alpha ,y}^{\otimes -1})$ o\`{u} $L_{\alpha ,y}$ est
maintenant un $B_{\alpha ,y}$-module inversible.  Notons
$W_{\alpha ,y}=\mathop{\rm Frac}(B_{\alpha ,y})
\otimes_{B_{\alpha ,y}}L_{\alpha ,y}$.

En rempla\c{c}ant successivement $Y$ par $Y_{1}$ et par $Y_{2}$, on
voit d'une part que l'ensemble $[{\cal N}(y)/Q(y)](k)_{\natural}$ est
isomorphe \`{a} l'ensemble des classes de couples $({\cal
W}_{1,y},{\cal W}_{2,y})$ o\`{u} ${\cal W}_{\alpha ,y}\subset W_{\alpha
,y}$ est un sous-$B_{\alpha ,y}$-module tel que $\mathop{\rm
Frac}(B_{\alpha ,y}){\cal W}_{\alpha ,y}=W_{\alpha ,y}$, modulo
l'action par homoth\'{e}ties de $(\mathop{\rm Frac}(B_{1,y})^{\times}/
B_{1,y}^{\times})\times (\mathop{\rm Frac}
(B_{2,y})^{\times}/B_{2,y}^{\times})$, et d'autre part que
$$\displaylines{
\qquad\mathop{\rm Aut}({\cal W}_{1,y},{\cal W}_{2,y})=\{(b_{1}\in \mathop{\rm
Frac} (B_{1,y})^{\times}\mid b_{1} {\cal W}_{1,y}={\cal W}_{1,y}\}/
B_{1,y}^{\times}
\hfill\cr\hfill
\times \{(b_{2}\in \mathop{\rm Frac} (B_{2,y})^{\times}\mid
b_{2}G_{2,y}=G_{2,y}\}/ B_{2,y}^{\times}.\qquad}
$$

Comme
$$
\mathop{\rm Frac}(B_{y})^{\times}=\mathop{\rm Frac}
(B_{1,y})^{\times} \times \mathop{\rm Frac}
(B_{2,y})^{\times},
$$
les $\mathop{\rm Frac}(B_{y})^{\times}$-modules $W_{y}$ et
$W_{1,y}\oplus W_{2,y}$ sont canoniquement isomorphes.  On a
la projection {\og}{non alg\'{e}brique}{\fg} introduite dans la
section 5 de [Ka-Lu] (voir aussi la section 6 de [La-Ra])
$$
{\cal W}_{y}\mapsto ({\cal W}_{1,y}={\cal W}_{y}\cap
W_{1,y},{\cal W}_{2,y}=
\mathop{\rm pr}\nolimits_{W_{2,y}}({\cal W}_{y}))
$$
dont on sait (cf.  loc.  cit.)  que toutes ses fibres ont 
$|\kappa (y)|^{r_{z}}$ \'{e}l\'{e}ments.  Cette projection 
est \'{e}quivariante pour les actions de 
$\mathop{\rm Frac}(B_{y})^{\times}/B_{y}^{\times}$ sur la
source et de $\mathop{\rm Frac} (B_{y})^{\times}/
(B_{1,y}^{\times} \times B_{2,y}^{\times})$ sur le but compte
tenu de l'\'{e}pimorphisme
$$
\mathop{\rm Frac}(B_{y})^{\times}/B_{y}^{\times}\twoheadrightarrow
\mathop{\rm Frac}(B_{y})^{\times}/(B_{1,y}^{\times}\times
B_{2,y}^{\times})
$$
qui est induit par l'inclusion $B_{y}^{\times}\subset
B_{1,y}^{\times} \times B_{2,y}^{\times}$.

La relation que l'on cherche \`{a} v\'{e}rifier se r\'{e}crit
$$
\sum_{{\cal W}_{y}\in }{1\over |\mathop{\rm Aut}({\cal W}_{y})|}=
|\kappa (y)|^{r_{y}} \left|\Bigl({B_{y}^{\times}\over
B_{1,y}^{\times} \times B_{2,y}^{\times}}\Bigr)\right|
\left(\sum_{({\cal W}_{1,y},{\cal W}_{2,y})}{1\over |\mathop{\rm
Aut}({\cal W}_{1,y},{\cal W}_{2,y})|}\right)
$$
et se d\'{e}montre maintenant par un simple comptage.
\hfill\hfill$\square$

\rem Remarque
\endrem
La formule pour $y\in Z-Z^{{\rm inerte}}$ ci-dessus est en fait un cas
particulier de la descente d'Harish-Chandra pour les int\'{e}grales
orbitales et la d\'{e}monstration que nous en avons donn\'{e}e est
une simple reformulation de l'argument usuel.
\hfill\hfill$\square$
\vskip 3mm

\section{4}{Le passage du local au global}
\vskip - 3mm

\subsection{4.1}{Retour \`{a} la situation locale}

Rappelons les donn\'{e}es locales du chapitre 1.

On s'est donn\'{e} un corps local non archim\'{e}dien $F$ d'anneau des
entiers ${\cal O}_{F}$, de corps r\'{e}siduel $k={\Bbb F}_{q}$ et on a
choisi une uniformisante $\varpi_{F}$ de $F$.  On s'est donn\'{e} une
extension quadratique non ramifi\'{e}e $F'$ de $F$ de corps
r\'{e}siduel $k'={\Bbb F}_{q^{2}}$, et on a not\'{e} $\tau$
l'\'{e}l\'{e}ment non trivial de $\mathop{\rm Gal}(F'/F)$.

On s'est aussi donn\'{e} une famille finie $(E_{i})_{i\in I}$
d'extensions finies s\'{e}parables $E_{i}$ de degr\'{e} $n_{i}$ de
$F$, toutes disjointes de $F'$.  On a not\'{e} $E_{i}'=E_{i}F'$ pour
$i\in I$ et encore par $\tau$ l'\'{e}l\'{e}ment non trivial de
$\mathop{\rm Gal}(E_{i}'/E_{i})$.

On s'est donn\'{e} enfin pour chaque $i\in I$ un \'{e}l\'{e}ment
$\gamma_{i}\in {\cal O}_{E_{i}'}$ hermitien, c'est-\`{a}-dire tel que
$\gamma_{i}^{\tau}+\gamma_{i}=0$.  On a suppos\'{e} que $\gamma_{i}$
engendre $E_{i}'$.

Le polyn\^{o}me minimal $P_{i}(T)$ de $\gamma_{i}$ sur $F'$ est un
polyn\^{o}me irr\'{e}ductible dans $F'[T]$.  On a suppos\'{e} que les
deux polyn\^{o}mes $P_{i}(T)$ sont deux \`{a} deux premiers entre eux.

Puisque les $\gamma_{i}$ sont entiers, les polyn\^{o}mes
$P_{i}(T)$ sont \`{a} coefficients dans ${\cal O}_{F'}$.
L'hypoth\`{e}se $\gamma_{i}$ hermitien implique que
$P_{i}^{\tau}(T)=(-1)^{n_{i}}P_{i}(-T)$.

Pour tout partie $J$ de $I$ on a pos\'{e} $E_{J}=\bigoplus_{i\in
J}E_{i}$, $P_{J}(T)=\prod_{i\in J}P_{J}(T)$, etc.

\subsection{4.2}{Un lemme d'approximation}

Rappelons quelques faits bien connus (voir les exercices, \S2, ch. II,
de [Ser]).
\vskip 2mm

Le discriminant $\mathop{\rm Disc}(P)\in F'$ d'un polyn\^{o}me unitaire
$P(T)\in F'[T]$ est le r\'{e}sultant des polyn\^{o}mes $P(T)$ et
${{\rm d}P\over {\rm d}T}(T)$.  C'est un polyn\^{o}me universel (\`{a}
coefficients entiers) en les coefficients de $P(T)$; il est donc dans
${\cal O}_{F'}$ si $P(T)\in {\cal O}_{F'}[T]$.  Le discriminant de $P(T)$
est non nul si et seulement si $P(T)$ est s\'{e}parable.

\thm LEMME 4.2.1
\enonce
Soit $J$ une partie de $I$ et $\widetilde{P}(T)$ un polyn\^{o}me
unitaire de degr\'{e} $n_{J}$ \`{a} coefficients dans ${\cal O}_{F'}$.
Alors, pour tout entier $m_{J}>n_{J}v_{F'}(\mathop{\rm Disc}(P))$, la
condition
$$
P_{J}(T)-\widetilde{P}(T)\in \varpi_{F}^{m_{J}}{\cal O}_{F'}[T]
$$
implique que
\vskip 1mm

\itemitem{-} $\widetilde{P}(T)$ est s\'{e}parable sur $F'$ et pour
chaque $i\in J$, $\widetilde{P}(T)$ admet une unique racine
$\widetilde{\gamma}_{i}$ dans $E_{i}'$ telle que
$$
v_{E_{i'}}(\gamma_{i}-\widetilde{\gamma}_{i})\geq e_{i}{m_{J}\over
n_{J}}
$$
o\`{u} on rappelle que $e_{i}$ est l'indice de ramification de
$E_{i}'$ sur $F'$; de plus, si on note $\widetilde{P}_{i}(T)$ le
polyn\^{o}me minimal de $\widetilde{\gamma}_{i}$ dans $F'[T]$, on a
$\widetilde{P}(T)= \prod_{i\in J}\widetilde{P}_{i}(T)$;
\vskip 1mm

\itemitem{-}  pour chaque $i\in J$, on a
$$
P_{i}(T)-\widetilde{P}_{i}(T)\in\varpi_{F'}^{n_{i}{m_{J}\over
n_{J}}}{\cal O}_{F'}[T]
$$
et les \'{e}l\'{e}ments ${{\rm d}P_{i}\over {\rm d}T}(\gamma_{i})
P_{J-\{i\}}(\gamma_{i})$ et ${{\rm d}\widetilde{P}_{i}\over {\rm
d}T}(\widetilde{\gamma}_{i})\widetilde{P}_{J-\{i\}}
(\widetilde{\gamma}_{i})$ ont m\^{e}me valuation dans $E_{i}'$;
\vskip 1mm

\itemitem{-} pour chaque partie $K$ de $J$, on a ${\cal O}_{F'}
[\gamma_{K}] ={\cal O}_{F'}[\widetilde{\gamma}_{K}]\subset E_{K}'$
o\`{u} bien entendu on a pos\'{e} $\widetilde{\gamma}_{K}=
\oplus_{i\in K}\widetilde{\gamma}_{i}\in E_{K}'$.

\endthm

\rem D\'{e}monstration
\endrem
Soient $\overline{F}{}'$ une cl\^{o}ture s\'{e}parable de $F'$ et
$\gamma_{i,1},\gamma_{i,2},\ldots ,\gamma_{i,n_{i}}$ les racines de
$P_{i}(T)$ dans $\overline{F}{}'$.  On note encore
$v_{F'}:\overline{F}{}'^{\,\times} \rightarrow {\Bbb Q}$ l'unique
valuation qui prolonge la valuation discr\`{e}te
$v_{F'}:F'^{\times}\rightarrow {\Bbb Z}$.  On rappelle que
$v_{F'}\circ\sigma =v_{F'}$ pour tout \'{e}l\'{e}ment $\sigma\in
\mathop{\rm Gal}(\overline{F}{}'/F')$ du groupe de Galois de
$\overline{F}{}'$ sur $F'$ et que $v_{E_{i}'}=e_{i}(v_{F'}\circ\sigma
)|E_{i}'^{\times}$ pour tout $F'$-plongement $\sigma
:E_{i}'\hookrightarrow \overline{F}{}'$.

Montrons dans un premier temps que $\widetilde{P}$ est
s\'{e}parable sur $F'$.  On a
$$
v_{F'}(\mathop{\rm Discr}(P_{J})-\mathop{\rm Discr}
(\widetilde{P}))\geq m_{J}\geq {m_{J}\over n_{J}}>
v_{F'}(\mathop{\rm Discr}(P_{J}))
$$
par hypoth\`{e}se.  Par suite, $\mathop{\rm Discr}(\widetilde{P})$
est non nul, et en fait de m\^{e}me valuation que $\mathop{\rm
Discr}(P_{J})$.

On notera $(\widetilde{P}_{\tilde{\imath}})_{\tilde{\imath}\in\widetilde{J}}$
l'ensemble des facteurs irr\'{e}ductibles (unitaires) de
$\widetilde{P}$ et pour chaque $\tilde{\imath}\in\widetilde{J}$,
$\widetilde{\gamma}_{\tilde{\imath},1},\ldots
,\widetilde{\gamma}_{\tilde{\imath},\tilde{n}_{\tilde{\imath}}}$ les racines
de $\widetilde{P}_{\tilde{\imath}}$ dans $\overline{F}{}'$.

Comme
$$
v_{F'}(\mathop{\rm Discr}(P_{J}))=\sum_{(i,j)\not=(i',j')}
v_{F'}(\gamma_{i,j}-\gamma_{i',j'})
$$
on a
$$
v_{F'}(\gamma_{i,j}-\gamma_{i',j'})<{m_{J}\over n_{J}}
$$
quels que soient $i,i'\in J$, $j\in \{1,\ldots ,n_{i}\}$ et $j'\in
\{1,\ldots ,n_{i'}\}$ tels que $(i,j)\not=(i',j')$.  En d'autres
termes, la {\og}{distance}{\fg}
$$
|\gamma_{i,j}-\gamma_{i',j'}|_{F'}=
q^{-v_{F'}(\gamma_{i,j}-\gamma_{i',j'})}
$$
entre deux racines distinctes de $P_{J}(T)$ est strictement plus
grande que $q^{-{m_{J}\over n_{J}}}$. De m\^{e}me, on a
$$
|\widetilde{\gamma}_{\tilde{\imath},\tilde{\jmath}}-
\widetilde{\gamma}_{\tilde{\imath}',\tilde{\jmath}'}|_{F'}>
q^{-{m_{J}\over n_{J}}}
$$
quels que soient $(\tilde{\imath},\tilde{\jmath})\not=
(\tilde{\imath}',\tilde{\jmath}')$.

Maintenant, soit $\widetilde{\gamma}$ une racine de
$\widetilde{P}(T)$ dans $\overline{F}{}'$.  On a
$$
\sum_{i\in J}\sum_{j=1}^{n_{i}}
v_{F'}(\widetilde{\gamma}-\gamma_{i,j})=
v_{F'}(P_{J}(\widetilde{\gamma}))=v_{F'}(P_{J}(\widetilde{\gamma})-
\widetilde{P}(\widetilde{\gamma}))\geq m_{J}.
$$
Il existe donc $i\in J$ et $j\in\{1,\ldots ,n_{i}\}$ tel que
$v_{F'}(\widetilde{\gamma}-\gamma_{i,j})\geq {m_{J}\over n_{J}}$, i.e.
$$
|\widetilde{\gamma}-\gamma_{i,j}|_{F'}\leq q^{-{m_{J}\over
n_{J}}},
$$
et le couple $(i,j)$ est unique puisque la distance entre deux
$\gamma_{i,j}$ distincts est $>q^{-{m_{J}\over
n_{J}}}$.

On construit ainsi une application $\widetilde{\gamma}\mapsto \gamma
(\widetilde{\gamma})$ de l'ensemble des racines de
$\widetilde{P}(T)$ dans $\overline{F}{}'$ dans l'ensemble des
racines de $P_{J}(T)$ dans $\overline{F}{}'$ qui est
caract\'{e}ris\'{e}e par la relation
$$
|\widetilde{\gamma}-\gamma (\widetilde{\gamma})|_{F'}\leq q^{-{m_{J}\over
n_{J}}}
$$
et qui est $\mathop{\rm Gal}(\overline{F}{}'/F')$-\'{e}quivariante
puisque la distance est $\mathop{\rm
Gal}(\overline{F}{}'/F')$-invariante. Cette application est
n\'{e}cessairement injective puisque la  distance entre deux
racines distinctes de $\widetilde{P}(T)$ est aussi $>q^{-{m_{J}\over
n_{J}}}$. Elle est donc bijective puisque $P_{J}(T)$ et
$\widetilde{P}(T)$ sont tous les deux s\'{e}parables de
degr\'{e} $n_{J}$. On en d\'{e}duit une bijection
$\tilde{\imath}\rightarrow i(\tilde{\imath})$ de $\widetilde{J}$ sur $J$ telle
que $\tilde{n}_{i(\tilde{\imath})}=n_{i}$ et
pour chaque $\tilde{\imath}$ une permutation $\tilde{\jmath}\rightarrow
j(\tilde{\jmath})$ de $\{1,\ldots ,n_{i(\tilde{\imath})}\}$ telle que
$$
\gamma (\widetilde{\gamma}_{\tilde{\imath},\tilde{\jmath}})=
\gamma_{i(\tilde{\imath}),j(\tilde{\jmath})}
$$
pour tout $\tilde{\imath}\in\widetilde{J}$ et tout $\tilde{\jmath}\in\{1,\ldots
,n_{i(\tilde{\imath})}\}$. Quitte \`{a} changer l'indexation des facteurs
de $\widetilde{P}$ et la num\'{e}rotation des racines des
$\widetilde{P}_{\tilde{\imath}}$, on peut supposer et on supposera dans la
suite que $\widetilde{J}=J$ et que les bijections $\tilde{\imath}\mapsto
i(\tilde{\imath})$ et $\tilde{\jmath}\mapsto j(\tilde{\jmath})$ 
sont les applications identiques.

La minimalit\'{e} de la distance entre $\widetilde{\gamma}$ et $\gamma
(\widetilde{\gamma})$ assure de plus qu'il ne peut pas exister
d'\'{e}l\'{e}ment de $\overline{F}{}'$ qui soit conjugu\'{e} \`{a}
$\gamma (\widetilde{\gamma})$ sur $F'[\widetilde{\gamma}]$ autre que
$\gamma (\widetilde{\gamma})$ lui-m\^{e}me, de sorte que $\gamma
(\widetilde{\gamma})\in F'[\widetilde{\gamma}]$ puisque
$F'[\widetilde{\gamma}]$ est s\'{e}parable sur $F'$. Par suite, on a
m\^{e}me $F'[\gamma_{i,j}]=F'[\widetilde{\gamma}_{i,j}]$ quels que
soient $i\in J$ et $j\in\{1,\ldots ,n_{i}\}$.

Les deux $F'$-alg\`{e}bres $F'[T]/(\widetilde{P}(T))$ et $E_{J}'$
sont donc $F'$-isomorphes. Plus pr\'{e}cis\'{e}ment on a l'isomorphisme
$$
E_{J}'=\bigoplus_{i\in J}F'[\gamma_{i}]
\buildrel\sim\over\longrightarrow\bigoplus_{i\in J}F'[\gamma_{i,1}]
=\bigoplus_{i\in J}F'[\widetilde{\gamma}_{i,1}]
\buildrel\sim\over\longrightarrow \bigoplus_{i\in
J}F'[T]/(\widetilde{P}_{i}(T))=F'[T]/(\widetilde{P}(T))
$$
o\`{u} $\gamma_{i}$ est envoy\'{e} sur $\gamma_{i,1}$ par le
premier isomorphisme et $\widetilde{\gamma}_{i,1}$ est envoy\'{e} sur
la classe de $T$ modulo $\widetilde{P}_{i}(T)$ par le second.

Notons $\widetilde{\gamma}=(\widetilde{\gamma}_{i})_{i\in J}\in
E_{J}'$ l'image inverse de la classe de $T$ modulo $\widetilde{P}$ par
l'isomorphisme ci-dessus.  Bien s\^{u}r, pour chaque $i\in J$,
l'\'{e}l\'{e}ment $\widetilde{\gamma}_{i}\in E_{i}'$ est entier sur
${\cal O}_{F'}$.  De plus, la valuation
$v_{E_{i}'}(\widetilde{\gamma}_{i}-\gamma_{i})$ est
$\geq e_{i}{m_{J}\over n_{J}}$.

On a
$$
v_{F'}(\mathop{\rm Discr}(P_{J}))=\sum_{i\in J}v_{F'}(\mathop{\rm
Discr}(P_{i}))+\sum_{i\not=i'\in J}v_{F'}(\mathop{\rm
Res}(P_{i},P_{i'}))
$$
avec $v_{F'}(\mathop{\rm Discr}(P_{i}))=n_{i}{v_{E'}\left({{\rm d}P_{i}
\over {\rm  d}T}(\gamma_{i})\right)\over e_{i}}=2\delta_{i}+n_{i}-{n_{i}
\over e_{i}}$ et $v_{F'}(\mathop{\rm Res}(P_{i},P_{i'}))=r_{i,i'}$,
de sorte que
$$
{m_{J}\over n_{J}}\geq v_{F'}(\mathop{\rm Discr}(P_{J}))=
\sum_{i\in J}\left({n_{i}a_{i}\over  e_{i}}+n_{i}-{n_{i}\over
e_{i}}\right)\geq {a_{i}\over e_{i}}
$$
o\`{u} on rappelle que l'on a pos\'{e}
$a_{i}={(2\delta_{i}+\sum_{i'\not=i}r_{i,i'})e_{i}\over n_{i}}$.

Comme ${\frak a}_{J}=\varpi_{E_{J}}^{\underline{a}_{J}}{\cal
O}_{E_{J}'}\subset {\cal O}_{F'}[\gamma_{J}]\subset {\cal O}_{E_{J}'}$
avec $\underline{a}_{J}=(a_{i})_{i\in J}$, la condition
$$
\widetilde{\gamma}_{i}-\gamma_{i}\in \varpi_{E_{i}}^{e_{i}{m_{J}\over
n_{J}}}{\cal O}_{E_{i}'}\subset \varpi_{E_{i}}^{a_{i}}{\cal
O}_{E_{i}'}
$$
quel que soit $i\in J$, implique que ${\cal
O}_{F'}[\widetilde{\gamma}_{J}]\subset {\cal O}_{F'}[\gamma_{J}]$.

Puisque les deux polyn\^{o}mes $P(T)$ et $\widetilde{P}(T)$ jouent des
r\^{o}les parfaitement sym\'{e}triques, on a aussi l'inclusion dans le
sens inverse, d'o\`{u} l'\'{e}galit\'{e} d\'{e}sir\'{e}e ${\cal
O}_{F'}[\widetilde{\gamma}_{J}]={\cal O}_{F'}[\gamma_{J}]$.
\hfill\hfill$\square$
\vskip 3mm

Une cons\'{e}quence de ce lemme est que les int\'{e}grales orbitales
qui figurent dans le lemme fondamental, dont on a vu qu'elles comptent
des id\'{e}aux fractionnaires dans les ${\cal O}_{F}$-alg\`{e}bres
${\cal O}_{F'}[T]/(P_{J}(T))$, ne d\'{e}pendent que de la
r\'{e}duction des $P_{i}(T)$ modulo une puissance assez grande de
$\varpi_{F}$.

\subsection{4.3}{Variantes du lemme d'approximation}

Il est commode de choisir un \'{e}l\'{e}ment $\zeta\in k'^{\times}$
tel que $\zeta^{\tau}+\zeta=0$.  Gr\^{a}ce \`{a} ce choix, les
polyn\^{o}mes $P_{i}(T)$ peuvent \^{e}tre \'{e}crits sous la forme
$P_{i}(T)=Q_{i}(\zeta T)$ o\`{u} $Q_{i}(T)\in {\cal O}_{F}[T]$.  Pour
toute partie $J$ de $I$ posons $Q_{J}(T)=\prod_{i\in J}Q_{i}(T)$; on
a alors $E_{J}=F[T]/Q_{J}(T)$.

La donn\'{e}e de $\gamma_{i}$ comme ci-dessus est \'{e}quivalente
\`{a} la donn\'{e}e d'un ${\cal O}_{F}$-sous-sch\'{e}ma fini plat
$$
Y_{i,{\cal O}_{F}}=\mathop{\rm Spec}({\cal O}_{F}[T]/Q_{i}(T))
$$
de la droite affine ${\Bbb A}_{{\cal O}_{F}}^{1}=\mathop{\rm
Spec}({\cal O}_{F}[T])$ au-dessus de ${\cal O}_{F}$.  L'anneau ${\cal
O}_{F}[T]/Q_{i}(T)$ est un anneau local int\`{e}gre complet, qui est
g\'{e}om\'{e}triquement unibranche.  Plus pr\'{e}cis\'{e}ment,
$Y_{i,{\cal O}_{F}}\otimes_{k}\overline{k}$ peut ne pas \^{e}tre
connexe, mais chacune de ses composantes connexes est unibranche.

Pour toute partie $J$ de $I$ notons $Y_{J,{\cal O}_{F}}=\mathop{\rm
Spec}({\cal O}_{F}[T]/Q_{J}(T))$ le diviseur de Cartier de ${\Bbb
A}_{{\cal O}_{F}}^{1}$ qui est la somme des $Y_{i,{\cal O}_{F}}$ pour
$i\in J$.

Pour tout entier $m\geq 0$, notons ${\cal O}_{F,m}$ la $k$-alg\`{e}bre
finie ${\cal O}_{F,m}={\cal O}_{F}/\varpi_{F}^{m}{\cal O}_{F}$ et
notons $(-)_{{\cal O}_{F,m}}$ la r\'{e}duction modulo $\varpi_{F}^{m}$
de $(-)_{{\cal O}_{F}}$.  Le lemme d'approximation de la section
pr\'{e}c\'{e}dente peut se reformuler comme suit.

\thm LEMME 4.3.1
\enonce
Soit $J$ une partie de $I$.  Il existe un entier $m_{J}\geq 0$ ayant
la propri\'{e}t\'{e} suivante: tout sous-${\cal O}_{F}$-sch\'{e}ma
fini et plat $\widetilde{Y}_{{\cal O}_{F}}={\cal
O}_{F}[T]/\widetilde{Q}(T)$ de ${\Bbb A}_{{\cal O}_{F}}^{1}$ tel que
$$
\widetilde{Y}_{{\cal O}_{F,m_{J}}}=Y_{J,{\cal O}_{F,m_{J}}}
$$
en tant que sous-sch\'{e}ma ferm\'{e} de $\mathop{\rm Spec}({\cal
O}_{F,m_{J}}[T])$, est isomorphe \`{a} $Y_{J,{\cal O}_{F}}$ en tant que
${\cal O}_{F}$-sch\'{e}ma non plong\'{e}.
\endthm

Notons $\overline{Q}_{J}(T)$ la r\'{e}duction modulo $\varpi_{F}$ de
$Q_{J}(T)$.  Pour tout entier $m\geq 0$, on peut voir
$\overline{Q}_{J}[T]$ comme un \'{e}l\'{e}ment de ${\cal O}_{F}[T]$ et
de ${\cal O}_{F,m}[T]$ via les plongements \'{e}vidents $k[T]\subset
{\cal O}_{F}[T]$ et $k[T]\subset {\cal O}_{F,m}[T]$.  Le sch\'{e}ma
$Y_{J,{\cal O}_{F,m}}$ est un sch\'{e}ma de longueur $mn_{J}$
support\'{e} par les z\'{e}ros de $\overline{Q}_{J}[T]$ dans la fibre
sp\'{e}ciale ${\Bbb A}_{k}^{1}$ de ${\Bbb A}_{\mathop{{\cal
O}_{F}}}^{1}$, de sorte que $Y_{{\cal O}_{F,m}}$ est
n\'{e}cessairement contenu dans le sch\'{e}ma fini
$Z_{J,m}=\mathop{\rm Spec} ({\cal
O}_{F,m}[T]/(\overline{Q}(T)^{mn_{J}}))$ qui ne d\'{e}pend que de
$\overline{Q}_{J}(T)$.

\thm LEMME 4.3.2
\enonce
Soient $J$ une partie de $I$ et $\widetilde{Q}(T)\in {\cal O}_{F}[T]$
un polyn\^{o}me unitaire de degr\'{e} $n_{J}$ ayant aussi pour
r\'{e}duction $\overline{Q}_{J}(T)\in k[T]$ modulo l'uniformisante
$\varpi_{F}$.  Notons $\widetilde{Y}_{{\cal O}_{F}}=\mathop{\rm
Spec}({\cal O}_{F}[T]/\widetilde{Q}(T))$.  Soit $m_{J}$ un entier
comme dans le lemme ci-dessus.  Supposons que
$$
Y_{J,{\cal O}_{F}}\cap Z_{J,m_{J}}=\widetilde{Y}_{{\cal O}_{F}}\cap
Z_{J,m_{J}}.
$$
Alors $\widetilde{Y}_{{\cal O}_{F}}$ est isomorphe \`{a} $Y_{J,{\cal
O}_{F}}$ en tant que ${\cal O}_{F}$-sch\'{e}ma non plong\'{e}.
\endthm

\rem D\'{e}monstration
\endrem
Comme on a vu plus haut, $Y_{J,{\cal O}_{F}} \cap Z_{J,m_{J}}=Y_{{\cal
O}_{F}}\cap {\Bbb A}_{{\cal O}_{F,m_{J}}}^{1}$.  De m\^{e}me, on voit que
$\widetilde{Y}_{{\cal O}_{F}} \cap Z_{J,m_{J}}=\widetilde{Y}_{{\cal
O}_{F}}\cap {\Bbb A}_{{\cal O}_{F,m_{J}}}^{1}$.  Ainsi, le lemme est une
cons\'{e}quence du lemme pr\'{e}c\'{e}dent.
\hfill\hfill$\square$
\vskip 3mm

Puisque toutes les int\'{e}grales orbitales que l'on consid\`{e}re ne
d\'{e}pendent que des classes d'iso\-mor\-phis\-me des $Y_{J,{\cal
O}_{F}}$, ces int\'{e}grales ne changent pas si on remplace les
polyn\^{o}mes $Q_{J}(T)$ par des polyn\^{o}mes proches
$\widetilde{Q}(T)$ o\`{u} la propri\'{e}t\'{e} d'\^{e}tre proche est
signifie que l'hypoth\`{e}se du lemme est satisfaite.

\subsection{4.4}{Th\'{e}or\`{e}me de Bertini-Poonen}

Une variante du th\'{e}or\`{e}me de Bertini, valable sur les corps
finis, a \'{e}t\'{e} d\'{e}montr\'{e} par Poonen ([Poo]).  Nous allons
rappeler pr\'{e}cis\'{e}ment son \'{e}nonc\'{e}.

\thm TH\'{E}OR\`{E}ME 4.4.1
\enonce
Soit $S$ un sous-sch\'{e}ma localement ferm\'{e} d'un espace projectif
standard de dimension finie ${\Bbb P}$ sur un corps fini $k={\Bbb
F}_q$.  Soit $Z$ un ferm\'{e} {\rm fini} de ${\Bbb P}$ tel que
$S-(S\cap Z)$ est lisse de dimension $d$.  Soit $z\in H^{0}(Z,{\cal
O}_{Z})$.  Alors, pour tout entier $N$ assez grand il existe un
polyn\^{o}me homog\`{e}ne $f\in H^{0}({\Bbb P},{\cal O}_{{\Bbb P}}(N))$
de degr\'{e} $N$ \`{a} coefficients dans $k$, tel que l'hypersurface
$H_{f}$ d\'{e}finie par $f$ v\'{e}rifie :
\vskip 1mm

\item{-} $H_{f}\cap S$ est lisse de dimension $d-1$ en dehors de
$Z\cap S$,
\vskip 1mm

\item{-} $H_{f}\cap Z$ est le sous-sch\'{e}ma ferm\'{e} de $Z$,
d\'{e}fini par l'\'{e}quation $z=0$.
\hfill\hfill$\square$

\endthm

En fait, Poonen a d\'{e}montr\'{e} un \'{e}nonc\'{e} plus pr\'{e}cis :
il a calcul\'{e} la densit\'{e} des $f$ v\'{e}rifiant ces
propri\'{e}t\'{e}s lorsque $N$ tend vers l'infini, et cette
densit\'{e} est un nombre strictement positif.  Pour notre besoin,
seule existence de $f$ suffit.

L'argument de Poonen se g\'{e}n\'{e}ralise sans difficult\'{e}s au cas
d'une famille finie de sous-sch\'{e}mas localement ferm\'{e}s
$(V_{i})_{i\in I}$ de
${\Bbb P}$.

\subsection{4.5}{Tracer des courbes sur la surface r\'{e}gl\'{e}e}

Soient $X$ une courbe propre, lisse et g\'{e}om\'{e}triquement connexe
sur $k={\Bbb F}_{q}$ et $X'\rightarrow X$ un rev\^{e}tement double
\'{e}tale qui est aussi g\'{e}om\'{e}triquement connexe.  Soient
$x_{0},x_{\infty} \in X(k)$ tels que $x_{0}$ est inerte et
$x_{\infty}$ est d\'{e}compos\'{e} dans le rev\^{e}tement
$X'\rightarrow X$.  On identifie le corps local $F$ avec le
compl\'{e}t\'{e} du corps des fonctions de $X$ en la place $x_{0}$ et
on identifie $F'$ \`{a} l'extension non ramifi\'{e}e de $F$
d\'{e}finie par $X'\rightarrow X$.

On se donne un diviseur $D$ de degr\'{e} $\geq g+1$ sur $X$ qui
\'{e}vite $x_{0}$.  Avec les notations du chapitre 2, on a le
fibr\'{e} en droites projectives $\Sigma ={\Bbb P}({\cal O}_{X}
\oplus ({\cal L}_{D})^{\otimes -1})\rightarrow X$ qui compl\`{e}te le
fibr\'{e} en droites vectorielles $\Sigma^{\circ}={\Bbb V}(({\cal
L}_{D})^{\otimes -1})\rightarrow X$ dont les sections sont celles de
${\cal L}_{D}$.

L'hypoth\`{e}se $\mathop{\rm deg}(D)\geq g+1$ implique que ${\cal
O}_{\Sigma}(1)\otimes p^{\ast}{\cal L}_{D}$ est tr\`{e}s ample sur
${\Bbb V}(({\cal L}_{D})^{\otimes -1})$.  En fait, consid\'{e}rons le
morphisme canonique
$$
\Sigma\rightarrow {\Bbb P}(H^{0}(\Sigma ,{\cal O}_{\Sigma}(1)\otimes
p^{\ast}{\cal L}_{D}))={\Bbb P}(H^{0}(X,{\cal L}_{D})\oplus k)={\Bbb
P}_{\Sigma}
$$
et son image $S$.  Notons $s_{\infty}\in {\Bbb P}_{\Sigma}(k)$ le
point d\'{e}fini par la projection sur le second facteur
$H^{0}(X,{\cal L}_{D})\oplus k\twoheadrightarrow k$.  Alors on a une
projection
$$
{\Bbb P}_{\Sigma}-\{s_{\infty}\}={\Bbb V}({\cal O}_{{\Bbb P}_{X}}(-1))
\rightarrow {\Bbb P}(H^{0}(X,{\cal L}_{D}))={\Bbb P}_{X}
$$
qui s'ins\`{e}re dans un carr\'{e} cart\'{e}sien
$$
\xymatrix{\Sigma^{\circ}\ar[d]_{p^{\circ}}\ar@{^{(}->}[r]&
{\Bbb P}_{\Sigma}-\{s_{\infty}\}\ar[d]\cr
X\ar@{^{(}->}[r] & {\Bbb P}_{X}}
$$
o\`{u} les fl\`{e}ches verticales sont des fibr\'{e}s en droites et
les fl\`{e}ches horizontales sont les fl\`{e}ches canoniques.  On a
$s_{\infty}\in S(k)$ et son image r\'{e}ciproque dans $\Sigma$ est le
diviseur \`{a} l'infini $\Sigma_{\infty}$ compl\'{e}mentaire de
$\Sigma^{\circ}$ dans $\Sigma$.  On a donc un isomorphisme
$\Sigma^{\circ}\buildrel\sim\over\longrightarrow S-\{s_{\infty}\}$.
Autrement dit, le morphisme $\Sigma\rightarrow S$ contracte exactement
le diviseur $\Sigma_{\infty}$ sur le point $s_{\infty}$.

Rappelons que pour tout entier $N\geq 1$ on a
$$
H^{0}(\Sigma ,({\cal O}_{\Sigma}(1)\otimes p^{\ast}{\cal
L}_D)^{\otimes N})= \bigoplus_{i=0}^{N}H^0(X,({\cal
L}_{D})^{\otimes i})
$$
qui est l'ensemble des $k$-points de l'espace de Hitchin ${\Bbb
A}_{{\rm U}(N)}$ pour le groupe unitaire ${\rm U}(N)$. Une section
$a=(a_{0},a_{1},\ldots ,a_{n})\in {\Bbb A}_{{\rm U}(N)}(k)$ o\`{u}
$a_{0}\not=0$, ne s'annule pas sur le diviseur $\Sigma_{\infty}$.

Prenons un voisinage de Zariski $U$ de $x_{0}$ dans $X$ et une
trivialisation de ${\cal L}_{D}$ sur $U$ qui identifie la
restriction de ${\Bbb V}(({\cal L}_{D})^{\otimes -1})$ \`{a} $U$ avec la
droite affine ${\Bbb A}_{U}^{1}$.  On en d\'{e}duit une
identification de ${\Bbb A}_{{\cal O}_{F}}^{1}$ avec ${\Bbb
V}(({\cal L}_{D})^{\otimes -1})\times_{X} X_{x_{0}}$ o\`{u} $X_{x_{0}}$
est le compl\'{e}t\'{e} de $X$ en $x_{0}$
\vskip 2mm

Pour $i\in I$ soit $Q_{i}(T)\in {\cal O}_{F}[T]$ un polyn\^{o}me comme
dans 4.3.  On note $Q_{J}(T)=\prod_{i\in J}^{}Q_{i}(T)$ pour toute
partie $J$ de $I$.  Pour $i\in I$ notons $\overline{Q}_{i}(T)\in k[T]$
la r\'{e}duction modulo $\varpi_{F}$ de $Q_{i}(T)$ et notons
$\overline{Q}_{J}(T)=\prod_{i\in J}\overline{Q}_{i}(T)$ pour toute
partie $J$ de $I$.  Soit $m=m_{I}$ un entier naturel comme dans les
lemmes 4.3.1 et 4.3.2, et
$$
Z_{m}=Z_{I,m}=\mathop{\rm Spec}({\cal O}_{F}[T]/\langle
\overline{Q}(t)^{mn},\varpi_{F}^{m}\rangle)
$$
le sous-sch\'{e}ma fini comme dans 4.3, vu comme un sous-sch\'{e}ma
ferm\'{e} de $\Sigma$.  Dans $Z_{m}$ on a le sous-sch\'{e}ma ferm\'{e}
$\mathop{\rm Spec}({\cal O}_{F}[T]/Q_{J}(T))\cap Z_{m}$ pour chaque
sous-ensemble $J\subset I$.  En prenant une partition $I=I_{1}\amalg
I_{2}$ comme dans la section (4.1), on a un sous-sch\'{e}ma ferm\'{e}
$$
\mathop{\rm Spec}({\cal O}_{F}[T]/Q_{I_\alpha}(T))\cap Z_{m}
$$
pour chaque $\alpha\in\{1,2\}$.

Au-dessus du point $x_{\infty}$ de $X$ choisissons arbitrairement
deux $k$-points distincts $y_{\infty,\alpha}$, $\alpha\in
\{1,2\}$.

\thm PROPOSITION 4.5.1
\enonce
Il existe des entiers naturels $N_{\alpha}\geq n_{\alpha}$ et des
points $a_{\alpha}\in {\Bbb A}_{{\rm U}(N_{\alpha})}(k)$ pour
$\alpha\in \{1,2\}$, tels que les courbes spectrales
$Y_{a_{\alpha}}\subset\Sigma$ correspondantes v\'{e}rifient les
propri\'{e}t\'{e}s suivantes:
\vskip 1mm

\item{-} $Y_{a_{\alpha}}\cap Z_{m}=\mathop{\rm Spec}({\cal
O}_{F}[T]/Q_{I_\alpha}(T))\cap Z_{m}$ quel que soit $\alpha\in
\{1,2\}$,
\vskip 1mm

\item{-} pour tout $\alpha\in\{1,2\}$, la courbe spectrale
$Y_{a_{\alpha}}$ est lisse en dehors de $Z_{m}\subset\Sigma$, elle est
\'{e}tale au-dessus de $x_\infty$ et elle passe par $y_{\infty
,\alpha}$,
\vskip 1mm

\item{-} les deux courbes spectrales $Y_{a_{1}}$ et $Y_{a_{2}}$ se
coupent transversalement en dehors de $Z_{m}$ et ne se coupent pas
au-dessus de $x_\infty$.

\endthm

\rem D\'{e}monstration
\endrem
En appliquant le th\'{e}or\`{e}me de Bertini-Poonen 4.4.1, on peut
d'abord trouver un entier $N_{1}\geq n_{1}$ et une section $a_{1}\in
H^{0}({\Bbb P},{\cal O}_{{\Bbb P}}(N_{1}))$ tels que l'hypersurface
correspondante $H_{a_{1}}\subset {\Bbb P}$ v\'{e}rifie les
propri\'{e}t\'{e}s suivantes:
\vskip 1mm

\item{-} $H_{a_{1}}$ ne passe pas ni par $s_{\infty}$, ni par les
$y_{\infty,2}$, mais elle passe par $y_{\infty,1}$,
\vskip 1mm

\item{-} $H_{a_{1}}$ coupe transversalement l'image de
$p^{-1}(x_{\infty})$ dans $S$,
\vskip 1mm

\item{-} $H_{a_{1}}\cap S$ est lisse en dehors de l'image de $Z_{m}$
dans $S$ et son intersection avec $Z_{m}$ est \'{e}gale \`{a}
$\mathop{\rm Spec}({\cal O}_{F}[T]/Q_{I_1}(T))\cap Z_{m}$.
\vskip 1mm

Bien entendu, le ferm\'{e} $Z\subset {\Bbb P}$ dans l'\'{e}nonc\'{e}
4.4.1 du th\'{e}or\`{e}me de Bertini-Poonen doit \^{e}tre la
r\'{e}union de notre $Z_{m}$ et des points $s_{\infty}$ et
$y_{\infty,i}$ avec $i\not=1$, ainsi que du premier voisinage
infinit\'{e}simal de $y_{\infty,1}$ dans ${\Bbb P}$.  De plus la
fonction $z\in H^{0}(Z,{\cal O}_{Z})$ doit \^{e}tre choisie pour que
$H_{a_{1}}$ ne passe pas par $s_{\infty}$, ni par $y_{\infty,i}$ avec
$i\not=1$, qu'elle passe par $y_{\infty,1}$, et enfin qu'elle coupe
$Z_{m}$ et le premier voisinage infinit\'{e}simal de $y_{\infty,1}$
comme on le souhaite.

On prend alors pour $Y_{a_{1}}$ l'image r\'{e}ciproque de
$H_{a_{1}}\cap S$ dans $\Sigma$.

En appliquant de nouveau le th\'{e}or\`{e}me de Bertini-Poonen 4.4.1,
on peut trouver un entier $N_{2}\geq n_{2}$ et une section $a_{2}\in
H^{0}({\Bbb P},{\cal O}_{{\Bbb P}}(N_{2}))$ telle que l'hypersurface
$H_{a_{2}}\subset {\Bbb P}$ correspondante v\'{e}rifie les
propri\'{e}t\'{e}s suivantes :
\vskip 1mm

\item{-} $H_{a_{2}}$ ne passe pas ni par $s_{\infty}$, ni par
$y_{\infty,1}$, mais elle passe par $y_{\infty,2}$,
\vskip 1mm

\item{-} $H_{a_{2}}$ coupe transversalement l'image de
$p^{-1}(x_{\infty})$ dans $S$,
\vskip 1mm

\item{-} $H_{a_{2}}\cap S$ est lisse en dehors de l'image de $Z_{m}$
dans $S$ et son intersection avec $Z_{m}$ est \'{e}gale \`{a}
$\mathop{\rm Spec}({\cal O}_{F}[T]/Q_{I_2}(T))\cap Z_{m}$,
\vskip 1mm

\item{-} en dehors de $Z_{m}$, $H_{a_{2}}$ coupe transversalement
$H_{a_{1}}$.
\vskip 1mm

On prend alors pour $Y_{a_{2}}$ l'image r\'{e}ciproque de
$H_{a_{2}}\cap S$ dans $\Sigma$.
\hfill\hfill$\square$

\thm LEMME 4.5.2
\enonce
Supposons que $\mathop{\rm deg}(D)>\mathop{\rm long}(Z_{m})$.  Alors
les courbes $Y_{a_{\alpha}}$ construites dans la proposition
pr\'{e}c\'{e}dente, ainsi que leurs rev\^{e}tements doubles \'{e}tales
$Y_{a_{\alpha}}'$ sont automatiquement g\'{e}om\'{e}triquement
irr\'{e}ductibles.  En particulier, le point g\'{e}om\'{e}trique
$a=(a_{1},a_{2})$ est un point elliptique au sens de la section {\rm
(2.8)}.
\endthm

\rem D\'{e}monstration
\endrem
Il suffit de d\'{e}montrer que $\overline{Y}{}_{a_{\alpha}}'=
Y_{a_{\alpha}}'\otimes_{k}\overline{k}$ est irr\'{e}ductible pour tout
$\alpha\in\{1,2\}$.  Raisonnons par l'absurde en supposant que
$\overline{Y}_{a_{\alpha}}'= \overline{Y}_{a_{\alpha},1}'\cup
\overline{Y}_{a_{\alpha},2}'$ o\`{u} $\overline{Y}_{a_{\alpha},j}'$
est une courbe finie et plate de degr\'{e} $N_{\alpha,j}\geq 1$
au-dessus de $\overline{X}'= X'\otimes_{k}\overline{k}$.  En
particulier, on a $N_{\alpha,1}+N_{\alpha,2}=N_{\alpha}$.  D'apr\`{e}s
le lemme 2.5.3, la composante $\overline{Y}_{a_{\alpha},j}'$ est le
diviseur d'une section de
$$
H^{0}(\Sigma'\otimes_{k}\overline{k},({\cal O}_{\Sigma'}(1)\otimes
p'^{\ast}{\cal O}_{X'}(2D))^{\otimes N_{\alpha,j}}).
$$
L'intersection $\overline{Y}_{a_{\alpha},1}'\cap
\overline{Y}_{a_{\alpha},2}'$ est alors un sous-sch\'{e}ma ferm\'{e}
fini de $\Sigma'\otimes_{k}\overline{k}$ de longueur
$$
2N_{\alpha,1}N_{\alpha,2}\mathop{\rm deg}(D)\geq
2\mathop{\rm deg}(D) >2\mathop{\rm long}(Z_{m}).
$$
Par cons\'{e}quent, cette intersection ne peut pas \^{e}tre
enti\`{e}rement contenue dans l'image r\'{e}ciproque de
$Z_{m}\otimes_{k}\overline{k}$ par $\pi_{\Sigma}:\Sigma'\rightarrow
\Sigma$ puisque cette image r\'{e}ciproque est de longueur
$2\mathop{\rm long}(Z_{m})$.  Mais, si $\overline{Y}_{a_{\alpha},1}'$
coupait $\overline{Y}_{a_{\alpha},2}'$ en dehors de
$\pi_{\Sigma}^{-1}(Z_{m}\otimes_{k}\overline{k})$, la r\'{e}union
$\overline{Y}_{a_{\alpha}}'$ de $\overline{Y}_{a_{\alpha},1}'$ et
$\overline{Y}_{a_{\alpha},2}'$ ne pourrait pas \^{e}tre lisse en dehors de
$\pi_{\Sigma}^{-1}(Z_{m}\otimes_{k}\overline{k})$, d'o\`{u} une
contradiction et le lemme est d\'{e}montr\'{e}.
\hfill\hfill$\square$

\subsection{4.6}{Fin de la d\'{e}monstration du lemme fondamental}

Consid\'{e}rons toujours la situation locale rappel\'{e}e en (4.1).
Fixons de plus une partition non triviale $I=I_{1}\amalg I_{2}$ comme
dans l'\'{e}nonc\'{e} du th\'{e}or\`{e}me 1.5.1.

En utilisant les r\'{e}sultats d\'{e}j\`{a} obtenus dans ce chapitre
on obtient la proposition suivante.

\thm PROPOSITION 4.6.1
\enonce
Pour $\mathop{\rm deg}(D)$ assez grand il existe
\vskip 1mm

\item{-} un entier $N\geq 1$ et une partition $N=N_{1}+N_{2}$ en
entiers $N_{1},N_{2}>0$,
\vskip 1mm

\item{-} un point $a=(a_{1},a_{2})\in {\Bbb A}_{H}(k)$ de l'espace de
Hitchin pour le groupe endoscopique $H={\rm U}(N_{1})\times {\rm U}(N_{2})$ du
groupe unitaire ${\rm U}(N)$,
\vskip 1mm

\item{-} un point ferm\'{e} $z_{0}$ de l'intersection $Z=Y_{a_{1}}\cap
Y_{a_{2}}$ des courbes spectrales $Y_{a_{1}}$ et $Y_{a_{2}}$ pour
${\rm U}(N_{1})$ et ${\rm U}(N_{2})$, composantes de la courbe spectrale
$Y_{a}=Y_{a_{1}}+Y_{a_{2}}\subset\Sigma$ pour ${\rm U}(N)$,
\vskip 1mm

\noindent ayant les propri\'{e}t\'{e}s suivantes:
\vskip 1mm

\item{-} pour $\alpha =1,2$, $Y_{a_{\alpha}}$ et son rev\^{e}tement
double \'{e}tale $Y_{a_{\alpha}}'$ sont g\'{e}om\'{e}triquement
irr\'{e}ductibles,
\vskip 1mm

\item{-} $Y_{a}$ est \'{e}tale au-dessus du point $x_{\infty}$ de $X$ et,
pour $\alpha=1,2$ il existe au moins un point rationnel $y_{\infty
,\alpha}$ de $Y_{a_{\alpha}}$ au-dessus de $x_{\infty}$,
\vskip 1mm

\item{-} en tout point ferm\'{e} $z$ de $Z-\{z_{0}\}$, les deux
courbes spectrales $Y_{a_{1}}$ et $Y_{a_{2}}$ se coupent
transversalement,
\vskip 1mm

\item{-} le corps r\'{e}siduel de $z_{0}$ est $k$, $z_{0}$ est inerte
dans $Z'$ et a fortiori, si on note $x_{0}$ son image par la
projection $p:\Sigma\rightarrow X$, $x_{0}$ est inerte dans $X'$,
\vskip 1mm

\item{-} si on note $x_{0}'$ et $z_{0}'$ les uniques points ferm\'{e}s
de $X'$ et de $Z'$ au-dessus de $x_{0}$ et $z_{0}$, et si on fixe une
identification du compl\'{e}t\'{e} de l'anneau local ${\cal
O}_{X',x_{0}'}$ \`{a} ${\cal O}_{F'}$ compatible aux involutions
$\tau$, le compl\'{e}t\'{e} de la ${\cal O}_{X',x_{0}'}$-alg\`{e}bre
finie et plate ${\cal O}_{Y_{a_{\alpha}}',z_{0}'}$ est isomorphe \`{a}
la ${\cal O}_{F'}$-alg\`{e}bre ${\cal O}_{F'}[T]/(P_{I_{\alpha}}(T))$,
et ce de mani\`{e}re compatible aux involutions $\tau$.

\hfill\hfill$\square$
\endthm

Il r\'{e}sulte alors du th\'{e}or\`{e}me 3.10.6 qu'avec
les notations de la section (3.10), on a l'\'{e}galit\'{e}
$$\displaylines{
\quad\prod_{z\in Z^{{\rm inerte}}}\left|
{(\widetilde{A}_{z}^{\times})^{\tau^{\ast}=(-)^{-1}}
\over (A_{z}^{\times})^{\tau^{\ast}=(-)^{-1}}}\right|
\left(\sum_{m_{z}\in [{\cal M}(z)/P(z)](k)_{\sharp}}{\kappa
(\mathop{\rm cl}\nolimits_{{\cal M}(z)}(m_{z}))\over
|\mathop{\rm Aut}(m_{z})|}\right)
\hfill\cr\hfill
=\prod_{z\in Z^{{\rm inerte}}}(-|\kappa (z)|)^{r_{z}}
\left|{(\widetilde{A}_{z}^{\times})^{\tau^{\ast} =(-)^{-1}}\over
(A_{1,z}^{\times}\times A_{2,z}^{\times})^{\tau^{\ast}
=(-)^{-1}}}\right|
\left(\sum_{n_{z}\in [{\cal N}(z)/Q(z)] (k)_{\sharp}}{1\over
|\mathop{\rm Aut}(n_{z})|}\right)\quad}
$$
pour chacun des deux caract\`{e}res endoscopiques $\kappa :({\Bbb
Z}/2{\Bbb Z})^{2}\rightarrow\{\pm 1\}$. Pour conclure, nous allons:
\vskip 1mm

\item{-} montrer par un calcul explicite que, pour chaque $z\in
Z^{{\rm inerte}}-\{z_{0}\}$, les facteurs index\'{e}s par $z$ des deux
produits ci-dessus sont \'{e}gaux;
\vskip 1mm

\item{-} identifier les facteurs index\'{e}s par $z_{0}$ aux deux 
membres du th\'{e}or\`{e}me 1.5.1 que l'on cherche \`{a} 
d\'{e}montrer.
\vskip 3mm

Commen\c{c}ons par donner une autre reformulation de chacun des
facteurs des produits de l'\'{e}galit\'{e} ci-dessus.

Fixons $z\in Z^{{\rm inerte}}$ et notons $z'$ l'unique point ferm\'{e}
de $Z'$ au-dessus de $z$.  Notons $s=[\kappa(z):k]$; on a alors
$[\kappa(z'):k]=2s$.  Puisque $\kappa (z')$ est un corps, $\mathop{\rm
Frob}\nolimits_{k}^{{\Bbb Z}}$ agit transitivement sur $\mathop{\rm
Hom}\nolimits_{k}(\kappa (z'),\overline{k})$.  De plus, l'involution
induite par l'action de $\tau$ sur $\kappa (z')$ co\"{\i}ncide avec
l'action de $\mathop{\rm Frob} \nolimits_{k}^{s}$ sur $\mathop{\rm
Hom}\nolimits_{k}(\kappa (z'),\overline{k})$.

Quitte \`{a} choisir un point base $\xi_{0}$ dans $\mathop{\rm
Hom}\nolimits_{k}(\kappa (z'),\overline{k})$, on
peut identifier de fa\c{c}on unique cet ensemble \`{a}
$$
\mathop{\rm Hom}\nolimits_{k}(\kappa(z'),\overline{k})\cong
\{0,1,\ldots,2s-1\}\cong {\Bbb Z}/2s{\Bbb Z}
$$
de telle sorte que, via cette bijection, $\mathop{\rm Frob}
\nolimits_{k}$ et $\tau$ agissent sur ${\Bbb Z}/2s{\Bbb Z}$
par $i\mapsto i+1$ et par $i\mapsto i+s$.  Pour tout
$i=0,1,\ldots,2s-1$, on note $\xi_{i}$
l'\'{e}l\'{e}ment correspondant de $\mathop{\rm
Hom}\nolimits_{k}(\kappa (z'),\overline{k})$.

Soit $A=A_{z}$ le compl\'{e}t\'{e} de l'anneau local de la courbe $Y'$
en le point $z'$; $A$ est muni d'une involution $\tau$.  Soit
$B=B_{z}$ le compl\'{e}t\'{e} de l'anneau local de $Y$ en le point
$z$; $B$ est alors le sous-anneau des \'{e}l\'{e}ments $a\in A$ tels
que $\tau (a)=a$.  La $\overline{k}$-alg\`{e}bre $\overline{A}=A
\widehat{\otimes}_{k} \overline{k}$ se d\'{e}compose en
un produit
$$
\overline{A}\buildrel\sim\over\longrightarrow
\prod_{i=0}^{2s-1}\overline{A}_{i}
$$
o\`{u} $\overline{A}_{i}=A\widehat{\otimes}_{\kappa (z'),
\xi_{i}}\overline{k}$.

L'automorphisme $A\widehat{\otimes}_{k} \mathop{\rm Frob}
\nolimits_{k}$ de $\overline{A}=A\widehat{\otimes}_{k}\overline{k}$
s'\'{e}crit
$$
A\widehat{\otimes}_{k}\mathop{\rm Frob}\nolimits_{k}
(\alpha_{0},\ldots,\alpha_{2s-1})=({\varphi_{2s-1}}(\alpha_{2s-1}),
\varphi_{0}(\alpha_{0}),\ldots,\varphi_{2s-2}(\alpha_{2s-2}))
$$
o\`{u} les $\varphi_{i}:\overline{A}_{i}\rightarrow\overline{A}_{i+1}$
sont des bijections $\mathop{\rm Frob}\nolimits_{k}$-lin\'{e}aires.
Pour tout entier $m$ notons $\varphi_{i}^{m}$ le compos\'{e}
$$
\varphi_{i}^{m}=\varphi_{i+m-1}\circ\cdots\circ\varphi_{i+1}\circ
\varphi_{i}:\overline{A}_{i}\rightarrow \overline{A}_{i+m}.
$$
En particulier, $\varphi_{0}^{2s}:\overline{A}_{0} \rightarrow
\overline{A}_{0}$ est l'automorphisme $\mathop{\rm Frob}
\nolimits_{k}^{2s}$-lin\'{e}aire $A\widehat{\otimes}_{\kappa
(z'),\xi_{0}}\mathop{\rm Frob}\nolimits_{k}^{2s}$, de sorte qu'on peut
identifier $A$ au sous-anneau de $\overline{A}_{0}$ des
\'{e}l\'{e}ments $a_{0}\in \overline{A}_{0}$ qui v\'{e}rifient
$\varphi_{0}^{2s}(a_{0})=a_{0}$.

L'automorphisme $\tau\widehat{\otimes}_{k}\overline{k}$ de
$\overline{A}$ se d\'{e}compose en le produit des isomorphismes
$\overline{k}$-lin\'{e}aires $\tau_{i}:\overline{A}_{i}\rightarrow
\overline{A}_{i+s}$ qui v\'{e}rifient les relations 
$\tau_{s+i}\circ\varphi_{i}^{s}=\varphi_{i+s}^{s}\circ \tau_{i}$.
En particulier, on a l'automorphisme $\mathop{\rm
Frob}_{k}^{s}$-lin\'{e}aire
$$
\overline{\tau}_{0}=\varphi_{s}^{s}\circ\tau_{0}=\tau_{s}\circ
\varphi_{0}^{s}:\overline{A}_{0}\rightarrow \overline{A}_{0}
$$
tel que $(\overline{\tau}_{0})^{2}=\varphi_{0}^{2s}$, de sorte qu'on
peut aussi identifier $B$ au sous-anneau de $\overline{A}_{0}$
form\'{e} des \'{e}l\'{e}ments $a_{0}\in \overline{A}_{0}$ tels que
$\overline{\tau}_{0}(a_{0})=a_{0}$.

L'anneau total des fractions ${\rm Frac}(B)$ est un produit de
corps
$$
\mathop{\rm Frac}(B)=E_{J}=\prod_{j\in J}E_{j}.
$$
Pour tout $j\in J$, notons $\kappa_{j}$ le corps r\'{e}siduel
de $E_{j}$; $\kappa_{j}$ contient naturellement $\kappa (z)$. On a
$$
\mathop{\rm Frac}(A)=E_{J}'=\prod_{j\in J}E_{j}'
$$
o\`{u}, pour chaque $j\in J$, $E_{j}'=E_{j}\widehat{\otimes}_{\kappa
(z)} \kappa (z')$ reste un corps de corps r\'{e}siduel
$\kappa_{j}'=\kappa_{j}\otimes_{\kappa (z)}\kappa (z')$.  En effet,
dans le cas $z=z_{0}$, c'est notre hypoth\`{e}se de d\'{e}part
rappel\'{e}e dans la section (4.1).  Dans le cas $z\not= z_{0}$, le
point $z$ est un point double de $Y$ \`{a} branches rationnelles sur
$\kappa (z)$, de sorte que l'ensemble $J$ a deux \'{e}l\'{e}ments $j$
et que pour chacun de ces deux \'{e}l\'{e}ments on a
$\kappa_{j}=\kappa (z)$.  Pour tout $j\in J$, notons
$s_{j}=\mathop{\rm deg}[\kappa_{j}:k]$; alors $s_{j}$ est divisible
par $s$ et $s_{j}/s$ est un nombre impair; autrement dit, on a
$$
s_{j}\equiv s\hbox{ (mod }2s).
$$

Pour chaque $j\in J$ on choisit un plongement $\zeta_{j,0}:
\kappa_{j}'\hookrightarrow\overline{k}$ qui prolonge le plongement
$\xi_{0}:\kappa (z')\hookrightarrow \overline{k}$. Ce choix induit
une bijection
$$
{\Bbb Z}/2s_{j}{\Bbb Z}\cong \mathop{\rm Hom}\nolimits_{k}
(\kappa_{j}',\overline{k}),~i_{j}\mapsto \zeta_{j,i_{j}}.
$$
qui rend commutatif le carr\'{e}
$$
\xymatrix{
\mathop{\rm Hom}\nolimits_{k}(\kappa_{j}',\overline{k}) 
\ar[r] \ar[d]& \mathop{\rm Hom}\nolimits_{k}(\kappa (z'),
\overline{k}) \ar[d] \cr
{\Bbb Z}/2s_{j}{\Bbb Z} \ar[r] & {\Bbb Z}/2s{\Bbb Z}}
$$
dont la fl\`{e}che horizontale du haut est la restriction de
$\kappa_{j}'$ \`{a} $\kappa (z')$.

Consid\'{e}rons l'ensemble
$$
\overline{J}=\{(j,i_{j})\mid j\in J,~i_{j}\in{\Bbb Z}/2s_{j}{\Bbb
Z}\}.
$$
On a alors
$$
\mathop{\rm Frac}(\overline{A})=\prod_{(j,i_{j})\in\overline{J}}
\overline{E}_{j,i_{j}}
$$
o\`{u} $\overline{E}_{j,i_{j}}=E_{j}'\widehat{\otimes}_{\kappa_{j}',
\zeta_{j,i_{j}}}\overline{k}$; chaque \'{e}l\'{e}ment $a\in 
\mathop{\rm Frac}(\overline{A})$ s'\'{e}crit donc sous
la forme $a=(a_{j,i_{j}})_{(j,i_{j})\in\overline{J}}$. Par ailleurs,
on a la d\'{e}composition
$$
\mathop{\rm Frac}(\overline{A}_{i})=\prod_{{\scriptstyle 
{j\in J}\atop\scriptstyle  i_{j}\equiv i~({\rm mod}~2s)}}
\overline{E}_{j,i_{j}},
$$
de sorte que chaque \'{e}l\'{e}ment $a\in \mathop{\rm Frac}
(\overline{A})$ s'\'{e}crit aussi sous la forme $a=(a_{i})_{i\in{\Bbb
Z}/2s{\Bbb Z}}$ avec $a_{i}\in\mathop{\rm Frac}(\overline{A}_{i})$,
quitte \`{a} regrouper les termes
$$
a_{i}=(a_{j,i_{j}})_{(j,i_{j})\in\overline{J},~i_{j}\equiv
i~({\rm mod}~2s)}.
$$

L'automorphisme de $\mathop{\rm Frac}(\overline{A})$ qui se d\'{e}duit 
de $A\widehat{\otimes}_{k}\mathop{\rm Frob}\nolimits_{k}:\overline{A}
\rightarrow\overline{A}$, se d\'{e}compose en le produit des 
homomorphismes $\mathop{\rm Frob}\nolimits_{k}$-lin\'{e}aires
$\varphi_{j,i_{j}}:\overline{E}_{j,i_{j}}\rightarrow {\overline
E}_{j,i_{j}+1}$. De m\^{e}me, $\tau:\mathop{\rm Frac}(\overline{A})
\rightarrow\mathop{\rm Frac}(\overline{A})$ se d\'{e}compose en le
produit des homomorphismes $\overline{k}$-lin\'{e}aires
$\tau_{j,i_{j}}:\overline{E}_{j,i_{j}}\rightarrow
\overline{E}_{j,i_{j}+s_{j}}$.
\vskip 2mm

Soit $\widetilde{P}(z)(\overline{k})\subset \mathop{\rm Frac}
(\overline{A})^{\times}$ le sous-groupe form\'{e} des $a\in
\mathop{\rm Frac}(\overline{A})^{\times}$ qui v\'{e}rifient $\tau
(a)a=1$.  Par d\'{e}finition $P(z)(\overline{k})$ est le quotient de 
$\widetilde{P}(z)(\overline{k})$ par son sous-groupe 
$\overline{A}{}^{\times}\cap \widetilde{P}(z)(\overline{k})$. 

On v\'{e}rifie que $\widetilde{P}(z)(\overline{k})\subset \mathop{\rm
Frac} (\overline{A})^{\times}$ est stable par $\mathop{\rm
Frob}\nolimits_{k}$.  Le groupe $\widetilde{P}(z)(\overline{k})$ agit
donc sur lui-m\^{e}me par la $\mathop{\rm Frob}
\nolimits_{k}$-conjugaison $(b,a)\rightarrow \mathop{\rm Frob}
\nolimits_{k}(b)ab^{-1}=\mathop{\rm Frob}\nolimits_{k}(b)b^{-1}a$
puisque $\widetilde{P}(z)(\overline{k})$ est commutatif.  Cette
$\mathop{\rm Frob}\nolimits_{k}$-conjugaison sur
$\widetilde{P}(z)(\overline{k})$ induit bien s\^{u}r par passage au
quotient la $\mathop{\rm Frob}\nolimits_{k}$-conjugaison sur
$P(z)(\overline{k})$.

\thm PROPOSITION 4.6.2
\enonce
L'homomorphisme
$$
\widetilde{P}(z)(\overline{k})\rightarrow ({\Bbb Z}/2{\Bbb Z})^{J}
$$
qui associe \`{a} $a=(a_{j,i_{j}})_{(j,i_{j})\in \overline{J}}$ 
l'\'{e}l\'{e}ment $(\lambda_{j})_{j\in J}$, o\`{u} $\lambda_{j}\in
{\Bbb Z}/2{\Bbb Z}$ est l'unique \'{e}l\'{e}ment tel que
$$
\sum_{i_{j}=0}^{s_{j}-1}\mathop{\rm val}
\nolimits_{\overline{E}_{j,i_{j}}}(a_{j,i_i})\equiv 
\lambda_{j}\hbox{ {\rm (}mod }2),
$$
induit un isomorphisme du groupe des classes de $\mathop{\rm
Frob}\nolimits_{k}$-conjugaison de $\widetilde{P}(z) (\overline{k})$
sur le groupe $({\Bbb Z}/2{\Bbb Z})^{J}$.  

Cet isomorphisme induit \`{a} son tour un isomorphisme du groupe des
classes de $\mathop{\rm Frob}\nolimits_{k}$-conjugaison de
$P(z)(\overline{k})$ sur le groupe $({\Bbb Z}/2{\Bbb Z})^{J}$.
\endthm

\rem D\'{e}monstration
\endrem
Soit $\overline {\widetilde A}$ le normalis\'{e} de $\overline{A}$
dans ${\rm Frac}(\overline{A})$ et soit $ P^{0}(z)(\overline{k})$ le
groupe des \'{e}l\'{e}ments $a\in {\overline{\widetilde A}}{}^\times$
tels que $\tau (a)a=1$.  Un th\'{e}or\`{e}me de Lang assure que
l'isog\'{e}nie $P^{0}(z)(\overline{k})\rightarrow
P^{0}(z)(\overline{k})$ d\'{e}finie par $b\mapsto \mathop{\rm Frob}
\nolimits_{k}(b)b^{-1}$ est surjective.  On en d\'{e}duit que
l'homomorphisme canonique du groupe des classes de $\mathop{\rm
Frob}\nolimits_{k}$-conjugaison de $\widetilde{P}(z)(\overline{k})$
sur le groupe des classes de $\mathop{\rm Frob}
\nolimits_{k}$-conjugaison de $\widetilde P(z)(\overline{k})/
P^{0}(z)(\overline{k})$, est un isomorphisme.

Le groupe $\widetilde{P}(z)(\overline{k})/P^{0}(z)(\overline{k})$ est
canoniquement isomorphe au groupe
$$
\{d=(d_{j,i_{j}})_{(j,i_{j})\in\overline{J}}\in{\Bbb Z}^{\overline{J}}
\mid d_{j,i_{j}}+d_{j,i_{j}+s_{j}}=0\}.
$$
On a donc un isomorphisme $\widetilde P(z)(\overline{k})/
P^{0}(z)(\overline{k})\cong {\Bbb Z}^{|\overline{J}|/2}$ en ne gardant
parmi les composantes de $d$ que les composantes $d_{j,i_{j}}$ o\`{u}
$i_{j}=0,1,\ldots,s_{j}-1$, les autres composantes $d_{j,i_{j}}$
o\`{u} $i_{j}=s_{j},s_{j}+1,\ldots,2s_{j}-1$ \'{e}tant uniquement
d\'{e}termin\'{e}es par les relations
$d_{j,i_{j}}+d_{j,i_{j}+s_{j}}=0$.

L'endomorphisme $\mathop{\rm Frob}\nolimits_{k}$ agit sur le groupe
$\widetilde{P}(z)(\overline{k})/P^{0}(z)(\overline{k})$ par
$\mathop{\rm Frob}\nolimits_{k}(d)=d'$ avec
$d'_{j,i_{j}}=d_{j,i_{j}-1}$.  Dans l'\'{e}criture ${\Bbb 
Z}^{|\overline{J}|/2}$ pour ce groupe, on a donc
$$
\mathop{\rm Frob}\nolimits_{k}((d_{j,0},d_{j,1},
\ldots,d_{j,s_{j}-1})_{j\in J})=
(-d_{j,s_{j}-1},d_{j,0},\ldots,d_{j,s_{j}-2})_{j\in J}.
$$
Par suite l'isog\'{e}nie de Lang 
$$
{\cal L}_{k}:{\Bbb Z}^{|\overline{J}|/2}
\rightarrow {\Bbb Z}^{|\overline{J}|/2},~d\mapsto 
\mathop{\rm Frob}\nolimits_{k}(d)d^{-1},
$$
est d\'{e}finie par
$$
{\cal L}_{k}((d_{j,0},d_{j,1},\ldots,d_{j,s_{j}-1})_{j\in J})=
(-d_{j,s_{j}-1}-d_{j,0},d_{j,0}-d_{j,1},\ldots,
d_{j,s_{j}-2}-d_{j,s_{j}-1})_{j\in J}.
$$
et a pour conoyau la fl\`{e}che ${\Bbb Z}^{|\overline{J}|/2} 
\twoheadrightarrow ({\Bbb Z}/2{\Bbb Z})^{J}$ qui envoie 
$(d_{j,i_{j}})_{j\in J,i_{j}=0,1,\ldots,s_{j}-1}$ sur 
$(\lambda_{j})_{j\in J}$ o\`{u} $\lambda_{j}\in {\Bbb Z}/2{\Bbb Z}$
est la classe modulo $2$ de la somme 
$\sum_{i_{j}=0}^{s_{j}-1}d_{j,i_{j}}$.

Puisque l'homomorphisme $(\overline{A}^{\times})^{\tau=-1}\rightarrow
(\overline{A}^{\times})^{\tau=-1}$ d\'{e}fini par $b\mapsto\mathop{\rm 
Fob}\nolimits_{k}(b)b^{-1}$ est surjectif, l'homomorphisme
$$
{\widetilde{P}}(z)(\overline{k})\rightarrow P(z)(\overline{k})=
{\widetilde{P}}(z)(\overline{k})/(\overline{A}^\times)^{\tau=-1}
$$
induit un isomorphisme entre les groupes des classes de $\mathop{\rm 
Fob}\nolimits_{k}$-conjugaisons et la proposition est d\'{e}montr\'{e}e.
\hfill\hfill $\square$
\vskip 3mm

Pour chaque $\lambda\in ({\Bbb Z}/2{\Bbb Z})^{J}$, on choisit un
repr\'{e}sentant $a(\lambda )\in P(z)(\overline{k})$ de la classe de
$\mathop{\rm Frob}\nolimits_{k}$-conjugaison correspondante comme
suit.  On choisit d'abord une fois pour toutes un repr\'{e}sentant
arbitraire $\widetilde\lambda \in {\Bbb Z}^{J}$ de la classe
$\lambda\in ({\Bbb Z}/2{\Bbb Z})^{J}$.  Puis, pour chaque $j\in J$, on
choisit une uniformisante $\varpi_{j}$ de $E_{j}$.  Cette
uniformisante d\'{e}finit une uniformisante $\varpi_{j}'=\varpi_{j}$
de l'extension non ramifi\'{e}e $E_{j}'$ de $E_{j}$ et une
uniformisante $\varpi_{j,i_{j}}=\varpi_{j}$ de l'extension non
ramifi\'{e}e $\overline{E}_{j,i_{j}}$ pour chaque
$j=0,1,\ldots,2s_{j}-1$.  Puisque $\varpi_{j,0}$ provient de $E_{j}$,
on a
$$
\overline{\tau}_{0}(\varpi_{j,0})=\varpi_{j,0}
$$
et a fortiori $\varphi_{0}^{2s_{j}}(\varpi_{j,0})=\varpi_{j,0}$. On 
d\'{e}finit alors l'\'{e}l\'{e}ment $a(\lambda )\in\widetilde{P}(z)
(\overline{k})$ de coordonn\'{e}es $a_{j,i_{j}}(\lambda)\in
\overline{E}_{j,i_{j}}^\times$ comme suit. Pour tout 
$(j,i_{j})\in\overline{J}$ 
on pose
$$
a_{j,i_{j}}(\lambda)=\openup\jot\cases{1 & si $i_{j}\not\equiv
0\hbox{ (mod }s)$,\cr
\varpi_{j,i_{j}}^{\widetilde\lambda_{j}} & si $i_{j}\equiv
0\hbox{ (mod }2s)$,\cr
\varpi_{j,i_{j}}^{-\widetilde\lambda_{j}} & si $i_{j}\equiv
s\hbox{ (mod }2s)$.\cr}
$$
Comme le nombre $s_{j}/s$ des classes de conjugaison modulo
$2s_{j}$ qui sont congrues \`{a} $0$ modulo $2s$ est impair, l'image
de de $a(\lambda )$ dans $({\Bbb Z}/2{\Bbb Z})^{J}$ est bien
$\lambda$.

Regroupant les termes $a_{j,i_{j}}(\lambda)$ suivant la classe modulo
$2s$ de $i_{j}$, on obtient
$$
a(\lambda )=(a_{0}(\lambda),a_{1}(\lambda),\ldots,
a_{2s-1}(\lambda))=(c_{\lambda},\overbrace{1,\ldots ,1}^{s-1}
,\tau_{0}(c_{\lambda})^{-1},\overbrace{1,\ldots ,1}^{s-1})
$$
o\`{u} $a_{i}(\lambda )\in\mathop{\rm Frac}(\overline{A}_{i})$ pour 
chaque $i=0,1,\ldots ,2s-1$ et
$$
c_{\lambda}=(\varpi_{j,i_{j}}^{\widetilde\lambda_{j}})_{j\in J,~i_{j}
\equiv 0~({\rm mod}~2s)}
$$
est un \'{e}l\'{e}ment de $\mathop{\rm
Frac}(\overline{A}_{0})^{\times}$ qui v\'{e}rifie
$\overline{\tau}_{0}(c_{\lambda})=c_{\lambda}$.  

Rappelons que $\mathop{\rm Frac}(B)$ s'identifie au sous-anneau de
$\mathop{\rm Frac}(\overline{A}_{0})$ form\'{e} des \'{e}l\'{e}ments
fixes sous l'auto\-mor\-phisme $\mathop{\rm Frob}
\nolimits_{k}^{s}$-lin\'{e}aire $\overline{\tau}_{0}$.  On peut donc
voir $c_{\lambda}$ comme un \'{e}l\'{e}ment de $\mathop{\rm Frac}(B)$ 
qui n'est autre que
$$
c_{\lambda}=(\varpi_{j}^{\widetilde{\lambda}_{j}})_{j\in J}
$$
dans la d\'{e}composition $\mathop{\rm Frac}(B)=\prod_{j\in J}E_{j}$.
\vskip 2mm

Les objets de la cat\'{e}gorie $[{\cal M}(z)/P(z)](k)$ sont les
couples $({\cal V},a)$ o\`{u} ${\cal V}\in {\cal M}(z)(\overline{k})$
et $a\in P(z)(\overline{k})$ satisfont \`{a} le relation $\mathop{\rm
Frob} \nolimits_{k}({\cal V})=a\cdot {\cal V}$.  \'{E}tant donn\'{e}s
deux objets $({\cal V},a)$ et $({\cal V}',a')$ de $[{\cal
M}(z)/P(z)](k)$ comme ci-dessus, un isomorphisme du premier dans le
second est un \'{e}l\'{e}ment $b\in P(z)(\overline{k})$ tel que ${\cal
V}'=b\cdot {\cal V}$ et que
$$
a'=\mathop{\rm Frob}\nolimits_{k}(b)ab^{-1}.
$$
Les automorphismes d'un objet $({\cal V},a)$ de $[{\cal
M}(z)/P(z)](k)$ sont donc les \'{e}l\'{e}ments $b\in P(z)(k)$ tels que
$b{\cal V}={\cal V}$.

La proposition pr\'{e}c\'{e}dente admet le corollaire suivant.

\thm COROLLAIRE 4.6.3
\enonce 
La cat\'{e}gorie $[{\cal M}(z)/P(z)](k)$ est \'{e}quivalente \`{a} sa
sous-cat\'{e}gorie pleine dont les objets sont les $({\cal V},a (\lambda))\in
\mathop{\rm ob}[{\cal M}(z)/P(z)](k)$ pour $\lambda\in ({\Bbb Z}/2{\Bbb 
Z})^{J}$ et dont l'ensemble des morphismes d'un objet $({\cal V},a(\lambda))$ 
vers un autre objet $({\cal V}',a(\lambda'))$ est
$$
\mathop{\rm Hom}(({\cal V},a(\lambda )),({\cal V}',a(\lambda')))=
\openup\jot\cases{\{b\in P(z)(k)\mid b\cdot {\cal V}={\cal V}'\} & 
si $\lambda=\lambda'$,\cr
\emptyset & sinon.\cr}
$$
\endthm

\rem D\'{e}monstration
\endrem
La proposition pr\'{e}c\'{e}dente implique que pour tout objet $({\cal
V},a)$ de $[{\cal M}(z)/P(z)](k)$, il existe $\lambda\in ({\Bbb
Z}/2{\Bbb Z})^{J}$ et $b\in P(z)(\overline{k})$ tel que
$a(\lambda)=\mathop{\rm Frob}\nolimits_{k}(b)ab^{-1}$, et que de plus
$\lambda$ est uniquement d\'{e}termin\'{e} par cette relation.
\hfill\hfill $\square$
\vskip 3mm

On se propose maintenant de d\'{e}crire la sous-cat\'{e}gorie pleine
de $[{\cal M}(z)/P(z)](k)$ des objets $({\cal V},a)$ avec
$a=a(\lambda)$ pour un $\lambda\in({\Bbb Z}/2{\Bbb Z})^J$ fix\'{e}.

Dans la section (3.10) on a pris l'objet de Kostant de ${\cal
N}(k)$  et son image dans ${\cal M}(k)$ comme points base pour
d\'{e}finir les variantes locales ${\cal M}(z)$ et ${\cal N}(z)$
de ${\cal M}$ et ${\cal N}$. L'ensemble $J$ des branches de $Y$ au
point ferm\'{e} $z$ se d\'{e}coupe en deux parties disjointes
$J=J_{1}\amalg J_{2}$ o\`{u} $J_{\alpha}$ est l'ensemble des branches de
$Y_{\alpha}$ au point ferm\'{e} $z$. On a donc les d\'{e}compositions
$$
E_{J}=E_{J_{1}}\times E_{J_{2}}
$$
o\`{u} $E_{J_{\alpha}}=\prod_{j\in J_{\alpha}}E_{j}$ pour tout
$\alpha=1,2$.  Notons $A_{\alpha}$ le compl\'{e}t\'{e} de l'anneau
local de $Y_{\alpha}'$ en $z'$; l'anneau total des fractions de
$A_{\alpha}$ est bien entendu $\mathop{\rm
Frac}(A_{\alpha})=E_{J_{\alpha}}$.

L'objet de Kostant de ${\cal N}$, restreint au compl\'{e}t\'{e} formel
de $Y'$ en $z'$, est la somme directe de deux couples $({\cal
K}_{\alpha},\iota_{{\cal K}_\alpha})$, un pour chaque $\alpha=1,2$,
o\`{u} ${\cal K}_\alpha$ est un $A_{\alpha}$-module libre de rang $1$ et
$\iota_{{\cal K}_\alpha}$ est une structure unitaire
$$
\iota_{{\cal K}_\alpha}:\tau^{\ast}{\cal K}_{\alpha}\buildrel\sim\over
\longrightarrow {\cal K}_{\alpha}^{\vee}={\rm Hom}_{A_{\alpha}}({\cal
K}_{\alpha},\omega_{\alpha}).
$$
On a not\'{e} ici $\omega_{\alpha}$ le compl\'{e}t\'{e} formel de la
fibre en $z'$ du dualisant relatif $\omega_{Y'_{\alpha}/X'}$ de
$Y'_{\alpha}/X'$; c'est un sous-$A_{\alpha}$-module de
$E_{J_{\alpha}}$, libre de rang $1$. Fixons un g\'{e}n\'{e}rateur 
$c_{J_{\alpha}}^{0}$ de $\omega_{\alpha}$ comme dans la section (1.4).

Le couple $({\cal K}_{\alpha}, \iota_{{\cal K}_{\alpha}})$ est
isomorphe au couple $(A_{\alpha},c_{J_{\alpha}}^{0})$ o\`{u} on note
encore $c_{J_{\alpha}}^{0}$ la structure unitaire
$$
\tau^{\ast}A_{\alpha}=A_{\alpha}\buildrel\sim\over\longrightarrow
{\rm Hom}_{A_{\alpha}}(A_{\alpha},\omega_{\alpha})=
\omega_{\alpha},~1\rightarrow c_{J_{\alpha}}^{0}.
$$
Fixons un tel isomorphisme ce qui nous permet de nous raccrocher au
langage commode de la section (1.4).  Soit $F'$ l'anneau local
compl\'{e}t\'{e} de $X'$ en le point ferm\'{e} $x'$ image de $z'$.
Consid\'{e}rons la forme hermitienne
$$
\Phi_{c_{J_{\alpha}}^{0}}:E'_{J_{\alpha}}\times E'_{J_{\alpha}}
\rightarrow F'
$$
d\'{e}finis par
$$
(x,y)\mapsto \Phi_{c_{J_{\alpha}}^{0}}(x,y)={\rm
Tr}_{E'_{J_{\alpha}}/F'}(c_{J_{\alpha}}^{0} \tau(x)y).
$$
D'apr\`{e}s le lemme 1.4.2, $A_{\alpha}$ est un r\'{e}seau auto-dual
pour cette forme hermitienne.  La fibre g\'{e}n\'{e}rique de $({\cal
K}_1\oplus {\cal K}_2,\iota_{{\cal K}_1}\oplus \iota_{{\cal K}_2})$
s'identifie \`{a} $E_{J}'=E_{J_{1}}'\times E_{J_{2}}'$ muni de la
forme unitaire $c_{\cal N}^{0}=c_{J_{1}}^{0}\times c_{J_{2}}^{0}$.

La forme hermitienne $c_{{\cal N}}^{0}$ induit forme hermitienne
$\overline{c}{}_{{\cal N}}^{\,0}$ sur
$$
E_{J}'\widehat{\otimes}_{k}\overline{k}
\buildrel\sim\over\longrightarrow\prod_{i=0}^{2s-1}
\mathop{\rm Frac}(\overline{A}_{i})
$$
qui \`{a} son tour induit une forme bilin\'{e}aire
$$
\overline{c}{}_{{\cal N},i}^{\,0}:\mathop{\rm Frac}(\overline{A}_{i})
\times \mathop{\rm Frac}(\overline{A}_{i+s})\rightarrow F
\widehat{\otimes}_{\kappa (x'),\xi_{i}}\overline{k}
$$
pour chaque $i=0,1,\ldots,2s-1$.  Les \'{e}l\'{e}ments de ${\cal
M}(z)(\overline{k})$ sont alors les collections ${\cal V}=({\cal
V}_{i})_{i}$ de sous-$\overline{A}_{i}$-modules ${\cal
V}_{i}\subset\mathop{\rm Frac}(\overline{A}_{i})$ tels que
$\mathop{\rm Frac} (\overline{A}_{i}){\cal V}_{i}=\mathop{\rm
Frac}(\overline{A}_{i})$ et que les r\'{e}seaux ${\cal V}_{i}$ et
${\cal V}_{i+s}$ sont orthogonaux par rapport \`{a} la forme
bilin\'{e}aire $\overline{c}{}_{{\cal N},i}^{\,0}$.

L'endomorphisme de Frobenius agit sur cet ensemble par
$$
\mathop{\rm Frob}\nolimits_{k}({\cal V}_{0},\ldots,{\cal
V}_{2s-1})=(\varphi_{2s-1}({\cal V}_{2s-1}),
\varphi_{0}({\cal V}_{0}),\ldots,\varphi_{2s-2}({\cal V}_{2s-2})).
$$

\thm PROPOSITION 4.6.4
\enonce 
Soit $\lambda\in({\Bbb Z}/2{\Bbb Z})^{J}$ et $a(\lambda)\in
P(z)(\overline{k})$ le repr\'{e}sentant de la classe de $\mathop{\rm
Frob}\nolimits_{k}$-conjugaison index\'{e}e par $\lambda$ construit
pr\'{e}c\'{e}demment.  La sous-cat\'{e}gorie pleine de $[{\cal
M}(z)/P(z)](k)$ des couples $({\cal V},a)$ avec $a=a(\lambda)$, est
\'{e}quivalente \`{a} la cat\'{e}gorie $[{\cal M}_{\lambda}(k)/P(z)(k)]$
dont:
\vskip 1mm

\item{-} les objets sont les sous-$A$-modules $M$ de $E_{J}'$ tels que
$E_{J}'M=E_{J}'$ qui sont auto-duaux par rapport \`{a} la forme
hermitienne
$$
\Phi_{c_{\lambda}c_{{\cal N}}^{0}}(x,y)={\rm Tr}_{E_{J}'/F'}
(c_{\lambda}c_{{\cal N}}^{0}\tau (x)y)
$$ 
o\`{u} $c_{\lambda}=(\varpi_{j}^{\widetilde{\lambda}_{j}})_{j\in J}$
et o\`{u} $\widetilde{\lambda}\in{\Bbb Z}^{J}$ est le rel\`{e}vement
de $\lambda$ qu'on a choisi pour d\'{e}finir $a(\lambda)$;
\vskip 1mm

\item{-} les morphismes d'un objet $M$ vers un objet $M'$ sont les
\'{e}l\'{e}ments $b\in P(z)(k)$ tel que $bM=M'$.

\endthm

\rem D\'{e}monstration
\endrem
Soit $({\cal V},a)$ un objet de $[{\cal M}(z)/P(z)](k)$ avec
$$
a=a(\lambda)=(c_{\lambda},\overbrace{1,\ldots ,1}^{s-1}
,\tau_{0}(c_{\lambda})^{-1},\overbrace{1,\ldots ,1}^{s-1})
$$
pour un $\lambda\in({\Bbb Z}/2{\Bbb Z})^{J}$.  La relation
$\mathop{\rm Frob}\nolimits_{k}({\cal V})=a\cdot {\cal V}$ est alors
\'{e}quivalente aux relations suivantes ${\cal
V}_{i}=\varphi_{0}^{i}({\cal V}_{0})$ pour tout $i=1,\ldots ,s-1$,
${\cal V}_{i}=\varphi_{s}^{i-s}(\tau_{0}(c_{\lambda}))
\varphi_{0}^{i}({\cal V}_{0})$ pour tout $i=s,\ldots ,2s-1$ et
$c_{\lambda}{\cal V}_{0}=\overline{\tau}_{0}(c_{\lambda})
\varphi_{0}^{2s}({\cal V}_{0})$.  Or, par construction, on a
$c_{\lambda}=\overline{\tau}_{0}(c_{\lambda})$, de sorte qu'on a
$$
{\cal V}_{0}=\varphi_{0}^{2s}({\cal V}_{0}).
$$
Soit $M$ l'ensemble des \'{e}l\'{e}ments de ${\cal V}_{0}$ qui sont
fix\'{e}s par $\varphi_{0}^{2s}$.  Alors $M$ est un sous-$A$-module de
$E_{J}'=\{a\in \mathop{\rm Frac}(\overline{A}_{0})\mid
\varphi_{0}^{2s}(a)=a\}$, sous-module qui v\'{e}rifie 
$E_{J}'M=E_{J}'$.

De plus, ${\cal V}_{0}$ et ${\cal V}_{s}=\overline{\tau}_{0}
(c_{\lambda})\varphi_{0}^{s}({\cal V}_{0})$ sont orthogonaux par
rapport \`{a} la forme bilin\'{e}aire
$$
\overline{c}{}_{{\cal N},i}^{\,0}:\mathop{\rm Frac}(\overline{A}_{i})
\times \mathop{\rm Frac}(\overline{A}_{i+s})\rightarrow 
F\widehat{\otimes}_{\kappa (x'),\xi_{i}}\overline{k}.
$$
si et seulement si $M$ est auto-dual par rapport \`{a} la forme
hermitienne $\Phi_{c_{\lambda}c_{{\cal N}}^{0}}$. La proposition est 
donc d\'{e}montr\'{e}e.
\hfill\hfill$\square$
\vskip 3mm

Traitons maintenant s\'{e}par\'{e}ment les cas $z\not=z_{0}$ et le 
cas $z=z_{0}$
\vskip 2mm

Dans le cas $z\not= z_{0}$, les deux courbes $Y_{a_{1}}$ et $Y_{a_{2}}$
sont lisses en le point $z$ et leur intersection y est transversale. 
On a $J=\{1,2\}$, $J_{1}=\{1\}$, $J_{2}=\{2\}$ et les uniformisantes 
$\varpi_{1}$ et $\varpi_{2}$ que l'on a fix\'{e} ci-dessus sont 
des uniformisante de $Y_{a_{1}}$ et $Y_{a_{2}}$ en $z$. En fait, on a 
$E_{\alpha}'=\kappa (z')((\varpi_{\alpha}))$,
$$
A=\kappa (z')[[\varpi_{1},\varpi_{2}]]/(\varpi_{1},\varpi_{2})\subset
\kappa (z')((\varpi_{1}))\times \kappa (z')((\varpi_{2})) =\mathop{\rm
Frac}(A)=E_{J}'
$$
et $\tau$ est le rel\`{e}vement canonique de l'\'{e}l\'{e}ment non
trivial de $\mathop{\rm Gal}(\kappa (z')/\kappa (z))$. On a sur
$E_{J}'$ la forme hermitienne $\Phi_{c_{{\cal N}}^{0}}$ pour
laquelle le normalis\'{e}
$$
\widetilde{A}=\kappa (z')[[\varpi_{1}]]\times\kappa (z')[[\varpi_{2}]]
\subset E_{J}'
$$
de $A$ dans $E_{J}'$ est un r\'{e}seau auto-dual.

Le groupe $({\Bbb Z}/2{\Bbb Z})^{J}$ a quatre \'{e}l\'{e}ments
$\lambda$ et on v\'{e}rifie que:
\vskip 1mm

\item{-} si $\lambda=(0,1)\hbox{ ou }(1,0)$, il n'y a aucun
sous-$A$-r\'{e}seau de $E_{J}'$ qui soit auto-dual par rapport \`{a}
la forme hermitienne $\Phi_{c_{\lambda}c_{{\cal N}}^{0}}$;
\vskip 1mm

\item{-} si $\lambda=(0,0)$, $\widetilde{A}$ est le seul
sous-$A$-r\'{e}seau de $E_{J}'$ qui est auto-dual par rapport \`{a} la
forme hermitienne $\Phi_{c_{\lambda}c_{{\cal N}}^{0}}$; de plus, pour tout 
$b\in P(z)(k)$, on a $b\widetilde{A}=\widetilde{A}$;
\vskip 1mm

\item{-} si $\lambda=(1,1)$, les sous-$A$-r\'{e}seaux de $E_{J}'$ qui
sont auto-duaux par rapport \`{a} la forme hermitienne
$\Phi_{c_{\lambda}c_{{\cal N}}^{0}}$ forment un espace principal
homog\`{e}ne sous l'action du groupe $P(z)(k)$.
\vskip 1mm

Comme le groupe $P(z)(k)$ a $|\kappa (z)|+1$ \'{e}l\'{e}ments et que 
$\kappa (0,0)=1$ et $\kappa (1,1)=-1$, on en d\'{e}duit que
$$
\sum_{m_{z}\in [{\cal M}(z)/P(z)](k)_{\sharp}}{\kappa
(\mathop{\rm cl}\nolimits_{{\cal M}(z)}(m_{z}))\over |\mathop{\rm
Aut}(m_{z})|} =-{|\kappa(z)|\over |\kappa(z)|+1}
$$
et que
$$
\left|{(\widetilde{A}_{z}^{\times})^{\tau^{\ast}=(-)^{-1}}
\over (A_{z}^{\times})^{\tau^{\ast}=(-)^{-1}}}\right|=|\kappa(z)|+1.
$$

On v\'{e}rifie facilement que $r_{z}=1$, que
$$
\left|{(\widetilde{A}_{z}^{\times})^{\tau^{\ast}
=(-)^{-1}}\over (A_{1,z}^{\times}\times
A_{2,z}^{\times})^{\tau^{\ast} =(-)^{-1}}}\right|=1
$$
et que
$$
\sum_{n_{z}\in [{\cal N}(z)/Q(z)] (k)_{\sharp}}{1\over
|\mathop{\rm Aut}(n_{z})|}=1.
$$

On a donc d\'{e}montr\'{e} que
$$\displaylines{
\qquad\left| {(\widetilde{A}_{z}^{\times})^{\tau^{\ast}=(-)^{-1}}
\over (A_{z}^{\times})^{\tau^{\ast}=(-)^{-1}}}\right|
\left(\sum_{m_{z}\in [{\cal M}(z)/P(z)](k)_{\sharp}}{\kappa
(\mathop{\rm cl}\nolimits_{{\cal M}(z)}(m_{z}))\over |\mathop{\rm
Aut}(m_{z})|}\right) \hfill\cr\hfill =(-|\kappa (z)|)^{r_{z}}
\left|{(\widetilde{A}_{z}^{\times})^{\tau^{\ast} =(-)^{-1}}\over
(A_{1,z}^{\times}\times A_{2,z}^{\times})^{\tau^{\ast}
=(-)^{-1}}}\right| \left(\sum_{n_{z}\in [{\cal N}(z)/Q(z)]
(k)_{\sharp}}{1\over |\mathop{\rm Aut}(n_{z})|}\right)\qquad}
$$
quel que soit $z\in Z^{\rm inerte}-\{z_0\}$.

La formule du produit, obtenue plus haut dans cette section, montre que
cette \'{e}galit\'{e} vaut aussi pour $z=z_{0}$.
\vskip 2mm

Enfin, dans le cas $z=z_{0}$, on a $J=I$, $J_{\alpha}=I_{\alpha}$ pour
$\alpha =1,2$, et on se retrouve dans la situation du chapitre 1.  Par
d\'{e}finition des int\'{e}grales orbitales (cf.  la section (1.5)),
on a
$$
\sum_{m_\lambda\in [{\cal M}_{\lambda}(k)/P(z)(k)]_{\sharp}}
{1\over \mathop{\rm Aut}(m_{\lambda})}={1\over |P(z)(k)|}
\mathop{\rm O}\nolimits_{\gamma_{I}}^{c_{\lambda}c_{{\cal
N}}^{0}},
$$
de sorte que
$$
\sum_{m\in [{\cal M}(z)/P(z)](k)_{\sharp}}{\kappa
(\mathop{\rm cl}\nolimits_{{\cal M}(z)}(m_{z}))\over 
|\mathop{\rm Aut}(m_{z})|}=\sum_{\lambda\in({\Bbb Z}/2{\Bbb Z})^{I}}
{\kappa (\lambda)\over |P(z)(k)|}
\mathop{\rm O}\nolimits_{\gamma}^{c_{\lambda}c_{{\cal N}}^{0}} 
={1\over |P(z)(k)|}\mathop{\rm O}\nolimits_{\gamma}^{\kappa}.
$$
Les m\^{e}mes arguments montrent que
$$
\sum_{n\in [{\cal N}(z)/Q(z)](k)_{\sharp}}{1\over |\mathop{\rm
Aut}(n)|}= {1\over |Q(z)(k)|}\mathop{\rm SO}\nolimits_{\gamma}^{H}.
$$
Comme
$$
{1\over |P(z)(k)|}
\left|{(\widetilde{A}_{z}^{\times})^{\tau^{\ast}=(-)^{-1}}\over
(A_{z}^{\times})^{\tau^{\ast}=(-)^{-1}}}\right|
= {1\over |Q(z)(k)|}
\left|{(\widetilde{A}_{z}^{\times})^{\tau^{\ast} =(-)^{-1}}\over
(A_{1,z}^{\times}\times A_{2,z}^{\times})^{\tau^{\ast}
=(-)^{-1}}}\right|
$$
on obtient l'\'{e}galit\'{e}
$$
\mathop{\rm O}\nolimits_{\gamma}^{\kappa}=(-|\kappa (z)|)^{r_{z}}
\mathop{\rm SO}\nolimits_{\gamma}^{H}.
$$
que l'on cherchait \`{a} d\'{e}montrer.
\hfill\hfill$\square$

\section{A}{Appendice}
\vskip - 3mm

\subsection{A.1}{Localisation \`{a} la Atiyah-Borel-Segal}

Soient $S$ un $k$-sch\'{e}ma de type fini et $f:X\rightarrow S$ un
$S$-sch\'{e}ma propre muni d'une action d'un $S$-sch\'{e}ma en groupes
commutatifs $P$ lisse et de type fini.  Soient de plus $T$ un
$S$-sch\'{e}ma en tores, $\Lambda$ un groupe commutatif de type fini
et $T\times\Lambda\rightarrow P$ un homomorphisme.  L'action de
$P$ induit des actions qui commutent de $T$ et $\Lambda$ sur le
$S$-sch\'{e}ma $X$.

Soient $g:Y=X^{T}\rightarrow S$ le $S$-sch\'{e}ma des points fixes pour
l'action de $T$ et $f^{T}:[X/T]\rightarrow S$ le $S$-champ quotient
pour l'action de $T$. Le groupe de type fini $\Lambda$ pr\'{e}serve
le ferm\'{e} $i:Y\hookrightarrow X$ et agit donc sur ce ferm\'{e} et son ouvert
compl\'{e}mentaire $j:U=X-Y\hookrightarrow X$; il agit aussi sur
$[X/T]$ en pr\'{e}servant le ferm\'{e} $[i]:[Z/T]\hookrightarrow [X/T]$
et l'ouvert compl\'{e}mentaire $[j]:[U/T]\hookrightarrow [X/T]$.

Consid\'{e}rons les faisceaux de cohomologie perverse
$$
{}^{\rm p}{\cal
H}^{n}(f_{\ast}{\Bbb Q}_{\ell}),~{}^{\rm p}{\cal
H}^{n}(g_{\ast}{\Bbb Q}_{\ell}),~{}^{\rm p}{\cal
H}^{n}(f_{\ast}^{T}{\Bbb Q}_{\ell}),~{}^{\rm p}{\cal
H}^{n}(f_{\ast}^{T}[i]_{\ast}{\Bbb Q}_{\ell})\hbox{ et }{}^{\rm p}{\cal
H}^{n}(f^{T}_{\ast}[j]_{!}{\Bbb Q}_{\ell})
$$
sur $S$.  L'action induite de $\Lambda$ sur ces faisceaux se factorise
\`{a} travers un quotient fini.  En effet, comme $P$ est de type fini
il existe un sous-groupe d'indice fini $\Lambda'\subset\Lambda$ tel
que, pour tout point g\'{e}om\'{e}trique $s$ dans $S$ et tout
$\lambda\in\Lambda'$, on ait $\lambda_{s}\in P_{s}^{0}$, et
d'apr\`{e}s le lemme d'homotopie 3.2.3, $\Lambda'$ agit trivialement
les faisceaux ci-dessus.

On a donc des d\'{e}compositions
$$
{}^{\rm p}{\cal H}^{n}(f_{\ast}{\Bbb Q}_{\ell})=\bigoplus_{\kappa}
{}^{\rm p}{\cal H}^{n}(f_{\ast}{\Bbb Q}_{\ell})_{\kappa},~{}^{\rm p}{\cal
H}^{n}(g_{\ast}{\Bbb Q}_{\ell})=\bigoplus_{\kappa}
{}^{\rm p}{\cal H}^{n}(f_{\ast}{\Bbb Q}_{\ell})_{\kappa},
$$
$$
{}^{\rm p}{\cal H}^{n}(f^{T}_{\ast}{\Bbb Q}_{\ell})=\bigoplus_{\kappa}
{}^{\rm p}{\cal H}^{n}(f^{T}_{\ast}{\Bbb Q}_{\ell})_{\kappa},
$$
$$
{}^{\rm p}{\cal H}^{n}(f^{T}_{\ast}[i]_{\ast}{\Bbb Q}_{\ell})=
\bigoplus_{\kappa}{}^{\rm p}{\cal H}^{n}(f^{T}_{\ast}[i]_{\ast}
{\Bbb Q}_{\ell})_{\kappa}\hbox{ et }
{}^{\rm p}{\cal H}^{n}(f^{T}_{\ast}[j]_{!}{\Bbb Q}_{\ell})=
\bigoplus_{\kappa}{}^{\rm p}{\cal H}^{n}(f^{T}_{\ast}[j]_{!}
{\Bbb Q}_{\ell})_{\kappa}.
$$
selon les caract\`{e}res d'ordre fini $\kappa$ de $\Lambda$ (tout du
moins si tous ces caract\`{e}res sont rationnels sur ${\Bbb Q}_{\ell}$;
sinon, il faut \'{e}tendre au pr\'{e}alable les scalaires \`{a} une
extension finie de ${\Bbb Q}_{\ell}$).

Consid\'{e}rons le diagramme commutatif \`{a} carr\'{e} cart\'{e}sien
$$
\xymatrix{
X\ar[d]_{f}\ar[r]^-{\,\pi_{X}} & [X/T] \ar[d]_{[f]}\ar[dr]^{f^{T}} &\cr
S\ar[r]_-{\pi} & [S/T]\ar[r]_-{\varepsilon} & S }
$$
o\`{u} $[f]$ est le morphisme qui se d\'{e}duit de $f$ par le passage
au quotient par l'action de $T$, o\`{u} $\pi :S\rightarrow [S/T]$
est le $T$-torseur universel et o\`{u} $\varepsilon :[S/T]\rightarrow
S$ est le morphisme structural.

D'apr\`{e}s le th\'{e}or\`{e}me de changement de base pour un morphisme
propre, on a $\pi^{\ast}[f]_{\ast}{\Bbb Q}_{\ell}=f_{\ast}{\Bbb
Q}_{\ell}$ et donc $\pi^{\ast}[f]_{\ast}{\Bbb
Q}_{\ell}=\pi^{\ast}\varepsilon^{\ast}f_{\ast}{\Bbb Q}_{\ell}$ puisque
$\varepsilon\circ\pi$ est l'identit\'{e} de $S$.  Comme $\pi$ est
lisse \`{a} fibres g\'{e}om\'{e}triquement connexes, on en d\'{e}duit
que
$$
{}^{\rm p}{\cal H}^{n}([f]_{\ast}{\Bbb Q}_{\ell})={}^{\rm p}{\cal 
H}^{n}(\varepsilon^{\ast}f_{\ast}{\Bbb Q}_{\ell})=
\varepsilon^{\ast}{}^{\rm p}{\cal H}^{n}(f_{\ast}{\Bbb Q}_{\ell})
$$
pour tout entier $n$ (cf. la Proposition 4.2.5 de [B-B-D]). 

Soit $X^{\ast}(T)$ le faisceau \'{e}tale sur $S$ dont les sections sur
un ouvert \'{e}tale $V$ de $S$ sont les homomorphismes
$V\times_{S}T\rightarrow {\Bbb G}_{{\rm m},V}$.  D'apr\`{e}s la
th\'{e}orie de Chern-Weil on a ${\cal H}^{n}(\varepsilon_{\ast}{\Bbb
Q}_{\ell})=(0)$ pour tout $n<0$ et pour tout $n$ impair, et pour tout
entier $n\geq 0$, on a un isomorphisme
$$
{\Bbb Q}_{\ell}[X^{\ast}(T)(-1)]^{2n}:=\mathop{\rm
Sym}\nolimits_{{\Bbb Q}_{\ell}}^{n}\left(X^{\ast}(T)\otimes_{{\Bbb
Z}_{S}}{\Bbb Q}_{\ell, S}(-1)\right)\buildrel\sim\over\longrightarrow
\bigoplus_{i\geq 0} {\cal H}^{2n}(\varepsilon_{\ast}{\Bbb Q}_{\ell})
$$
induit par la fl\`{e}che
$$
X^{\ast}(T)\rightarrow {\cal H}^{2}(\varepsilon_{\ast}{\Bbb Q}_{\ell})(1)
$$
qui envoie un caract\`{e}re $\chi$ sur la classe de Chern du fibr\'{e}
inversible sur $[S/T]$ obtenu en poussant le $T$-torseur tautologique
par $\chi :T\rightarrow {\Bbb G}_{{\rm m},S}$.

La structure multiplicative et le morphisme $[f]^{\ast}$ munissent la
somme formelle de faisceaux pervers
$$
\bigoplus_{n}{}^{\rm p}{\cal H}^{n}(f^{T}_{\ast}{\Bbb Q}_{\ell}),
$$
ainsi que la somme formelle de faisceaux pervers
$$
\bigoplus_{n}{}^{\rm p}{\cal H}^{n}(f^{T}_{\ast}{\Bbb
Q}_{\ell})_{\kappa},
$$
pour tout $\kappa$, de structures de ${\Bbb Q}_{\ell}
[X^{\ast}(T)(-1)]$-modules gradu\'{e}s o\`{u}
$$
{\Bbb Q}_{\ell}[X^{\ast}(T)(-1)]:=\bigoplus_{n}{\Bbb Q}_{\ell}
[X^{\ast}(T)(-1)]^{2n}
$$
est vu comme un faisceau en ${\Bbb Q}_{\ell}$-alg\`{e}bres
gradu\'{e}es.  Ces sommes infinies de faisceaux pervers n'ont a priori
pas de sens. En fait il s'agit juste d'une fa\c{c}on commode
d'\'{e}crire les choses: on pourrait tr\`{e}s bien travailler degr\'{e}
par degr\'{e}.

Supposons de plus que, pour un caract\`{e}re $\kappa_{0}$ particulier
de $\overline{\Lambda}$, les faisceaux pervers ${}^{\rm p}{\cal
H}^{n}(f_{\ast}{\Bbb Q}_{\ell})_{\kappa_{0}}$ sont purs de poids $n$
pour tous $n\in {\Bbb Z}$.

\thm PROPOSITION A.1.1
\enonce
Sous cette hypoth\`{e}se, il existe un isomorphisme
$$
\bigoplus_{n}{}^{\rm p}{\cal H}^{n}(f^{T}_{\ast}{\Bbb
Q}_{\ell})_{\kappa_{0}}\cong \left(\bigoplus_{n}{}^{\rm p}{\cal H}^{n}
(f_{\ast}{\Bbb Q}_{\ell})_{\kappa_{0}}\right)\otimes_{{\Bbb Q}_{\ell ,S}}
{\Bbb Q}_{\ell}[X^{\ast}(T)(-1)].
$$
En particulier, cette somme directe est un module libre de type fini
sur ${\Bbb Q}_{\ell}[X^{\ast}(T)(-1)]$.
\endthm

\rem D\'{e}monstration
\endrem
On a la suite spectrale de Leray
$$
E_{2}^{pq}=\,{}^{\rm p}{\cal H}^{p}(\varepsilon_{\ast}{}^{\rm p}{\cal
H}^{q}([f]_{\ast}{\Bbb Q}_{\ell}))\Rightarrow {}^{\rm p}{\cal H}^{p+q}
(f^{T}_{\ast}{\Bbb Q}_{\ell}).
$$
compatible \`{a} l'action de $\Lambda$.  Or on a vu que ${}^{\rm p}{\cal
H}^{q}([f]_{\ast}{\Bbb Q}_{\ell})=\varepsilon^{\ast}{}^{\rm p}{\cal 
H}^{q}(f_{\ast}{\Bbb Q}_{\ell})$ pour tout entier $q$, de sorte que 
la formule de projection permet de r\'{e}crire le terme $E_{2}^{pq}$ de 
cette suite spectrale sous la forme
$$
E_{2}^{pq}=\,{}^{\rm p}{\cal H}^{p}(\varepsilon_{\ast}{\Bbb Q}_{\ell}
\otimes_{{\Bbb Q}_{\ell ,S}}{}^{\rm p}{\cal H}^{q}(f_{\ast}{\Bbb
Q}_{\ell})).
$$
Comme $\varepsilon_{\ast}{\Bbb Q}_{\ell}$ est \`{a} cohomologie 
ordinaire lisse sur $S$, on a encore
$$
E_{2}^{pq}=\,{\cal H}^{p}(\varepsilon_{\ast}{\Bbb Q}_{\ell})
\otimes_{{\Bbb Q}_{\ell ,S}}{}^{\rm p}{\cal H}^{q}(f_{\ast}{\Bbb
Q}_{\ell}).
$$
On en d\'{e}duit la suite spectrale
$$
E_{2}^{pq}=\,{\cal H}^{p}(\varepsilon_{\ast}{\Bbb
Q}_{\ell})\otimes_{{\Bbb Q}_{\ell ,S}} {}^{\rm p}{\cal
H}^{q}(f_{\ast}{\Bbb Q}_{\ell})_{\kappa_{0}}\Rightarrow {}^{\rm
p}{\cal H}^{p+q} (f^{T}_{\ast}{\Bbb Q}_{\ell})_{\kappa_{0}}.
$$
Le ${\Bbb Q}_{\ell}$-faisceau lisse ${\cal H}^{p}(\varepsilon_{\ast}{\Bbb
Q}_{\ell})$ \'{e}tant pur de poids $p$, l'hypoth\`{e}se de puret\'{e}
du faisceau pervers ${}^{{\rm p}}{\cal H}^{q}(f_{\ast}T)_{\kappa_{0}}$
implique que cette derni\`{e}re suite spectrale d\'{e}g\'{e}n\`{e}re
en $E_{2}$ puisque toutes ses fl\`{e}ches $d_{r}^{pq}$ relient des
termes de poids diff\'{e}rents et sont donc n\'{e}cessairement nulles.

\hfill\hfill$\square$
\vskip 3mm

Le triangle distingu\'{e} habituel
$$
[j]_{!}{\Bbb Q}_{\ell,[U/T]}\rightarrow {\Bbb Q}_{\ell,[X/T]}
\rightarrow [i]_{\ast}{\Bbb Q}_{\ell,[Y/T]}\rightarrow
$$
induit un triangle distingu\'{e}
$$
f^{T}_{\ast}[j]_{!}{\Bbb Q}_{\ell}\rightarrow f^{T}_{\ast}{\Bbb
Q}_{\ell}\rightarrow f^{T}_{\ast}[i]_{\ast}{\Bbb Q}_{\ell}\rightarrow
$$
qui induit \`{a} son tour une suite exacte longue de faisceaux
cohomologie perverse sur $S$
$$\eqalign{
\cdots &\rightarrow {}^{\rm p}{\cal H}^{n}(f^{T}_{\ast}[j]_{!}
{\Bbb Q}_{\ell}) \rightarrow {}^{\rm p}{\cal H}^{n} (f^{T}_{\ast}{\Bbb
Q}_{\ell})\rightarrow {}^{\rm p}{\cal H}^{n}
(f^{T}_{\ast}[i]_{\ast}{\Bbb Q}_{\ell})\cr
&\rightarrow {}^{\rm p}{\cal H}^{n+1}(f^{T}_{\ast}[j]_{!}{\Bbb
Q}_{\ell})\rightarrow \cdots .\cr}
$$
Cette suite exacte est \'{e}quivariante par rapport \`{a} l'action du
groupe de type fini $\Lambda$ et induit donc une suite exacte
$$\eqalign{
\cdots & \rightarrow \,{}^{\rm p}{\cal H}^{n}(f^{T}_{\ast}[j]_{!}
{\Bbb Q}_{\ell})_{\kappa_{0}}\rightarrow {}^{\rm p}{\cal H}^{n}
(f^{T}_{\ast}{\Bbb Q}_{\ell})_{\kappa_{0}}\rightarrow
{}^{\rm p}{\cal H}^{n}(f^{T}_{\ast}[i]_{\ast}{\Bbb
Q}_{\ell})_{\kappa_{0}}\cr
&\rightarrow {}^{\rm p}{\cal H}^{n+1}(f^{T}_{\ast}[j]_{!}
{\Bbb Q}_{\ell})_{\kappa_{0}}\rightarrow \cdots .\cr}
$$
On en d\'{e}duit une fl\`{e}che
$$
\bigoplus_{n}{}^{\rm p}{\cal H}^{n}(f^{T}_{\ast}[j]_{!}
{\Bbb Q}_{\ell})_{\kappa_{0}}\rightarrow\bigoplus_{n}
{}^{\rm p}{\cal H}^{n}(f^{T}_{\ast}{\Bbb Q}_{\ell})_{\kappa_{0}}
$$
et une fl\`{e}che de restriction
$$
\bigoplus_{n}{}^{\rm p}{\cal H}^{n}(f^{T}_{\ast}{\Bbb Q}_{\ell})_{\kappa_{0}}
\rightarrow \bigoplus_{n}{}^{\rm p}{\cal H}^{n}
(f^{T}_{\ast}[i]_{\ast}{\Bbb Q}_{\ell})_{\kappa_{0}}
$$
qui sont ${\Bbb Q}_{\ell}[X^{\ast}(T)(-1)]$-lin\'{e}aires
gradu\'{e}es.  Ici, on a
$$
\bigoplus_{n}{}^{\rm p}{\cal H}^{n}(f^{T}_{\ast}[i]_{\ast}{\Bbb
Q}_{\ell})_{\kappa_{0}}= \left(\bigoplus_{n}{}^{\rm p}{\cal H}^{n}
(g_{\ast}{\Bbb Q}_{\ell})_{\kappa_{0}}\right)
\otimes_{{\Bbb Q}_{\ell ,S}}{\Bbb Q}_{\ell}[X^{\ast}(T)(-1)]
$$
puisque $T$ agit trivialement sur $Y$.

\thm PROPOSITION A.1.2
\enonce
Le ${\Bbb Q}_{\ell}[X^{\ast}(T)(-1)]$-module gradu\'{e} $\bigoplus_{n}
{}^{\rm p}{\cal H}^{n}(f^{T}_{\ast}[j]_{!}{\Bbb Q}_{\ell})$ est
de torsion.
\endthm

\rem D\'{e}monstration
\endrem
Notons $h=f\circ j:U\rightarrow S$ et $[h]:[U/T]\rightarrow {\rm
B}(T/S)$ le morphisme quotient de $h$ par l'action de $T$.  On a
$$
f^{T}_{\ast}[j]_{!}{\Bbb Q}_{\ell}=\varepsilon_{\ast}[h]_{!}{\Bbb
Q}_{\ell}
$$
o\`{u} $\varepsilon :{\rm B}(T/S)\rightarrow S$ est le morphisme
structural.

On va en fait d\'{e}montrer plus g\'{e}n\'{e}ralement que pour tout
$S$-sch\'{e}ma $h:U\rightarrow S$ de type fini, muni d'une action de
$T$ qui rel\`{e}ve l'action triviale de $T$ sur $S$ et qui est sans
point fixe, le faisceau pervers
$$
\bigoplus_{n} {}^{\rm p}{\cal H}^{n}(\varepsilon_{\ast} [h]_{!}{\Bbb
Q}_{\ell})
$$
sur $S$ en modules gradu\'{e}s sur le faisceau de ${\Bbb
Q}_{\ell}$-alg\`{e}bres ${\Bbb Q}_{\ell}[X^{\ast}(T)(-1)]$ est de
torsion.
\vskip 2mm

L'\'{e}nonc\'{e} est local pour la topologie \'{e}tale sur $S$.  On
peut donc supposer que $T={\Bbb G}_{{\rm m},S}^{r}$ pour un entier
$r\geq 1$, ou ce qui revient au m\^{e}me que l'on a une action de
$T={\Bbb G}_{{\rm m},k}^{r}$ sur $X$ au-dessus de l'action triviale de
$T$ sur $S$.
\vskip 2mm

Il existe une filtration finie de $U$ par des parties ferm\'{e}es
$T$-invariantes
$$
U=U_{0}\supset U_{1}\supset\cdots\supset U_{p}\supset
U_{p+1}\supset\cdots
$$
et pour chaque $p$, un sous-tore $T_{p}\subset T$ tel que, quel que
soit le point g\'{e}om\'{e}trique $x$ de $U^{p}:=U_{p}-U_{p+1}$, la
composante neutre du fixateur $T_{x}\subset \kappa(x)\otimes_{k}T$ de
$x$ n'est autre que $\kappa (x)\otimes_{k}T_{p}$.

Alors, si on note $h^{p}:U^{p}\rightarrow S$ la restriction de $h$
\`{a} la partie localement ferm\'{e}e $U^{p}=U_{p}-U_{p+1}$ de $U$, on
a la suite spectrale
$$
E_{1}^{pq}={}^{\rm p}{\cal H}^{p+q}(\varepsilon_{\ast}
[h^{p}]_{!}{\Bbb Q}_{\ell})\Rightarrow {}^{\rm p}{\cal
H}^{p+q}(\varepsilon_{\ast}[h]_{!}{\Bbb Q}_{\ell}).
$$
Par suite, il suffit de v\'{e}rifier que chaque faisceau pervers en
${\Bbb Q}_{\ell}[X^{\ast}(T)(-1)]$-modules gradu\'{e}s $\bigoplus_{n}
{}^{\rm p}{\cal H}^{n}(\varepsilon_{\ast} [h^{p}]_{!}{\Bbb Q}_{\ell})$
est de torsion.

On peut donc supposer dans notre probl\`{e}me initial qu'il existe une
sous-tore $T'\subset T$ tel que, quel que soit le point
g\'{e}om\'{e}trique $x$ de $U$, la composante neutre du fixateur
$T_{x}\subset \kappa(x)\otimes_{k}T$ de $x$ est \'{e}gale \`{a}
$\kappa (x)\otimes_{k}T'$.
\vskip 2mm

Fixons arbitrairement un sous-tore $T''\subset T$ tel que la
fl\`{e}che produit $T'\times_{k}T''\rightarrow T$ est surjective et
\`{a} noyau fini.  On remarque que $T''$ agit sans point fixe sur $U$
puisque $T$ agit sans point fixe sur $U$.

On a un diagramme commutatif \`{a} carr\'{e} cart\'{e}sien
$$
\xymatrix{[U/(T'\times_{k}T'')]\ar[d]_{[h]'}^{}\ar[rr]_{}^{} & &
[U/T]\ar[d]_{}^{[h]}\cr
[S/(T'\times_{k}T'')]\ar[rd]_{\varepsilon'}^{}\ar[rr]_{}^{} & & [S/T]
\ar[ld]_{}^{\varepsilon}\cr
 & S & }
$$
o\`{u} les fl\`{e}ches horizontales sont induites par
l'\'{e}pimorphisme $T'\times_{k}T''\rightarrow T$ et sont des gerbes
finies, o\`{u} les fl\`{e}ches verticales sont induites par $h$ et
o\`{u} les fl\`{e}ches obliques sont les morphismes structuraux.

On a donc un isomorphisme de restriction
$$
{}^{\rm p}{\cal H}^{p+q}(\varepsilon_{\ast} [h]_{!}{\Bbb
Q}_{\ell})\buildrel\sim\over\longrightarrow {}^{\rm p}{\cal
H}^{p+q}(\varepsilon_{\ast}' [h]_{!}'{\Bbb Q}_{\ell})
$$
de faisceaux pervers en ${\Bbb Q}_{\ell}[X^{\ast}(T)(-1)]$-modules
gradu\'{e}s, o\`{u} bien entendu ${\Bbb Q}_{\ell}[X^{\ast}(T)(-1)]$
agit sur le but via l'inclusion $X^{\ast}(T)\hookrightarrow
X^{\ast}(T'\times_{k}T'')$.  Comme le conoyau de cette inclusion est de
torsion, pour d\'{e}montrer notre assertion il suffit donc de
v\'{e}rifier que chaque faisceau pervers en ${\Bbb
Q}_{\ell}[X^{\ast}(T'\times_{k}T'')(-1)]$-modules gradu\'{e}s
$\bigoplus_{n} {}^{\rm p}{\cal H}^{n}(\varepsilon_{\ast}'
[h]_{!}'{\Bbb Q}_{\ell})$ est de torsion.

On peut donc supposer de plus dans notre probl\`{e}me initial que
$T=T'\times_{k}T''$.
\vskip 2mm

Comme $[h]:[U/T]\rightarrow [S/T]$ est alors le produit du
$T'$-torseur universel $S\rightarrow [S/T']$ par le morphisme $[U/T'']
\rightarrow [S/T'']$ induit par $h$, on voit par application de la
formule de K\"{u}nneth que l'on peut supposer en outre dans notre
probl\`{e}me initial que $T'=(1)$, c'est-\`{a}-dire que $T$ agit
librement sur $U$.
\vskip 2mm

On a donc un $S$-sch\'{e}ma de type fini $h:U\rightarrow S$ sur lequel
le tore $T={\Bbb G}_{{\rm m},k}^{r}$ agit librement.  Le $S$-champ
alg\'{e}brique quotient $[U/T]$ est alors un $S$-espace
alg\'{e}brique.  Il admet donc un recouvrement ouvert pour la
topologie \'{e}tale (par des sch\'{e}mas affines si l'on veut) qui
trivialise le $T$-torseur $U\rightarrow [U/T]$.  Par image inverse sur
$U$, on obtient alors un recouvrement ouvert $T$-\'{e}quivariant
$(U_{\alpha})_{\alpha\in A}$ de $U$ pour la topologie \'{e}tale et,
pour chaque $\alpha\in A$, un morphisme de $k$-sch\'{e}mas
$\varphi_{\alpha}:U_{\alpha}\rightarrow T$ tel que
$$
\varphi_{\alpha}(t\cdot x)=t\varphi_{\alpha}(x).
$$
Pour toute partie $B$ de $A$ on note $U_{B}=\prod_{\alpha\in
B}(U_{\alpha}/U)$ {\og}{l'intersection}{\fg} des ouverts $U_{\alpha}$
pour $\alpha\in B$ et $h_{B}:U_{B}\rightarrow U\rightarrow S$ le
morphisme canonique.

La suite spectrale
$$
E_{1}^{pq}=\bigoplus_{|B|=1-p>0}{}^{\rm p}{\cal H}^{p+q}(\varepsilon_{\ast}
[h_{B}]_{!}{\Bbb Q}_{\ell})\otimes_{{\Bbb Z}}\bigwedge^{|B|}{\Bbb Z}^{B}
\Rightarrow {}^{\rm p}{\cal
H}^{p+q}(\varepsilon_{\ast}[h]_{!}{\Bbb Q}_{\ell}).
$$
pour le recouvrement $([U_{\alpha}/T])_{\alpha\in A}$ de $[U/T]$
montre qu'il suffit de v\'{e}rifier que les ${\Bbb
Q}_{\ell}[X^{\ast}(T)(-1)]$-modules gradu\'{e}s $\bigoplus_{n}{}^{\rm
p}{\cal H}^{n}(\varepsilon_{\ast}[h]_{B,!}{\Bbb Q}_{\ell})$ pour $B$
parcourant l'ensemble des parties de $A$, sont tous de torsion.

On peut donc supposer dans notre probl\`{e}me initial qu'en plus de
toutes les hypoth\`{e}ses d\'{e}j\`{a} formul\'{e}es, il existe un
morphisme de $k$-sch\'{e}mas $\varphi :U\rightarrow T$ tel que
$$
\varphi (t\cdot x)=t\varphi (x).
$$
Mais alors le morphisme $[h]$ se factorise en
$$
[h]:[U/T]\,\smash{\mathop{\hbox to 8mm{\rightarrowfill}}
\limits^{\scriptstyle [\varphi]}}\, [S\times_{k}T/T]=S
\,\smash{\mathop{\hbox to 8mm{\rightarrowfill}}
\limits^{\scriptstyle [\mathop{\rm pr}\nolimits_{T}]}}\, [S/T]
$$
o\`{u} $T$ agit par translation sur lui-m\^{e}me, o\`{u} $\mathop{\rm
pr}\nolimits_{T}:S\times_{k}T\rightarrow S$ est la projection
canonique et o\`{u} $[\mathop{\rm pr}\nolimits_{T}]$ est donc le
$T$-torseur universel.  Par suite le morphisme de restriction
$[h]^{\ast}:\varepsilon_{\ast}{\Bbb Q}_{\ell}\rightarrow
\varepsilon_{\ast}[h]_{\ast}{\Bbb Q}_{\ell}$ se factorise en
$$
[h]^{\ast}:\varepsilon_{\ast}{\Bbb Q}_{\ell}\,\smash{\mathop{\hbox to
8mm{\rightarrowfill}} \limits^{\scriptstyle [\mathop{\rm pr}
\nolimits_{T}]^{\ast}}}\, {\Bbb Q}_{\ell}\, \smash{\mathop{\hbox to
8mm{\rightarrowfill}} \limits^{\scriptstyle [\varphi ]^{\ast}}}
\,\varepsilon_{\ast}[h]_{\ast}{\Bbb Q}_{\ell}
$$
puisque $\varepsilon\circ [\mathop{\rm pr} \nolimits_{T}]$ est
l'identit\'{e} de $S$.

Dans le cas particulier o\`{u} l'est s'est ramen\'{e}, on a finalement
d\'{e}montr\'{e} que l'id\'{e}al d'augmentation de ${\Bbb
Q}_{\ell}[X^{\ast}(T)(-1)]$ tout entier annule $\bigoplus_{n}{}^{\rm
p}{\cal H}^{n}([h]_{!}{\Bbb Q}_{\ell})$.
\hfill\hfill$\square$
\vskip 3mm

On d\'{e}duit de la proposition l'\'{e}nonc\'{e} suivant qui est une
variante du th\'{e}or\`{e}me de localisation d'Atiyah-Borel-Segal.

\thm COROLLAIRE A.1.3
\enonce
Sous l'hypoth\`{e}se de puret\'{e}, la fl\`{e}che
$$
\bigoplus_{n}{}^{\rm p}{\cal H}^{n}(f^{T}_{\ast}[j]_{!}
{\Bbb Q}_{\ell})_{\kappa_{0}}\rightarrow\bigoplus_{n}{}^{\rm p}
{\cal H}^{n}(f^{T}_{\ast}{\Bbb Q}_{\ell})_{\kappa_{0}}
$$
ci-dessus est nulle et la fl\`{e}che de restriction
$$
\bigoplus_{n}{}^{\rm p}{\cal H}^{n} (f^{T}_{\ast}{\Bbb Q}_{\ell})_{\kappa_{0}}
\rightarrow \left(\bigoplus_{p}{}^{\rm p}{\cal H}^{p}(g_{\ast}
{\Bbb Q}_{\ell})_{\kappa_{0}}\right)\otimes_{{\Bbb Q}_{\ell ,S}}
{\Bbb Q}_{\ell}[X^{\ast}(T)(-1)]
$$
est injective.
\endthm

\rem D\'{e}monstration
\endrem
Une fl\`{e}che d'un module de torsion dans un module libre sur ${\Bbb
Q}_{\ell}[X^{\ast}(T)(-1)]$ est n\'{e}cessairement nulle.  La
nullit\'{e} de la premi\`{e}re fl\`{e}che implique l'injectivit\'{e} de
la seconde d'apr\`{e}s la suite exacte longue.
\hfill\hfill$\square$
\vskip 3mm

\subsection{A.2}{Cohomologie \'{e}quivariante de la droite projective
pinc\'{e}e}

Soient $Y$ un $k$-sch\'{e}ma et ${\cal L}$ un ${\cal O}_{Y}$-Module
inversible.  On note $p:P={\Bbb P}({\cal L}\oplus {\cal
O}_{Y})\rightarrow Y$ le fibr\'{e} en droites projectives quotients de
${\cal L}\oplus {\cal O}_{Y}$ et on note $\sigma_{0}:Y\rightarrow
{\Bbb V}({\cal L})\subset P$ et $\sigma_{\infty}:Y
\buildrel\sim\over\longrightarrow P-{\Bbb V}({\cal L})\subset P$ ses
sections nulle et infinie.

On fait agir ${\Bbb G}_{{\rm m},Y}$ par homoth\'{e}ties sur le
fibr\'{e} en droites ${\Bbb V}({\cal L})\rightarrow Y$ et on prolonge
cette action au fibr\'{e} projectif $P\rightarrow Y$.  On note
$[p]:[P/{\Bbb G}_{{\rm m},Y}]\rightarrow [Y/{\Bbb G}_{{\rm m},Y}]={\rm
B}({\Bbb G}_{{\rm m},Y}/Y)$ le morphisme de champs alg\'{e}briques
quotient de $p$ par cette action de ${\Bbb G}_{{\rm m},Y}$ sur $P$
au-dessus de l'action triviale sur $Y$.  Alors $[p]$ est un fibr\'{e}
en droites projectives et s'\'{e}crit sous la forme ${\Bbb
P}(\widetilde{{\cal L}}\oplus {\cal O}_{{\rm B}({\Bbb G}_{{\rm
m},Y}/Y)})$ pour ${\cal O}_{{\rm B}({\Bbb G}_{{\rm m},Y}/Y)}$-Module
inversible $\widetilde{{\cal L}}$. On a
le diagramme :
$$\xymatrix{
P \ar[r] \ar[d]_{p}& [P/{\Bbb G}_m] \ar[d]^{[p]}\cr
Y \ar[r] & B({\Bbb G}_{m,Y}/Y) }
$$

On peut aussi faire agir trivialement ${\Bbb G}_{{\rm m},Y}$ sur $P$
et $Y$, et $p$ passe au quotient pour cette action en un morphisme de
champs alg\'{e}briques ${\rm B}(p):{\rm B}({\Bbb G}_{{\rm m},P}/P)
\rightarrow {\rm B}({\Bbb G}_{{\rm m},Y}/Y)$ qui n'est autre que le
fibr\'{e} en droites projectives ${\Bbb P}(\varepsilon_{Y}^{\ast}{\cal
L}\oplus {\cal O}_{{\rm B}({\Bbb G}_{{\rm m},Y}/Y)})$ o\`{u}
$\varepsilon_{Y}:{\rm B}({\Bbb G}_{{\rm m},Y}/Y)\rightarrow Y$ est le
morphisme structural.  On a le diagramme :
$$
\xymatrix{
P \ar[d]_{p} \ar@<0.6ex>[r] & B({\Bbb G}_{{\rm m},P}/P) \ar[d]^{B(p)}
\ar@<0.6ex>[l]\cr
Y \ar@<0.6ex>[r] & B({\Bbb G}_{{\rm m},Y}/Y)\ar@<0.6ex>[l]^-{\epsilon_{Y}}}
$$

\thm LEMME A.2.1
\enonce
Soit ${\cal T}$ le ${\cal O}_{{\rm B}({\Bbb G}_{{\rm m},Y}/Y)}$-Module
inversible universel sur le champ classifiant ${\rm B}({\Bbb G}_{{\rm
m},Y}/Y)$.  Alors on a un isomorphisme canonique
$$
\widetilde{{\cal L}}={\cal T}\otimes \varepsilon_{Y}^{\ast}{\cal L}
$$
de ${\cal O}_{{\rm B}({\Bbb G}_{{\rm m},Y}/Y)}$-Modules inversibles.
\hfill\hfill$\square$
\endthm

Compte tenu de ce lemme, on a comme pour tout fibr\'{e} en droites
projectives:

\thm PROPOSITION A.2.2
\enonce
La somme
$$
[p]_{\ast}{\Bbb Q}_{\ell ,[P/{\Bbb G}_{{\rm m},Y}]}\rightarrow {\Bbb
Q}_{\ell ,[Y/{\Bbb G}_{{\rm m},Y}]}\oplus {\Bbb Q}_{\ell ,[Y/{\Bbb
G}_{{\rm m},Y}]}
$$
des fl\`{e}ches de restriction induites par $\sigma_{0}$ et
$\sigma_{\infty}$ est canoniquement isomorphe \`{a} la fl\`{e}che
$$\eqalign{
{\Bbb Q}_{\ell ,[Y/{\Bbb G}_{{\rm m},Y}]}\oplus {\Bbb Q}_{\ell
,[Y/{\Bbb G}_{{\rm m},Y}]}[-2](-1)&\rightarrow {\Bbb Q}_{\ell ,[Y/{\Bbb
G}_{{\rm m},Y}]}\oplus {\Bbb Q}_{\ell ,[Y/{\Bbb G}_{{\rm m},Y}]}\cr
a\oplus b&\mapsto (a+\widetilde{c}(b))\oplus (a-\widetilde{c}(b))\cr}
$$
o\`{u} la classe de Chern $\widetilde{c}=c_{1}(\widetilde{{\cal
L}}):{\Bbb Q}_{\ell ,[Y/{\Bbb G}_{{\rm m},Y}]}[-2](-1)\rightarrow
{\Bbb Q}_{\ell ,[Y/{\Bbb G}_{{\rm m},Y}]}$ de $\widetilde{{\cal L}}$
est la somme
$$
\widetilde{c}=t+\varepsilon_{Y}^{\ast}(c)
$$
de la classe de Chern $t=c_{1}({\cal T})$ de ${\cal T}$ et de l'image
r\'{e}ciproque par $\varepsilon_{Y}^{\ast}$ de la classe de Chern
$c=c_{1}({\cal L})$ de ${\cal L}$.
\hfill\hfill$\square$
\endthm

\thm COROLLAIRE A.2.3
\enonce
La fl\`{e}che de restriction
$$
\varepsilon_{Y,\ast}[p]_{\ast}{\Bbb Q}_{\ell ,[P/{\Bbb G}_{{\rm
m},Y}]}\rightarrow \varepsilon_{Y,\ast}{\Bbb Q}_{\ell ,{\rm B}({\Bbb
G}_{{\rm m},Y})}\oplus \varepsilon_{Y,\ast}{\Bbb Q}_{\ell ,{\rm
B}({\Bbb G}_{{\rm m},Y})}
$$
induite par $\sigma_{0}$ et $\sigma_{\infty}$ est
canoniquement isomorphe \`{a} la fl\`{e}che
$$\displaylines{
\qquad\bigoplus_{n\geq 0}{\Bbb Q}_{\ell ,Y}[-2n](-n)\oplus
\bigoplus_{n\geq 0}{\Bbb Q}_{\ell ,Y}[-2n-2](-n-1)
\hfill\cr\hfill
\rightarrow \bigoplus_{n\geq 0}{\Bbb Q}_{\ell ,Y}[-2n](-n)\oplus
\bigoplus_{n\geq 0}{\Bbb Q}_{\ell ,Y}[-2n](-n)\qquad}
$$
d\'{e}finie par
$$
\oplus_{n}(a_{n}\oplus b_{n})\mapsto
\oplus_{n}\left((a_{n}+b_{n}+c_{n+1}(b_{n+1}))\oplus
(a_{n}-b_{n}-c_{n+1}(b_{n+1}))\right)
$$
o\`{u} $c_{n+1}:{\Bbb Q}_{\ell ,Y}[-2n-2](-n-1)\rightarrow {\Bbb Q}_{\ell
,Y}[-2n](-n)$ est induit par la premi\`{e}re classe de Chern
$c$ de ${\cal L}$.
\hfill\hfill$\square$
\endthm

Notons $q:P-(\sigma_{0}(Y)\cup \sigma_{\infty}(Y))\rightarrow Y$ le
${\Bbb G}_{{\rm m},Y}$-torseur compl\'{e}mentaire des sections nulle
et infini et $[q]:Y=[P-(\sigma_{0}(Y)\cup \sigma_{\infty}(Y))/{\Bbb
G}_{{\rm m},Y}]\rightarrow [Y/{\Bbb G}_{{\rm m},Y}]$ le quotient de
$q$ par l'action de ${\Bbb G}_{{\rm m},Y}$ de la structure de torseur.
La fl\`{e}che de restriction ci-dessus s'ins\`{e}re dans le triangle
distingu\'{e}
$$
\underbrace{\varepsilon_{Y,\ast}[q]_{!}{\Bbb Q}_{\ell ,Y}}_{A}
\rightarrow \underbrace{\varepsilon_{Y,\ast}[p]_{\ast}{\Bbb Q}_{\ell,
[P/{\Bbb G}_{{\rm m},Y}]}}_{B}\rightarrow
\underbrace{\varepsilon_{Y,\ast}{\Bbb Q}_{\ell ,{\rm B}({\Bbb G}_{{\rm
m},Y})}\oplus \varepsilon_{Y,\ast}{\Bbb Q}_{\ell ,{\rm B}({\Bbb G}_{{\rm
m},Y})}}_{C}\rightarrow A[1]
$$
dans la cat\'{e}gorie d\'{e}riv\'{e}e des faisceaux $\ell$-adiques
sur $Y$.
\vskip 2mm

Soient maintenant $S$ un $k$-sch\'{e}ma, $g:Y\rightarrow S$ un
morphisme propre de $k$-sch\'{e}mas et $p:P={\Bbb P}({\cal
L}\oplus{\cal O}_Y)\rightarrow Y$ un fibr\'{e} en droites
projectives comme pr\'{e}c\'{e}demment.

On note $\varepsilon_{S}:{\rm B}({\Bbb G}_{{\rm m},S})\rightarrow S$
le morphisme structural et $t$ la premi\`{e}re classe de Chern du
${\cal O}_{{\rm B}({\Bbb G}_{{\rm m},S})}$-Module inversible
universel, vue comme une section de ${\cal H}^{2}(\varepsilon_{S,\ast}{\Bbb
Q}_{\ell ,{\rm B}({\Bbb G}_{{\rm m},S})})(-1)$.  On a l'anneau gradu\'{e}
$$
\bigoplus_{n\geq 0}{\cal H}^{2n}(\varepsilon_{S,\ast}{\Bbb Q}_{\ell ,{\rm
B}({\Bbb G}_{{\rm m},S})})(-n)={\Bbb Q}_{\ell ,S}[t]
$$
et tous les autres ${\cal H}^{n}(\varepsilon_{S,\ast}{\Bbb Q}_{S})$ sont
nuls.

La premi\`{e}re classe de Chern $c:{\Bbb
Q}_{\ell,Y}\rightarrow {\Bbb Q}_{\ell ,Y}[2](1)$ de ${\cal L}$ induit
une fl\`{e}che
$$
g_{\ast}(c):g_{\ast}{\Bbb Q}_{\ell ,Y}\rightarrow g_{\ast}{\Bbb
Q}_{\ell ,Y}[2](1)
$$
et donc une fl\`{e}che
$$
c_{S}=\bigoplus_{n}{}^{\rm p}{\cal H}^{n}(g_{\ast}(c)):
\bigoplus_{n}{}^{\rm p}{\cal H}^{n}(g_{\ast}{\Bbb Q}_{\ell
,Y})\rightarrow \bigoplus_{n}{}^{\rm p}{\cal H}^{n+2}(g_{\ast}{\Bbb
Q}_{\ell ,Y})(1).
$$

\thm COROLLAIRE A.2.4
\enonce
Le triangle distingu\'{e}
$$
g_{\ast}A\rightarrow g_{\ast}B\rightarrow g_{\ast}C\rightarrow
g_{\ast}A[1]
$$
induit une suite exacte courte
$$
0\rightarrow \bigoplus_{n}{}^{\rm p}{\cal H}^{n}(B) \rightarrow
\bigoplus_{n}{}^{\rm p}{\cal H}^{n}(C) \rightarrow \bigoplus_{n}{}^{\rm
p}{\cal H}^{n+1}(A)\rightarrow 0
$$
en ${\Bbb Q}_{\ell ,S}[t]$-modules gradu\'{e}s en faisceaux pervers
$\ell$-adiques sur $S$.  De plus, cette suite exacte courte peut
\^{e}tre d\'{e}crite comme suit.  Si on note $D=\bigoplus_{n}{}^{\rm
p}{\cal H}^{n}(g_{\ast}{\Bbb Q}_{\ell ,Y})$ gradu\'{e} par $n$, et
$D[t]=D\otimes_{{\Bbb Q}_{\ell ,S}}{\Bbb Q}_{\ell ,S}[t]$ avec la
graduation diagonale, alors la suite ci-dessus est canoniquement
isomorphe \`{a} la suite exacte
$$
0\rightarrow D[t]\oplus tD[t]\,\smash{\mathop{\hbox to
6mm{\rightarrowfill}} \limits^{\scriptstyle \alpha}}\,D[t]\oplus D[t]
\,\smash{\mathop{\hbox to 6mm{\rightarrowfill}} \limits^{\scriptstyle
\beta}}\,D \rightarrow 0
$$
o\`{u} les fl\`{e}ches $\alpha$ et $\beta$ sont donn\'{e}es par
$$
\alpha :d_{0}(t)\oplus td_{1}(t)\mapsto
(d_{0}(t)+(t+c_{S})d_{1}(t))\oplus (d_{0}(t)-(t+c_{S})d_{1}(t))
$$
et
$$
\beta:d_{0}(t)\oplus d_{\infty}(t)\mapsto
d_{0}(-c_{S})-d_{\infty}(-c_{S})
$$
o\`{u} $c_{S}d(t)=\sum_{n}c_{S}(d_{n})t^{n}$ et
$d(c_{S})=\sum_{n}c_{S}^{n}(d_{n})$ pour tout
$d(t)=\sum_{i}d_{n}t^{n}\in D[t]$.
\endthm

\rem D\'{e}monstration
\endrem
La partie $\varepsilon_{Y}^{\ast}(c)$ de $\widetilde{c}$ \'{e}tant une
fl\`{e}che entre cohomologies de degr\'{e}s diff\'{e}rents, elle ne
compte pas.
\hfill\hfill$\square$
\vskip 3mm

Comme ci-dessus soient $S$ un sch\'{e}ma, $g:Y\rightarrow S$ un
$S$-sch\'{e}ma propre et $p:P\rightarrow Y$ un fibr\'{e} en droites
projectives, muni des deux sections $\sigma_{0},\sigma_{\infty}:
Y\hookrightarrow P$ et de l'action par homoth\'{e}ties de ${\Bbb
G}_{{\rm m},Y}$.

On consid\`{e}re alors un {\og}{pincement ${\Bbb G}_{{\rm
m},S}$-\'{e}quivariant de $P$ le long de ces deux sections}{\fg},
c'est-\`{a}-dire un $S$-sch\'{e}ma $f:X\rightarrow S$ muni d'une
action de ${\Bbb G}_{{\rm m},S}$, d'un morphisme ${\Bbb G}_{{\rm
m},S}$-\'{e}quivariant de $S$-sch\'{e}mas $\rho :P \rightarrow X$ qui
s'ins\`{e}re dans le diagramme
$$
\xymatrix{
& P\ar[dl]^p \ar[dd]^{h}\ar[dr]^{\rho} &\cr
Y \ar[dr]_{g}\ar@<4pt>@/^1pc/[ur]^{\sigma_{0},\sigma_{\infty}} 
\ar@/^1pc/[ur] & & X
\ar[dl]^{f}\cr
&S& }
$$
\vskip 1mm

\item{-}$\rho$ est fini,
\vskip 1mm

\item{-}$\rho$ induit un isomorphisme de chacune des images de
$\sigma_{0}$ et $\sigma_{\infty}$ sur le lieu des points fixes
$X^{{\Bbb G}_{{\rm m}}}$ de ${\Bbb G}_{{\rm m}}$ dans $X$,
\vskip 1mm

\item{-} $\rho$ induit un isomorphisme de $P-(\sigma_{0}(Y)\cup
\sigma_{\infty}(Y))$ sur $X-X^{{\Bbb G}_m}$,
\vskip 1mm

\noindent Notons qu'il s'agit d'un pincement au-dessus de $S$ et
non au-dessus de $Y$ de sorte qu'on ne se donne pas a priori un
morphisme $X\rightarrow Y$.

On se donne plut\^{o}t une $S$-immersion ferm\'{e}e
$i:Y\hookrightarrow X$ d'image $X^{{\Bbb G}_m}$ le lieu des points
fixes de ${\Bbb G}_m$ agissant sur $X$. Cette immersion ferm\'{e}e
est trivialement ${\Bbb G}_m$-\'{e}quivariante. La donn\'{e}e de $i$
induit la donn\'{e}e de deux automorphismes $\iota_{0}$ et
$\iota_{\infty}$ de $Y$ tels que les deux carr\'{e}s
$$
\xymatrix{Y\ar[d]_{\iota_{0}}\ar@{^{(}->}[r]^{\sigma_{0}} &
P\ar[d]^{\rho}\cr
Y\ar@{^{(}->}[r]_{i} & X}\hskip 20mm
\xymatrix{Y\ar[d]_{\iota_{\infty}}\ar@{^{(}->}[r]^{\sigma_{\infty}} &
P\ar[d]^{\rho}\cr
Y\ar@{^{(}->}[r]_{i} & X}
$$
soient commutatifs.
\vskip 2mm

On munit en plus les $S$-sch\'{e}mas $Y$, $P$ et $X$ d'une action de
${\Bbb Z}^{I}$ pour un ensemble fini $I$, qui commutent aux actions de
${\Bbb G}_{{\rm m},S}$ et pour laquelle les morphismes $f$, $g$, $p$
et $i$ sont \'{e}quivariants.

On fait alors les hypoth\`{e}ses suivantes sur la cohomologie
$\ell$-adique:
\vskip 1mm

\item{-} l'action de ${\Bbb Z}^{I}$ sur les faisceaux pervers
$\ell$-adiques ${}^{{\rm p}}{\cal H}^{n}(f_{\ast}{\Bbb Q}_{\ell ,X})$
et ${}^{{\rm p}}{\cal H}^{n}(g_{\ast}{\Bbb Q}_{\ell ,Y})$ se
factorise \`{a} travers le quotient ${\Bbb
Z}^{I}\twoheadrightarrow ({\Bbb Z}/2{\Bbb Z})^{I}$,
\vskip 1mm

\item{-} l'action de l'automorphisme $\iota_{0}$ (resp. $\iota_{\infty}$)
sur chaque ${}^{{\rm p}}{\cal H}^{n}(g_{\ast}{\Bbb Q}_{\ell ,Y})$
co\"{\i}ncide avec l'action d'un \'{e}l\'{e}ment $e_{0}$ (resp.
$e_{\infty}$) de $({\Bbb Z}/2{\Bbb Z})^{I}$.

\thm PROPOSITION A.2.5
\enonce
Soit $\kappa:({\Bbb Z}/2{\Bbb Z})^{I}\rightarrow {\Bbb
Q}_{\ell}^{\times}$ un caract\`{e}re tel que $\kappa
(e_{0})\not=\kappa (e_{\infty})$.  Alors la fl\`{e}che de restriction
$i^\ast$ en cohomologie \'{e}quivariante induit sur la partie
$\kappa$-isotypique une fl\`{e}che injective
$$
\bigoplus_{n}{}^{{\rm p}}{\cal H}^{n}(f^{{\Bbb G}_{{\rm
m},S}}_\ast{\Bbb Q}_{\ell ,X})_\kappa \rightarrow
\Bigl(\bigoplus_{m}{}^{{\rm p}}{\cal H}^{m}(g_\ast{\Bbb Q}_{\ell
,Y})_\kappa\Bigr)[t]=\Bigl(\bigoplus_{m}{}^{{\rm p}}{\cal
H}^{m}(g_\ast{\Bbb Q}_{\ell ,Y})_\kappa\Bigr)\otimes_{{\Bbb Q}_{\ell
,S}}{\Bbb Q}_{\ell ,S}[t]
$$
d'image
$$
(t+c_{S})\Bigl(\bigoplus_{m}{}^{{\rm p}}{\cal H}^{m}(g_\ast{\Bbb Q}_{\ell
,Y})_\kappa\Bigr)[t]\subset \Bigl(\bigoplus_{m}{}^{{\rm p}}{\cal
H}^{m}(g_\ast{\Bbb Q}_{\ell ,Y})_\kappa\Bigr)[t]
$$
o\`{u} $c_{S}$ est la fl\`{e}che d\'{e}finie avant le corollaire {\rm A.2.4}.
\endthm

\rem D\'{e}monstration
\endrem
La fl\`{e}che d'adjonction ${\Bbb Q}_{\ell,
X}\rightarrow\rho_{\ast}{\Bbb Q}_{\ell ,P}$ induit par restriction au
ferm\'{e} $i:Y\hookrightarrow X$ la fl\`{e}che
$$\eqalign{
{\Bbb Q}_{\ell ,Y}=i^{\ast}{\Bbb Q}_{\ell ,X}\rightarrow
i^{\ast}\rho_{\ast}{\Bbb Q}_{\ell
,P}&=i^{\ast}\rho_{\ast}\sigma_{0,\ast}{\Bbb Q}_{\ell ,Y}\oplus
i^{\ast}\rho_{\ast}\sigma_{\infty,\ast}{\Bbb Q}_{\ell
,Y}\cr
&=i^{\ast}i_{\ast}\iota_{0,\ast}{\Bbb Q}_{\ell ,Y}\oplus
i^{\ast}i_{\ast}\iota_{\infty,\ast}{\Bbb Q}_{\ell
,Y}\cr
&=\iota_{0,\ast}{\Bbb Q}_{\ell ,Y}\oplus\iota_{\infty ,\ast}{\Bbb
Q}_{\ell ,Y}}
$$
qui est compos\'{e}e de la diagonale ${\Bbb Q}_{\ell ,Y}\rightarrow {\Bbb
Q}_{\ell ,Y}\oplus{\Bbb Q}_{\ell ,Y}$ et de la somme directe des
fl\`{e}ches d'adjonction ${\Bbb Q}_{\ell
,Y}\rightarrow\iota_{0,\ast}\iota_{0}^{\ast}{\Bbb Q}_{\ell
,Y}=\iota_{0,\ast}{\Bbb Q}_{\ell ,Y}$ et ${\Bbb Q}_{\ell
,Y}\rightarrow\iota_{\infty ,\ast}\iota_{\infty}^{\ast}{\Bbb Q}_{\ell
,Y}=\iota_{\infty,\ast}{\Bbb Q}_{\ell ,Y}$. On a donc un morphisme de
triangles distingu\'{e}s
$$
\xymatrix{f_{\ast}^{{\Bbb G}_{{\rm m},S}}{\Bbb Q}_{\ell ,
X}\ar[d]_{i^{\ast}}\ar[r] & h_{\ast}^{{\Bbb G}_{{\rm m},S}}{\Bbb
Q}_{\ell ,P}\ar[d]_{\sigma_{0}^{\ast}\oplus\sigma_{\infty}^{\ast}}\ar[r]
& g_{\ast}^{{\Bbb G}_{{\rm m},S}}{\Bbb Q}_{\ell , Y}\ar@{=}[d]\ar[r] &\cr
g_{\ast}^{{\Bbb G}_{{\rm m},S}}{\Bbb Q}_{\ell , Y}\ar[r]_-{u} &
g_{\ast}^{{\Bbb G}_{{\rm m},S}}{\Bbb Q}_{\ell , Y}\oplus
g_{\ast}^{{\Bbb G}_{{\rm m},S}}{\Bbb Q}_{\ell , Y}\ar[r]_-{v}&
g_{\ast}^{{\Bbb G}_{{\rm m},S}}{\Bbb Q}_{\ell , Y}\ar[r] &}
$$
o\`{u} $h=f\circ\rho=g\circ p$, $u=\pmatrix{\iota_{0}^{\ast}\cr
\iota_{\infty}^{\ast}\cr}$ et $v=(\iota_{\infty}^{\ast},
-\iota_{0}^{\ast})$, d'o\`{u} en passant \`{a} la cohomologie un
digramme \`{a} lignes exactes
$$
\xymatrix{E\ar[d]_{\alpha'}\ar[r]^-{u'} & D[t]\oplus tD[t]
\ar[d]_{\alpha}\ar[r]^-{v'} &
D[t]\ar@{=}[d]\cr
D[t]\ar[r]_-{u}^{} & D[t]\oplus D[t]\ar[r]_-{v}^{} & D[t] }
$$
o\`{u} on a pos\'{e} $E=\bigoplus_{n}{}^{{\rm p}}{\cal H}^{n}(f^{{\Bbb
G}_{{\rm m},S}}_\ast{\Bbb Q}_{\ell ,X})$ et $D=\bigoplus_{n}{}^{{\rm
p}}{\cal H}^{n}(g_\ast{\Bbb Q}_{\ell ,X})$, o\`{u} d'apr\`{e}s le
corollaire pr\'{e}c\'{e}dent, on a $D[t]\oplus
tD[t]=\bigoplus_{n}{}^{{\rm p}}{\cal H}^{n}(h^{{\Bbb G}_{{\rm
m},S}}_\ast{\Bbb Q}_{\ell ,X})$ et on a une formule pour la fl\`{e}che
$\alpha$ qui montre en particulier que $\alpha$ est injective.

Prenons la partie $\kappa$ de ce diagramme
$$
\xymatrix{E_{\kappa}\ar[d]_{\alpha'}\ar[r]^-{u'} & D_{\kappa}[t]
\oplus tD_{\kappa}[t] \ar[d]_{\alpha}\ar[r]^-{v'} & D_{\kappa}[t]
\ar@{=}[d]\cr
D_{\kappa}[t]\ar[r]_-{u}^{} & D_{\kappa}[t]\oplus D_{\kappa}[t]
\ar[r]_-{v}^{} & D_{\kappa}[t]}
$$
Par hypoth\`{e}se, $u$ a maintenant pour composantes les
multiplications par $\kappa (e_{0})$ et $\kappa (e_{\infty})$, c'est-\`{a}-dire
$u$ est au signe pr\`{e}s l'anti-diagonale, et $v$ est donc la
fl\`{e}che $(-\kappa (e_{\infty}),\kappa (e_{0}))$, c'est-\`{a}-dire au signe
pr\`{e}s la somme.  On en d\'{e}duit que la fl\`{e}che compos\'{e}e
$v'=v\circ\alpha$ envoie $d_{0}(t)\oplus d_{1}(t)$ sur $\pm 2d_{0}(t)$
et est par cons\'{e}quent surjective. En consid\'{e}rant la suite
exacte longue de cohomologie perverse du triangle distingu\'{e} de la
ligne du haut de diagramme ci-dessus, on obtient que $u'$ est injective. Ceci
montre que la fl\`{e}che de restriction $\alpha'$ de l'\'{e}nonc\'{e}
est injective. En outre, son image est l'image r\'{e}ciproque par $u$
de l'image de $\alpha$, c'est-\`{a}-dire
$$
\{d(t)\in D[t]\mid \kappa (e_{0})d(t)-\kappa (e_{\infty})d(t)\in
(t+c_{S})D[t]\}=(t+c_{S})D[t],
$$
ce que l'on voulait d\'{e}montrer.
\hfill\hfill$\square$

\subsection{A.3}{Une formule des points fixes}

Soit $X$ un sch\'{e}ma propre sur $k={\Bbb F}_{q}$ et $G$ un
$k$-sch\'{e}ma en groupes commutatifs lisse de type fini qui agit sur
le $k$-sch\'{e}ma $X$.  Notons $G^{0}$ la composante neutre de $G$ et
$\pi_{0}(G)=G/G^{0}$ son $k$-sch\'{e}ma en groupes des composantes
connexes.  On suppose que le $k$-sch\'{e}ma en groupes fini \'{e}tale
$\pi_{0}(G)$ est d\'{e}ploy\'{e} sur $k$.  Autrement dit on suppose
que chaque composante connexe de $G$ admet un point rationnel sur $k$.
On a alors
$$
\pi_{0}(G)(k)=\pi_{0}(G)(\overline{k})=G(\overline{k})/G^{0}(\overline{k})
$$
et la suite exacte
$$
0\rightarrow G^{0}(k)\rightarrow G(k)\rightarrow \pi_{0}(G)(k)\rightarrow
0.
$$
Le th\'{e}or\`{e}me de Lang assure que
$$
G^{0}(\overline{k})=\{{\cal L}_{q}(g):=\mathop{\rm Frob}\nolimits_{q}
(g)g^{-1}\mid g\in G(\overline{k})\}.
$$

On consid\`{e}re le $k$-champ alg\'{e}brique quotient $[X/G]$, la
cat\'{e}gorie $[X/G](\overline{k})$ des $\overline{k}$-points de ce
champ et sa cat\'{e}gorie $[X/G](k)$ des $k$-points.  La cat\'{e}gorie
$[X/G](\overline{k})$ est la cat\'{e}gorie dont les objets sont les
$x\in X(\overline{k})$ et dont les fl\`{e}ches de $x$ vers $x'$ sont
les $g\in G(\overline{k})$ tels que $g\cdot x=x'$.  La cat\'{e}gorie
$[X/G](k)$ a pour objets les couples $(x,g)\in
X(\overline{k})\times G(\overline{k})$ tels que $\mathop{\rm
Frob}\nolimits_{q}(x)=g\cdot x$, et pour morphismes de $(x,g)$ dans un
autre objet $(x',g')$ les $h\in G(\overline{k})$ tels que $h\cdot
x=x'$ et que ${\cal L}_{q}(h)g=g'$.

Le groupe des automorphismes d'un objet $x$ de $[X/G](\overline{k})$
est \'{e}gal \`{a}
$$
\mathop{\rm Aut}(x)=G_{x}(\overline{k}).
$$
Si $(x,g)\in [X/G](k)$, comme $G$ est commutatif le fixateur
$G_{x}\subset G$ de $x$ est stable par $\mathop{\rm
Frob}\nolimits_{q}: G\rightarrow G$ et est donc d\'{e}fini sur $k$; le
groupe des automorphismes de l'objet $(x,g)$ dans $[X/G](k)$ est le
groupe fini
$$
\mathop{\rm Aut}(x,g)=\mathop{\rm Aut}(x)^{\mathop{\rm Frob}
\nolimits_{q}}=G_{x}(k).
$$

L'ensemble des classes d'isomorphie d'objets de la cat\'{e}gorie
$[X/G](k)$ est l'ensemble quotient
$$
[X/G](k)_{\sharp}=\{(x,g)\in X(\overline{k})\times G(\overline{k})
\mid\mathop{\rm Frob}\nolimits_{q}(x)=g\cdot x\}/G(\overline{k})
$$
o\`{u} l'action de $h\in G(\overline{k})$ est donn\'{e}e par $h\cdot
(x,g)=(h\cdot x,{\cal L}_{q}(h)g)$.  On note aussi $[X/G](k)_{\sharp}$
un syst\`{e}me de repr\'{e}sentants de ces classes d'isomorphie.  On a
une application
$$
\mathop{\rm cl}:[X/G](k)_{\sharp}\rightarrow\pi_{0}(G)(\overline{k})=
\pi_{0}(G)(k)
$$
qui envoie la classe d'isomorphie de $(x,g)$ sur la composante connexe
de $G$ contenant $g$. Cette application est bien d\'{e}finie car ${\cal
L}_{q}(h)\in G_{0}(\overline{k})$ pour tout $h\in G(\overline{k})$.

Le groupe fini $G(k)$ agit sur $X$ et donc sur les groupes de
cohomologie $\ell$-adique $H^{n}(\overline{k} \otimes_{k}X,{\Bbb
Q}_{\ell})$, cette derni\`{e}re action commutant \`{a} l'action par
transport de structure de $\mathop{\rm Gal}(\overline{k}/k)$. Pour
chaque caract\`{e}re $\chi :\pi_{0}(G)(k) \rightarrow {\Bbb
Q}_{\ell}^{\times}$ on peut donc consid\'{e}rer la partie
$\chi$-isotypique $H^{n}(\overline{k} \otimes_{k}X,{\Bbb
Q}_{\ell})_{\chi}$, munie de l'action de l'\'{e}l\'{e}ment de
Frobenius g\'{e}om\'{e}trique $\mathop{\rm Frob}\nolimits_{q}\in
\mathop{\rm Gal}(\overline{k}/k)$.

\thm PROPOSITION A.3.1
\enonce
Pour tout caract\`{e}re $\chi :\pi_{0}(G)(k)\rightarrow {\Bbb
Q}_{\ell}^{\times}$, on a la formule des points fixes
$$
\sum_{n}(-1)^{n}\mathop{\rm Tr}(\mathop{\rm Frob}\nolimits_{q},
H^{n}(\overline{k}\otimes_{k}X,{\Bbb Q}_{\ell})_{\chi})=
|G^{0}(k)|\sum_{(x,g)\in [X/G](k)_{\sharp}}{\chi (\mathop{\rm cl}
(x,g))\over |\mathop{\rm Aut}(x,g)|}.
$$
\endthm

\rem Remarque
\endrem
L'ensemble $[X/G](k)_{\sharp}$ et les cardinaux $|\mathop{\rm
Aut}(x,g)|$ ne d\'{e}pendent que de la cat\'{e}gorie
$[X/G](\overline{k})$ et du foncteur $\mathop{\rm Frob}\nolimits_{q}:
[X/G](\overline{k})\rightarrow [X/G](\overline{k})$.
\hfill\hfill$\square$
\vskip 3mm

Pour d\'{e}montrer cette proposition nous aurons besoin du lemme suivant.

\thm LEMME A.3.2
\enonce
{\rm (i)} Pour tout caract\`{e}re $\chi : G(k)\rightarrow {\Bbb
Q}_{\ell}^{\times}$, on a
$$
\sum_{n}(-1)^{n}\mathop{\rm Tr}(\mathop{\rm Frob}\nolimits_{q},
H^{n}(\overline{k}\otimes_{k}X,{\Bbb Q}_{\ell})_{\chi}) ={1\over
|G(k)|}\sum_{g\in G(k)}\chi (g)|X^{\mathop{\rm Frob}
\nolimits_{q}\circ g^{-1}}|
$$
o\`{u}
$$
X^{\mathop{\rm Frob}\nolimits_{q}\circ g^{-1}}= \{x\in
X(\overline{k})\mid \mathop{\rm Frob}\nolimits_{q}(x)= g\cdot x\}.
$$

\decale{\rm (ii)} La fonction $g\rightarrow |X^{\mathop{\rm
Frob}\nolimits_{q}\circ g^{-1}}|$ est constante sur chacune des
composantes connexes de $G$.

\decale{\rm (iii)} La trace altern\'{e}e ci-dessus est non nulle
seulement si $\chi$ se factorise \`{a} travers $G(k)\twoheadrightarrow
\pi_{0}(G)(k)$ et si c'est le cas on a encore
$$
\sum_{n}(-1)^{n}\mathop{\rm Tr}(\mathop{\rm Frob}\nolimits_{q},
H^{n}(\overline{k}\otimes_{k}X,{\Bbb Q}_{\ell})_{\chi}) ={1\over
|\pi_{0}(G)(k)|}\sum_{\gamma\in \pi_{0}(G)(k)}\chi (\dot\gamma)|
X^{\mathop{\rm Frob} \nolimits_{q}\circ \dot\gamma^{-1}}|
$$
o\`{u}, pour chaque $\gamma\in\pi_{0}(G)(k)$, $\dot\gamma\in G(k)$ est
n'importe quel \'{e}l\'{e}ment de la composante connexe $\gamma$.
\endthm

\rem D\'{e}monstration
\endrem
On a par d\'{e}finition de la partie $\chi$-isotypique
$$\displaylines{
\qquad\sum_{n}(-1)^{n}\mathop{\rm Tr}(\mathop{\rm Frob}
\nolimits_{q},H^{n}(\overline{k}\otimes_{k}X,{\Bbb Q}_{\ell})_{\chi})
\hfill\cr\hfill
={1\over |G(k)|}\sum_{g\in G(k)}\chi (g)\mathop{\rm Tr}
(\mathop{\rm Frob}\nolimits_{q}\circ g^{-1},H^{n}(\overline{k}
\otimes_{k}X,{\Bbb Q}_{\ell})).\qquad}
$$
Or $g$ est d'ordre fini puisque $G(k)$ est un groupe fini.  On peut
donc appliquer la variante due \`{a} Deligne et Lusztig de la formule
des points fixes de Grothendieck (cf.  [De-Lu]) et l'assertion (i) du
lemme est d\'{e}montr\'{e}.

Pour d\'{e}montrer l'assertion (ii) il suffit de remarquer que, pour
chaque $h\in G(\overline{k})$ tel que ${\cal L}_{q}(h)\in G^{0}(k)$,
on a une bijection
$$
X^{\mathop{\rm Frob}\nolimits_{q}\circ g^{-1}}
\buildrel\sim\over\longrightarrow
X^{\mathop{\rm Frob}\nolimits_{q}\circ ({\cal
L}_{q}(h)g^{-1})},~x\mapsto h\cdot x.
$$

D'apr\`{e}s le lemme d'homotopie 3.2.3, $G^{0}(k)$ agit trivialement
sur chaque groupe de cohomologie $\ell$-adique
$H^{n}(\overline{k}\otimes_{k}X,{\Bbb Q}_{\ell})$.  Par suite, la
composante $\chi$-isotypique est non nulle seulement si $\chi$ est
trivial sur $G^{0}(k)\subset G(k)$, c'est-\`{a}-dire se factorise par
$G(k)\twoheadrightarrow \pi_{0}(G)(k)$, d'o\`{u} l'assertion (iii).
\hfill\hfill$\square$
\vskip 3mm

\rem D\'{e}monstration de proposition
\endrem
Soit $\dot\gamma\in G(k)$ dans la composante connexe
$\gamma\in\pi_{0}(G)(k)$.  Compte tenu du lemme il s'agit de montrer
que
$$
|X^{\mathop{\rm Frob}\nolimits_{q}\circ \dot\gamma^{-1}}|
=\sum_{{\scriptstyle (x,g)\in [X/G](k)_{\sharp}\atop\scriptstyle
\mathop{\rm cl}(x,g)=\gamma}}{|G(k)|\over |G_{x}(k)|}.
$$

Or l'application
$$
\varphi_{\dot\gamma}:X^{\mathop{\rm Frob}\nolimits_{q}\circ
\dot\gamma^{-1}}\rightarrow \mathop{\rm cl}
\nolimits^{-1}(\gamma)\subset [X/G](k)_{\sharp}
$$
qui associe \`{a} $x$ la classe d'isomorphie de $(x,\dot\gamma )$, est
surjective puisque, pour tout $(x,g)\in \mathop{\rm cl}
\nolimits^{-1}(\gamma)$, il existe $h\in G^{0}(\overline{k})$ tel que
$g={\cal L}_{q}(h)\dot\gamma$ et on a
$$
(x,g)=h\cdot (h^{-1}\cdot x,\dot\gamma )
$$
o\`{u} $h^{-1}\cdot x\in X^{\mathop{\rm Frob}\nolimits_{q}\circ
\dot\gamma^{-1}}$.

Il ne reste donc plus qu'\`{a} v\'{e}rifier que le cardinal de la
fibre passant par $x$ de l'application $\varphi_{\dot\gamma}$ est
pr\'{e}cis\'{e}ment ${|G(k)|\over |G_{x}(k)|}$.  Mais cette fibre
$$
\{x'\in X^{\mathop{\rm Frob}\nolimits_{q}\circ \dot\gamma^{-1}}\mid
\exists h\in G(\overline{k})\hbox{ tel que }x'=h\cdot x\hbox{ et
}{\cal L}_{q}(h)=0\}
$$
est isomorphe \`{a} $G(k)/G_{x}(k)$, d'o\`{u} la conclusion.
\hfill\hfill$\square$

\vskip 10mm

\centerline{{\bf Bibliographie}}
\vskip 5mm

\newtoks\ref \newtoks\auteur \newtoks\titre \newtoks\annee
\newtoks\revue \newtoks\tome \newtoks\pages \newtoks\reste

\def\bibitem#1{\parindent=20pt\itemitem{#1}\parindent=24pt}

\def\article{\bibitem{[\the\ref]}%
\the\auteur~-- \the\titre, {\sl\the\revue} {\bf\the\tome},
({\the\annee}), \the\pages.\smallskip\filbreak}

\def\autre{\bibitem{[\the\ref]}%
\the\auteur~-- \the\reste.\smallskip\filbreak}

\ref={A-I-K}
\auteur={A. {\pc ALTMAN}, A. {\pc IARROBINO}, S. {\pc KLEIMAN}}
\reste={Irreducibility of the Compactified Jacobian, dans ``Real and
complex singularities. Proceedings, Oslo 1976, P. Holm (ed.)'',
Sijthoff \& Nordhoff, (1977), 1-12}
\autre

\ref={Al-Kl}
\auteur={A.B. {\pc ALTMAN}, S.L. {\pc KLEIMAN}}
\reste={The Presentation Functor and the Compactified Jacobian, dans
{\it The Grothendieck Festschrift, Volume I}, Birkh\"{a}user, (1990),
15-32}
\autre

\ref={B-B-D}
\auteur={A.A. {\pc BEILINSON}, J. {\pc BERNSTEIN}, P. {\pc DELIGNE}}
\reste={Faisceaux pervers, {\it Analyse et topologie sur les
espaces singuliers}, Ast\'{e}risque {\bf 100}, (1982)}
\autre

\ref={B-N-R}
\auteur={A. {\pc BEAUVILLE}, M.S. {\pc NARASIMHAN}, S. {\pc
RAMANAN}}
\titre={Spectral curve and the generalised theta divisor}
\revue={J. reine angew. Math.}
\tome={398}
\annee={1989}
\pages={169-179}
\article

\ref={Del}
\auteur={P. {\pc DELIGNE}}
\titre={La conjecture de Weil. II}
\revue={Publications Math\'{e}\-ma\-ti\-ques de
l'I.H.\'{E}.S.}
\tome={52}
\annee={1980}
\pages={313-428}
\article

\ref={De-Lu}
\auteur={P. {\pc DELIGNE}, G. {\pc LUSZTIG}}
\titre={Representations of reductive groups over finite fields}
\revue={Ann. of Math.}
\tome={103}
\annee={1976}
\pages={103-161}
\article

\ref={Est}
\auteur={E. {\pc ESTEVES}}
\titre={Compactifying the relative Jacobian over families of
reduced curves}
\revue={Trans.  Amer.  Math.  Soc.}
\tome={353}
\annee={2001}
\pages={3045-3095}
\article

\ref={Fal}
\auteur={G. {\pc FALTINGS}}
\titre={Stable $G$-bundles and projective connections}
\revue={J. Alg. Geom.}
\tome={2}
\annee={1993}
\pages={507-568}
\article

\ref={F-G-S}
\auteur={D. {\pc FANTECHI}, L. {\pc G\"{O}TTSCHE}, D. van {\pc STRATEN}}
\titre={Euler Number of the Compactified Jacobian and Multiplicity of 
Rational Curves}
\revue={J. Algebraic Geometry}
\tome={8}
\annee={1999}
\pages={115-133}
\article

\ref={Gab}
\auteur={O. {\pc GABBER}}
\titre={On space filling curves and Albanese varieties}
\revue={Geom.  Funct.  Anal.}
\tome={11}
\annee={2001}
\pages={1192Ð1200}
\article

\ref={G-K-M}
\auteur={M. {\pc GORESKY}, R. {\pc KOTTWITZ}, R. {\pc MACPHERSON}}
\reste={Homology of affine Springer fibers in the unramified case,
http://arxiv.org/abs/math.RT/0305144, (2003)}
\autre

\ref={Hal}
\auteur={T. {\pc HALES}}
\reste={A simple definition of the transfer factors for unramified
groups, dans {\it Representation theory of groups and algebras},
Contemp.  Math. {\bf 145}, (1993), 109-134}
\autre

\ref={Hit}
\auteur={N. {\pc HITCHIN}}
\titre={Stable bundles and integrable systems}
\revue={Duke Math. J.}
\tome={54}
\annee={1987}
\pages={91-114}
\article

\ref={Ka-Lu}
\auteur={D. {\pc KAZHDAN}, G. {\pc LUSZTIG}}
\titre={Fixed Point Varieties on Affine Flag Manifolds}
\revue={Israel J. of Math.}
\tome={62}
\annee={1988}
\pages={129-168}
\article

\ref={Kot 1}
\auteur={R.E. {\pc KOTTWITZ}}
\titre={Stable formula: Elliptic Singular Terms}
\revue={Math. Ann.}
\tome={275}
\annee={1986}
\pages={365-399}
\article

\ref={Kot 2}
\auteur={R.E. {\pc KOTTWITZ}}
\titre={Transfer factors for Lie algebras}
\revue={Represent. Theory}
\tome={3}
\annee={1999}
\pages={127-138}
\article

\ref={Lan}
\auteur={R.P. {\pc LANGLANDS}}
\reste={Les d\'{e}buts d'une formule des traces stable, Publications
math\'{e}matiques de l'Universit\'{e} Paris VII, (1979)}
\autre

\ref={La-Sh}
\auteur={R.P. {\pc LANGLANDS}, D. {\pc SHELSTAD}}
\titre={On the definition of the transfer factors}
\revue={Math. Ann.}
\tome={278}
\annee={1987}
\pages={219-271}
\article

\ref={Lau}
\auteur={G. {\pc LAUMON}}
\reste={Sur le lemme fondamental pour les groupes
unitaires, math.AG/ 0212245, (2002)}
\autre

\ref={La-Ra}
\auteur={G. {\pc LAUMON}, M. {\pc RAPOPORT}}
\reste={A geometric approach to the fundamental lemma for unitary
groups, math.AG/9711021, (1997)}
\autre

\ref={Ngo}
\auteur={B. C. {\pc NG\^{O}}}
\reste={Fibration de Hitchin et endoscopie, en pr\'{e}paration}
\autre

\ref={Poo}
\auteur={B. {\pc POONEN}}
\reste={Bertini theorems over finite fields,
math.AG/0204002, (2002)}
\autre

\ref={Reg}
\auteur={C.J. {\pc REGO}}
\titre={The Compactified Jacobian}
\revue={Ann. Scient. \`{E}c. Norm. Sup.}
\tome={13}
\annee={1980}
\pages={211-223}
\article

\ref={Ser}
\auteur={J.-P. {\pc SERRE}}
\reste={Corps locaux, Hermann, (1968)}
\autre

\ref={Wal 1}
\auteur={J.-L. {\pc WALDSPURGER}}
\reste={Comparaison d'int\'{e}grales orbitales pour des groupes
$p$-adiques, dans {\it Proceedings of the International Congress of
Mathematicians (Z\"{u}rich, 1994), Vol.  1}, Birkh\"{a}user, (1995),
807-816}
\autre

\ref={Wal 2}
\auteur={J.-L. {\pc WALDSPURGER}}
\titre={Homog\'{e}n\'{e}it\'{e} de certaines distributions sur les groupes
$p$-adiques}
\revue={Publ. Math. Inst. Hautes \'{E}tudes Sci.}
\tome={81}
\annee={1995}
\pages={25-72}
\article

\ref={Wal 3}
\auteur={J.-L. {\pc WALDSPURGER}}
\reste={Endoscopie et changement de caract\'{e}ristiques,
pr\'{e}\-publi\-ca\-tion, (2004)}
\autre

\ref={Wal 4}
\auteur={J.-L. {\pc WALDSPURGER}}
\reste={Int\'{e}grales orbitales nilpotentes et endoscopie pour les
groupes classiques non ramifi\'{e}s, Ast\'{e}risque {\bf 269}, (2001)}
\autre

\ref={EGA II}
\auteur={A. {\pc GROTHENDIECK}}
\reste={\'{E}l\'{e}ments de g\'{e}om\'{e}trie alg\'{e}brique
(r\'{e}dig\'{e}s avec la collaboration de Jean Dieudonn\'{e}): II.
\'{E}tude globale \'{e}l\'{e}mentaire de quelques classes de
morphismes, {\it Publications Math\'{e}\-ma\-ti\-ques de
l'I.H.\'{E}.S.} {\bf 8}, (1961)}
\autre

\ref={EGA IV}
\auteur={A. {\pc GROTHENDIECK}}
\reste={\'{E}l\'{e}ments de g\'{e}om\'{e}trie alg\'{e}brique
(r\'{e}dig\'{e}s avec la collaboration de Jean Dieudonn\'{e}): IV.
\'{E}tude locale des sch\'{e}mas et des morphismes de sch\'{e}mas,
Troisi\`{e}me partie, {\it Publications Math\'{e}\-ma\-ti\-ques de
l'I.H.\'{E}.S.} {\bf 28}, (1966)}
\autre

\ref={SGA 3}
\auteur={A. {\pc DEMAZURE}, A. {\pc GROTHENDIECK}}
\reste={Sch\'{e}mas en groupes, Groupes de Type Multiplicatif et
Structure des Sch\'{e}mas en Groupes G\'{e}n\'{e}raux, Lecture Notes
in Math. {\bf 152}, (1970)}
\autre

\bye